7/21/2014

# Compactness theorems for SL(2; ℂ) generalizations of the 4-dimensional anti-self dual equations

Clifford Henry Taubes†

Department of Mathematics
Harvard University
Cambridge, MA 02138

chtaubes@math.harvard.edu

Uhlenbeck's compactness theorem can be used to analyze sequences of connections with anti-self dual curvature on principal SU(2) bundles over oriented 4-dimensional manifolds. The theorems in this paper give an extension of Uhlenbeck's theorem for sequences of solutions of certain SL(2;ℂ) analogs of the anti-self dual equations.

† Supported in part by the National Science Foundation

## 1. Introduction

The upcoming Theorems 1.1 and 1.2 state the promised analog of Uhlenbeck's compactness theorem for sequences of solutions to various $SL(2;\mathbb{C})$ generalizations of the self-dual equations on a given oriented, Riemannian 4-manifold. Theorems 1.3, 1.4 and 1.5 say more about the singular sets of the limits of these sequences.

### a) The equations

The equations of interest are defined on an oriented 4-manifold with a given Riemannian metric. This manifold is denoted in what follows by X. Use the metric on X to define the Hodge star operator on $\wedge^*T^*X$ and use this Hodge star to write the bundle of 2-forms as the direct sum $\wedge^2T^*X = \Lambda^+ \oplus \Lambda^-$ with $\Lambda^+$ denoting the bundle of self-dual 2-forms and with $\Lambda^-$ denoting the bundle of anti self-dual 2 forms. If $\omega$ denotes a given 2 form, then its respective self-dual and anti self-dual parts are denoted by $\omega^+$ and $\omega^-$.

The definition of the equations require the specification of a principal SO(3) bundle, this denoted by P. Supposing that A is a connection on P, then $F_A$ is used below to denote its curvature 2-form. This is a 2-form on X with values in the bundle associated to P with fiber the Lie algebra of SO(3). The Lie algebra of SO(3) is denoted by $\mathfrak{S}$ and this associated bundle by $P \times_{SO(3)} \mathfrak{S}$. The connection A defines an exterior covariant derivative on sections of $\wedge^*T^*X \times (P \times_{SO(3)} \mathfrak{S})$, this denoted by $d_A$.

The equations ask a pair $(A, \mathfrak{a})$ of connection on P and section of $T^*X \otimes (P \times_{SO(3)} \mathfrak{S})$ obey

$$(F_A - \mathfrak{a} \wedge \mathfrak{a})^+ = 0 \quad \textit{and} \quad (d_A\mathfrak{a})^- = 0 \quad \textit{and} \quad d_A{*}\mathfrak{a} = 0. \tag{1.1}$$

Witten [W1]-[W4], [GW], [KW] and also [Hay] proposed that certain linear combinations of the equations in (1.1) and the version with the self and anti-self dual forms interchanged should also be considered. The latter are parametrized by $\tau \in [0,1]$ and can be written as

- $\tau(F_A - \mathfrak{a} \wedge \mathfrak{a})^+ = (1-\tau)(d_A\mathfrak{a})^+$ ,
- $(1-\tau)(F_A - \mathfrak{a} \wedge \mathfrak{a})^- = -\tau(d_A\mathfrak{a})^-$ ,
- $d_A{*}\mathfrak{a} = 0$ .

(1.2)

The $\tau = 0$ version of (1.2) is (1.1) and the $\tau = 1$ version of (1.2) is the version of (1.1) that is defined on X with its same metric but with its orientation reversed. In the case when X is compact, Kapustin and Witten [KW] prove that the solutions to (1.2) with $\tau \in (0,1)$ exist only in the case when $P \times_{SO(3)} \mathfrak{S}$ has zero first Pontrjagin number; and that in this case, the solutions are such that $A + i\,\mathfrak{a}$ defines a flat $PSL(2;\mathbb{C})$ connection. A nice discussion of these equations can be found in [GU].

A salient special case of (1.1) and (1.2) arises when X is the product of an interval in $\mathbb{R}$ and a compact, oriented 3-manifold. To say more about the product case, introduce by way of notation M to denote the 3-manifold and let $\mathbb{P}$ denote for the moment a principal $PSl(2;\mathbb{C})$ bundle over M. Such a bundle is isomorphic to an associated principal bundle to a principal SO(3) bundle, the latter to be denoted by P. A connection $\mathbb{A}$ on P can be written using this isomorphism as $A + i\,\mathfrak{a}$ with A being a connection on P and with $\mathfrak{a}$ being a section of $T^*M \otimes (P \times_{SO(3)} \mathfrak{S})$. There is a Chern Simons function that maps the space of connections on $\mathbb{P}$ to $\mathbb{C}$, this denoted by $\mathfrak{cs}$. To define $\mathfrak{cs}$, fix a fiducal connection, $A_0$, on P. Let $\mathbb{A}$ denote a given connection on $\mathbb{P}$. Write $\mathbb{A}$ as $A_0 + \hat{a}_{\mathbb{A}}$ with $\hat{a}_{\mathbb{A}}$ being a section over M of the complexification of the bundle $(P \times_{SO(3)} \mathfrak{S}) \otimes T^*M$. The value of function $\mathfrak{cs}$ on $\mathbb{A}$ is

$$\mathfrak{cs}(\mathbb{A}) = \tfrac{1}{2} \int_M \operatorname{trace}(2\hat{a}_{\mathbb{A}} \wedge F_{A_0} + \hat{a}_{\mathbb{A}} \wedge d_{A_0}\hat{a}_{\mathbb{A}} + \tfrac{2}{3}\hat{a}_{\mathbb{A}} \wedge \hat{a}_{\mathbb{A}} \wedge \hat{a}_{\mathbb{A}})\,. \tag{1.3}$$

The gradient flow equations for the real part of $\mathfrak{cs}$ ask that a parametrized path of connection on P and section of $T^*M \otimes (P \times_{SO(3)} \mathfrak{S})$ obey

- $\frac{\partial}{\partial s} A + *(F_A - \mathfrak{a} \wedge \mathfrak{a}) = 0$ ,
- $\frac{\partial}{\partial s} \mathfrak{a} - *d_A \mathfrak{a} = 0$ ,

(1.4)

with s denoting here the Euclidean parameter along the path and with $*$ denoting here the Hodge $*$ operator on M. Let $\mathbb{I}$ denote the interval on which this path is defined. The projection $\mathbb{I} \times M \to M$ pulls P back as a principal SO(3) bundle over M and this pull-back allows $(A, \mathfrak{a})$ to be viewed as a pair of connection on P over $\mathbb{I} \times M$ and 1-form over $\mathbb{I} \times M$ with values in corresponding $P \times_{SO(3)} \mathfrak{S}$. With $(A, \mathfrak{a})$ viewed in this way, then the equations in (1.4) are identical to the left most and center equation in (1.1) for the case when X is the interior of $\mathbb{I} \times M$ with its Riemannian metric being the product of the metric $ds \otimes ds$ on $\mathbb{I}$ and the given metric on M.

The gradient flow equations for the imaginary part of $\mathfrak{cs}$ ask that $(A, \mathfrak{a})$

- $\frac{\partial}{\partial s} A - *d_A \mathfrak{a} = 0.$
- $\frac{\partial}{\partial s} \mathfrak{a} - *(F_A - \mathfrak{a} \wedge \mathfrak{a}) = 0.$

(1.5)

These equations are the $\tau = \frac{1}{2}$ version of the equations in the first two bullets of (1.2) for the case when X is the interior of $\mathbb{I} \times M$ if $(A, \mathfrak{a})$ is again viewed as a pair of connection on P's pull-back to $\mathbb{I} \times M$ and 1-form with values $P \times_{SO(3)} \mathfrak{S}$. By way of a parenthetical remark, the equations in (1.5) have a second interpretation: They are proportional to the flow equations for the Hamiltonian vector field that is defined from the differential of the real part of $\mathfrak{cs}$ using a certain canonical symplectic form on the space connections on the principal PSL(2; $\mathbb{C}$) bundle $\mathbb{P}$.

Fix $\tau \in [0, 1]$. The formal downward gradient flow equations for the real part of $\frac{(\tau + i(1-\tau))^2}{\tau^2 + (1-\tau)^2}\,\mathfrak{cs}$ ask a path $(A, \mathfrak{a})$ to obey

- $\frac{\partial}{\partial s} A + \frac{\tau^2 - (1-\tau)^2}{\tau^2 + (1-\tau)^2} *(F_A - \mathfrak{a} \wedge \mathfrak{a}) - \frac{2\tau(1-\tau)}{\tau^2 + (1-\tau)^2} *d_A\mathfrak{a} = 0$ .
- $\frac{\partial}{\partial s} \mathfrak{a} - \frac{\tau^2 - (1-\tau)^2}{\tau^2 + (1-\tau)^2} *d_A\mathfrak{a} - \frac{2\tau(1-\tau)}{\tau^2 + (1-\tau)^2} *(F_A - \mathfrak{a} \wedge \mathfrak{a}) = 0$ .

(1.6)

Some straightforward algebraic manipulations identify the equations in (1.6) with the corresponding $\tau \in [0, 1]$ version of the equations in the first two bullets of (1.2) for the case when X is again the interior of $\mathbb{I} \times M$ with its product metric.

Witten augments (1.6) with the constraint that $d_A{*}\mathfrak{a} = 0$. As Witten points out, the section $*d_A{*}\mathfrak{a}$ of $P \times_G \mathfrak{S}$ is, whether zero or not, constant along the flow given by (1.4). With the added $d_A{*}\mathfrak{a} = 0$ condition, the equations in (1.4) can be written as the corresponding $\tau \in [0,1]$ version of (1.2) in the case when X is the interior of $\mathbb{I} \times M$.

The interest in what follows are the maps that solve (1.6) on the interior of $\mathbb{I}$ and are such that the function $|\frac{\partial}{\partial s} A|^2 + |\frac{\partial}{\partial s} \mathfrak{a}|^2 + |d_A\mathfrak{a}|^2 + |F_A - \mathfrak{a} \wedge \mathfrak{a}|^2$ is integrable on $\mathbb{I} \times M$. A map with these properties is said to be an *instanton*. As explained in Section 2d, a solution to (1.6) is an instanton solution if and only if the real part of the function $\frac{(\tau + i(1-\tau))^2}{\tau^2 + (1-\tau)^2}\,\mathfrak{cs}$ has finite variation across $\mathbb{I}$.

**b) Uhlenbeck's theorem**

Let X denote a given compact, oriented Riemannian 4-manifold and let $P \to X$ denote a given principal SO(3) bundle. Uhlenbeck's compactness theorem in [U] can be invoked to describe the limiting behavior of sequences of connections on P with anti-self dual curvature. These are solutions to (1.1) with $\mathfrak{a} = 0$. Uhlenbeck's theorem asserts the following: There exists a data set $(P_\Delta, \Theta, A_\Delta)$ with $P_\Delta \to X$ being a principal SO(3) bundle, $\Theta \to X$ being a finite set of points and $A_\Delta$ being connection on $P_\Delta$ with anti-self dual curvature. There is also a sequence of isomorphisms $\{g_i\}_{i \in \{1,2,\ldots\}}$ from $P_\Delta|_{X-\Theta} \to P|_{X-\Theta}$ and a subsequence $\Xi \subset \{1, 2, \ldots\}$ such that $\{g_i{*}A_i\}_{i \in \Xi}$ converges to $A_\Delta$ in the $C^\infty$ topology on compact subsets to $X-\Theta$. The number of elements in the set $\Theta$ is a priori bounded by a universal multiple of the first Pontrjagin class of the bundle $P \times_{SO(3)} \mathfrak{S}$.

Let M denote a given compact, oriented Riemannian 3-manifold and let P now denote a principal SO(3) bundle over M. Uhlenbeck's compactness theorem can also be used to describe the limiting behavior of sequences of parametrized paths in the space of connections on P that obey $\mathfrak{a} = 0$ version of (1.4). To say more, let $\{A_i\}_{i\in\{1,2,...\}}$ denote a sequence of such paths with each path having the same parametrizing interval. Denote this interval by $\mathbb{I}$. Uhlenbeck's theorem can be invoked if certain conditions are satisfied by the elements of the sequence $\{A_i\}_{i\in\{1,2,...\}}$ at any boundary points of $\mathbb{I}$ and if certain conditions are satisfied by the limits of the elements in $\{A_i\}_{i\in\{1,2,...\}}$ along any non-compact rays in $\mathbb{I}$. The appeal to Uhlenbeck's theorem yields a finite set $\Theta \subset \mathbb{I}\times M$, a principal SO(3) bundle $P_\Delta$ over $\mathbb{I} \times M$, and a connection $A_\Delta$ over $\mathbb{I} \times M$ with anti-self dual curvature with the following significance: View the sequence of paths $\{A_i\}_{i\in\{1,2,...\}}$ as a sequence of connections with anti-self dual curvature on the pull-back of the bundle P to $\mathbb{I} \times M$ via the projection to M. There exists a subsequence $\Xi \subset \{1, 2, \ldots\}$ and a sequence of isomorphisms $\{h_n\}_{n\in\Xi}$ from $P_\Delta$ to P over the complement of $\Theta$ such that $\{h_n{}^*A_i\}_{i\in\Xi}$ converges to $A_\Delta$ in the $C^\infty$ topology on compact subsets in $(\mathbb{I} \times M)-\Theta$.

The number of elements in the $\mathbb{I} \times M$ incarnation of $\Theta$ can be bounded a priori given certain information about the behaviour of the elements of the sequence $\{A_i\}_{i\in\Xi}$ on the boundary points of $\mathbb{I}$ and the limits of the elements along non-compact rays in $\mathbb{I}$.

**c) Notions and terminology**

This subsection introduces some notions and terminology that are used in the upcoming Theorems 1.1 and 1.2. Assume in what follows that X is a given 4-dimensional Riemannian manifold.

*Part 1*: Let $P \to X$ denote a principal SO(3) bundle and let A denote a connection on P. The connection A with the Levi-Civita connection from the metric defines a covariant derivative on $(P \times_{SO(3)} \mathfrak{G})$-valued tensor fields on X. This covariant derivative is denoted by $\nabla_A$. Meanwhile, the symbol $\nabla$ is used for the Levi-Civita connection's covariant derivative on tensor fields with values in $\mathbb{R}$. The covariant derivative $\nabla_A$ is metric compatible if the Riemannian metric is used in conjunction with an ad-invariant metric on $\mathfrak{G}$ to define a fiber metric for tensor valued sections of $P \times_{SO(3)} \mathfrak{G}$. The ad-invariant metric that is used here is obtained by viewing $\mathfrak{G}$ as the vector space of skew-Hermitian, $2 \times 2$ complex matrices with zero trace. Granted this identification, take the inner product that assigns the number $-\frac{1}{2}\,\mathrm{trace}(\sigma\sigma')$ to matrices $\sigma$ and $\sigma'$. This inner product is denoted by $\langle\sigma\sigma'\rangle$.

*Part 2*: Let $\hat{U} \subset X$ denote an open set. A *real line* bundle over $\hat{U}$ is the associated $\mathbb{R}$ bundle to a principal $\{1, -1\}$ bundle over $\hat{U}$. Let $\mathcal{I} \to \hat{U}$ denote such a bundle. The notion of the exterior derivative of an $\mathcal{I}$ valued form on $\hat{U}$ is defined a priori with no need to introduce an auxilliary connection. An $\mathcal{I}$ valued form is said to be *harmonic* when the exterior derivative of the form and the exterior derivative of its Hodge dual are both equal to zero. Keep in mind that harmonic $\mathcal{I}$ valued 1-forms are a priori smooth.

Let V a direct sum of tensor products of TX and T*X. The Riemannian metric on TX defines a metric on $V \otimes \mathcal{I}$ and the Levi-Civita connection defines a corresponding, metric compatible covariant derivative for sections of $V \otimes \mathcal{I}$. The latter is denoted in what follows by $\nabla$.

Let $P_\Delta \to \hat{U}$ denote a given principal SO(3) bundle and let $A_\Delta$ denote a given connection on $P_\Delta$. A homomorphism $\sigma_\Delta$: $\mathcal{I} \to P_\Delta \times_{SO(3)} \mathfrak{S}$ is isometric if it is norm preserving, and it is $A_\Delta$-covariantly constant if it intertwines the Levi-Civita covariant derivative $\nabla$ with the covariant derivative $\nabla_{A_\Delta}$.

*Part 3*: Suppose again that $\hat{U} \subset X$ is an open set and that $P_\Delta \to \hat{U}$ is a given principal SO(3) bundle and that $A_\Delta$ is a given connection on $P_\Delta$. The curvature 2-form of the given connection $A_\Delta$ on $P_\Delta$ is said to be *harmonic* if the Hodge star of its curvature is annihilated by $A_\Delta$'s exterior covariant derivative.

Let $\{A_n\}_{n\in\{1,2,\ldots\}}$ denote a sequence of connections on $P_\Delta$. This sequence is said to converge to $A_\Delta$ in the $C^\infty$ topology on compact subsets of $\hat{U}$ when the following occurs: Write each $n \in \{1, 2, \ldots\}$ version of $A_n$ as on $\hat{U}$ as $A_\Delta + a_{A_n}$ with $a_{A_n}$ being a section over $\hat{U}$ of the bundle $(P\times_{SO(3)} \mathfrak{S}) \otimes T^*X$. Let $V \subset \hat{U}$ denote any given open set in $\hat{U}$ with compact closure. Then $\{ a_{A_n} \}_{n=\{1,2,\ldots\}}$ converges to zero on V in each $k \in \{1, 2, \ldots\}$ version of the $C^k$ topology on sections over V of the bundle $(P\times_{SO(3)} \mathfrak{S}) \otimes T^*X$. The sequence $\{A_n\}_{n\in\{1,2,\ldots\}}$ is said to converge to $A_\Delta$ in the $L^2_1$ topology on compact subsets of $\hat{U}$ to $A_\Delta$ if

$$\lim_{n\to\infty} \int_V (|\nabla_{A_\Delta} a_{A_n}|^2 + |a_{A_n}|^2) = 0 \tag{1.8}$$

whenever $V \subset \hat{U}$ is an open set with compact closure.

Let $a$ denote a given section of $(P_\Delta\times_{SO(3)} \mathfrak{S}) \otimes T^*X$ over $\hat{U}$ and let $\{a_n\}_{n\in\{1,2,\ldots\}}$ denote a sequence of sections of $(P_\Delta\times_{SO(3)} \mathfrak{S}) \otimes T^*X$. The sequence $\{a_n\}$ converges to $a$ in the $C^0$ topology on compact subsets in $\hat{U}$ when the following is true: Let $V \subset \hat{U}$ denote any given open set with compact closure. Then the sequence $\{a - a_n\}_{n\in\{1,2,\ldots\}}$ converges to zero in the $C^0$ topology on the space of sections over V of the bundle $(P_\Delta\times_{SO(3)} \mathfrak{S}) \otimes T^*X$. The sequence $\{a_n\}_{i\in\{1,2,\ldots\}}$ is said to converge to $a$ in the $C^\infty$ topology on compact subsets of $\hat{U}$ when the following occurs: If V is any given open set in $\hat{U}$ with compact closure, the sequence $\{a - a_n\}_{n\in\{1,2,\ldots\}}$ converges to zero in the $C^\infty$ topology on the space of smooth sections over V of $(P_\Delta\times_{SO(3)} \mathfrak{S}) \otimes T^*X$. The sequence $\{a_n\}_{n\in\{1,2,\ldots\}}$ is said to converge to $a$ in the $L^2_1$ topology on compact subsets of $\hat{U}$ if

$$\lim_{n\to\infty} \int_V (|\nabla_{A_\Delta}(a - a_n)|^2 + |a - a_n|^2) = 0$$

(1.9)

whenever $V \subset \hat{U}$ is an open set with compact closure.

*Part 4*: A function $f$: $X \to \mathbb{R}$ is said to be an $L^2$ function if $f^2$ is integrable on X. It is said to be an $L^2_1$ function if both $f^2$ and $|\nabla f|^2$ are integrable on X.

**d) The extensions of Uhlenbeck's theorem**

What follows directly is the promised analog of Uhlenbecks theorem for sequence of solutions to (1.2). Keep in mind Kapustin and Witten's theorem [KW] to the effect that $\tau \notin \{0,1\}$ solutions exists on X when X is compact only if the first Pontrjagin class of $P\times_{SO(3)}\mathfrak{S}$ is zero; and that if $(A, \mathfrak{a})$ is such a solution, then $A + i\,\mathfrak{a}$ defines a flat, $PSL(2;\mathbb{C})$ connection. It is because of this theorem that Theorem 1.1 talks only about sequences of solutions to (1.2) in the cases when $\tau = 0$ or 1, or when the first Pontrjagin class of $P\times_{SO(3)}\mathfrak{S}$ is zero and the sequence is a sequence of flat $PSL(2;\mathbb{C})$ connections.

**Theorem 1.1**: *Let* X *denote a compact and oriented Riemannian 4-manifold; and let* P *denote a principal* SO(3) *bundle over* X. *If* $P\times_{SO(3)}\mathfrak{G}$ *has positive first Pontrjagin class, set* $\tau = 0$ *and if this Pontrjagin class is negative, set* $\tau = 1$. *Let* $\{(A_n, \mathfrak{a}_n)\}_{n=\{1,2,\ldots\}}$ *in these cases denote a sequence with each member being a pair of connection on the bundle* P *and section of* $T^*X\otimes(P\times_{SO(3)}\mathfrak{S}))$ *that obeys the relevant* $\tau \in \{0, 1\}$ *version of (1.2). In the case when the first Pontrjagin class of* $P\times_{SO(3)}\mathfrak{G}$ *is zero, assume that* $\{A_n + i\,\mathfrak{a}_n)\}_{n\in\{1,2,\ldots\}}$ *defines a sequence of flat* $PSL(2;\mathbb{C})$ *connections. Given* $n \in \{1, 2, \ldots\}$, *set* $r_n$ *to be the maximum of 1 and the* $L^2$ *norm of* $\mathfrak{a}_n$.

- *Assume that the sequence* $\{r_n\}_{n=1,2,\ldots}$ *has a bounded subsequence. There exists a principal* SO(3) *bundle* $P_\Delta \to X$ *and a pair* $(A_\Delta, \mathfrak{a}_\Delta)$ *with* $A_\Delta$ *being a connection on* $P_\Delta$ *and* $\mathfrak{a}_\Delta$ *being a section* $T^*X \otimes (P_\Delta \times_{SO(3)}\mathfrak{S})$ *that obeys the relevant* $\tau \in \{0,1\}$ *version of (1.2) when the first Pontrjagin class of* $P\times_{SO(3)}\mathfrak{G}$ *is non-zero, and is such that* $A_\Delta + i\,\mathfrak{a}_\Delta$ *is a flat* $PSL(2;\mathbb{C})$ *connection when this Pontrjagin class vanishes. Moreover, there exists a finite set* $\Theta \subset X$, *a subsequence* $\Xi \subset \{1,2, \ldots\}$, *and a sequence,* $\{g_n\}_{n\in\Xi}$, *of automorphisms of* $P_\Delta|_{X-\Theta}$ *such that* $\{(g_n{}^*A_n, g_n{}^*\mathfrak{a}_n)\}_{n\in\Xi}$ *converges to* $(A_\Delta, \mathfrak{a}_\Delta)$ *in the* $C^\infty$ *topology on compact subsets in* $X-\Theta$.
- *Assume that the sequence* $\{r_n\}_{n=1,2,\ldots}$ *has no bounded subsequence. There exists in this case the following data*:
  1) *A finite set* $\Theta \subset X$ *and a closed, nowhere dense set* $Z \subset X-\Theta$.
  2) *A real line bundle* $\mathcal{I} \to X-(Z\cup\Theta)$.

3) *A harmonic,* $\mathcal{I}$ *valued 1-form* $\nu$ *on* $X-(Z\cup\Theta)$. *The norm of* $\nu$ *extends over the whole of* $X$ *as a bounded,* $L^2_1$ *function. In addition,*

   a) *The extension of* $|\nu|$ *is a continuous on* $X-\Theta$ *and its zero locus is the set* $Z$.

   b) *Let* $U$ *denote an open set in* $X-\Theta$ *with compact closure. The function* $|\nu|$ *is Hölder continuous on* $U$ *with Hölder exponent that is independent of* $U$ *and of the original sequence* $\{(A_n, \mathfrak{a}_n)\}_{n=\{1,2,\ldots\}}$.

   c) *If* $p$ *is any given point in* $X$, *then the function* $\mathrm{dist}(\cdot, p)^{-1}\,|\nabla\nu|$ *extends to the whole of* $X$ *as an* $L^2$ *function.*

4) *A principal* $SO(3)$ *bundle* $P_\Delta \to X-(Z\cup\Theta)$ *and a connection* $A_\Delta$ *on* $P_\Delta$ *with harmonic curvature.*

5) *An isometric* $A_\Delta$ *covariantly constant homorphism* $\sigma_\Delta\colon \mathcal{I} \to P_\Delta\times_{SO(3)}\mathfrak{S}$.

*In addition, there exists a subsequence* $\Xi \subset \{1, 2, \ldots\}$ *and a sequence of isomorphisms* $\{g_n\colon P_\Delta \to P|_{X-(Z\cup\Theta)}\}_{n\in\Xi}$ *such that*

i) *The sequence* $\{g_n{}^*A_n\}_{n\in\Xi}$ *converges in the* $L^2_1$ *topology on compact subsets of* $X-(Z\cup\Theta)$ *to* $A_\Delta$.

ii) *The sequence* $\{r_n^{-1}\, g_n{}^*\mathfrak{a}_n\}_{n\in\Xi}$ *converges in the* $L^2_1$ *topology on compact subsets of* $X-(Z\cup\Theta)$ *to* $\nu\otimes\sigma_\Delta$. *This convergence is also in the* $C^0$ *topology on compact subsets of* $X-\Theta$. *Moreover, the sequence* $\{r_n^{-1}|\mathfrak{a}_n|\}_{n\in\Xi}$ *is pointwise bounded and it converges in the* $L^2_1$ *topology on* $X$ *to* $|\nu|$.

As was the case for Uhlenbeck's theorem, the number of elements in $\Theta$ is a priori bounded by a universal multiple of the absolute value of the first Pontrjagin class of $P\times_{SO(3)}\mathfrak{S}$. In particular, $\Theta$ is empty if the first Pontrjagin class of $P\times_{SO(3)}\mathfrak{S}$ is zero. More is said about the set $Z$ in the upcoming Theorems 1.3 and 1.4. Meanwhile, Theorem 1.5 in [T1] says more about $|\nu|$ near $Z$.

The upcoming Theorem 1.2 states an analog of Uhlenbeck's theorem for sequences of instanton solutions to the equations in (1.6). As noted previously, the equations in (1.6) on the interior of the interval $\mathbb{I}$ with the constraint that $d_A{*}\mathfrak{a} = 0$ are the same as those in (1.2) when $(A, \mathfrak{a})$ is viewed as a pair of connection $P \to \mathbb{I}\times M$ and section over $\mathbb{I}\times M$ of $(P\times_{SO(3)}\mathfrak{S})\otimes T^*(\mathbb{I}\times M)$. Theorem 1.2 views such a pair in this way as it considers the interior of $\mathbb{I}\times M$ as a 4-manifold and refers only to solutions to (1.2) on this 4-manifold with no a priori assumption that the solution is the incarnation of a solution of (1.4). The statement of the theorem takes $\mathbb{I}$ to have finite length. This is for convenience only and of no real concern because any given interval is a union of the latter sort.

**Theorem 1.2**: *Let* $\mathbb{I} \subset \mathbb{R}$ *denote a closed, bounded interval of length at least 4 and let* M *denote a compact, oriented Riemannian 3-manifold. Let* P *denote a principal* SO(3) *bundle over* $\mathbb{I} \times M$. *Fix* $\mathrm{E} \geq 0$ *and an open interval* $I \subset \mathbb{I}$ *with compact closure in the interior of* $\mathbb{I}$. *Suppose that* $\{(A_n, \mathfrak{a}_n)\}_{n\in\{1,2,\ldots\}}$ *is sequence of pairs consisting of a connection on* P *and a section of* $(P\times_{SO(3)} \mathfrak{G}) \otimes T^*(\mathbb{I} \times M)$. *Let* $\{\tau_n\}_{n=1,2,\ldots}$ *denote a convergent sequence in* $[0,1]$ *and assume that each* $n \in \{1, 2, \ldots\}$ *version of* $(A_n, \mathfrak{a}_n)$ *obeys the* $\tau = \tau_n$ *version of* (1.2) *on the interior of* $\mathbb{I} \times M$. *Assume in addition that*

$$\int_{\mathbb{I}\times M} (|d_{A_n}\mathfrak{a}_n|^2 + |F_{A_n} - \mathfrak{a}_n \wedge \mathfrak{a}_n|^2) \leq \mathrm{E}\,.$$

*Given* $n \in \{1, 2, \ldots\}$, *set* $r_n$ *to be the maximum of* 1 *and the* $L^2$ *norm of* $\mathfrak{a}_n$ *on* $I\times M$.

- *Assume that the sequence* $\{r_n\}_{n=1,2,\ldots}$ *has a bounded subsequence. There exists a principal* SO(3) *bundle* $P_\Delta \to I\times M$ *and a pair* $(A_\Delta, \mathfrak{a}_\Delta)$ *of connection on* $P_\Delta$ *and section* $T^*(I\times M) \otimes (P_\Delta \times_{SO(3)} \mathfrak{G})$ *that obeys some* $\tau \in [0,1]$ *version of* (1.2). *Moreover, there exists a finite set* $\Theta \subset I\times M$, *a subsequence* $\Xi \subset \{1,2, \ldots\}$, *and a sequence,* $\{g_n\}_{n\in\Xi}$, *of automorphisms of* $P_\Delta|_{(I\times M)-\Theta}$ *such that* $\{(g_n{}^*A_n, g_n{}^*\mathfrak{a}_n)\}_{n\in\Xi}$ *converges to* $(A_\Delta, \mathfrak{a}_\Delta)$ *in the* $C^\infty$ *topology on compact subsets in* $(I\times M)-\Theta$.
- *Assume that the sequence* $\{r_n\}_{n=1,2,\ldots}$ *lacks bounded subsequence. If there exists* $m > 1$ *such that*

$$\int_{\mathbb{I}\times M} |\mathfrak{a}_n|^2 < m\, r_n{}^2$$

*for all* $n \in \{1, 2, \ldots\}$, *then the conclusions of the second bullet of Theorem 1.1 hold with* $X = I\times M$.

By way of a parenthetical remark, the maximum number of elements in Theorem 1.2's set $\Theta$ can be bounded a priori given $\mathrm{E}$, $\mathbb{I}$, I and $m$.

What follows are remarks about Theorem 1.2's assumptions when the sequence $\{(A_n, \mathfrak{a}_n)\}_{n\in\{1,2,\ldots\}}$ in Theorem 1.2 is a sequence of solutions to (1.6). Theorem 1.2 first asks that the integral over $\mathbb{I}\times M$ of the function $|\frac{\partial}{\partial s} A_n|^2 + |\frac{\partial}{\partial s}\mathfrak{a}_n|^2 + |d_{A_n}\mathfrak{a}_n|^2 + |F_{A_n} - \mathfrak{a}_n \wedge \mathfrak{a}_n|^2$ be finite for each n. Thus the theorem asks that each pair in $\{(A_n, \mathfrak{a}_n)\}_{n\in\{1,2,\ldots\}}$ be an instanton solution. The upcoming Lemma 2.3 asserts that a solution to (1.6) is an instanton if and only if the map from $\mathbb{I}$ to $\mathbb{R}$ given by the real part of $\frac{(\tau + i(1-\tau))^2}{\tau^2 + (1-\tau)^2}\, \mathfrak{cs}(\mathbb{A} + i\,\mathfrak{a})$ has bounded variation on $\mathbb{I}$. Keep in mind that this function is in any event decreasing. Theorem 1.2 also asks for an a priori upper bound for the corresponding sequence of integrals.

Note that the upcoming Section 2e proves that there exists in all cases a number *m* for Theorem 1.2's second bullet when each member of the theorem's sequence $\{(A_n, \mathfrak{a}_n)\}_{n\in\{1,2,\ldots\}}$ is an instanton solution to (1.6).

**e) The structure of the set Z**

Theorems 1.1 and 1.2 say only that their versions of the set Z are closed and nowhere dense. The upcoming Theorems 1.3 and 1.4 gives more information about Z. Section 1h explains why these two theorems can be viewed as special cases of the theorems in [T1].

Theorem 1.3 introduces the notion of a *point of discontinuity for* $\mathcal{I}$. A given point $p \in Z$ is a point of discontinuity for $\mathcal{I}$ if the bundle $\mathcal{I}$ is not isomorphic to the product bundle on the complement of Z in any neighborhood of p. Supposing that p is not a point of discontinuity for $\mathcal{I}$, then there is a ball $B \subset U-\Theta$ centered at p such that $\mathcal{I}$ extends from $B-(B\cap Z)$ to B as a real line bundle. This extension is necessarily isomorphic to the product $\mathbb{R}$ bundle over B. Such an isomorphism identifies $\nu$ on B as an $\mathbb{R}$-valued harmonic 1-form. This being the case, there is nothing really novel about $\nu$ near p. The mysterious points in Z are those that are points of discontinuity for $\mathcal{I}$.

The upcoming theorem also talks about 2-dimensional Lipshitz submanifolds in X and 2-dimensional differentiable submanifolds in X. To define the first of these terms, write $\mathbb{R}^4$ as $\mathbb{R}^2\times\mathbb{R}^2$. Let $\Delta$ denote a neighborhood of the origin in $\mathbb{R}^d$. A 2-dimensional *Lipshitz graph* through the origin in $\mathbb{R}^4$ is the image of a map from $\Delta$ to $\mathbb{R}^4$ that has the form $z \to (\varphi(z), z)$ with $\varphi$ being a Lipshitz map from $\Delta$ to $\mathbb{R}^2$ sending the origin to the origin. By way of a reminder, a map $\varphi$ from a domain in $\mathbb{R}^2$ to $\mathbb{R}^2$ is said to be Lipshitz when there is a positive number $\kappa$ such that $\mathrm{dist}(\varphi(z),\varphi(z')) \le \kappa|z - z'|$ for all pairs of points z and z´ in the domain. For the purposes of this paper, a subset in X is said to be a 2-dimensional Lipshitz submanifold if it intersects a ball about each of its points as the image via a coordinate chart map of a Lipshitz graph through the origin in $\mathbb{R}^4$.

A 2-dimensional Lipshitz submanifold is said to be a differentiable submanifold if it has a tangent plane at each point and if these tangent planes vary continuously along the submanifold.

**Theorem 1.3**: *Let* $\Theta$ *and* Z *denote the sets that are described in Item 1) of the second bullet of Theorem 1.1 either in the case when* X *is compact or in the context of Theorem 1.2 when* X *refers to* $I\times M$.

- *The points of discontinuity for* $\mathcal{I}$ *are the points in the closure of an open subset of* Z *that has the structure of a 2-dimensional, differentiable submanifold in* $X-\Theta$.
- *The set* Z *has an open dense set that is contained in a countable union of 2-dimensional Lipshitz submanifolds in* $X-\Theta$.

The preceding theorem talks about open dense sets in Z. The next theorem gives something of a picture of the whole of Z. It says that Z looks to a first approximation like a stratified set with 2-dimensional top strata, then strata of dimension 1 and dimension 0. This upcoming theorem introduces the notions of *geodesic arcs* and *geodesic disks* through a given point. Supposing that p is the point in question, a geodesic arc through p is the image in X via a Gaussian coordinate chart map centered at p of the small |x| part of a line through the origin in $\mathbb{R}^4$; and a geodesic disk through p is the corresponding image of the small |x| part of plane through the origin in $\mathbb{R}^4$. The theorem also introduces the notion of a *union of geodesic half-disks with common boundary arc* through p. The definition of this notion requires first the definition of a half-plane through the origin in $\mathbb{R}^4$: Supposing that $\Pi$ is a plane through the origin and $f$ is a non-trivial linear function on $\Pi$, then the $f \geq 0$ part of $\Pi$ is said to be a half-plane through the origin. The line where $f = 0$ is said to be the boundary line of the half plane. A union of geodesic half-disks with common boundary arc is the image via a Gaussian coordinate chart map centered at p of the small |x| part of the union of two or more half-planes all with common boundary line. A geodesic arc or disk or union of half-disks with common boundary arc in a given open ball is said to be *proper* if the set in question is a relatively closed, path connected subset of the ball.

**Theorem 1.4**: *Let* $\Theta$ *and* Z *denote the sets that are described in Item 1) of the second bullet of Theorem 1.1 either in the case when* X *is compact or in the context of Theorem 1.2 when* X *refers to* $I \times M$. *The set* Z *has Hausdorff dimension at most 2. Moreover,* Z *looks to a first approximation like a stratified 2 dimensional set in the following sense: Given* $\theta \in (0,1)$ *and* $z > 1$, *there are finite collections* $\mathfrak{U}$ *and* $\mathfrak{V}$ *of the following sort:*

- *Both* $\mathfrak{U}$ *and* $\mathfrak{V}$ *are composed of balls of radius at most* $z^{-1}$.
- *The balls from* $\mathfrak{U}$ *are pairwise disjoint and* $\sum_{B\in\mathfrak{U}}(\text{radius}(B))^{\theta} < z^{-1}$.
- *The balls in the set* $\mathfrak{V}$ *cover the part of* Z *in* $X - (\cup_{B\in\mathfrak{U}} B)$.
- *Supposing that* B *is a ball from* $\mathfrak{V}$, *let* p *denote its center point and* r *its radius. The part of* Z *in* B *is contained in the radius* $z^{-1}r$ *tubular neighborhood of either a proper geodesic arc through* p *or a proper geodesic disk through* p *or a proper union of half-disks through* p *with common boundary arc.*

By way of a parenthetical remark, it is likely that the whole of Z is a contained in a finite union of finite 2-dimensional, rectifiable subsets. This was a claimed in an earlier version of this article, but the proof of the latter assertion had an error. This error was found by Thomas Walpuski.

**f) The example of flat PSL(2; $\mathbb{C}$) connections**

If X is compact, then the only solutions to the $\tau \neq 0$ or 1 versions of (1.2) are pairs $(A, \mathfrak{a})$ with $F_A - \mathfrak{a} \wedge \mathfrak{a} = 0$ and $d_A\mathfrak{a} = 0$. That this is so follows from the upcoming (2.17). The 2-forms $F_A - \mathfrak{a} \wedge \mathfrak{a}$ and $i\, d_A\mathfrak{a}$ are the respective anti-Hermitian and Hermitian parts of the PSL(2; $\mathbb{C}$) connection $\mathbb{A} = A + i\mathfrak{a}$ and so their vanishing says that $\mathbb{A}$ is a flat PSL(2; $\mathbb{C}$) connection on X. The additional $d_A{*}\mathfrak{a} = 0$ constraint can be viewed as a *stability* condition on the space of flat, PSL(2; $\mathbb{C}$) connections. Let $\mathcal{R}$ denote the quotient of the space of solutions to the $\tau \in (0, 1)$ version of (1.2) by the action of the group of automorphisms of the bundle P. It follows from what was just said about flat connections that $\mathcal{R}$ can be viewed as a subspace of the quotient of the space of representations of $\pi_1(X)$ into PSL(2; $\mathbb{C}$) by the conjugation action of SO(3). Theorem 1.1 defines a sort of compactification of this space $\mathcal{R}$ as it associates a data set $(Z, \mathcal{I}, \nu, P_\Delta, A_\Delta, \sigma_\Delta)$ to sequences in $\mathcal{R}$ with no convergent subsequences.

The data supplied by Theorem 1.1 for non-convergent sequences in $\mathcal{R}$ is related to Morgan and Shalen's compactification [MS1]-[MS4] of the quotient of the space irreducible representations of $\pi_1(X)$ into PSL(2; $\mathbb{C}$) by the conjugation action of the non-compact group PSL(2; $\mathbb{C}$). The next paragraphs gives a brief account of how this relationship comes about. This account borrows from Part 1 in Section 1c of [T2].

The Morgan-Shalen compactification of the space PSL(2; $\mathbb{C}$) conjugacy classes of representations of $\pi_1(X)$ into PSL(2; $\mathbb{C}$) involves $\pi_1(X)$-equivariant maps from the universal cover of X to $\mathbb{R}$-trees. An $\mathbb{R}$-tree is a contractible metric space with the property that there is a unique path between any two points. Daskoloupoulus, Dostoglu and Wentworth [DDW] use work of Korevaar and Schoen [KS] to prove that the Morgan-Shalen maps can be taken to be harmonic. Theorem 1.1's limit 1-form $\nu$ and the set Z have the following interpretation in the context of [DDW1]-[DDW3]: Use $\mathbb{T}$ to denote the [DDW1] limit $\mathbb{R}$ tree and $u$ their limit harmonic map from X's universal cover to $\mathbb{T}$. Gromov and Schoen [GS] and [S]) proved that $u$ can be viewed as an honest harmonic function on small balls in the complement of a set with Hausdorff dimension at most 2. The latter set appears in the context of Theorem 1.1 as the inverse image in X's universal cover of the set Z. Meanwhile, the differential of $u$ where it is an honest harmonic function is the pull-back of $\nu$ to X's universal cover.

A tetrad of closely related notions in 3-manifold topology are singular measured foliations, measured laminations, weighted branched surfaces and maps to $\mathbb{R}$-trees. See for example [B], [O], [Hat], [MS1]-[MS4], [GO] and [HO]. An elegant account of the 2-dimensional story can found in Calegari's book [C]. Whether talking about flat PSL(2; $\mathbb{C}$) connections or not, singular measured foliations appear directly in the context of Theorem 1.2 because the kernel of the $\mathfrak{I}$ valued 1-form $\nu$ defines an integrable distribution on $X-(Z\cup\Theta)$. The set Z plays the role of the branching locus of the corresponding foliation and $|\nu|$ where non-zero gives the transverse measure.

**g) Outline of the proofs of Theorems 1.1 and 1.2**

The remaining sections of this paper prove Theorems 1.1 and 1.2. As already noted, the proofs of Theorems 1.3 and 1.4 are in the sequel to this paper [T1]. An outline of the proof of Theorems 1.1 and 1.2 is given momentarily. What follows directly is a list of the headings for the subsequent sections.

SECTION 2: *$L^2_1$ and $L^\infty$ convergence*

SECTION 3: *The functions $r_{c_\wedge}$ and $r_{cF}$*

SECTION 4: *Unexpectedly small curvature*

SECTION 5: *Monotonicity for the $(A, a)$ version of* N

SECTION 6: *Convergence in the $C^0$ topology*

SECTION 7: *The data Z, $\mathfrak{I}$ and $\nu$, the connection $A_\Delta$ and the homomorphism $\sigma_\Delta$.*

SECTION 8: *Holder continuity of $|\nu|$ along* Z.

The first bullet of Theorem 1.1 when X is compact is restated as the first bullet of Proposition 2.1; and the statement of the first bullet of Theorem 1.1 in the context of Theorem 1.2 when $X = I\times M$ is restated as the first bullet of Proposition 2.2. These propositions are proved in Section 2. The second bullets of Propositions 2.1 and 2.2 concern the case in Theorems 1.1 and 1.2 when the sequence $\{r_n\}_{n\in\{1,2,\ldots\}}$ has no convergent subsequences. These bullets assert in part that there is a subsequence of positive integers, this denoted by $\Lambda$, such that $\{r_n^{-1}|\mathfrak{a}_n|\}_{n\in\Lambda}$ converges in the $L^2_1$ topology and pointwise to a bounded, $L^2_1$ function on X. This function is denoted by $|\hat{a}_\diamond|$ and it is defined at each point of X by the rule $|\hat{a}_\diamond|(p) = \limsup_{n\in\Lambda} r_n^{-1}|\mathfrak{a}_n|(p)$. The function $|\hat{a}_\diamond|$ will turn out to be the norm of the 1-form $|\nu|$.

Proposition 6.1 makes a formal statement to the effect that $|\hat{a}_\diamond|$ is continuous on the complement of a finite set in X. This finite set is $\Theta$. Proposition 6.1 also states that the sequence $|r_n^{-1}\mathfrak{a}_n|$ converges to $|\hat{a}_\diamond|$ in the $C^0$ topology on compact subsets in $X-\Theta$. The set $\Theta$ is a version of one of the sets from Lemma 6.2. The criteria for membership is such that a point p is *not* in $\Theta$ when

$$\lim_{r\to 0}\Big(\limsup_{n\in\Lambda}\int_{\mathrm{dist}(\cdot,p)<r}|F_{A_n} - \mathfrak{a}_n\wedge\mathfrak{a}_n|^2\Big) \tag{1.11}$$

is small in a suitable sense. The fact that (1.11) is small for all but finitely many points is seen to be a consequence of the integral identities in the upcoming equations (2.17) and (2.23). The identity in (2.17) holds when X is compact and (2.23) holds in the context of Theorem 1.2. These integral identities play a crucial role in the proof of Theorems 1.1 and 1.2 because they imply that $\Theta$ as defined in Section 6 is finite.

The set Z is defined in Section 7 as the locus in $X-\Theta$ where $|\hat{a}_{\diamond}|$ is zero. This is a closed set in $X-\Theta$ because $|\hat{a}_{\diamond}|$ is continuous. Lemma 7.1 proves that Z is nowhere dense. This is to say that Z has empty interior. Except for the assertion of Item 3b) of Theorem 1.1's second bullet, all remaining assertions in the second bullet of Theorem 1.1 follow directly from what is said in Proposition 7.2 given the preceding definition of Z. The assertion in Item 3b) of the second bullet of Theorem 1.1 follows from Proposition 8.1.

The appearance in Proposition 7.2 and Theorem 1.1 of the 1-form $\nu$ and the homomorphism $\sigma_{\Delta}$ can be traced to the existence of another identity, this being the Bochner-Weitzenbock formula in the upcoming equation (2.5). In particular, the integral forms of the ($A = A_n$, $a = r_n^{-1}\mathfrak{a}_n$) and $r = r_n$ versions of (2.5) for $n \in \{1, 2, \ldots\}$ suggest that a 'reasonable' limit of the sequence $\{r_n^{-1}\mathfrak{a}_n\}_{n\in\Lambda}$, call it $q$, must be such that $q\wedge q = 0$. This can happen only if $q$ is decomposable in the sense that it can be written near any point where it is non-zero as $\nu \otimes \sigma_{\Delta}$ with $\nu$ being real valued 1-form and $\sigma_{\Delta}$ being a unit length element in $\mathfrak{S}$. If the corresponding subsequence of $\{A_n\}_{n\in\{1,2,\ldots\}}$ converges and the limit, call it $A_{\Delta}$, is such that the $A_{\Delta}$-exterior derivatives of both $q$ and $*q$ are zero, then $\nu$ must be harmonic and $\sigma_{\Delta}$ must be $A_{\Delta}$-covariantly constant. The identification of $\nu$ as an $\mathbb{R}$-valued 1-form can be done locally but not globally because $\sigma_{\Delta}$ can be defined at best up to a sign. The line bundle $\mathcal{I}$ accounts for this sign ambiguity.

Sections 3, 4 and 5 supply the tools that are used in Section 6 to prove that $|\hat{a}_{\diamond}|$ is continuous on the complement of a finite set in X. Section 3 associates two length scales to a solution to (1.2) and any given point in X, these denoted by $r_{c\wedge}$ and $r_{cF}$. Letting p denote the given point, the definition of $r_{c\wedge}$ is such that the integral of $|\mathfrak{a}\wedge\mathfrak{a}|^2$ over a ball centered p of radius less than $r_{c\wedge}$ is bounded by a certain small number that is ultimately determined by the geometry of X. By the same token, the definition of $r_{cF}$ is such that the integral of $|F_A|^2$ over a ball of radius less than $r_{cF}$ centered at p is smaller than a number that is determined by the geometry of X. Section 3 states and then proves three propositions that give a priori estimates for norms of $F_A$ and the covariant derivative $\nabla_A\mathfrak{a}$ on a ball centered p of radius less than the minimum of $r_{c\wedge}$ and $r_{cF}$. These a priori bounds depend on a function of the radius of the ball that is defined in (3.5). This function is constructed using the pair (A, $\mathfrak{a}$); it is an analog of Almgren's frequency function [Al], [HHL], [Han], [DF].

Section 4 states and then proves a proposition asserting that the size of $|\mathfrak{a} \wedge \mathfrak{a}|$ on a ball centered at p of radius less than the minimum of $r_{c_\wedge}$ and $r_{c_F}$ is much smaller than one might guess from the definition of $r_{c_\wedge}$. This surprisingly small upper bound holds on the ball if the frequency function in (3.5) is small. The proof of the proposition in Section 4 exploits properties of the linearized versions of the equations in (1.2) on $\mathbb{R}^4$.

Section 5 states and then proves a monotonicity formula for the frequency function in (3.5). This monotonic behavior of the frequency function is used in Section 6 to relate the behavior of $(A, \mathfrak{a})$ on balls centered at p of radius less than the minimum of $r_{c_\wedge}$ and $r_{c_F}$ to the behavior of $(A, \mathfrak{a})$ on balls of radius $\mathcal{O}(1)$ centered at p. The use in Section 6 of the monotonicity formula from Section 5 requires that p come from the compliment of the set $\Theta$. The use in Section 6 of what is said in Section 4 also requires that p not be in $\Theta$.

**h) The proofs of Theorems 1.3 and 1.4**

The first bullet of Theorem 1.3 is a special case of Theorem 1.2 in [T1] and the second bullet in Theorem 1.3 is a special case of Theorem 1.4 in [T1]. Meanwhile, Theorem 1.4 here is a special case of Theorem 1.2 in [T1]. The input from this paper that is needed to invoke Theorems 1.2, 1.3 and 1.4 in [T1] consists of following data:

- *Let* $U \subset X - \Theta$ *denote a given open set with compact closure. There is a continuous non-negative function on* U *to be denoted by* $f$ *obeying*
  - i) $f > 0$ *somewhere.*
  - ii) *There exists* $\varepsilon > 0$ *such that if* $p \in U$ *and* $f(p) = 0$, *and if* r *is small and positive, then* $\int_{\mathrm{dist}(p,\cdot)<r} f^2 \le r^{4+\varepsilon}$.
- *Let* Z *denote* $f^{-1}(0)$. *A real line bundle* $\mathcal{I} \to U - Z$.
- *A section* $\nu$ *of* $T^*X \otimes \mathcal{I}$ *over* $U - Z$ *that obeys*
  - i) $d\nu = 0$ *and* $d{*}\nu = 0$.
  - ii) $|\nu| = f$.
  - iii) *The function* $|\nabla\nu|^2$ *is integrable on* $U - Z$.

(1.12)

The data set $(Z, \mathcal{I}, \nu)$ from the second bullet of either Theorem 1.1 or Theorem 1.2 with the function $f = |\nu|$ obeys the conditions in (1.12). The condition in Item ii) of the first bullet is a consequence of the assertion in Item 3b) of the second bullet of Theorem 1.1 or its analog for Theorem 1.2 to the effect that $|\nu|$ is Hölder continuous on $Z \cap U$ with an exponent that is bounded away from zero. The condition in Item ii) of the first bullet also follows from the first and fourth bullets of Lemma 8.2.

### i) Conventions

This article uses $c_0$ to denote a number that is greater than 1. Unless said otherwise, the value of $c_0$ depends only on the manifold X, its metric, the bundle P and on occasion a specified open set X with compact closure to be denoted by U. A second convention has $\chi$ denoting non-increasing function on $\mathbb{R}$ chosen once and for ever that is equal to 1 on $(-\infty, \frac{1}{4}]$ and equal to 0 on $[\frac{3}{4}, \infty)$.

## 2. $L^2_1$ and $L^\infty$ convergence

The two propositions that follow momentarily constitute the starting point for the subsequent analysis of sequences of solutions to (1.2) or (1.4). The first of these propositions talks about solutions to (1.2) when X is compact. The manifold X and the bundle P constitute what is said by the proposition to be the *geometric data*. This first proposition introduces one additional bit of notation, this being a quadratic, fiber preserving map from $T^*X \otimes (P \times_G \mathfrak{S})$ to $T^*X \otimes T^*X$ that is defined with the help of the G-invariant inner product on $(P \times_G \mathfrak{S})$. To define this map, at a given point $p \in X$, suppose that $q$ is a 1-form at p with values in $(P \times_G \mathfrak{S})|_p$. Write $q$ as $q_\alpha \hat{e}^\alpha$ with $\{\hat{e}^\alpha\}_{\alpha=1,2,3,4}$ being an orthonormal basis for $T^*X$ at p and with $\{q_\alpha\}_{\alpha\in 1,2,3,4}$ being in $(P \times_G \mathfrak{S})|_p$. Here and in what follows, summation over repeated indices is implicit. The quadratic map sends $q$ to what is denoted by $\langle q \otimes q \rangle$, this being $\langle q_\alpha q_\beta \rangle\, \hat{e}^\alpha \otimes \hat{e}^\beta$. What is denoted by Ric is the Ricci tensor of the Riemannian metric, this being a section of $T^*X \otimes T^*X$ also. The metric pairing between Ric and $\langle q \otimes q \rangle$ is denoted by $\mathrm{Ric}(\langle q \otimes q \rangle)$.

**Proposition 2.1**: *Let* X *denote a compact and oriented Riemannian 4-manifold; and let* P *denote a principal* SO(3) *bundle over* X. *If* $P \times_{SO(3)} \mathfrak{G}$ *has positive first Pontrjagin class, set* $\tau = 0$ *and if this Pontrjagin class is negative, set* $\tau = 1$. *Let* $\{(A_n, \mathfrak{a}_n)\}_{n=\{1,2,\ldots\}}$ *in these cases denote a sequence with each member being a pair of connection on the bundle* P *and section of* $T^*X \otimes (P \times_{SO(3)} \mathfrak{S}))$ *that obeys the relevant* $\tau \in \{0, 1\}$ *version of (1.2). In the case when the first Pontrjagin class of* $P \times_{SO(3)} \mathfrak{G}$ *is zero, assume that* $\{A_n + i\,\mathfrak{a}_n)\}_{n\in\{1,2,\ldots\}}$ *defines a sequence of flat* $PSL(2;\mathbb{C})$ *connections. Given* $n \in \{1, 2, \ldots\}$, *set* $r_n$ *to be the maximum of 1 and the* $L^2$ *norm of* $\mathfrak{a}_n$.

- *Assume that the sequence* $\{r_n\}_{n=\{1,2,\ldots\}}$ *has a bounded subsequence. There exists a principal* SO(3) *bundle* $P_\Delta \to X$ *and a pair* $(A_\Delta, \mathfrak{a}_\Delta)$ *with* $A_\Delta$ *being a connection on* $P_\Delta$ *and* $\mathfrak{a}_\Delta$ *being a section* $T^*X \otimes (P_\Delta \times_{SO(3)} \mathfrak{S})$ *that obeys the relevant* $\tau \in \{0,1\}$ *version of (1.2) when the first Pontrjagin class of* $P \times_{SO(3)} \mathfrak{G}$ *is non-zero, and is such that* $A_\Delta + i\,\mathfrak{a}_\Delta$ *is a flat* $PSL(2;\mathbb{C})$ *connection when this Pontrjagin class vanishes. There is, in addition, a finite set* $\Theta \subset X$, *a subsequence* $\Xi \subset \{1,2, \ldots\}$, *and a sequence,* $\{g_n\}_{n\in\Xi}$, *of automorphisms of* $P_\Delta|_{X-\Theta}$ *such that* $\{(g_n{}^*A_n, g_n{}^*\mathfrak{a}_n)\}_{n\in\Xi}$ *converges to* $(A_\Delta, \mathfrak{a}_\Delta)$ *in the*

*$C^\infty$ topology on compact subsets in* $X-\Theta$. *The set* $\Theta$ *is empty and* $P_\Delta = P$ *if the first Pontrjagin class of* $P\times_{SO(3)}\mathfrak{G}$ *is zero*.

- *Assume that the sequence* $\{r_n\}_{\{1,2,\ldots\}}$ *has no bounded subsequence. There exists a subsequence* $\Lambda \subset \{1, 2, \ldots\}$ *such that Items a)-e) listed below hold. This list uses* $a_n$ *when* $n \in \Lambda$ *to denote* $r_n^{-1}\mathfrak{a}_n$.

  a) *The sequences* $\{ \int_X (|d|a_n||^2 + |a_n|^2) \}_{n\in\Lambda}$ *and* $\{\sup_X |a_n|\}_{n\in\Lambda}$ *are bounded by a number that depends only on the geometric data*.

  b) *The sequence* $\{|a_n|\}_{n\in\Lambda}$ *converges weakly in the* $L^2{}_1$ *topology and strongly in all* $p < \infty$ *versions of the* $L^p$ *topology. The limit function is denoted by* $|\hat{a}_\Diamond|$, *this being an* $L^\infty$ *function with* $L^2$ *norm equal to 1 The function* $|a_\Diamond|$ *is defined at each point in* $X$ *by the rule whereby* $|\hat{a}_\Diamond|(p) = \lim\sup_{n\in\Lambda} |a_n|(p)$ *for each* $p \in X$.

  c) *The sequence* $\{\langle a_n \otimes a_n\rangle\}_{n\in\Lambda}$ *converges strongly in any* $q < \infty$ *version of the* $L^q$ *topology on the space of sections of* $T^*X \otimes T^*X$. *The limit section is denoted by* $\langle \hat{a}_\Diamond \otimes \hat{a}_\Diamond\rangle$ *and its trace is the function* $|\hat{a}_\Diamond|^2$.

  d) *Use* $f$ *to denote a given* $C^0$ *function*.

     i) *The sequences* $\{ \int_X f\, |F_{A_n} - r_n^2 a_n \wedge a_n|^2 \}_{n\in\Lambda}$ *and* $\{ \int_X f\, |d_{A_n} a_n|^2 \}_{n\in\Lambda}$ *converge.*

     ii) *The sequences* $\{ \int_X f\, |\nabla_{A_n} a_n|^2 \}_{n\in\Lambda}$ *and* $\{2\, r_n^2 \int_X f\, |a_n \wedge a_n|^2 \}_{n\in\Lambda}$ *converge. The limit of the first sequence is denoted by* $Q_{\nabla,f}$ *and that of the second by* $Q_{\wedge,f}$. *These are such that if* $f$ *is a* $C^2$ *function, then*

     $$\tfrac{1}{2} \int_X d^* df\, |\hat{a}_\Diamond|^2 + Q_{\nabla,f} + Q_{\wedge,f} + \int_X f\, \mathrm{Ric}(\langle \hat{a}_\Diamond \otimes \hat{a}_\Diamond\rangle) = 0\ ,$$

  e) *Fix* $p \in X$ *and let* $G_p$ *denote the Green's function with pole at* $p$ *for the operator* $d^\dagger d + 1$. *The sequence that is indexed by* $\Lambda$ *with n'th term being the integral of* $G_p(|\nabla_{A_n} a_n|^2 + 2r_n^2 |a_n \wedge a_n|^2)$ *is bounded by a number that depends only on the geometric data. In either case, let* $Q_{\Diamond,p}$ *denote the lim-inf of this sequence. The function* $|\hat{a}_\Diamond|^2$ *obeys the equation*

     $$\tfrac{1}{2} |\hat{a}_\Diamond|^2(p) + Q_{\Diamond,p} = - \int_X G_p\, (\tfrac{1}{2} |\hat{a}_\Diamond|^2 - \mathrm{Ric}(\langle \hat{a}_\Diamond \otimes \hat{a}_\Diamond\rangle))\ .$$

To set the stage for Proposition 2.2, suppose that $\mathbb{I} \subset \mathbb{R}$ is a closed interval. The projection map from $\mathbb{I} \times M$ to $M$ pulls back $P$, $P\times_{SO(3)}\mathfrak{G}$ and $T^*M$ as bundles over $\mathbb{I} \times M$. Suppose that $A$ is a path of connections on $P$ parametrized by $\mathbb{I}$. Use what was just said about pull-backs to view $A$ as a connection on the pull-back of $P$ over $\mathbb{I} \times M$. Let $\mathfrak{a}$ denote a path in $C^\infty(M; (P\times_{SO(3)}\mathfrak{G}) \otimes T^*M)$ parametrized by $\mathbb{I}$. What was just said about pull-backs is likewise used to view $\mathfrak{a}$ as a section over $\mathbb{I} \times M$ of $(P\times_{SO(3)} \mathfrak{S}) \otimes T^*(\mathbb{I} \times M)$. The incarnations A and $\mathfrak{a}$ as a connection over $\mathbb{I} \times M$ and section of a vector bundle over $\mathbb{I} \times M$ are still denoted by A and $\mathfrak{a}$.

As noted previously, the equations in (1.6) on the interior of $\mathbb{I}$ with the constraint that $d_A{*}\mathfrak{a} = 0$ are the same as those in (1.2) when $(A, \mathfrak{a})$ are viewed in their incarnation as a pair of connection $P \to \mathbb{I} \times M$ and section over $\mathbb{I} \times M$ of $T^*(P\times_{SO(3)} \mathfrak{S}) \otimes (\mathbb{I} \times M)$. Proposition 2.2 views such a pair in this way as it considers the interior of $\mathbb{I} \times M$ as a 4-manifold and refers only to solutions to (1.2) on this 4-manifold with no assumption that the solution is the incarnation of a solution of (1.6).

**Proposition 2.2**: *Let* $\mathbb{I} \subset \mathbb{R}$ *denote a closed, bounded interval of length at least 4 and let* M *denote a compact, oriented Riemannian 3-manifold. Let* P *denote a principal* SO(3) *bundle over* $\mathbb{I} \times M$. *Fix* $E \geq 0$ *and an open interval* $I \subset \mathbb{I}$ *with compact closure in the interior of* $\mathbb{I}$. *Suppose that* $\{(A_n, \mathfrak{a}_n)\}_{n=1,2,\ldots}$ *is a sequence of pairs of connection on* P *and section of* $(P\times_{SO(3)} \mathfrak{S}) \otimes T^*(\mathbb{I} \times M)$. *Let* $\{\tau_n\}_{n=1,2,\ldots}$ *denote a convergent sequence in* $[0,1]$ *and assume that each* $n \in \{1, 2, \ldots\}$ *version of* $(A_n, \mathfrak{a}_n)$ *obeys the* $\tau = \tau_n$ *version of* (1.2) *on the interior of* $\mathbb{I} \times M$. *Assume in addition that*

$$\int_{\mathbb{I}\times M} (|d_{A_n} \mathfrak{a}_n|^2 + |F_{A_n} - \mathfrak{a}_n \wedge \mathfrak{a}_n|^2) \leq E .$$

*Given* $n \in \{1, 2, \ldots\}$, *set* $r_n$ *to be the maximum of* 1 *and the* $L^2$ *norm of* $\mathfrak{a}_n$ *on* $I \times M$.

- *Assume that the sequence* $\{r_n\}_{n=1,2,\ldots}$ *has a bounded subsequence. There is a principal* G *bundle* $P_\Delta \to I \times M$ *and a pair* $(A_\Delta, \mathfrak{a}_\Delta)$ *of connection on* $P_\Delta$ *and section of the bundle* $(P_\Delta\times_{SO(3)} \mathfrak{S}) \otimes T^*(\mathbb{I} \times M))$ *that obeys some* $\tau \in [0,1]$ *version of* (1.2). *Moreover, there is a finite set* $\Theta \subset I \times M$; *a subsequence* $\Xi \subset \{1, 2, \ldots\}$ *and a sequence* $\{g_n\}_{n\in\Xi}$ *of isomorphisms from* $P_\Delta|_{(I\times M)-\Theta}$ *to* $P|_{(I\times M)-\Theta}$ *with the following property: The sequence* $\{(g_n{}^*A_n, g_n{}^*\mathfrak{a}_n)\}_{n\in\Xi}$ *converges to* $(A_\Delta, \mathfrak{a}_\Delta)$ *in the* $C^\infty$ *topology on open subsets in* $(I \times M)-\Theta$ *with compact closure. The set* $\Theta$ *is empty and* $P_\Delta = P$ *if* $\tau \in (0,1)$.
- *Assume that the sequence* $\{r_n\}_{n=1,2,\ldots}$ *lacks bounded subsequence. If there exists* $m > 1$ *such that*

  $$\int_{\mathbb{I}\times M} |\mathfrak{a}_n|^2 < m r_n{}^2$$

  *for all* $n \in \{1, 2, \ldots\}$, *then there exists a subsequence* $\Lambda \subset \{1, 2, \ldots\}$ *such that Items a)-e) listed below hold. This list uses* $a_n$ *when* $n \in \Lambda$ *to denote* $r_n{}^{-1}\mathfrak{a}_n$.

  a) *The sequences* $\{ \int_{I\times M} (|d|a_n||^2 + |a_n|^2) \}_{n\in\Lambda}$ *and* $\{\sup_{I\times M} |a_n|\}_{n\in\Lambda}$ *are bounded by a*

*number that depends only on the geometric data and the boundary limits of the sequence* $\{(A_n, a_n)\}_{n\in\Lambda}$.

b) *The sequence* $\{|a_n|\}_{n\in\Lambda}$ *converges weakly in the* $L^2_1$ *topology on* $I \times M$ *and strongly in all* $p < \infty$ *versions of the* $L^p$ *topology on* $I \times M$. *The limit function is denoted by* $|\hat{a}_\Diamond|$, *this being an* $L^\infty$ *function with* $L^2$ *norm equal to 1. The function* $|a_\Diamond|$ *is defined at each point of* $I \times M$ *by the rule* $p \to |\hat{a}_\Diamond|(p)$ *with* $|\hat{a}_\Diamond|(p) = \lim\sup_{n\in\Lambda} |a_n|(p)$.

c) *The sequence* $\{\langle a_n \otimes a_n\rangle\}_{n\in\Lambda}$ *converges strongly in any* $q < \infty$ *version of the* $L^q$ *topology on the space of sections of* $T^*(I\times M) \otimes T^*(I\times M)$. *The limit section is denoted by* $\langle \hat{a}_\Diamond \otimes \hat{a}_\Diamond\rangle$ *and its trace is the function* $|\hat{a}_\Diamond|^2$.

d) *Use* $f$ *to denote a given* $C^0$ *function.*

i) *The sequences* $\{ \int_{I\times M} f\,|F_{A_n} - r_n^2\, a_n \wedge a_n|^2 \}_{n\in\Lambda}$ *and* $\{ \int_{I\times M} f\,|d_{A_n} a_n|^2 \}_{n\in\Lambda}$ *converge.*

ii) *The sequences* $\{ \int_{I\times M} f\,|\nabla_{A_n} a_n|^2 \}_{n\in\Lambda}$ *and* $\{2\, r_n^2 \int_{I\times M} f\,|a_n \wedge a_n|^2 \}_{n\in\Lambda}$ *converge. The limit of the first sequence is denoted by* $Q_{\nabla,f}$ *and that of the second by* $Q_{\wedge,f}$. *If* $f$ *is a* $C^2$ *function with compact support in* $I\times M$, *then*

$$\frac{1}{2}\int_{I\times M} d^*df\,|\hat{a}_\Diamond|^2 + Q_{\nabla,f} + Q_{\wedge,f} + \int_{I\times M} f\,\mathrm{Ric}(\langle \hat{a}_\Diamond \otimes \hat{a}_\Diamond\rangle) = 0 ,$$

e) *Fix a point* $p \in I \times M$ *and let* $G_p$ *denote the Green's function on* $\mathbb{I} \times M$ *for the operator* $-\frac{\partial^2}{\partial s^2} + d^\dagger d + 1$ *that vanishes on any boundary component of* $\mathbb{I} \times M$ *and has its pole at* $p$. *The sequence that is indexed by* $\Lambda$ *with n'th term being the integral of* $G_p(|\nabla_{A_n} a_n|^2 + 2r_n^2 |a_n \wedge a_n|^2)$ *over* $I \times M$ *is bounded by a number that depends solely on the geometric data and the boundary limits of the sequence* $\{(A_n, a_n)\}_{n\in\Lambda}$. *Let* $Q_{\Diamond,p}$ *denote the lim-inf of the aforementioned sequence of integrals. The function* $|\hat{a}_\Diamond|^2$ *obeys*

$$\frac{1}{2}|\hat{a}_\Diamond|^2(p) + Q_{\Diamond,p} = -\int_{I\times M} G_p\left(\frac{1}{2}|\hat{a}_\Diamond|^2 - \mathrm{Ric}(\langle \hat{a}_\Diamond \otimes \hat{a}_\Diamond\rangle)\right) + Q_p ,$$

*where* $Q_{(\cdot)}$ *is a smooth function on* $I \times M$.

Lemma 2.4 implies in part that the requirement in the second bullet on the two integrals of each $\mathfrak{a} \in \{\mathfrak{a}_n\}_{n\in\{1,2,\ldots\}}$ is always satisfied when $\{(A_n, \mathfrak{a}_n)\}_{n\in\{1,2,\ldots\}}$ is an instanton solution to (1.6) on $\mathbb{I} \times M$. More is said about this in Section 2e.

By way of a parenthetical remark, the conclusions of Proposition 2.2 hold if the $\{(A_n, \mathfrak{a}_n)\}_{n\in\{1,2,\ldots\}}$ is an instanton solution to (1.6) without the assumption about the vanishing of each $n \in \{1, 2, \ldots\}$ version of $d_{A_n} * \mathfrak{a}_n$. It is sufficient to assume only a number $m$ such that $\sup_{\mathbb{I}\times M} |d_{A_n} * \mathfrak{a}_n| \le m$ for all n.

The proofs of Propositions 2.1 and 2.2 are in Section 2f. Sections 2a-2e supply auxilliary results that are used in the proof and in subsequent sections of this paper. Section 2g has some parenthetical remarks about a gradient flow interpretation of (1.6).

**a) Sobolev and Hardy inequalities**

The subsequent arguments in this sections and in later sections invoke a dimension 4 Sobolev inequality and Hardy's inequality. To say what these are, let X denote for the moment a given Riemannian 4-manifold and let $f$ denote an $L^2_1$ function on a ball in X. Let p denote the center point of this ball. Use r to denote the radius of this ball and $B_r$ to denote the ball. If r is such that this ball is well inside a Gaussian coordinate chart, then both $f^2$ and $\mathrm{dist}(p,\cdot)^{-1} f$ are square integrable on the ball and their integrals are such that the following is true:

- $\left(\int_{B_r} f^4\right)^{1/2} < c_0 \int_{B_r} (|df|^2 + r^{-2}\,|f|^2)$ .
- $\int_{B_r} \frac{1}{\mathrm{dist}(p,\cdot)^2} f^2 < c_0 \int_{B_r} (|df|^2 + r^{-2}\,|f|^2)$ .
- *If f has compact support in* $B_r$ *or if* $\int_{B_r} f = 0$*, then both of the preceding inequalities hold with the right hand side being* $c_0 \int_{B_r} |df|^2$ .

(2.1)

It follows from the first bullet in (2.1) that $L^2_1$ functions are in $L^k$ for $k \le 4$ and the resulting map from $L^2_1$ to $L^k$ is bounded. This map is also compact when $k < 4$. The top bullet in (2.1) is a Sobolev inequality and the second bullet is Hardy's inequality.

Note that if $\mathfrak{a}$ is a section of a vector bundle with fiber metric and $\nabla$ is a metric compatible connection, then $|d|\mathfrak{a}|| \le |\nabla\mathfrak{a}|$. This being the case, if $|\nabla\mathfrak{a}|$ and $|\mathfrak{a}|$ are square integrable, then $|\mathfrak{a}|$ is an $L^2_1$ function. This fact is used without comment in what follows.

**b) Identities for solutions to (1.1)**

Suppose that $r \ge 1$ and that (A, $a$) is a pair consisting of a connection on P and a section of $T^*X \otimes (P \times_{SO(3)} \mathfrak{S})$ that obeys the equations

$$(F_A - r^2\, a \wedge a)^+ = 0 \quad \textit{and} \quad (d_A a)^- = 0 \;\; \textit{and} \;\; d_A{*}a = 0\ .$$

(2.2)

This subsection derives some basic identities that are obeyed by solutions to (2.2).

The first identity is an integral identity that holds in the case when X is compact. The statement of the identity uses $p_1(P \times_{SO(3)} \mathfrak{S})$ to denote the first Pontrjagin number of the bundle $P \times_{SO(3)} \mathfrak{S}$. As explained in the next paragraph, if X is compact, if $r \geq 1$ and if $(A, a)$ obeys (2.2), then

$$\int_X |(F_A - r^2 a \wedge a)^-|^2 + r^2 \int_X |(d_A a)^+|^2 = -32\pi^2 p_1(P \times_{SO(3)} \mathfrak{S}) \ . \tag{2.3}$$

This identity implies, among other things, that there are no solutions to (1.1) in the case when X is compact if the first Pontrjagin number of $P \times_{SO(3)} \mathfrak{S}$ is positive. It also implies that $(A, a)$ solves (2.2)when X is compact and $p_1(P \times_{SO(3)} \mathfrak{S}) = 0$ if and only if $A + i\, r a$ defines a flat, $PSl(2; \mathbb{C})$ connection on X.

To derive (2.3), introduce the principal $PSl(2; \mathbb{C})$ bundle $\mathbb{P} = P \times_{SO(3)} PSl(2; \mathbb{C})$. Let $\mathbb{A}$ for the moment denote any given connection on $\mathbb{P}$. The group $PSl(2; \mathbb{C})$ has its adjoint action on the Lie algebra of $Sl(2; \mathbb{C})$ this being $\mathfrak{S}_{\mathbb{C}} = \mathfrak{S} \otimes_{\mathbb{R}} \mathbb{C}$. The second Chern number of the associated, rank 3 complex vector bundle is the first Pontrjagin number of the bundle $P \times_{SO(3)} \mathfrak{S}$. With this in mind, let $\mathbb{A}$ denote any given connection on $\mathbb{P}$. The Chern-Weil theory identifies this first Pontrjagin number with

$$\frac{1}{32\pi^2} \int_X \langle F_{\mathbb{A}} \wedge F_{\mathbb{A}} \rangle \ . \tag{2.4}$$

The first Pontrjagin number is real so only the real part of $\langle F_{\mathbb{A}} \wedge F_{\mathbb{A}} \rangle$ contributes to this integral. This understood, then (2.3) follows from the $\mathbb{A} = A + i\, r a$ version of (2.4) because the real part of $\langle F_{\mathbb{A}} \wedge F_{\mathbb{A}} \rangle$ is $-|(F_A - r^2 a \wedge a)^-|^2 - |(d_A a)^+|^2$ if $(A, a)$ obeys (2.2).

The next identity asserts that $a$ obeys the second order differential equation in the upcoming (2.5). This identity holds whether or not X is compact. To set the notation, suppose that $p \in X$. If an orthonormal frame for TX is chosen at p, then the pairing between $\mathfrak{a}$ and the frame vectors defines a set $\{a_\alpha\}_{\alpha=1,2,3,4}$ of elements in $(P \times_{SO(3)} \mathfrak{S})|_x$. These define an endomorphism of the fiber at x of any given $(P \times_{SO(3)} \mathfrak{S})|_x$ valued bundle of tensors using the rule that sends a tensor $q$ to $[a_\alpha, [q, a_\alpha]]$ with it understood that there is a summation over repeated indices. This endomorphism does not depend on the chosen orthonormal basis. The upcoming equation (2.5) also denotes the endomorphism of $T^*X \otimes (P \times_{SO(3)} \mathfrak{S})$ that is defined by the Ricci tensor by $\mathrm{Ric}(\cdot)$.

The promised second order equation asserts that

$$\nabla_A^\dagger \nabla_A a + r^2 [a_\alpha, [a, a_\alpha]] + \mathrm{Ric}(a) = 0 \tag{2.5}$$

when (A, $a$) obeys (2.2). This equation follows from the first order equations in (2.2) using the Bochner-Weitzenboch identity for the operator on sections of $T^*X \otimes (P \times_{SO(3)} \mathfrak{S})$ that acts by the rule $q \to *2d_A(d_A q)^- - d_A*(d_A*q)$. The equation in (2.5) leads in turn to a differential equation for $|a|^2$ that reads

$$\tfrac{1}{2}\, d^\dagger d\,|a|^2 + |\nabla_A a|^2 + 2\,r^2|a\wedge a|^2 + \mathrm{Ric}(\langle a\otimes a\rangle) = 0\ , \tag{2.6}$$

this derived by taking the inner product of both sides of (2.5) with $a$.

The identity in (2.6) is used to derive various useful integral identities in the case when X is compact. The first of these is obtained by integrating both sides of (2.5) over X; it asserts that

$$\int_X (|\nabla_A a|^2 + 2r^2\ |a\wedge a|^2 + \mathrm{Ric}(\langle a\otimes a\rangle))\ = 0\ . \tag{2.7}$$

This identity leads directly to the a priori bound

$$\int_X (|\nabla_A a|^2 + 2r^2|a\wedge a|^2)\ \le c_0 \int_X |a|^2\ . \tag{2.8}$$

The second integral identity gives an a priori bound for the pointwise norm of $|a|$ when X is compact. To state this second identity, fix for the moment $p \in X$ and let $G_p$ denote the Green's function for the operator $d^\dagger d + 1$ on X with pole at p. The latter obeys the bounds

$$|G_p(\cdot)| < c_0 \frac{1}{\mathrm{dist}(p,\cdot)^2} \quad \textit{and} \quad |dG_p(\cdot)| \le c_0 \frac{1}{\mathrm{dist}(p,\cdot)^3}\ . \tag{2.9}$$

Multiply both sides of (2.6) by $G_p$ and integrate over X. Integrate by parts to see that

$$\tfrac{1}{2}\,|a|^2(p) + \int_X G_p(|\nabla_A a|^2 + r^2\ |a\wedge a|^2) = \int_X G_p(\tfrac{1}{2}\,|a|^2 - \mathrm{Ric}(\langle a\otimes a\rangle))\ . \tag{2.10}$$

The middle bullet in (2.1) with (2.8) can be used to bound the absolute value of the right hand side of (2.10) by $c_0 \int_X |a|^2$ ; and doing so gives a $c_0 \left(\int_X |a|^2\right)^{1/2}$ bound for $\sup_X |a|$.

**c) Identities for solutions to (1.2)**

Suppose that $(A, \mathfrak{a})$ is a solution to some $\tau \in [0, 1]$ version of (1.2) on a Riemannian 4-manifold X. Supposing that $r \geq 1$ has been specified, set $a = r^{-1}\mathfrak{a}$. The equations in (1.2) in terms of A and $a$ are

- $\tau(F_A - r^2 a \wedge a)^+ = (1-\tau)\, r\, (d_A a)^+$ ,
- $(1-\tau)(F_A - r^2 a \wedge a)^- = -\tau\, r\, (d_A a)^-$ ,
- $d_A{*}a = 0$ .

(2.11)

As explained momentarily, the pair (A, $a$) also obey (2.5). Granted that this is so, then this pair obeys (2.6), (2.7), (2.8) and (2.9) when X is compact; and it is also the case that $|a|$ is nowhere greater than $c_0$ times the $L^2$ norm of $a$ when X is compact.

To see about (2.5), note first that the $\tau = 1$ version of (2.11) is (2.2) and so what is said in the previous subsection applies when $\tau = 1$. What is said in this same subsection applies to the $\tau = 0$ version of (2.11) because the $\tau = 0$ version of (2.11) is a version of (2.2) on the orientation reversed manifold. Granted these observations, assume henceforth that $\tau$ is neither 0 nor 1. Write $(1-\tau)\tau^{-1}$ as $e^\theta$ with $\theta \in \mathbb{R}$. The equations in (2.11) are the respective self-dual and anti-self dual projections of the single equation

$$F_A - r^2 a \wedge a = \sinh\theta\, r\, d_A a + \cosh\theta\, r\, {*d_A a} \,. \tag{2.12}$$

Take the covariant exterior derivative of the latter and use the Bianchi identity to see that

$$-r^2 d_A a \wedge a + r^2 a \wedge d_A a = \sinh\theta\, r\, (F_A \wedge a - a \wedge F_A) + \cosh\theta\, r\, d_A{*}d_A a \,. \tag{2.13}$$

Using (2.12) to write $F_A$ leads from (2.13) to the identity

$$-\cosh\theta\, r^2 (d_A a \wedge a - a \wedge d_A a) = \sinh\theta\, r^2 ({*d_A a} \wedge a - a \wedge d_A a) + r\, d_A{*}d_A a \,. \tag{2.14}$$

To go further, act on both sides of (2.12) by the Hodge star and use the resulting identity to write (2.14) as

$$d_A{*}d_A a + r^2 (a \wedge {*(a \wedge a)} - {*(a \wedge a)} \wedge a) + {*F_A} \wedge a - a \wedge {*F_A} = 0 \,. \tag{2.15}$$

Take (2.15) and use the Weitzenboch formula for the operator $d_A{*}d_A$ with the fact that $d_A{*}a = 0$ to see that $\nabla_A^\dagger\nabla_A a + {*r^2}(a \wedge {*(a \wedge a)} - {*(a \wedge a)} \wedge a) + \mathrm{Ric}(a) = 0$. The latter equation is (2.5).

Suppose for the moment that X is compact. The identity in (2.3) holds only in the case that $\tau = 0$ but there is an analog for the $\tau = 1$ case. The latter is the $\tau = 1$ version of an identity that holds when $r$, $\tau$ and $(A, a)$) obey (2.11), this being

$$\int_X |F_A - r^2 a \wedge a|^2 + r^2 \int_X |d_A a|^2 = \frac{(1-2\tau)}{\tau^2 + (1-\tau)^2} 32\pi^2 p_1(P \times_{SO(3)} \mathfrak{S}) \ . \tag{2.16}$$

(This with the equation $d_A {*} a = 0$ gives (3.33) in [KW].) The derivation of (2.16) starts by using (2.11) to recast the real part of the product of $(\tau + i(1-\tau))^2$ times the $\mathbb{A} = A + i\, r a$ version of (2.4) as the identity

$$\int_X |\tau(F_A - r^2 a \wedge a)^- - (1-\tau) r (d_A a)^-|^2 + \int_X |(1-\tau)(F_A - r^2 a \wedge a)^+ + \tau r (d_A a)^+|^2 = (1-2\tau)\, 32\pi^2 p_1(P \times_{SO(3)} \mathfrak{S}) \ . \tag{2.17}$$

The equations in (2.11) and the identity in (2.16) lead directly to (2.17).

The assertion in [KW] to the effect that $r$ and $(A, a)$ solve any given $\tau \in (0,1)$ version of (2.11) only if $A + i\, r a$ is a flat $PSL(2;\mathbb{C})$ connection follows from (2.16) given the assertion that the $\tau \in (0, 1)$ versions of (2.11) have solutions only in the case when $p_1(P \times_{SO(3)} \mathfrak{S})$ is zero. As explained by Kapustin and Witten, this assertion about the vanishing of the Pontrjagin class follows from two facts, the first being that the 4-form $\langle (\tau(F_A - r^2 a \wedge a) - (1-\tau) r (d_A a)) \wedge ((1-\tau)(F_A - r^2 a \wedge a) + \tau r (d_A a)) \rangle$ is zero when (2.11) is obeyed; and the second being that this form in any event differs from $\tau(1-\tau) \langle F_A \wedge F_A \rangle$ by the exact form $d\langle r a \wedge ((\tau^2 - (1-\tau)^2)(F_A - \frac{1}{3} r^2 a \wedge a) - \tau(1-\tau) r d_A a) \rangle$.

**d) Identities for solutions to (1.6)**

Let $\mathbb{I} \subset \mathbb{R}$ denote a closed interval of length greater than 4. Fix $\tau \in [0,1]$ and suppose that $(A, \mathfrak{a})$ is a solution to the corresponding version of (1.10) on $\mathbb{I} \times M$. Assume in what follows that $(A, \mathfrak{a})$ is boundary convergent. Analogs of the bounds in Section 2c are derived momentarily for $(A, \mathfrak{a})$.

As in the preceding subsections, it proves convenient to fix $r \geq 1$ and set $a = r^{-1} \mathfrak{a}$ so as to write the equations in (1.6) as

- $\frac{\partial}{\partial s} A + \frac{\tau^2 - (1-\tau)^2}{\tau^2 + (1-\tau)^2} {*}(F_A - r^2 a \wedge a) - \frac{2\tau(1-\tau)}{\tau^2 + (1-\tau)^2} r {*} d_A a = 0 \ .$
- $\frac{\partial}{\partial s} a - \frac{\tau^2 - (1-\tau)^2}{\tau^2 + (1-\tau)^2} r {*} d_A a - \frac{2\tau(1-\tau)}{\tau^2 + (1-\tau)^2} {*}(F_A - r^2 a \wedge a) = 0 \ .$

(2.18)

Keep in mind that the Hodge star here is the 3-dimensional Hodge star on M and that the covariant exterior derivatives are derivatives on M. There are three parts to the subsequent derivation of the identites that are obeyed by the solutions to (2.18).

*Part 1*: The first Bochner-Weitzenboch formula for solutions to (2.18) follows from the fact that the equations are the formal gradient flow equations of the real part of $\frac{(\tau - i(1-\tau))^2}{\tau^2 + (1-\tau)^2}\,\mathfrak{cs}$, this being the identity

$$\int_M \left(|\tfrac{\partial}{\partial s}A|^2 + r^2|\tfrac{\partial}{\partial s}a|^2 + r^2|d_A a|^2 + |F_A - r^2 a \wedge a|^2\right) = -\frac{d}{ds}\,\mathrm{re}\left(\frac{(\tau + i(1-\tau))^2}{\tau^2 + (1-\tau)^2}\,\mathfrak{cs}\right) \tag{2.19}$$

that holds at any given $s \in \mathbb{I}$. The equation in (2.19) is derived by taking the $L^2$ norm of the two equations in (2.18) and then integrating by parts. The identity in (2.19) leads directly to the observation stated below for reference as Lemma 2.3.

**Lemma 2.3**: *Suppose that* $\mathbb{I} \subset \mathbb{R}$ *is an interval and suppose that* $(A, a)$ *is a solution to (2.18) on* $\mathbb{I} \times M$. *The real part of the function* $\frac{(\tau + i(1-\tau))^2}{\tau^2 + (1-\tau)^2}\,\mathfrak{cs}(A + i\,r a)$ *is non-increasing on* $\mathbb{I}$ *and* $\int_{\mathbb{I}\times M} (|\frac{\partial}{\partial s}A|^2 + r^2|\frac{\partial}{\partial s}a|^2 + r^2|d_A a|^2 + |F_A - r^2 a \wedge a|^2)$ *is equal to the the total decrease of its real part along* $\mathbb{I}$.

Lemma 2.3 implies that $(A, \mathfrak{a} = r a)$ is an instanton solution to (1.6) if and only if the total decrease along $\mathbb{I}$ of the real part of $\frac{(\tau + i(1-\tau))^2}{\tau^2 + (1-\tau)^2}\,\mathfrak{cs}(A + i\,\mathfrak{a})$ is finite.

Lemma 2.3 leads to a bound on the $L^2$ norm of $\frac{\partial}{\partial s} a$ that depends only on the total decrease of the real part of $\frac{(\tau + i(1-\tau))^2}{\tau^2 + (1-\tau)^2}\,\mathfrak{cs}(A + i r a)$. The next lemma exploits this fact.

**Lemma 2.4**: *Suppose that* $\mathbb{I} \subset \mathbb{R}$ *is an interval and suppose that* $(A, a)$ *is a solution to (2.18) on* $\mathbb{I}\times M$. *Assume that the total decrease along* $\mathbb{I}\times M$ *of the real part of the function* $\frac{(\tau + i(1-\tau))^2}{\tau^2 + (1-\tau)^2}\,\mathfrak{cs}(A + i r a)$ *is finite and denote this number by* $c_{\mathfrak{cs}}$. *Fix* $\delta \in (0, \infty)$. *Supposing that* $s, s' \in \mathbb{I}$*, then*

$$\int_{\{s\}\times M} |a|^2 \le (1+\delta) \int_{\{s'\}\times M} |a|^2 + (1+\delta^{-1})\, r^{-2} c_{\mathfrak{cs}} |s - s'| ,$$

***Proof of Lemma 2.4***: Let $\mathrm{I} \subset \mathbb{I}$ denote the interval between s and s´. Use the fundamental theorem of calculus to see that respective $L^2$ norms of $a$ along $\{s\}\times M$ and $\{s´\}\times M$ differ by at most $|s - s´|^{1/2}$ times the $L^2$ norm of $\frac{\partial}{\partial s} a$ on $\mathrm{I} \times M$. The latter bound with Lemma 2.3 implies the lemma's bound.

A second Bochner-Weitzenboch formula for solutions to (1.10) is obtained from (2.19) by adding and subtracting from the left hand side the square of the $L^2$ norm on M of $d_A{*}a$ at the relevant $s \in \mathbb{I}$. Having done this rewrite the exterior covariant derivatives using integration by parts to obtain the identity

$$\int_M \big(r^{-2}\, |\tfrac{\partial}{\partial s} A|^2 + |\tfrac{\partial}{\partial s} a|^2 + |\nabla_A a|^2 + r^{-2}\, |F_A|^2 + r^2\, |a \wedge a|^2 + \mathrm{Ric}(\langle a \otimes a\rangle)\big) = \int_M |d_A * a|^2 - r^{-2} \frac{d}{ds}\, \mathrm{re}\Big(\frac{(\tau + i(1-\tau))^2}{\tau^2 + (1-\tau)^2}\, \mathfrak{cs}\Big). \qquad (2.20)$$

The subsequent applications of (2.20) exploit the following observation of Witten about solutions to (2.18):

*Let $\mathbb{I} \subset \mathbb{R}$ denote a closed interval of length at least 4 and let $r$ and $(A, a)$ denote a solution to any given $\tau \in [0, 1]$ version of (2.11) on $\mathbb{I} \times M$. Then $\frac{\partial}{\partial s}(d_A{*}a) = 0$.*

(2.21)

This assertion is proved by first writing $\frac{\partial}{\partial s}(d_A{*}a)$ as $d_A{*}(\frac{\partial}{\partial s} a) + \frac{\partial}{\partial s} A \wedge a - a \wedge \frac{\partial}{\partial s} A$ and then using (2.18) to see that the contributions from $\frac{\partial}{\partial s} a$ and $\frac{\partial}{\partial s} A$ cancel.

Granted (2.21), the Bochner formula in (2.20) leads to an a priori bound for the $L^2$ norm of $\nabla_A a$ on domains of the form $\mathrm{I} \times M$ with $\mathrm{I} \subset \mathbb{I}$ being a bounded interval if the total drop along $\mathbb{I}$ of $\frac{(\tau + i(1-\tau))^2}{\tau^2 + (1-\tau)^2}\, \mathfrak{cs}(A + i r a)$ is finite. Assume henceforth that such is the case and denote this drop by $c_{\mathfrak{cs}}$. Given that the function $\frac{d}{ds}\, \mathrm{re}(\frac{(\tau + i(1-\tau))^2}{\tau^2 + (1-\tau)^2}\, \mathfrak{cs})$ on $\mathbb{I}$ is non-positive, integration of (2.20) over I finds

$$\int_{\mathrm{I}\times M} \big(r^{-2} |\tfrac{\partial}{\partial s} A|^2 + |\tfrac{\partial}{\partial s} a|^2 + |\nabla_A a|^2 + r^{-2} |F_A|^2 + r^2 |a \wedge a|^2\big) \le c_0 \int_{\mathrm{I}\times M} |a|^2 + r^{-2}\big(c_{\mathfrak{cs}} + c_{\mathrm{div}}\, \mathrm{length}(\mathrm{I})\big), \qquad (2.22)$$

with $c_{\mathrm{div}}$ equal to the square of the $L^2$ norm on M of $(r d_A{*}a)|_s$ with s being any given point in $\mathbb{I}$. Since $\frac{\partial}{\partial s}(\mathrm{div}_A a) = 0$, the precise choice for s is of no consequence.

The inequality in (2.22) is invoked in the proof of Proposition 2.2 in the following circumstance: Supposing that $(A, \mathfrak{a})$ obeys (1.4), take $r$ to be the maximum of 1 and the $L^2$ norm of $\mathfrak{a}$ on $I \times M$. Set $a = r^{-1}\mathfrak{a}$. Let $I' \subset \mathbb{I}$ denote a second interval of length 1. The inequality in Lemma 2.4 with the $I'$ version of (2.22) lead directly to the bound

$$\int_{I' \times M} \left(|\tfrac{\partial}{\partial s} a|^2 + |\nabla_A a|^2 + |a|^2\right) \le c_0\left(1 + r^{-2}(c_{\mathfrak{cs}} + c_{\mathrm{div}})\right)\left(1 + \mathrm{dist}(I', I)\right) \tag{2.23}$$

By way of a parenthetical remark, this last bound is the analog for solutions to (2.18) of the bound in (2.8).

*Part 2*: View A as a connection on P's pull-back over $\mathbb{I} \times M$ and view $\mathfrak{a}$ as a 1-form on $\mathbb{I} \times M$ with values in the pull-back of $P \times_{SO(3)} \mathfrak{S}$. In this guise, the equations in (2.18) for $(A, a)$ assert that $(A, a)$ obey the equations in the first two bullets of (2.11). Keep in mind in this regard that the notation in (2.11) differs from that in (2.18). In particular, the exterior derivative of $a$ and the curvature 2-form $F_A$ that appear in (2.11) have components that do not annihilate the tangent vector to the $\mathbb{I}$ factor of $\mathbb{I} \times M$.

Granted what was said in the preceding paragraph about (2.11), then the arguments from the preceding subsection that lead from (2.11) to (2.5) can be repeated in the context at hand to see that $a$ on $I \times M$ obeys the equation

$$-\tfrac{\partial^2}{\partial s^2} a + \nabla_A{}^{\dagger}\nabla_A a + r^2{*}(a \wedge {*}(a \wedge a) - {*}(a \wedge a) \wedge a) + d_A({*}d_A{*}a)) + \mathrm{Ric}\, a = 0\ , \tag{2.24}$$

where the notation is as follows: The symbol $*$ denotes the Hodge star on $\wedge^* T^*M$ and Ric denotes M's Ricci curvature tensor. The covariant derivative as defined by A along the constant s slices of $I \times M$ is denoted by $\nabla_A$; and $d_A$ denotes the exterior covariant derivative along these same slices. This last equation in turn leads to an analog of (2.6) that reads as follows:

$$\tfrac{1}{2}\left(-\tfrac{\partial^2}{\partial s^2}|a|^2 + d^{\dagger}d|a|^2\right) + |\tfrac{\partial}{\partial s} a|^2 + |\nabla_A a|^2 + 2r^2|a \wedge a|^2 + \langle a \wedge {*}\, d_A({*}d_A{*}a)\rangle + \mathrm{Ric}(\langle a \otimes a\rangle) = 0. \tag{2.25}$$

Fix a bounded, open interval $I \subset \mathbb{I}$ of length 1 with compact closure in the interior of $\mathbb{I}$. Let $r$ again denote the maximum of 1 and the $L^2$ norm of $\mathfrak{a}$ on $I \times M$. Set $a = r^{-1}\mathfrak{a}$ as before. As explained directly, the equation in (2.24) leads to an a priori bound for the pointwise norm of $a$ on $I \times M$ that depend only on the numbers $c_{\mathfrak{cs}}$ and $c_{\mathrm{div}}$. To derive this bound, let p denote a given point in I and let $G_p$ denote the Green's function on $\mathbb{I} \times M$ for the operator $-\frac{\partial^2}{\partial s^2} + d^{\dagger}d + 1$ with pole at p. This Green's function should also vanish on the boundary or boundaries of $\mathbb{I} \times M$ if $\mathbb{I}$ is not the whole of $\mathbb{R}$.

An analog of (2.10) in the present context reads

$$\tfrac{1}{2}|a|^2(\mathrm{p}) + \int_{\mathrm{I}\times\mathrm{M}} \mathrm{G_p}\big(|\tfrac{\partial}{\partial s}a|^2 + |\nabla_A a|^2 + 2r^2\,|a\wedge a|^2\big) \le \mathrm{c_I} + \mathrm{c_0}r^{-2}(c_{\mathfrak{cs}} + c_{\mathrm{div}*}))\,.$$

(2.26)

with $\mathrm{c_I} \le \mathrm{c_0}(1 + \mathrm{dist}(\mathrm{I},\partial\mathbb{I})^{-4})$ and with $c_{\mathrm{div}*}$ denoting the square of the sup norm of $|r\mathrm{d}_A{*}a|$.

The derivation of (2.26) exploits bounds for $\mathrm{G_p}$ and $\mathrm{dG_p}$ that are the $\mathbb{I}\times\mathrm{M}$ analogs of those in (2.9), these being

$$|\mathrm{G_p}| \le \mathrm{c_0}\big(1 + \frac{1}{\mathrm{dist(p,\cdot)}^2}\big)\,e^{-\mathrm{dist}(\cdot,\mathrm{p})/\mathrm{c_0}} \quad \textit{and} \quad |\mathrm{dG_p}| \le \mathrm{c_0}\big(1 + \frac{1}{\mathrm{dist(p,\cdot)}^3}\big)\,e^{-\mathrm{dist}(\cdot,\mathrm{p})/\mathrm{c_0}}\quad .$$

(2.27)

To derive (2.26), let $\chi$ denote now and henceforth a favorite, nonincreasing smooth function on $\mathbb{R}$ that equals 1 on $(-\infty, \frac{1}{4}]$ and equals zero on $[\frac{3}{4}, \infty)$. Use $\chi$ to construct a smooth, nonnegative function on $\mathbb{I}$ that equals 1 on I and where the distance to I is greater than the minimum of 2 and $\frac{1}{2}$ of the distance from I to $\partial\mathbb{I}$. This function can and should be constructed so that the norm of its derivative is bounded by $\mathrm{c_0}(1 + \mathrm{dist}(\mathrm{I},\partial\mathbb{I})^{-1})$ and so that the norm of its second derivative is bounded by $\mathrm{c_0}(1+\mathrm{dist}(\mathrm{I},\partial\mathbb{I})^{-2})$. Denote this function by $\chi_{\mathrm{I}}$.

Take the inner product of (2.25) with $\chi_{\mathrm{I}}\mathrm{G_p}a$ and integrate the result over $\mathbb{I}\times\mathrm{M}$. A subsequent integration by parts leads to a bound on the expression that appears on the left hand side of (2.26) by

$$\int_{\mathbb{I}\times\mathrm{M}} \chi_{\mathrm{I}}\mathrm{G_p}\big(\tfrac{1}{2}|a|^2 - \mathrm{Ric}(\langle a\otimes a\rangle)\big) + \tfrac{1}{2}\int_{\mathbb{I}\times\mathrm{M}} \big(2\tfrac{\partial}{\partial s}\mathrm{G_p}\tfrac{\partial}{\partial s}\chi_{\mathrm{I}} - \mathrm{G_p}\tfrac{\partial^2}{\partial s^2}\chi_{\mathrm{I}}\big)|a|^2$$

$$+ \int_{\mathbb{I}\times\mathrm{M}} \chi_{\mathrm{I}}\mathrm{G_p}|\mathrm{d}_A{*}a|^2 - \int_{\mathbb{I}\times\mathrm{M}} \chi_{\mathrm{I}}\,\mathrm{dG_p}\wedge{*}\langle a\,({*}\mathrm{d}_A{*}a)\rangle\ .$$

(2.28)

To bound the integrals in (2.28) by the expression that appears on the right hand side of (2.26), it proves useful to consider the four integrals in (2.28) in turn starting with the left most integral. To bound the left most integral, use (2.27) with the bound in (2.23) and the second bullet in (2.1) to bound the norm of the left most integral in (2.28) by $\mathrm{c_0}(1+r^{-2}(c_{\mathfrak{cs}} + c_{\mathrm{div}}))$. Consider next the integral that is second from the left in (2.28). To this end, note first that the support of the integrand of this integral has distance no less than the minimum of 1 and $\mathrm{c_0}^{-1}\mathrm{dist}(\mathrm{I},\partial\mathbb{I})$ from p. This understood, it follows from (2.27) with (2.23) that the norm of this integral is at most than $\mathrm{c_I}(1+r^{-2}(c_{\mathfrak{cs}} + c_{\mathrm{div}}))$ with $\mathrm{c_I}$ being less than the maximum of $\mathrm{c_0}$ and $\mathrm{c_0}\,\mathrm{dist}(\mathrm{I},\partial\mathbb{I})^{-4}$. Meanwhile, it follows from (2.27) and that the integral in (2.28) that is third from the left is no greater than $\mathrm{c_0}$ times the square of the supremum norm of $d_A{*}a$, this being $c_0\, r^{-2}\, c_{div*}$. The four steps that follow deal with the right most integral in (2.28).

<u>Step</u> <u>1</u>: The right most integral in (2.28) is no greater than $r^{-1}(c_{div*})^{1/2}$ times the integral of function $\chi_I |dG_p||a|$. This function in turn is no greater than $c_0\,\mathrm{dist}(p,\cdot)^{-3}\chi_I|a|$. The integral of the latter function is no greater than $c_0$ times the $L^6$ norm of $\chi_I|a|$. The subsequent steps bound this $L^6$ norm by $c_{I*}(1+r^{-2}(c_{\mathfrak{cs}}+c_{div*}))^{1/2}$ with $c_{I*}$ being a number that is less than $c_0(1 + \mathrm{dist}(I,\partial\mathbb{I})^{-1})$.

<u>Step</u> <u>2</u>: To see about the $L^6$ norm of $\chi_I|a|$, use the function $\chi$ to define a non-negative function on $\mathbb{I}$ that is equal to 1 on the support of $\chi_I$ and equal to zero where the distance to I is greater than the minimum of 3 and $\frac{3}{4}$ times the distance from I to $\partial\mathbb{I}$. This function can and should be constructed so that the norm of its derivative is no greater than $c_0(1 + \mathrm{dist}(I,\partial\mathbb{I})^{-1})$ and so that the norm of its second derivative is no greater than $c_0(1+\mathrm{dist}(I,\partial\mathbb{I})^{-2})$. Denote this function by $\chi_{I*}$.

<u>Step</u> <u>3</u>: Multiply both sides of (2.25) by $\chi_{I*}|a|$ and then integrate the resulting identity. Let $u = |a|^{3/2}$. Keeping in mind that $\chi_{I*} = 1$ on the support of $\chi_I$, an integration by parts leads to an inequality that has the form

$$\int_{\mathbb{I}\times M} \chi_{I*}\,|du|^2 \le c_0 \int_{\mathbb{I}\times M} |\tfrac{\partial}{\partial s}\chi_{I*}||a|^2|\nabla_A a| + c_0 \int_{\mathbb{I}\times M} \chi_{I*}\,|a|^3 + c_0 \int_{\mathbb{I}\times M} \chi_I\,|a||d_A{*}a|^2 + c_0 \int_{\mathbb{I}\times M} \chi_{I*}\,|a||\nabla_A a||d_A{*}a|\,. \tag{2.29}$$

The terms on the left hand side of this equation can be bounded using (2.23) and the supremum norm bound on $|r d_A{*}a|$ with the result being a bound on the integral of the function $\chi_{I*}|du|^2$ by $c_{I\diamond}(1+r^{-2}(c_{\mathfrak{cs}}+c_{div*}))^{3/2}$ with $c_{I\diamond} \le c_0(1+\mathrm{dist}(I,\partial\mathbb{I})^{-1})$.

<u>Step</u> <u>4</u>: Use Step 3's bound for the left hand side of (2.29) with the inequality in the top bullet of (2.1) to bound the $L^4$ norm of $|a|^{3/2}$ on the support of $\chi_I$. This bound leads to a $c_{I*}(1+r^{-2}(c_{\mathfrak{cs}}+c_{div*}))^{1/2}$ bound for the $L^6$ norm $|a|$ on the support of $\chi_I$ with $c_{I*}$ denoting a number less than $c_0(1+\mathrm{dist}(I,\partial\mathbb{I})^{-1})$.

**e) The integrals of $|a|^2$ over I and over $\mathbb{I}$**

Take $\mathbb{I}$ in this section to be a closed, bounded interval of length greater than 4 and take $\mathrm{I} \subset \mathbb{I}$ to be an open interval with compact closure in the interior of $\mathbb{I}$. Suppose for the moment that $(\mathrm{A}, \mathfrak{a})$ is a solution to (1.2) on $\mathbb{I}$. Fix $m > 1$. With the second bullets of Theorem 1.2 and Proposition 2.2 in mind, this subsection talks about the inequality

$$\int_{\mathbb{I}\times \mathrm{M}} |\mathfrak{a}|^2 < m \int_{\mathrm{I}\times \mathrm{M}} |\mathfrak{a}|^2 . \tag{2.30}$$

It proves convenient to take $r$ in this section to be the maximum of 1 and the $\mathrm{L}^2$ norm of $\mathfrak{a}$ on I. Define $a$ to be $r^{-1}\mathfrak{a}$ so that (2.30) asserts the bound

$$\int_{\mathbb{I}\times \mathrm{M}} |a|^2 \le m . \tag{2.31}$$

If $(\mathrm{A}, a)$ is an instanton solution to (2.18), then (2.31) holds for $m = 2$ and thus for $m > 2$ if $r$ is larger than $c_0(1 + c_{\mathfrak{cs}} + c_{\mathrm{div}})\,\mathrm{length}(\mathbb{I})$. That this is so follows from (2.22). If $(\mathrm{A}, a)$ does not obey (2.18), then (2.31) need not hold for any a priori choice for $m$ even in the case when $|a|$ is an $\mathrm{L}^2$ function on $\mathbb{I}\times\mathrm{M}$.

Suppose that $(\mathrm{A}, a)$ obeys (2.2) on $\mathbb{I}\times\mathrm{M}$ and that $m > 1$ is such that (2.31) holds. Equation (2.5) holds on $\mathbb{I}\times\mathrm{M}$ as does (2.6). Much the same integration by parts arguments that led from (2.6) to (2.8) can be repeated using the function $\chi_{\mathrm{I}}$ with (2.31) to restrict the integration domain. The result is a bound of the form

$$\int_{\mathrm{I}\times \mathrm{M}} (|\nabla_{\mathrm{A}} a|^2 + 2r^2 |a\wedge a|^2) \le c_{\mathrm{I}}\, m \tag{2.32}$$

with $c_{\mathrm{I}} \le c_0(1 + \mathrm{dist}(\mathrm{I}, \partial\mathbb{I})^{-2})$. The notation here is such that the covariant derivative is the full covariant derivative on the 4-manifold $\mathrm{X} = \mathrm{I}\times\mathrm{M}$. This bound is superseded by the bound in (2.22) when $(\mathrm{A}, a)$ is an instanton solution to (2.18) and $r$ is sufficiently large.

By the same token, the integration by parts that leads from (2.6) to (2.10) can likewise be repeated using the Green's function in (2.27) and $\chi_{\mathrm{I}}$ with the result this time being a bound that has the form

$$\tfrac{1}{2}|a|^2(\mathrm{p}) + \int_{\mathrm{I}\times \mathrm{M}} \mathrm{G}_{\mathrm{p}}\big(|\tfrac{\partial}{\partial s} a|^2 + |\nabla_{\mathrm{A}} a|^2 + 2r^2\, |a\wedge a|^2\big) \le c_{\mathrm{I}}(1+m). \tag{2.33}$$

with $c_I \le c_0(1+\mathrm{dist}(I,\partial\mathbb{I})^{-4})$. This bound is identical to that in (2.26) when $(A, a)$ is an instanton solution to (2.18) and $r$ is large.

**f) Proof of Proposition 2.1 and Proposition 2.2**

The proof of Proposition 2.1 is contained in Parts 1-4 of what follows. The proof of Proposition 2.2 is nearly word for word identical to that of Proposition 2.1. Part 5 discusses only the salient differences.

*Part 1*: This part proves what is asserted by the first bullet of the proposition. To start, let $\Lambda \subset \{1, 2, \ldots\}$ denote a subsequence chosen so that $\{r_n\}_{n\in\Lambda}$ is bounded. Let $r_*$ denote the upper bound. Fix $n \in \Lambda$ and invoke the $(A = A_n, a = \mathfrak{a}_n)$ versions of (2.5) and (2.8) to see that the sequence $\{A_n\}_{n\in\Lambda}$ is a sequence of connections on P with an a priori $L^2$ bound on its curvature.

Granted the a priori $L^2$ bound, Uhlenbeck's theorem [U] (see also Section 5c) supplies the following data set: A finite set $\Theta \subset X$, a principal G bundle $P_\Delta \to X$ and an $L^2_1$ connection, $A_\Delta$, on $P_\Delta$, a subsequence $\Lambda' \subset \Lambda$ and a corresponding set $\{g_n\}_{n\in\Lambda'}$ of isomorphisms from $P_\Delta|_{X-\Theta}$ to $P|_{X-\Theta}$. This data has the following property: The sequence $\{g_n{}^*A_n\}_{n\in\Lambda'}$ converges weakly in in the $L^2_1$ topology on compact subsets of $X-\Theta$ to the connection $A_\Delta$.

With $A_\Delta$ in hand, invoke the bound in (2.8) to find a subsequence $\Xi \subset \Lambda'$ such that $\{g_n{}^*\mathfrak{a}_n\}_{n\in\Xi}$ converges weakly in the $L^2_1$ topology on compact subsets of $X-\Theta$ to an $L^2_1$ section over X of $(P_\Delta \times_{SO(3)} \mathfrak{G}) \otimes T^*X$. Denote the latter by $\mathfrak{a}_\Delta$. The weak convergence of $\{(g_n{}^*A_n, g_n{}^*\mathfrak{a}_n)\}_{n\in\Lambda}$ to $(A_\Delta, \mathfrak{a}_\Delta)$ implies among other things that $(A_\Delta, \mathfrak{a}_\Delta)$ obeys the equations in (1.1) or (1.2) as the case may be. If it is (1.2), then the relevant version of $\tau$ is a limit of $\{\tau_n\}_{n\in\Lambda}$. Since the equations in (1.1) and (1.2) are elliptic, straightforward bootstrapping arguments can be used to prove that $\{(g_n{}^*A_n, g_n{}^*\mathfrak{a}_n)\}_{n\in\Lambda}$ converges to $(A_\Delta, \mathfrak{a}_\Delta)$ in the $C^\infty$ topology on compact subsets of $X-\Theta$. The fact that $\Theta$ is empty when $\tau$ is not 0 or 1 follows from (2.10) because the latter when $\tau \in (0, 1)$ with the equations in (2.11) leads to an a priori bound on any given $p \in X$ version of the sequence whose n'th term is the integral over X of the function $\mathrm{dist}(\cdot,p)^{-2}|F_{A_n}|^2$ with the bound being independent of the chosen point p.

*Part 2*: This part of the proof addresses Item a), the first two assertions in Item b) and Item c) of the second bullet of the proposition. Item a) holds with the subsequence $\Lambda$ being the whole of $\{1, 2, \ldots\}$. This follows from the various $(A, a) \in \{(A_n, a_n)\}_{n=1,2,\ldots}$ versions of (2.8) and the fact that the corresponding versions of (2.10) with (2.8) and middle bullet of (2.1) lead to a $c_0$ bound on the sequence $\{\sup_X |a_n|\}_{n=1,2,\ldots}$.

The a priori bound for the sequence whose n'th term is the $L^2_1$ norm of $|a_n|$ leads directly to the following conclusion: There is a subsequence $\Lambda' \subset \{1, 2, \ldots\}$ and an $L^2_1$ function, $|\hat{a}_\diamond|$, on X such that $\{|a_n|\}_{n\in\Lambda'}$ converges weakly in the $L^2_1$ topology to $|\hat{a}_\diamond|$. The a priori bound for $\sup_X |a_n|$ implies that $|\hat{a}_\diamond|$ is in $L^\infty$ with its $L^\infty$ norm bounded by a purely geometric constant. The weak $L^2_1$ convergence of $\{|a_n|\}_{n\in\Lambda'}$ to $|\hat{a}_\diamond|$ implies strong $L^p$ convergence for $p < 4$. This and the fact that $\{\sup_X |a_n|\}_{n\in\Lambda'}$ is bounded implies that $\{|a_n|\}_{n\in\Lambda'}$ converges strongly to $|\hat{a}_\diamond|$ in all $p < \infty$ version of the $L^p$ topology. The part of Item b) that concerns the pointwise definition of $|\hat{a}_\diamond|$ is addressed in Part 4.

Consider now Item c). To this end, keep in mind that the elements of the sequence $\{a_n\}_{n\in\Lambda'}$ enjoy an a priori $L^\infty$ bound. This being the case, there is a subsequence $\Lambda'' \subset \Lambda'$ such that $\{\langle a_n \otimes a_n\rangle\}_{n\in\Lambda''}$ converges weakly in the $L^2_1$ topology to an $L^2_1$ section of $T^*X \otimes T^*X$. The convergence of $\{\langle a_n \otimes a_n\rangle\}_{n\in\Lambda''}$ to its limit is strong $L^2$ convergence. The a priori $L^\infty$ bound for $\{a_n\}_{n\in\Lambda''}$ implies that the $\{\langle a_n \otimes a_n\rangle\}_{n\in\Lambda''}$ converges strongly in all $p < \infty$ versions of the $L^p$ topology. Since the trace of any given $n \in \{1, 2, \ldots\}$ version of $\langle a_n \otimes a_n\rangle$ is $|a_n|^2$, it follows that the trace of the weak $L^2_1$ limit of $\{\langle a_n \otimes a_n\rangle\}_{n\in\Lambda''}$ is $|\hat{a}_\diamond|^2$.

*Part 3*: This part of the proof addresses Item d). It follows from (2.17) that the sequence $\{\int_X f\,|F_{A_n} - r_n^2 a_n \wedge a_n|^2\}_{n\in\{1,2,\ldots\}}$ is bounded if $f$ is bounded. The existence of a subsequence $\Lambda \subset \Lambda''$ such that any $f \in C^0(X;\mathbb{R})$ version of this sequence converges follows from using two facts, the first being that the assignment of a bounded function $f$ to $\{\int_X f\,|F_{A_n} - r_n^2 a_n \wedge a_n|^2\}_{n\in\{1,2,\ldots\}}$ defines a bounded sequence of linear functionals on $C^0(X;\mathbb{R})$. The second fact is that $C^0(X;\mathbb{R})$ has a countable dense set. Keeping in mind the a priori bounds in (2.8), much the same argument proves that $\Lambda$ can be chosen so that the sequences $\{\int_X f\,|\nabla_{A_n} a_n|^2\}_{n\in\Lambda}$ and $\{r_n^2 \int_X f\,|a_n \wedge a_n|^2\}_{n\in\Lambda}$ converge when $f \in C^0(X;\mathbb{R})$.

The equation in Item d) for the limits $Q_{\nabla,f}$ and $Q_{\wedge,f}$ follows by invoking the various $(A, a) \in \{(A_n, a_n)\}_{n\in\Lambda}$ versions of (2.6).

*Part 4*: This last part of the proof addresses the assertions in Item e) of the second bullet and the part of Item b) of the second bullet that defines $|\hat{a}_\diamond|$ at each point in X.

The existence of an a priori bound for the integrals of the various $n \in \Lambda$ versions of $G_p(|\nabla_{A_n} a_n|^2 + r_n^2 |Q(a_n)|^2)$ follows directly from (2.10).

The assertion in Item e) and the part of Item b) that gives a pointwise definition for $|\hat{a}_\diamond|$ are addressed simultaneously in three steps. With regards to the definition of $|\hat{a}_\diamond|$ at each point, keep in mind that an $L^\infty$ function is a priori defined on the complement of a measure zero set and thus any pointwise definition can be changed at will on any measure zero set in X.

Step 1: Fix $p \in \Lambda$. The first observation concerns the limit of the sequence that is indexed by $\Lambda$ with n'th term being

$$\mathrm{I}_n(p) = -\int_X G_p\left(\tfrac{1}{2}|a_n|^2 - \mathrm{Ric}(\langle a_n \otimes a_n\rangle)\right). \tag{2.34}$$

The inequalities in (2.9) with what is said by Item a) and Item c) of the second bullet about strong $L^p$ convergence for $p < \infty$ imply that the sequence of functions $\{\mathrm{I}_n(\cdot)\}_{n\in\Lambda}$ converges in the $C^1$ topology to the function on X whose value at any given point p is equal to $-\int_X G_p\left(\tfrac{1}{2}|\hat{a}_\Diamond|^2 - \mathrm{Ric}(\langle \hat{a}_\Diamond \otimes \hat{a}_\Diamond\rangle)\right)$. Use $\mathrm{I}_\Diamond$ to denote this limit function.

Step 2: Let $\mathrm{v}$ denote the function given by the rule $p \to \mathrm{v}(p) = \limsup_{n\to\infty} |a_n|(p)$. Item e) of the proposition's second bullet defines each $p \in X$ version of $Q_{\Diamond,p}$. As explained directly,

$$\tfrac{1}{2}\mathrm{v}^2(p) + Q_{\Diamond,p} = \mathrm{I}(p) \quad \textit{for each } p \in X. \tag{2.35}$$

To derive this identity, it proves useful to introduce the following notation: Given $n \in \Lambda$, let $Q_{n,p}$ denote the integral of $G_p(|\nabla_{A_n} a_n|^2 + 2r_n^{\,2}|a_n \wedge a_n|^2)$. Choose a sequence $\Xi \subset \Lambda$ such that $\lim_{n\in\Xi} Q_{n,p} = Q_{\Diamond,p}$. Fix $\varepsilon > 0$. If $n \in \Xi$ is sufficiently large, then $Q_{n,p}$ will be less than $Q_{\Diamond,p} + \varepsilon$. It is also the case that $|a_n|(p)$ will be less than $\limsup_{n'\to\infty} |a_{n'}|(p) + \varepsilon$ when n is large because of what is meant by lim sup. Since $\tfrac{1}{2}|a_n|^2(p) + Q_{n,p} = \mathrm{I}_n(p)$ and since $\{\mathrm{I}_n\}_{n\in\Lambda}$ converges uniformly to $\mathrm{I}_\Diamond$, it follows as a consequence that $\tfrac{1}{2}\mathrm{v}(p)^2 + Q_{\Diamond,p}$ is no greater than $\mathrm{I}_\Diamond(p) + c_0\varepsilon$. Since $\varepsilon$ can be as small as desired, this implies that $\tfrac{1}{2}\mathrm{v}(p)^2 + Q_{\Diamond,p} \le \mathrm{I}_\Diamond(p)$.

Choose a second sequence $\Xi' \subset \Lambda$ such that $\lim_{n\in\Xi'} |a_n|(p) = \mathrm{v}(p)$. If $n \in \Xi'$ is sufficiently large, then $|a_n|(p) > \mathrm{v}(p) - \varepsilon$. Meanwhile, $Q_{n,p}$ will be greater than $Q_{\Diamond,p} - \varepsilon$ when n is sufficiently large because of what is meant by lim inf. Since $\{\mathrm{I}_n\}_{n\in\Lambda}$ converges to $\mathrm{I}_\Diamond$, these last observations imply that $\tfrac{1}{2}\mathrm{v}(p)^2 + Q_{\Diamond,p} \ge \mathrm{I}_\Diamond(p) - c_0\varepsilon$. Since $\varepsilon$ can be as small as desired, this last inequality implies that $\tfrac{1}{2}\mathrm{v}(p)^2 + Q_{\Diamond,p} \ge \mathrm{I}_\Diamond(p)$.

Step 3: The function $\mathrm{v}(\cdot)^2$ is a measurable function (see [Rudin, Theorem 11.17], and since it is bounded, Fatou's lemma (see [Rudin, Theorem 11.31]) can be invoked to see that $\mathrm{v} \ge |\hat{a}_\Diamond|$ almost everywhere on X in the sense that this inequality must hold on a set of full measure. The proof that $|\hat{a}_\Diamond| \geq \mathrm{v}$ on a set of full measure mimics what is done in Section 2c of [T2] to prove an analogous assertion in the context of $PSL(2;\mathbb{C})$ connections on 3-dimensional manifolds. In particular, the arguments that follow mimic those in Step 8 in the proof of the sixth bullet of Proposition 2.2 in [T2].

To start the proof, fix $\varepsilon \in (0, 1)$. Given $p \in X$, define functions $\delta_{p,\varepsilon}$ and $f_{p,\varepsilon}$ as follows: Let $\upsilon_{p,\varepsilon}$ denote the function that is equal to 1 on the radius $\varepsilon(1-\varepsilon^3)$ ball centered at p, equal to $(1-\varepsilon^{-1}\mathrm{dist}(p,\cdot))/\varepsilon^3$ on the spherical annulus centered at p with inner radius $\varepsilon(1-\varepsilon^3)$ and outer radius $\varepsilon$, and equal to zero on the complement of the radius $\varepsilon$ ball centered at p. Write the integral of $\upsilon_{p,\varepsilon}$ as $A_{p,\varepsilon}\varepsilon^4$. The number $A_{p,\varepsilon}$ differs from $\frac{1}{2}\pi^2$ by at most $c_0\varepsilon^2$. Define $\delta_{p,\varepsilon}$ to equal $A_{p,\varepsilon}{}^{-1}\varepsilon^{-4}\upsilon_{p,\varepsilon}$. The function $f_{p,\varepsilon}$ is the solution to the equation $d^\dagger d f + f = \delta_{p,\varepsilon}$. This function $f_{p,\varepsilon}$ turns out to be a $C^2$ function on X.

Given that $\varepsilon < c_0{}^{-1}$, the function $f_{p,\varepsilon}$ can be written on the complement of the radius $\varepsilon$ ball centered at p as

$$f_{p,\varepsilon} = G_p + \mathfrak{t}_{p,\varepsilon} \quad \textit{with } |\mathfrak{t}_{p,\varepsilon}| \leq c_0\,\varepsilon\, G_p\,, \tag{2.36}$$

Meanwhile, $f_{p,\varepsilon}$ on the radius $\varepsilon$ ball centered at p can be written as

$$f_{p,\varepsilon} = \tfrac{1}{8} A_{p,\varepsilon}{}^{-1}\varepsilon^{-2}(2 - \varepsilon^{-2}\mathrm{dist}(\cdot,p)^2) + \mathfrak{w}_{p,\varepsilon} \quad \textit{with } |\mathfrak{w}_{p,\varepsilon}| \leq c_0\,\varepsilon^{-1}\,. \tag{2.37}$$

The sequence $\{f_{p,\varepsilon}\}_{\varepsilon\in(0,1)}$ converges to $G_p$ as $\varepsilon \to 0$ in the $C^\infty$ topology on compact subsets in X−p. The following is also true

*Given* $\rho \in (0, c_0{}^{-1})$*, there exists* $c_\rho > 1$ *which is independent of* $p \in X$ *and such that if* $\varepsilon < c_\rho{}^{-1}$*, then* $f_{p,\varepsilon}$ *obeys* $f_{p,\varepsilon} < G_p + \mathfrak{t}_{p,\varepsilon,\rho}$ *on the whole of* X *where* $\mathfrak{t}_{p,\varepsilon,\rho}$ *is a continuous function with* $|\mathfrak{t}_{p,\varepsilon,\rho}| < \rho$.

(2.38)

This follows from the comparison principle given that $G_p > f_{p,\varepsilon'}$ where the distance to p is less than $\varepsilon$. The proof that $G_p > f_{p,\varepsilon'}$ where $\mathrm{dist}(\cdot,p) < \varepsilon$ when $\varepsilon$ is small follows from (2.38) using the fact that $G_p$ differs from $\frac{1}{8}A_{p,\varepsilon}{}^{-1}\mathrm{dist}(\cdot,p)^{-2}$ by at most $c_0\,\mathrm{dist}(\cdot,p)^{-1}$ on the ball of radius $c_0{}^{-1}$ centered at p.

<u>Step 4</u>: Fix $\rho \in (0, c_0{}^{-1})$ and $\varepsilon < c_\rho{}^{-1}$. Given $n \in \{1, 2, \ldots\}$ and $p \in X$, multiply both sides of (2.6) by $f_{p,\varepsilon}$, and then integrate the result over X. Having done so, invoke (2.36) to see that

$$\tfrac{1}{2}\int_X \delta_{p,\varepsilon}|a_n|^2 + Q_{n,p} = I_\cdot(p) - e_{n,\varepsilon,\rho} \tag{2.39}$$

where $e_{n,\varepsilon,\rho}$ is such that $\lim_{n\in\Lambda} \sup_{p\in M} |e_{n,\varepsilon,\rho}| \le c_0\,\rho$. Meanwhile, the fact that $\{|a_n|\}_{n\in\Lambda}$ converges strongly to $|\hat{a}_\Diamond|$ in the $L^2$ topology has the following implication: Given $\varepsilon < c_\rho^{-1}$ and $\varepsilon' \in (0, 1]$, there exists $n_{\varepsilon,\varepsilon'} \ge 1$ which is independent of $p \in X$ and such that

$$\sup_{p\in X} \Big| \int_X \delta_{p,\varepsilon} |\hat{a}_\Diamond|^2 - \int_X \delta_{p,\varepsilon} |a_n|^2 \Big| < \varepsilon' \quad \textit{when } n > n_{\varepsilon,\varepsilon'} \,. \tag{2.40}$$

Fix $p \in X$, and take $\varepsilon' = \rho$ and then $n > n_{\varepsilon,\varepsilon'=\rho}$ with two additional constraints. The first requires that $|Q_{n,p} - Q_\Diamond(p)| < \rho$. This constraint is obeyed if n is chosen from a subsequence in $\Lambda$ whose versions of $Q_{n,p}$ converge to $Q_{\Diamond,p}$ as the index n increases. The second constraint asks that n be such that $\sup_{p\in M} |e_{n,\varepsilon,\rho}| \le c_0\,\rho$ with $e_{n,\varepsilon,\rho}$ as in (2.36). This constraint is obeyed if n is sufficiently large. Having chosen n to obey all of these constraints, then the $\varepsilon' = \rho$ version of (2.40) and (2.39) imply that

$$\tfrac{1}{2} \int_X \delta_{p,\varepsilon} |\hat{a}_\Diamond|^2 + Q_\Diamond(p) \ge I_\Diamond(p) - c_0\rho \,. \tag{2.41}$$

This equation holds for each $p \in X$. Since $\rho$ and then $\varepsilon$ can be chosen as small as desired, a comparison between (2.41) and (2.35) finds $|\hat{a}_\Diamond| \ge v$ on the complement of a measure zero set in X.

*Part 5*: This last part of section concerns the proof of Proposition 2.2. The proof in the case when P is viewed as the pull-back of a namesake on M and $\{(A_n, \mathfrak{a}_n)\}_{n=1,2,\ldots}$ is the $\mathbb{I} \times M$ incarnation of a sequence of maps from $\mathbb{I}$ to $\mathrm{Conn}(P|_M) \times C^\infty((P \times_G \mathfrak{g})|_M \otimes T^*M)$ is word for word identical to the proof given in Parts 1-4 of Proposition 2.1 but for the replacement of referrals to equations in Sections 2a and 2b with their counterparts in Section 2c. In particular, (2.23) and (2.26) replace (2.8) and (2.10). Note in particular that there is no need to assume that $d_{A_n}\mathfrak{a}_n = 0$ for each n if there is an n-independent bound on $|d_{A_n}\mathfrak{a}_n|$.

Suppose now that P is not a priori given as the pull-back of a namesake on M, and that $\{(A_n, a_n)\}_{n=1,2,\ldots}$ is not a priori the incarnation of a corresponding sequence of maps from $\mathbb{I}$ to $\mathrm{Conn}(P|_M) \times C^\infty(M; (P \times_G \mathfrak{g})|_M \otimes T^*M)$. The proof of Proposition 2.2 in this case is also a nearly verbatim repetition of what is said in Parts 1-4 with (2.32) and (2.33) playing the respective roles of (2.8) and (2.10). This being the case, the remaining paragraphs of Part 5 give a very brief summary of facts that are needed to derive the requisite versions of (2.23) and (2.26).

**g) The equations in (1.2) on $\mathbb{I} \times M$ as gradient flow equations**

This subsection is parenthetical in nature; it has no bearing on what is done in subsequent sections. What follows explains how the equations in (1.2) can be seen as (formal) gradient flow equations when the solution of (1.2) does not also obey (1.6).

To start, suppose that A is a given connection on P. Fix a point $s_0 \in \mathbb{I}$ and use parallel transport along the constant M slices of $\mathbb{I} \times M$ from $P|_{s_0 \times M}$ to define an isomorphism between P and $P|_{s_0 \times M}$. Doing so writes A as a map from $\mathbb{I}$ to $\mathrm{Conn}(P|_{s_0 \times M})$. Now suppose that $\mathfrak{a}$ is a section over $\mathbb{I} \times M$ of $T^*(\mathbb{I} \times M) \otimes (P \times_{SO(3)} \mathfrak{S})$. Write this section $\mathfrak{a}$ as $\varphi\, ds + \mathfrak{w}$ with $\mathfrak{w}$ annihilating the vector field $\frac{\partial}{\partial s}$. The isomorphism between P and $P|_{s_0 \times M}$ allows $\mathfrak{w}$ to be viewed as a map from $\mathbb{I}$ to $T^*M \otimes (P|_{s_0 \times M} \times_{SO(3)} \mathfrak{S})$ and $\varphi$ to be viewed as a map from $\mathbb{I}$ to $P|_{s_0 \times M} \times_{SO(3)} \mathfrak{S}$.

If $(A, \mathfrak{a} = \varphi\, ds + \mathfrak{w})$ obeys (1.2), then the maps A, $\mathfrak{w}$ and $\varphi$ obey

- $\frac{\partial}{\partial s} A - [\varphi, \mathfrak{w}] + \frac{\tau^2 - (1-\tau)^2}{\tau^2 + (1-\tau)^2} *(F_A - \mathfrak{w} \wedge \mathfrak{w}) - \frac{2\tau(1-\tau)}{\tau^2 + (1-\tau)^2} *d_A\mathfrak{w} = 0$ .
- $\frac{\partial}{\partial s} \mathfrak{w} - d_A\varphi - \frac{\tau^2 - (1-\tau)^2}{\tau^2 + (1-\tau)^2} *d_A\mathfrak{w} - \frac{2\tau(1-\tau)}{\tau^2 + (1-\tau)^2} *(F_A - \mathfrak{w} \wedge \mathfrak{w}) = 0$ .
- $\frac{\partial}{\partial s} \varphi + *d_A{*}\mathfrak{w} = 0$ .

(2.42)

These equations are, up to a sign, the gradient flow equations of a function on the space $\mathrm{Conn}(P|_{s_0 \times M}) \times C^\infty(M; T^*M \otimes (P|_{s_0 \times M} \times_G \mathfrak{S})) \times C^\infty(M; P|_{s_0 \times M} \times_G \mathfrak{S})$ that is given up to a sign by the rule that assigns a given triple $(A, \mathfrak{w}, \varphi)$ to

$$\mathrm{re}\Big(\frac{(\tau + i(1-\tau))^2}{\tau^2 + (1-\tau)^2}\, \mathfrak{cs}(A, \mathfrak{w})\Big) + \int_{s \times M} \langle \mathfrak{w} \wedge *d_A\varphi \rangle .$$

(2.43)

Note that if $(A, \mathfrak{a})$ obeys (1.6) then $\varphi$ obeys the equation

$$\frac{\partial^2}{\partial s^2} \varphi + d_A{}^\dagger d_A\varphi + *(\mathfrak{w} \wedge [\varphi, *\mathfrak{w}] + [\varphi, *\mathfrak{w}] \wedge \mathfrak{w}) = 0.$$

(2.44)

The latter leads to the identity

$$\tfrac{1}{2}\Big(-\frac{\partial^2}{\partial s^2} |\varphi|^2 + d^\dagger d\, |\varphi|^2\Big) + \Big|\frac{\partial}{\partial s} \varphi\Big|^2 + |\nabla_A\varphi|^2 + 2\,|[\varphi, \mathfrak{w}]|^2 = 0 .$$

(2.45)

This last equation with the maximum principle has the following implication: If $\sup_M |\varphi|$ has limit zero in the non-compact directions on $\mathbb{I}$ and on approach to any boundary point of $\mathbb{I}$, then $\varphi = 0$ on the whole of $\mathbb{I} \times M$ and the solution in question is a solution to (1.6).

## 3. The numbers $r_{c_\wedge}$ and $r_{cF}$

To set the stage for what is to come, suppose in what follows that X is an oriented, Riemannian 4-manifold and that U is a specified open set in X with compact closure. Use P to denote a given principal SO(3) bundle over X. The term *geometric data* is used in what follows to denote the data set consisting of the manifold X, the subset U and the bundle P. A second convention in what follows has $c_0$ denoting a number that is greater than 1 and depends only on this geometric data. The precise value of $c_0$ can be assumed to increase between incarnations.

Fix $r > 1$ and let (A, $a$) denote a pair of connection on P and section of the bundle $T^*X \otimes (P \times_{SO(3)} \mathfrak{S})$ that obey a given $\tau \in [0, 1]$ version of (2.11) on X. These equations imply among other things that (A, $a$) obey (2.5) and (2.6). The numbers $r_{c_\wedge}$ and $r_{cF}$ that are defined momentarily measure the respective size of $r^4\, a \wedge a$ and $F_A$ near a given point in X. The definitions of $r_{c_\wedge}$ and $r_{cF}$ require the a priori specification of a positive number to be denoted by $c$. The convention takes $c > 100$. With $c$ chosen, fix $p \in U$ and define the number $r_{c_\wedge}$ to be the largest of the numbers $r \in (0, c_0^{-1}]$ for which

$$r^4 \int_{B_r} |a \wedge a|^2 \le c^{-2} \,. \tag{3.1}$$

By the same token, define the number $r_{cF}$ to be the largest of those $r \in (0, c_0^{-1}]$ for which

$$\int_{B_r} |F_A|^2 \le c^{-2} \,. \tag{3.2}$$

The numbers $r_{c_\wedge}$ and $r_{cF}$ depend on the chosen point p but since p is fixed in most of what follows, this dependence is not indicated by the notation. The subsequent subsections of this section and Section 4 state and prove various properties of (A, $a$) on the ball centered at p with radius less than the smaller of $r_{c_\wedge}$ and $r_{cF}$.

### a) The functions K and N

Keep in mind in what follows that a ball of radius less than $c_0^{-1}$ about any given point in U will be well inside a Gaussian coordinate chart for a neighborhood of the point in the manifold X. Suppose again that $r > 1$ and (A, $a$) are a pair of connection on P and section of the bundle $T^*X \otimes (P \times_{SO(3)} \mathfrak{S})$ that obey a given $\tau \in [0, 1]$ version of (2.11) on X. The data consisting of $r$, (A, $a$) and a point $p \in U$ is used momentarily to define two positive functions on $[0, c_0^{1}]$, these denoted by K and N.

The definition of $\mathrm{K}$ uses two auxilliary functions on $[0, c_0^{-1}]$. The first of the auxilliary functions is denoted by h and it is defined by the rule

$$r \to h(r) = \int_{\partial B_r} |a|^2 \quad . \tag{3.3}$$

With regards to the notation, $\partial B_r$ denotes the boundary of the closure of $B_r$. Arguments much like those used by Aronszajn [Ar] can be used with (2.6) to prove that $a$ can not vanish identically on $\partial B_r$. This implies in particular that h is positive on $(0, c_0^{-1}]$.

The second function is denoted by $\mathfrak{d}$. The definition of $\mathfrak{d}$ is given below in (3.4). What is denoted by $\mathrm{M}$ in this equation is defined by writing the trace of the second fundamental form of $\partial B_r$ for any given $r \in (0, c_0^{-1}]$ as $\frac{3}{r} + \mathrm{M}$. The function $\mathfrak{d}$ is defined by the rule

$$r \to \mathfrak{d}(r) = \int_0^r \Big(\frac{1}{h(s)}\Big(\int_{B_s} \mathrm{Ric}(\langle a \otimes a \rangle) + \tfrac{1}{2} \int_{\partial B_s} \mathrm{M}\, |a|^2\Big)\Big)\, ds \quad . \tag{3.4}$$

As explained momentarily, the absolute value of the function $\mathfrak{d}$ is bounded by $c_0 r^2$.

With $\mathfrak{d}$ in hand, define $\mathrm{K}$ to be the positive square root of $e^{-2\mathfrak{d}}\, r^{-3} h$. With $\mathrm{K}$ in hand, the function $\mathrm{N}$ is then defined by the rule

$$r \to \mathrm{N}(r) = \frac{1}{r^2 \mathrm{K}(r)^2} \int_{B_r} \big(|\nabla_A a|^2 + 2r^2 |a \wedge a|^2\big) \quad . \tag{3.5}$$

The definition of $\mathfrak{d}$ is designed so as to relate $\mathrm{N}$ to $\mathrm{K}$ via the identity

$$\frac{d}{dr} \mathrm{K} = \frac{1}{r} \mathrm{N}\, \mathrm{K} \quad . \tag{3.6}$$

This identity is proved by first using the definition of $\mathrm{M}$ to write the derivative of h as

$$\frac{d}{dr} h = \frac{3}{r} h + \int_{\partial B_r} \frac{\partial}{\partial r} |a|^2 + r^{-3} \int_{\partial B_r} 2\mathrm{M}\, |a|^2 \quad . \tag{3.7}$$

Stoke's theorem with (2.6) is then used to identify the integral over $\partial B_r$ of $\frac{\partial}{\partial r} |a|^2$ with

$$2 \int_{B_r} \big(|\nabla_A a|^2 + 2r^2 |a \wedge a|^2 + \mathrm{Ric}(\langle a \otimes a \rangle)\big) \quad . \tag{3.8}$$

The identity in (3.6) follows directly from (3.4), (3.5) and (3.8) because $\mathrm{K}^2 = e^{2\mathfrak{d}} r^{-3} h$.

The preceding paragraph asserted the bound $|\mathfrak{d}| \leq c_0 r^2$. The proof that this is so invokes the identity in the upcoming (3.9) for $L^2_1$ functions on $B_r$. To set the stage for (3.9), assume that $r \in (0, c_0^{-1}]$ and that $f$ is a given Sobolev class $L^2_1$ function on $B_r$. The identity writes

$$\int_{B_r} f^2 = \tfrac{1}{4} r \int_{\partial B_r} f^2 + \mathfrak{e}_f(r) \tag{3.9}$$

with $\mathfrak{e}_f$ denoting a function whose absolute value obeys

$$|\mathfrak{e}_f(r)| \leq c_0 r^2 \Big( \int_{B_r} f^2 \Big) + c_0\, r \Big( \int_{B_r} f^2 \Big)^{1/2} \Big( \int_{B_r} |df|^2 \Big)^{1/2} . \tag{3.10}$$

To prove that $|\mathfrak{d}| \leq c_0 r^2$, invoke (3.9) and (3.10) with $f = |a|$ to see that what is written in (3.8) is no smaller than $-c_0 r h$. As $\mathrm{M}$ is no smaller than $-c_0 r$, it then follows by integrating (3.7) that the function h obeys

$$h(r) \geq \big(\tfrac{r}{s}\big)^3 e^{-c_0 (r^2 - s^2)} h(s) \tag{3.11}$$

when $r \geq s$ and both are from $(0, c_0^{-1})$. It follows in turn from (3.11) that

$$\int_{B_r} |a|^2 \leq 4\, e^{c_0 r^2}\, r\, h(r) \tag{3.12}$$

when $r \in [0, c_0^{-1}]$. To see the significance of (3.12), use (3.5) to write $\mathfrak{d}(r)$ as an integral from 0 to r of a function of the variable s. Denote this function by $\mathfrak{q}$. It follows from (3.12) and the fact that $|\mathrm{M}|$ on $\partial B_s$ is bounded by $c_0 s$ that $|\mathfrak{q}(s)|$ is also bounded $c_0 s$. The latter bound implies directly the asserted $c_0 r^2$ bound for $|\mathfrak{d}|$.

**b) The pair $(A, a)$ on a ball of radius less than $r_{c_F}$ and $r_{c_\wedge}$**

Fix $p \in U$. The three propositions that follows state a priori bounds that hold on balls centered at p with radii less than the minimum of the versions of $r_{c_F}$ and $r_{c_\wedge}$. With a number $r_c \in (0, c_0^{-1}]$ specified, the upcoming propositions introduce by way of notation $r_\ddagger$ to denote $\kappa(r_c)\, r$ and $a_\ddagger$ to denote $\kappa^{-1}(r_c)\, a$.

**Proposition 3.1**: *There exists $\kappa > 1$ that depends only on the geometric data and given $\mu \in (0, \frac{1}{4}]$, there exists $\kappa_\mu > \kappa$ that depends only on $\mu$ and the geometric data; these numbers having the following significance: Fix $r > 1$ and suppose that $(A, a)$ is a* *connection on* P *and section of* $T^*X \otimes (P \times_{SO(3)} \mathfrak{S})$ *that obeys some* $\tau \in [0,1]$ *version of (2.11). Fix* $c \geq \kappa$ *and* $p \in U$ *so as to define* $r_{cF}$ *and* $r_{c\wedge}$ *using* $(A, a)$. *If* $r_c \in (0, \kappa^{-1}]$, *then the first bullet below holds. If, in addition,* $r_c$ *is less than the minimum of* $r_{cF}$ *and* $r_{c\wedge}$, *then the second bullet holds as well.*

- $|a_\ddagger| \leq \kappa_\mu$ *on* $B_{(1-\mu)r_c}$ .
- $\int_{B_{(1-\mu)r_c}} (|\nabla_A(\nabla_A a_\ddagger)|^2 + r_\ddagger^2 |\nabla_A(a_\ddagger \wedge a_\ddagger)|^2) \leq \kappa_\mu (c^{-2} + r_c^4) + \kappa_\mu N(r_c).$

This proposition is proved in Section 3c.

The upcoming Proposition 3.2 uses $e$ to denote both a unit length vector in $TX|_p$ and the vector field on $B_{r_c}$ that is obtained by the latter's parallel transport along the radial geodesic rays from p using the Levi-Civita connection. The notation uses $a_\ddagger(e)$ to denote the section of $P \times_{SO(3)} \mathfrak{S}$ over $B_{r_c}$ is obtained by taking the metric pairing between $a_\ddagger$ and the vector field $e$.

**Proposition 3.2**: *There exists* $\kappa > 1$ *that depends only on the geometric data and given* $\mu \in (0, \frac{1}{4}]$, *there exists* $\kappa_\mu > \kappa$ *that depends only on* $\mu$ *and the geometric data; and these numbers having the following significance: Fix* $r > 1$ *and suppose that* $(A, a)$ *is a connection on* P *and section of* $T^*X \otimes (P \times_{SO(3)} \mathfrak{S})$ *that obeys some* $\tau \in [0,1]$ *version of (2.11). Fix* $c \geq \kappa$ *and a point* $p \in U$ *so as to define* $r_{cF}$ *and* $r_{c\wedge}$ *using* $(A, a)$. *Suppose that* $r_c$ *is positive but no greater than the minimum of* $r_{cF}$, $r_{c\wedge}$ *and* $c^{-1}$. *If* $N(r_c)$ *is less than* $c^{-2}$, *then*

- *If* $e \in TX|_p$, *then* $\max_{B_{(1-\mu)r_c}} |a_\ddagger(e)|^2 - \min_{B_{(1-\mu)r_c}} |a_\ddagger(e)|^2 \leq \kappa_\mu c^{-1} |e|^2$
- $\left| |a_\ddagger| - \frac{1}{\sqrt{2}\,\pi} \right| < \kappa_\mu c^{-1}$ *on* $B_{(1-\mu)r_c}$ .

The proof of Proposition 3.2 is in Section 3d.

**Proposition 3.3**: *There exists* $\kappa > 1$ *that depends only on the geometric data and given* $\mu \in (0, \frac{1}{4}]$, *there exists* $\kappa_\mu > \kappa$ *that depends only on* $\mu$ *and the geometric data; these numbers having the following significance: Fix* $r > 1$ *and suppose that* $(A, a)$ *is a connection on* P *and section of* $T^*X \otimes (P \times_{SO(3)} \mathfrak{S})$ *that obeys some* $\tau \in [0,1]$ *version of (2.11). Fix* $c \geq \kappa$ *and a point* $p \in U$ *so as to define* $r_{cF}$ *and* $r_{c\wedge}$ *using* $(A, a)$. *Suppose that* $r_c$ *is positive but no greater than the minimum of* $r_{cF}$, $r_{c\wedge}$ *and* $c^{-1}$. *If* $N(r_c) \leq \kappa_\mu^{-1} c^{-2}$, *then*

- $\int_{B_{(1-\mu)r_c}} |\nabla_A F_A|^2 \leq \kappa_\mu r_c^{-2} (1 + r_c^2 r_\ddagger^2) c^{-2}$ .

- $\sup_{B_{(1-\mu)r_c}} |\nabla_A a_‡| \le \kappa_\mu r_c^{-1}(1 + r_c r_‡)\, c^{-1}$ .

Proposition 3.3 is proved in Sections 3e.

By way of notation, the proofs that follow use $c_\mu$ to denote a number that is greater than 1 and depends only on $\mu$ and the geometric data. It can be assumed to increase between successive appearances.

**c) Proof of Proposition 3.1**

The pair $(A, a_‡)$ obey (2.5) and (2.6) with $r$ replaced by $r_‡$. Keep this fact in mind. Part 1 of the proof proves the first bullet and Part 2 proves the second bullet.

*Part 1*: The proof of Proposition 3.1's top bullet is given in three steps.

Step 1: Use the function $\chi$ to construct a smooth, nonnegative function on $B_{r_c}$ that equals 0 where $\text{dist}(\cdot, p) \ge (1 - \frac{1}{4}\mu) r_c$, equals 1 where $\text{dist}(\cdot, p) \le (1 - \frac{3}{4}\mu) r_c$ and with exterior derivative having norm less than $c_0 \mu^{-1} r_c^{-1}$. Denote this function by $\chi_\mu$. Multiply both sides of the $(r_‡, (A, a_‡))$ version of (2.6) by $\chi_\mu^2$. Having done so, integrate by parts twice to see that

$$\int_{B_{(1-\mu)r_c}} (|\nabla_A a_‡|^2 + r_‡^2 |a_‡ \wedge a_‡|^2) \le c_0 \mu^{-2} r_c^{-2} \int_{B_{r_c}} |a_‡|^2 \ . \tag{3.13}$$

Invoke (3.12) to bound the right hand side of (3.13) by $c_0 \mu^{-2} r_c^2$.

Step 2: Fix a point $q \in B_{(1-\mu)r_c}$ and introduce by way of notation $G_q$ to denote the Dirichlet Green's function for the operator $d^\dagger d$ on $B_{r_c}$ . Multiply both sides of the $(r_‡, (A, a_‡))$ version of (2.6) by $\chi_\mu^2 G_q$ and then integrate the result over $B_{r_c}$ . Integrate by parts twice obtain an inequality that reads

$$\begin{aligned} \tfrac{1}{2}|a_‡|^2 + \int_{B_{r_c}} G_q \chi_\mu^2 (|\nabla_A a_‡|^2 + r_‡^2 |a_‡ \wedge a_‡|^2) &= - \int_{B_{r_c}} G_q \chi_\mu^2 \text{Ric}(\langle a_‡ \otimes a_‡ \rangle)) \\ &\quad - \int_{B_{r_c}} ((\chi_\mu d^\dagger d\chi_\mu - 2|d\chi_\mu|^2) G_q + 2\mathfrak{m}(\chi_\mu d\chi_\mu, dG_q)) |a_‡|^2 \end{aligned} \tag{3.14}$$

with $\mathfrak{m}(\cdot, \cdot)$ denoting the metric inner product. The right hand side of (3.14) is at most

$$c_0 \int_{B_{(1-\mu/4)r_c}} \frac{1}{\mathrm{dist}(\cdot,q)^2}|a_\ddagger|^2 + c_0 r_c^{-4}\mu^{-2} \int_{B_{r_c}} |a_\ddagger|^2 \tag{3.15}$$

because $G_q \le c_0 \frac{1}{\mathrm{dist}(\cdot,q)^2}$ and $|dG_q| \le c_0 \frac{1}{\mathrm{dist}(\cdot,q)^3}$ and because the support of $d\chi_\mu$ has distance no less than $c_0^{-1}\mu\, r_c$ from q.

<u>Step 3</u>: Use the second bullet in (2.1) and what is said in Step 1 with $\mu$ replaced by $\mu/4$ to bound the left most term in (3.15) by $c_0\mu^{-2} r_c^2$. Use the inequality in (3.12) meanwhile to bound the right most term in (3.15) by $c_0\mu^{-2}$. Use these bounds in (3.14) to bound $|a_\ddagger|^2$ by $c_0\mu^{-2}$.

*Part 2*: The proof of the second bullet of Proposition 3.1 is given in five steps.

<u>Step 1</u>: Fix an orthonormal frame for T*X on the radius $c_0^{-1}$ ball centered at p. What follows uses $\{a_{\ddagger\alpha}\}_{\alpha\in\{1,2,3,4\}}$ to denote the components of $a$ when written using this basis. Use $\{\nabla_{A,\beta}a_\ddagger\}_{\beta\in\{1,2,3,4\}}$ to denote the components of the four directional covariant derivatives of $a$ at p along the dual vectors to the chosen basis. Each of these four directional covariant derivatives is a section of $T^*X \otimes (P\times_{SO(3)} \mathfrak{G})$ over $B_{1/c_0}$. For each frame label $\beta \in \{1, 2, 3, 4\}$, the next equation uses $\{(\nabla_{A,\beta}a_\ddagger)_\alpha\}_{\alpha\in\{1,2,3,4\}}$ to denote the four components of $\nabla_{A,\beta}a_\ddagger$ when the latter is written using the chosen orthonormal basis. Meanwhile, $\{\nabla_{A,\beta}F_A\}_{\beta\in\{1,2,3,4]}$ denotes the 4 directional covariant derivatives of $F_A$ and for each $\beta \in \{1, 2, 3, 4\}$, the set $\{(\nabla_{A,\beta}F_A)_{\mu\nu}\}_{\mu,\nu\in\{1,2,3,4\}}$ denote the components of $\nabla_{A,\beta}F_A$. Differentiate (2.5) and commute covariant derivatives to obtain an equation for $\nabla_{A,\beta}a$ that has the schematic form

$$\begin{aligned}(\nabla_A{}^\dagger\nabla_A(\nabla_{A,\beta}a_\ddagger))_\alpha + 2\,[F_{A,\mu\beta},(\nabla_{A,\mu}a_\ddagger)_\alpha] + [(\nabla_{A,\mu}F_A)_{\mu\beta},a_{\ddagger\alpha}] + r_\ddagger{}^2[(\nabla_{A,\beta}a_\ddagger)_\mu,[a_{\ddagger\alpha},a_{\ddagger\mu}]] \\ + r_\ddagger{}^2[a_{\ddagger\mu},[(\nabla_{A,\beta}a_\ddagger)_\alpha,a_{\ddagger\mu}]] + r_\ddagger{}^2[a_{\ddagger\mu},[a_{\ddagger\alpha},(\nabla_{A,\beta}a_\ddagger)_\mu]] + \mathcal{R}^1{}_{\alpha\beta}(\nabla_A a_\ddagger) + \mathcal{R}^0{}_{\alpha\beta}(a_\ddagger) = 0\end{aligned} \tag{3.16}$$

where $\mathcal{R}^1$ and $\mathcal{R}^0$ are endomorphisms that are defined by the Riemann curvature tensor and its covariant derivatives. In particular, both have absolute value bounded by $c_0$. The convention in (3.16) and subsequently with respect to repeated indices is the same as in (2.5), this being the summation over repeated indices is implicit.

<u>Step 2</u>: Take the inner product of both sides of (3.16) with $\nabla_A a_\ddagger$ to obtain

$$\tfrac{1}{2}\, d^\dagger d\,|\nabla_A a_\ddagger|^2 + |\nabla_A(\nabla_A a_\ddagger)|^2 + \tfrac{1}{2}\, r_\ddagger{}^2|\nabla_A(a_\ddagger \wedge a_\ddagger)|^2 + 2\langle(\nabla_{A,\beta}a_\ddagger)_\alpha\,[F_{\mu\beta},(\nabla_{A,\mu}a_\ddagger)_\alpha]\rangle +$$

$$\langle(\nabla_{A,\mu}F_A)_{\mu\beta}, [a_{\ddagger\alpha},(\nabla_{A,\beta}a_\ddagger)_\alpha]\rangle + r_\ddagger^{\,2}\langle[(\nabla_{A,\beta}a_\ddagger)_\alpha,(\nabla_{A,\beta}a_\ddagger)_\mu]\,[a_{\ddagger\alpha}, a_{\ddagger\beta}]\rangle \le c_0\,(|\nabla_A a_\ddagger|^2 + |a_\ddagger|^2)\ . \tag{3.17}$$

Multiply both sides of (3.17) by $\chi_\mu^{\,2}$ and then integrate over $B_{r_c}$. After doing so, an integration by parts leads to an inequality that can be written schematically as

$$\int_{B_{r_c}} (|\nabla_A(\chi_\mu\nabla_A a_\ddagger)|^2 + r_\ddagger^{\,2}\chi_\mu^{\,2}|\nabla_A(a_\ddagger\wedge a_\ddagger)|^2) + \mathfrak{e}_Q + \mathfrak{e}_\nabla + \mathfrak{e}_{d*} \le c_\mu\, r_c^{-2}\int_{B_{r_c}} (|\nabla_A a_\ddagger|^2 + c_0\int_{B_{r_c}} |a_\ddagger|^2 \tag{3.18}$$

with $\{\mathfrak{e}_Q, \mathfrak{e}_\nabla, \mathfrak{e}_{d*}\}$ described in the next three steps.

<u>Step</u> <u>3</u>: The $\mathfrak{e}_Q$ term is the contribution from the integral of the product of $\chi_\mu^{\,2}$ with $r^2\langle[(\nabla_{A,\beta}a_\ddagger)_\alpha,(\nabla_{A,\beta}a_\ddagger)_\mu]\,[a_{\ddagger\alpha}, a_{\ddagger\beta}]\rangle$. This being the case, its absolute value obeys

$$|\mathfrak{e}_Q| \le c_0\, c^{-1}\,(\int_{B_{r_c}} |\chi_\mu\nabla_A a_\ddagger|^4\,)^{1/2} \tag{3.19}$$

because $r_\ddagger^{\,2}|a_\ddagger\wedge a_\ddagger| = r^2|a\wedge a|$ and the $L^2$ norm of $r^2|a\wedge a|$ on $B_{r_c}$ is less than $c^1$. The term denoted by $\mathfrak{e}_\nabla$ is the contribution from the integral of $2\chi_\mu^{\,2}\langle(\nabla_{A,\beta}a_\ddagger)_\alpha\,[F_{A,\mu\beta},(\nabla_{A,\mu}a_\ddagger)_\alpha]\rangle$ and its norm is also bounded by what is written on the right hand side of (3.19) because the $L^2$ norm of $F_A$ on $B_{r_c}$ is also less than $c^{-1}$. Note in this regard that these $c^{-1}$ bounds for the $L^2$ norms of $r^2|a\wedge a|$ and $F_A$ exist because that $r_c$ is no greater $r_{cF}$ and no greater than $r_{c\wedge}$. The term $\mathfrak{e}_{d*}$ is the contribution from the integral of $\chi_\mu^{\,2}\langle\nabla_{A,\mu}F_{A,\mu\beta}, [a_{\ddagger\alpha},\nabla_{A,\beta}a_{\ddagger\alpha}]\rangle$. The next step says more about this term.

<u>Step</u> <u>4</u>: An integration by parts writes $\mathfrak{e}_{d*}$ as $\frac{1}{2}\,(\mathfrak{e}_F + \mathfrak{e}_\mu - \mathfrak{e}_\nabla)$ with $\mathfrak{e}_F$ and $\mathfrak{e}_\mu$ as follows: What is denoted by $\mathfrak{e}_F$ is $-c_0$ times the integral of $\chi_\mu^{\,2}\,|[F_A, a_\ddagger]|^2$. Use the inequality in the first bullet with $\mu$ replaced by $\mu/4$ to bound the norm of $\mathfrak{e}_F$ by $c_\mu\, c^{-2}$. What is denoted by $\mathfrak{e}_\mu$ has the contribution to the integration by parts identity from the derivative of $\chi_\mu$. In particular, the $\mathfrak{e}_\mu$ term has norm no greater than the sum of $c_0$ times the square of the $L^2$ norm of $[F_A, a_\ddagger]$ and $c_0\mu^{-2}r_c^{-2}$ times that of $\nabla_A a_\ddagger$. It follows as a consequence that $\mathfrak{e}_\mu$ has norm bounded by $c_\mu\,(c^{-2} + N(r_c))$ and thus $|\mathfrak{e}_{d*}| \le c_\mu\,(c^{-2} + N(r_c))$ also. Note that this uses the fact that $r_c^{-2}\int_{B_{r_c}} (|\nabla_A a_\ddagger|^2$ is by definition equal to $N(r_c)$.

<u>Step</u> <u>5</u>: Use the bounds from the preceding steps in (3.18) with the bound in (3.12) and the definition of $N(r_c)$ conclude that

$$\int_{B_{r_c}} (|\nabla_A(\chi_\mu \nabla_A a_\ddagger)|^2 + r_\ddagger^2 \chi_\mu^2 |\nabla_A(a_\ddagger \wedge a_\ddagger)|^2) \le c_0 c^{-1} \Big( \int_{B_{r_c}} |\chi_\mu \nabla_A a|^4 \Big)^{1/2} + c_\mu(c^{-2} + r_c^4) + c_\mu N(r_c).$$

(3.20)

To finish the proof, use the fact that the norm of $d|\chi_\mu \nabla_A a|$ is at most $|\nabla_A(\chi_\mu \nabla_A a)|$ and the Sobolev inequalities from the first and third bullets in (2.1) to bound the square of the $L^4$ norm of $|\chi_\mu \nabla_A a|$ that appears on the right hand side (3.20) by $c_0$ times the left hand side of (3.20). Insert the latter bound to obtain the bound that is claimed by the second bullet of Proposition 3.1 bound when $c > c_0$.

**d) Proof of Proposition 3.2**:

The proof has four steps; the first three prove the proposition's first bullet and the fourth step proves the second bullet. The vector $e$ in the proof that follows of the first bullet is assumed to have norm equal to 1.

Step 1: The $(r_\ddagger, (A, a_\ddagger))$ version of (2.5) implies an equation of the form

$$\nabla_A^\dagger \nabla_A a_\ddagger(e) + r_\ddagger^2 [a_{\ddagger\mu}, [a_\ddagger(e), a_{\ddagger\mu}]] + \mathcal{R}^1_e(\nabla_A a_\ddagger) + \mathcal{R}^0_e(a_\ddagger) = 0$$

(3.21)

with $\mathcal{R}^1_e$ and $\mathcal{R}^0_e$ being endomorphisms with norms that are bound respectively by $c_0 r_c$ and $c_0$. Take the inner product of (3.21) with $a_\ddagger(e)$ to obtain an identity that reads

$$\tfrac{1}{2} d^\dagger d |a_\ddagger(e)|^2 + |\nabla_A a_\ddagger(e)|^2 + r_\ddagger^2 |[a_\ddagger, a_\ddagger(e)]|^2 + \mathfrak{R}^1 + \mathfrak{R}^0 = 0$$

(3.22)

with $\mathfrak{R}^1$ and $\mathfrak{R}^0$ denoting the the inner products between $a_\ddagger(e)$ and $\mathcal{R}^1_e(\nabla_A a_\ddagger)$ and $\mathcal{R}^0_e(a_\ddagger)$.

Step 2: Fix $q \in B_{(1-\mu)r_c}$ and reintroduce $G_q$ to denote the Dirichlet Green's function for the operator $d^\dagger d$ on $B_{r_c}$. Multiply both (3.22) by $\chi_\mu^2 G_q$ and then integrate the result over $B_{r_c}$. Integrate by parts twice to obtain an analog of (3.14) that has the form:

$$\tfrac{1}{2} |a_\ddagger(e)|^2(q) + \int_{B_{r_c}} ((\chi_\mu d^\dagger d\chi_\mu - 2|d\chi_\mu|^2) G_q + 2\mathfrak{m}(\chi_\mu d\chi_\mu, dG_q)) |a_\ddagger(e)|^2 + \mathfrak{r}_a + \mathfrak{r}_m = 0$$

(3.23)

with $\mathfrak{r}_a$ and $\mathfrak{r}_m$ as described in the next two paragraphs.

What is denoted by $\mathfrak{r}_a$ is the integral of the sum of two terms, the first being $\chi_\mu^2 G_q |\nabla_A a_\ddagger(e)|^2$ and the second being $r_\ddagger^2 \chi_\mu^2 G_q |[a_\ddagger(e), a_\ddagger]|^2$. Keeping in mind that $G_q$ is

non-negative and bounded by $c_0 \frac{1}{\mathrm{dist}(\cdot,q)^2}$, the second bullet of Proposition 3.1 with $\mu$ replaced by $\mu/4$ and the second bullet of (2.1) imply the bound $|\mathfrak{r}_a| \le c_\mu (c^{-2} + r_c^4 + N(r_c))$.

What is denoted by $\mathfrak{r}_m$ has the contribution from integral of $\chi_\mu^2 G_q (\mathfrak{R}^1 + \mathfrak{R}^0)$. This being the case, it follows that $|\mathfrak{r}_m|$ is in any case no greater than the integral over $B_{(1-\mu/4)r_c}$ of $c_0 \frac{1}{\mathrm{dist}(\cdot,q)^2} (r_c^2 |\nabla_A a_\ddagger|^2 + |a_\ddagger|^2)$. Granted that this is so, then the second bullet of Proposition 3.1 with $\mu$ replaced by $\mu/4$ and the second bullet of (2.1) can be invoked again to see that $|\mathfrak{r}_m| \le c_\mu r_c^2 (c^{-2} + r_c^4 + N(r_c))$.

Step 3: Let $V_c$ denote the volume of the $B_{(1-\mu/4)r_c}$ and let $I_e$ denote for the moment the integral of $V_c^{-1}|a_\ddagger(e)|$ over this same ball. Since $I_e$ is constant, the equation in (3.23) leads directly to the equation:

$$\tfrac{1}{2} (|a_\ddagger(e)|^2(q) - I_e^2) + \int_{B_{r_c}} ((\chi_\mu d^\dagger d\chi_\mu - 2\,|d\chi_\mu|^2) G_q + 2m(\chi_\mu d\chi_\mu, dG_q))(|a_\ddagger(e)|^2 - I_e^2) + \mathfrak{r}_a + \mathfrak{r}_m = 0. \tag{3.24}$$

Use (3.24) with the $c_0 \frac{1}{\mathrm{dist}(\cdot,q)^2}$ bound for $G_q$ and what Step 2 said about $\mathfrak{r}_a$ and $\mathfrak{r}_m$ to conclude that

$$\big|\,|a_\ddagger(e)|^2(q) - I_e^2\big| \le c_0\, \mu^{-2} r_c^{-4} \int_{B_{(1-\mu/4)r_c}} \big|\,|a_\ddagger(e)|^2 - I_e^2\big| + c_\mu (c^{-2} + r_c^4 + N(r_c)). \tag{3.25}$$

The integral on the right hand side of (3.25) is bounded

$$c_\mu r_c^2 \sup_{B_{(1-\mu/4)r_c}} |a_\ddagger(e)| \Big( \int_{B_{(1-\mu/4)r_c}} \big|\,|a_\ddagger(e)| - I_e\big|^2 \Big)^{1/2}. \tag{3.26}$$

Meanwhile, a standard Dirichlet inequality bounds the $L^2$ norm of $|a_\ddagger(e)| - I_e$ on $B_{(1-\mu/4)r_c}$ by $c_\mu r_c$ times the $L^2$ norm $|\nabla_A a_\ddagger(e)|$ on $B_{(1-\mu/4)r_c}$ and the latter $L^2$ norm is at most $c_0\, r_c\, N(r_c)^{1/2}$. Since the top bullet of Proposition 3.1 bounds $|a_\ddagger(e)|$ by $c_\mu$, this last bound leads to the bound $\sup_{B_{(1-\mu)r_c}} \big|\,|a_\ddagger(e)|^2 - I_e^2\big| \le c_\mu (N(r_c)^{1/2} + c^{-2} + r_c^4 + N(r_c))$. This gives the bound in the top bullet of Proposition 3.1 when $N(r_c) \le c^{-2}$.

Step 4: This step proves the second bullet of Proposition 3.2. To do this, start with (3.15) and but for notational changes, the arguments used in the previous step that lead from (3.24) to (3.27) can be repeated to obtain the bound

$$\sup_{B_{(1-\mu)r_c}} \big||a_\ddagger|^2 - I^2\big| \le c_\mu\,(c^{-2} + r_c^{\,4} + N(r_c)^{1/2} + N(r_c)) \tag{3.27}$$

with $I$ denoting the integral of $V_c^{-1}\,|a_\ddagger|$ over $B_{r_c}$.

To obtain an estimate for $I$, fix $r \in (\frac{1}{2}r_c, (1-\mu)r_c)$ and invoke (3.27) to see that the integral of $r^{-3}|a_\ddagger|^2$ over the boundary of the radius $r$ ball centered at $p$ differs by at most $c_\mu\,(c^{-2}+r_c^{\,4}+N(r_c))$ from that of $r^{-3}I^2$. Since $I$ is constant, the latter integral can be written as $2\pi^2 I^2(1+\mathfrak{e})$ with $\mathfrak{e}$ being a number with absolute value at most $c_0 r^2$.

The integral of $r^{-3}|a_\ddagger|^2$ over $\partial B_r$ can be written as $(1+\mathfrak{e}(r))\,K(r)^2/K(r_c)^2$ with the absolute value of $\mathfrak{e}$ being at most $c_0 r_c^{\,2}$. This is a direct consequence of the definitions of $K$ and $a_\ddagger$. With the preceding in mind, use an integration by parts with the fact that $K$ is increasing to see that the difference between the integral of $r^{-3}|a_\ddagger|^2$ over $\partial B_r$ and its integral over $\partial B_{r_c}$ differ by at most $c_0 N(r_c)^{1/2}$.

The conclusions of the preceding two paragraphs imply that $I^2$ and $\frac{1}{2\pi^2}$ can differ by at most $c_\mu(c^{-2}+r_c^{\,4}+N(r_c)^{1/2}+N(r_c))$. This bound implies Proposition 3.2's second bullet.

## e) Proof of Proposition 3.3

The proof has two parts, the first part proves the top bullet in the proposition and the second part proves the lower bullet.

*Part 1*: The proof of the first bullet in Proposition 3.3 has two steps. These steps assume that $c$, $r_c$ and $N(r_c)$ obey the bounds that are needed to invoke the $\mu/16$ version of Propositions 3.1 and 3.2.

<u>Step 1</u>: The proof starts with the formula given below for the covariant divergence of $F_A$. This formula writes the components of $a_\ddagger$ with respect to a chosen orthonormal frame for $T^*X$ at any given point in X as $\{a_{\ddagger\alpha}\}_{\alpha\in\{1,2,3,4\}}$, these being in the fiber of of $P\times_{SO(3}\mathfrak{G}$ at the given point. The formula writes the components of $\nabla_A a_\ddagger$ using this same orthonormal frame as $\{(\nabla_A a_\ddagger)_\alpha\}_{\alpha\in\{1,2,3,4\}}$. Using this notation, the covariant divergence of $F_A$ can be written at the point in question as

$$*d_A{*}F_A = r_\ddagger^{\,2}\,[a_{\ddagger\alpha}, (\nabla_A a_\ddagger)_\alpha] \tag{3.28}$$

with it understood that repeated indices are summed. The proof of (3.28) starts by taking the covariant divergence of both sides of (2.12). The formula in (3.28) is obtained from the resulting expression after commuting certain derivatives and then invoking the identities in (2.11) and (2.5).

The identity in (3.28) implies a second order differential equation for $F_A$ that can be written as

$$\nabla_A^{\dagger}\nabla_A F_A + r_{\ddagger}^{2}[a_{\ddagger\alpha}, [F_A, a_{\ddagger\alpha}]] + \{F_A, F_A\} = -2\, r_{\ddagger}^{2}\, \nabla_A a_{\ddagger\alpha} \wedge \nabla_A a_{\ddagger\alpha} + r_{\ddagger}^{2}\, \mathcal{R}([a, a]) \tag{3.29}$$

with $\{\cdot, \cdot\}$ denoting a certain homomorphism from $\otimes_2((P\times_{SO(3)} \mathfrak{su}(2)) \otimes (\wedge^2 T^*X))$ to $(P\times_{SO(3)} \mathfrak{su}(2)) \otimes (\wedge^2 T^*X)$. This homomorphism factors as the composition of first the commutator map from $\otimes_2((P\times_{SO(3)} \mathfrak{su}(2)) \otimes (\wedge^2 T^*X))$ to $(P\times_{SO(3)} \mathfrak{su}(2)) \wedge^2 ((\wedge^2 T^*X)$ and then a linear map from the latter space to $(P\times_{SO(3)} \mathfrak{su}(2)) \otimes (\wedge^2 T^*X)$. Though not relevant directly, the homorphism $\{\cdot, \cdot\}$ also respects the writing of $\wedge^2 T^*X$ as the direct sum of the bundles of self dual and anti-self dual forms. What is denoted by $\mathcal{R}$ in (3.29) comes from the Riemann curvature tensor.

<u>Step 2</u>: Take the inner product of both sides of (3.29) with $\chi_\mu^2 F_A$ and integrate the resulting equality. An integration by parts leads directly from the latter to the following:

$$\int_{B_{r_c}} (|\nabla_A(\chi_\mu F_A)|^2 + r_{\ddagger}^{2}\chi_\mu^2 |[F_A, a_{\ddagger}]|^2) \le c_0 \mu^{-2} r_c^{-2} \int_{B_{r_c}} |F_A|^2 + c_0\, c^{-2}$$
$$+\, c_0 \Big(\int_{B_{r_c}} |F_A|^2\Big)^{1/2} \Big(\int_{B_{r_c}} \chi_\mu^4 |F_A|^4\Big)^{1/2} - 2 r_{\ddagger}^{2} \int_{B_{r_c}} \chi_\mu^2 \big\langle *F_A \wedge (\nabla_A a_{\ddagger\alpha} \wedge \nabla_A a_{\ddagger\alpha})\big\rangle . \tag{3.30}$$

The left most term on the right hand side of (3.30) is no greater than $c_0 \mu^{-2} r_c^{-2} c^{-2}$, this because $r_c \le r_{cF}$. For the same reason, the middle term on the right hand side of (3.30) is no greater than

$$c_0\, c^{-1} \Big(\int_{B_{r_c}} \chi_\mu^4 |F_A|^4\Big)^{1/2} ; \tag{3.31}$$

and it follows from the first and third bullets of (2.1) that the latter is at most $c_0 c^{-1}$ times the right hand side of (3.31). The right most integral in (3.30) is $-r_{\ddagger}^{2}$ times the right hand side of (3.19). This being the case, it follows from the first bullet of Proposition 3.1 and the first bullet of (2.1) that the absolute value of the right most integral in (3.30) is at most $c_\mu r_{\ddagger}^{2} (c^{-2} + r_c^4 + \mathrm{N}(r_c))$.

The preceding bounds for the integrals on (3.30)'s right hand side lead directly to the bound asserted by the proposition's first bullet when $c > c_0$.

*Part 2*: The proof of the second bullet in Proposition 3.3 has four steps. These steps assume that $c$ and $r_c$ obey the bounds that are needed to invoke the $\mu/16$ version of Proposition 3.1 and that $c$, $r_c$ and $N(r_c)$ obey the bounds that are needed to invoke the $\mu/16$ versions of Proposition 3.2 and the first bullet in Proposition 3.3.

<u>Step 1</u>: The inequality in (3.17) implies directly the coarser inequality

$$\tfrac{1}{2}\, d^{\dagger}d\,|\nabla_A a_{\ddagger}|^2 \le c_0\,|F_A|\,|\nabla_A a_{\ddagger}|^2 + c_0\, r_{\ddagger}^{\,2}\,|a_{\ddagger}|^2\,|\nabla_A a_{\ddagger}|^2 + c_0\,(|\nabla_A a_{\ddagger}|^2 + |a_{\ddagger}|^2)\ . \tag{3.32}$$

By way of an explanation, the appearance of $|F_A|\,|\nabla_A a_{\ddagger}|^2$ accounts for the term $2\langle(\nabla_{A,\beta} a_{\ddagger})_\alpha\,[F_{A,\mu\beta}, (\nabla_{A,\mu} a_{\ddagger})_\alpha]\rangle$ on (3.17)'s left hand side. The appearance of $r_{\ddagger}^{\,2}|a_{\ddagger}|^2|\nabla_A a_{\ddagger}|^2$ accounts for the term $\langle(\nabla_{A,\mu}F_A)_{\mu\beta}, [a_{\ddagger\alpha}, (\nabla_{A,\beta} a_{\ddagger})_\alpha\,]\rangle$ because $\{(\nabla_{A,\mu}F_A)_{\mu\beta}\}_{\beta\in\{1.2.3.4\}}$ are the components of $*d_A{*}F_A$ when the latter is written using the chosen orthonormal frame for $T^*X$, and because $*d_A{*}F_A$ is given by (3.28).

<u>Step 2</u>: Fix $q \in B_{r_c}$ and introduce again $G_q$ to denote the Dirichlet Green's function with pole at q for the Laplacian $d^\dagger d$ on $B_{r_c}$ . Multiply both sides of (3.32) by $\chi_\mu^{\,2} G_q$, integrate the resulting inequality and then integrate by parts to see that

$$\begin{aligned}\chi_\mu^{\,2}|\nabla_A a_{\ddagger}|^2(q) \le c_0\,\mu^{-2}\, r_c^{-4} \int_{B_{r_c}} |\nabla_A a_{\ddagger}|^2 \\ + c_0 \int_{B_{r_c}} \chi_\mu^2 G_q\, |F_A|\,|\nabla_A a_{\ddagger}|^2 + c_0\,(r_{\ddagger}^{\,2} + 1) \int_{B_{r_c}} \chi_\mu^2 G_q |\nabla_A a_{\ddagger}|^2 + c_0 \int_{B_{r_c}} \chi_\mu^2 G_q\, |a_{\ddagger}|^2\end{aligned} \tag{3.33}$$

The subsequent analysis of (3.33) implicitly uses the bound $|G_q| \le c_0\, \frac{1}{\mathrm{dist}(\cdot,q)^2}$ .

<u>Step 3</u>: Use the second bullet of (2.1) with the definition of N and with (3.12) to bound the right most integral in (3.33) by $c_0\, r_c^{\,2}(N(r_c) + 1)$. The first bullet in Proposition 3.1 with the middle bullet in (2.1) supplies a $c_\mu(c^{-2} + r_c^{\,4} + N(r_c))$ bound for $\int_{B_{r_c}} \chi_\mu^2 G_q |\nabla_A a_{\ddagger}|^2$ .

To see about the integral of $\chi_\mu^{\,2} G_q|F_A|\,|\nabla_A a_{\ddagger}|^2$, first invoke the bound

$$\int_{B_{r_c}} \chi_\mu^2 G_q\, |F_A|\,|\nabla_A a_{\ddagger}|^2 \le c_0\,(\sup_{B_{r_c}} \chi_\mu|\nabla_A a_{\ddagger}|)\,\Big(\int_{B_{r_c}} \chi_\mu G_q\, |F_A|^2\Big)^{1/2}\Big(\int_{B_{r_c}} \chi_\mu G_q |\nabla_A a_{\ddagger}|^2\Big)^{1/2} . \tag{3.34}$$

With the preceding understood, invoke the second bullet of (2.1) with the $\mu/4$ versions of the first bullets in Propositions 3.1 and 3.3 to bound the right hand side of (3.34) by

$$c_\mu \left( \sup_{B_{r_c}} \chi_\mu |\nabla_A a_\ddagger| \right) \left( r_c^{-1} + r_\ddagger (c^{-2} + r_c^2 + N(r_c))^{1/2} \right) . \tag{3.35}$$

The latter is at most the sum of two terms, the first being $\frac{1}{100} \sup_{B_{r_c}} \chi_\mu^2 |\nabla_A a_\ddagger|^2$ and the second being $c_\mu r_c^{-2}((1 + r_c^2 r_\ddagger^2)(c^{-2} + r_c^2 + N(r_c))$.

Step 4: The inequality in (3.33) with the bounds in Step 3 lead directly to the assertion of Proposition 3.3's second bullet.

## 4. Unexpectedly small curvature

The geometric data in this section is the same as in Section 3, this being the oriented Riemannian 4-manifold X, the open set $U \subset X$ with compact closure and the principal SO(3) bundle $P \to X$. Fix $r > 1$ and let $(A, a)$ denote a pair of connection on P and section of $T^*X \otimes (P \times_{SO(3)} \mathfrak{S})$ that obeys some $\tau \in [0, 1]$ version of (2.11). Fix a large positive number $c$ and a point $p \in U$ so as to define the numbers $r_{cF}$ and $r_{c\wedge}$ as done in (3.1) and (3.2). The upcoming Proposition 4.1 asserts that the integral of $r^2 a \wedge a$ is smaller than might be expected on balls of radius less than a given $r_c \in (0, \min(r_{cF}, r_{c\wedge}))$ when $r_c$ and the value of (3.5)'s function N at $r_c$ are small.

**Proposition 4.1**: *Fix $\mu \in (0, \frac{1}{2}]$ and there exists $\kappa_\mu > 1$ that depends only on $\mu$ and the geometric data, and is such that if $c > \kappa_\mu$, then the subsequent assertion is true. Given $\varepsilon \in (0, 1]$, there exists $\kappa_{*\varepsilon} > 1$ that depends on $\mu$ and $\varepsilon$ and whose significance is as follows: Fix $r > 1$ and assume that $(A, a)$ is a pair of connection on P and section of $T^*X \otimes (P \times_{SO(3)} \mathfrak{S})$ that obey a given $\tau \in [0,1]$ version of (2.11). Fix $p \in U$ to define $r_{cF}$ and $r_{c\wedge}$ using $(A, a)$, and suppose that $r_c$ is positive but less than the minimum of $r_{c\wedge}$, $r_{cF}$ and $c^{-1}$. Suppose in addition that $N(r_c)$ is less than the minimum of $\kappa_{*\varepsilon}^{-1}$ and $c^{-2}$. Then*

$$r^4 \int_{B_{(1-\mu) r_c}} |a \wedge a|^2 < \varepsilon .$$

Section 4a proves this proposition when $r_c \kappa(r_c) r$ has a specified a priori upper bound; and Section 4d proves the proposition when $r_c \kappa(r_c) r$ is larger than a number that depends on $\mu$, $\varepsilon$ and the geometric data. The intervening sections supply input for the arguments in Section 4d. The notation here as in the previous section uses $r_\ddagger$ to denote $\kappa(r_c) r$ and it uses $a_\ddagger$ to denote $\kappa(r_c)^{-1} a$. The notation also uses $c_\mu$ to denote a number that is greater than 1 and that depends only on the geometric data and the chosen value of $\mu$. As in Section 3, the value of $c_\mu$ can be assumed to increase between its appearances.

**a) The proof of Proposition 4.1 when $r_c r_\ddagger$ is small**

Fix for the moment $z > 1$. What follows directly proves that the integral over $B_{r_c}$ of $r^4|a \wedge a|^4$ is less that ε when $r_c r_\ddagger$ is bounded by $z$ and when $\mathrm{N}(r_c)$ is smaller than a number that depends $z$ and ε: Use the definition of $\mathrm{N}$ to bound the integral of $r^4|a \wedge a|^2$ over $B_{r_c}$ by $c_0 r_c^2 \mathrm{K}(r_c)^2 r^2 \mathrm{N}(r_c)$. This in turn is less than ε when $\mathrm{N}(r_c)$ is less than $c_0^{-1} z^{-2} \varepsilon$.

**b) Decomposing $a_\ddagger$ when $r_c\, r_\ddagger$ is large**

A convention used in what follows assumes that the constant $c$ is large enough so as to invoke Proposition 3.2 with $\mu$ replaced by $\mu/16$. In particular, $c$ is chosen so that what follows is true.

- $\big||a_\ddagger| - \frac{1}{\sqrt{2}\,\pi}\big| < \frac{1}{1000}$ *on* $B_{(1-\mu/16)r_c}$ .
- *If* $e \in TX|_p$, *then* $\max_{B_{(1-\mu/16)\, r_c}} |a_\ddagger(e)| - \min_{B_{(1-\mu/16)\, r_c}} |a_\ddagger(e)| < \frac{1}{1000}|e|$ .

(4.1)

By way of notation, the convention in (4.1) as in Proposition 3.2 is to use $e$ to denote both a vector in $TX|_p$ and a vector field on $B_{1/c_0}$ with the latter incarnation obtained from the former by parallel transport along the radial geodesics from p by the Levi-Civita connection. This convention is used also in the subsequent parts of the proof.

A positive definite, Hermitian endomorphism of $P \times_{SO(3)} \mathfrak{S}$ over $B_{(1-\mu/16)r_c}$ is defined directly using $a_\ddagger$ and the Riemannian metric. The definition uses a chosen orthonormal frame for T*X on $B_{(1-\mu/16)r_c}$ and writes the components of $a_\ddagger$ with respect to this frame as $\{a_{\ddagger\alpha}\}_{\alpha \in \{1,2,3,4\}}$. These are sections over $B_{(1-\mu/16)r_c}$ of the bundle $P \times_{SO(3)} \mathfrak{S}$. With this notation understood, define $\mathbb{T}_\ddagger$ by the rule

$$\sigma \to \mathbb{T}_\ddagger(\sigma) = a_{\ddagger\alpha} \langle a_{\ddagger\alpha} \sigma \rangle \,,$$

(4.2)

with it understood that there is an implicit summation over the repeated indices. Note that $\mathbb{T}$ does not depend on the chosen orthonormal frame.

The following two assertions are respective consequences of the first and second bullets of (4.1) when $c > c_\mu$:

- *The endomorphism* $\mathbb{T}_\ddagger$ *has at least one eigenvalue that is greater than* $\frac{1}{80}$ *at each point of* $\mathrm{B}_{(1\text{-}\mu/16)\mathrm{r}_c}$.
- *Suppose that* $\mathbb{T}_\ddagger$ *has only one eigenvalue greater than the minimum of* $c^{-2}$ *and* $\frac{1}{750}$ *at some point in* $\mathrm{B}_{(1\text{-}\mu/16)\mathrm{r}_c}$ *and that the corresponding eigenspace has multiplicity 1. Then there is only one eigenvector with eigenvalue greater than the minimum of* $\mathrm{c}_\mu c^{-1}$ *and* $\frac{1}{500}$ *at all points of* $\mathrm{B}_{(1\text{-}\mu/16)\mathrm{r}_c}$ *and the eigenspace has multiplicity one also.*

(4.3)

Conditions that guarantee the existence of only one $\mathcal{O}(1)$ eigenvalue of $\mathbb{T}_\ddagger$ at all points in $\mathrm{B}_{(1\text{-}\mu/16)\mathrm{r}_c}$ are given in the upcoming Lemma 4.2. In any event, suppose that it is the case that $\mathbb{T}_\ddagger$ has just one eigenvalue that is greater than $\frac{1}{500}$ at all points of $\mathrm{B}_{(1\text{-}\mu/16)\mathrm{r}_c}$ and that the corresponding eigenspace is 1-dimensional. There is in this case section of $\mathrm{P}\times_{\mathrm{SO}(3)}\mathfrak{su}(2)$ over $\mathrm{B}_{(1\text{-}\mu/16)\mathrm{r}_c}$ whose restriction to any given point is an eigenvector of $\mathbb{T}_+$ with the largest eigenvalue. This eigensection is denoted in what follows by $\sigma_\ddagger$. Use $\lambda_\ddagger$ to denote the corresponding eigenvalue. It follows from Proposition 3.2 that $\lambda_\ddagger$ differs from $\frac{1}{2\pi}$ by at most $\mathrm{c}_\mu c^{-1}$, this difference being less than $c^{-1/2}$ when $c \ge \mathrm{c}_\mu$.

The $\mathrm{P}\times_{\mathrm{SO}(3)}\mathfrak{su}(2)$ valued 1-form $a_\ddagger$ is written over $\mathrm{B}_{(1\text{-}\mu/16)\mathrm{r}_c}$ using $\sigma_\ddagger$ as

$$a_\ddagger = \nu\sigma_\ddagger + \mathfrak{a}$$

(4.4)

with $\nu$ being an $\mathbb{R}$-valued 1-form on $\mathrm{B}_{(1\text{-}\mu/16)\mathrm{r}_c}$ and with $\mathfrak{a}$ being a section of the bundle $\mathrm{T}^*X\otimes(\mathrm{P}\times_{\mathrm{SO}(3)}\mathfrak{S})$ over $\mathrm{B}_{(1\text{-}\mu/16)\mathrm{r}_c}$ that obeys $\langle\sigma_\ddagger\,\mathfrak{a}\rangle = 0$ at each point. Moreover,

- $\lambda_\ddagger = |\nu|^2$.
- $|\mathfrak{a}| \le \mathrm{c}_\mu\, c^{-1/2}$.
- *The metric pairing between* $\nu$ *and* $\mathfrak{a}$ *is zero.*
- $|a_\ddagger \wedge a_\ddagger|^2 = 4|\nu|^2|\mathfrak{a}|^2 + |\mathfrak{a}\wedge\mathfrak{a}|^2$.
- *Let* $\nu_\mathrm{p}$ *denote the section of* $\mathrm{T}^*\mathrm{B}_{(1\text{-}\mu/16)\mathrm{r}_c}$ *that is defined by the parallel transport of* $\nu|_\mathrm{p}$ *along the radial geodesics from* p. *Then* $|\nu - \nu_\mathrm{p}| \le \mathrm{c}_\mu\, c^{-1}$ *on* $\mathrm{B}_{(1\text{-}\mu/16)\mathrm{r}_c}$.

(4.5)

The first three bullets of (4.5) are proved by writing $\mathbb{T}_\ddagger$ using (4.4) as the endomorphism

$$\sigma \to |\nu|^2\,\sigma_\ddagger\langle\sigma_\ddagger\sigma\rangle + \sigma_\ddagger\nu_\alpha\langle\mathfrak{a}_\alpha\sigma\rangle + \mathfrak{a}_\alpha\nu_\alpha\langle\sigma_\ddagger\sigma\rangle + \mathfrak{a}_\alpha\langle\mathfrak{a}_\alpha\sigma\rangle\,.$$

(4.6)

The first and third bullets of (4.5) follow from the $\sigma = \sigma_\ddagger$ version of (4.6) and the second bullet follows from the second bullet of (4.3) by taking $\sigma$ in (4.6) to be pointwise orthogonal to $\sigma_{\ddagger}$. The fourth bullet in (4.5) can be seen by first writing $a_{\ddagger} \wedge a_{\ddagger}$ using (4.6) and then invoking the third bullet of (4.5) and the fact that $\langle \mathfrak{a}\, \sigma_{\ddagger} \rangle = 0$. Note in this regard that these two fact imply that $|\nu \wedge [\sigma_{\ddagger}, \mathfrak{a}]| = 2\,|\nu|\,|\mathfrak{a}|$ and the fact that $\langle \mathfrak{a}\, \sigma_{\ddagger} \rangle = 0$ implies by itself that $\mathfrak{a} \wedge \mathfrak{a}$ is everywhere the tensor product of an $\mathbb{R}$-valued 2-form and $\sigma_{\ddagger}$. The fifth bullet follows from the top bullet of Proposition 3.2.

The derivatives of $\nu$, $\sigma_{\ddagger}$ and $\mathfrak{a}$ are such that the following are true when $c > c_{\mu}$:

- $(\nabla_{\beta}\nu)_{\alpha} = \langle \sigma_{\ddagger}(\nabla_{A,\beta} a_{\ddagger})_{\alpha} \rangle + \langle \nabla_{A,\beta}\sigma_{\ddagger}\, \mathfrak{a}_{\alpha} \rangle$ .
- $|\nu|^2 [\sigma_{\ddagger}, \nabla_{A,\beta}\sigma_{\ddagger}] - \langle \nabla_{A,\beta}\sigma_{\ddagger}\, \mathfrak{a}_{\alpha} \rangle [\sigma_{\ddagger}, \mathfrak{a}_{\alpha}] = [\sigma_{\ddagger}, \nu_{\alpha}(\nabla_{A,\beta} a_{\ddagger})_{\alpha}] + \langle \sigma_{\ddagger}\, (\nabla_{A,\beta} a_{\ddagger})_{\alpha} \rangle [\sigma_{\ddagger}, a_{\ddagger\alpha}]$ .
- $(\nabla_{A,\beta}\, \mathfrak{a})_{\alpha} = -\frac{1}{4}\,([\nabla_{A,\beta}\sigma_{\ddagger}, [a_{\ddagger\alpha}, \sigma_{\ddagger}]] + [\sigma_{\ddagger}, [a_{\ddagger\alpha}, \nabla_{A,\beta}\sigma_{\ddagger}]] + [\sigma_{\ddagger}, [(\nabla_{A,\beta} a_{\ddagger})_{\alpha}, \sigma_{\ddagger}]]$.
- $|\nabla\nu| + |\nabla_A\sigma_{\ddagger}| + |\nabla_A\mathfrak{a}| \le c_0\, |\nabla_A a|$ .

(4.7)

The first bullet is proved by differentiating both sides of (4.4) and then using the fact that $\langle \sigma_{\ddagger}\, \nabla_A\sigma_{\ddagger} \rangle = 0$ and $\langle \sigma_{\ddagger}\, \mathfrak{a} \rangle$ both vanish. The second bullet's proof starts with the identity $[\sigma_{\ddagger}, \nu_{\alpha}\nabla_A a_{\alpha}] = |\nu|^2 [\sigma_{\ddagger}, \nabla_A\sigma_{\ddagger}] - \nabla\nu_{\alpha}\, [\sigma_{\ddagger}, \mathfrak{a}]$ which is a consequence the vanishing of $\nu_{\alpha}\, \mathfrak{a}_{\alpha}$. This rewriting of $[\sigma_{\ddagger}, \nu_{\alpha}\nabla_A a_{\alpha}]$ and (4.7)'s first bullet lead directly to the identity in the second bullet of (4.7). The third bullet of (4.7) is obtained by differentiating the identity $\mathfrak{a} = -\frac{1}{4}\,[\sigma_{\ddagger}, [a, \sigma_{\ddagger}]]$. The fourth bullet follows from the first three when $c > c_0$. To prove that this is so, invoke the first and third bullets of (4.7) to bound $|\nabla\nu|$ and $|\nabla_A\mathfrak{a}|$ by $c_0\,(|\nabla_A a_{\ddagger}| + |\nabla_A\sigma_{\ddagger}|)$. These bounds imply the fourth bullet if $|\nabla_A\sigma_{\ddagger}| \le c_0\, |\nabla_A a_{\ddagger}|$. Such a bound follows from the second bullet of (4.7) when $c$ is large because the latter bullet with no assumption on $c$ implies a bound of the form $(1 - c_0|\mathfrak{a}|^2)\, |\nabla_A\sigma_{\ddagger}| \le c_0\, |\nabla_A a_{\ddagger}|$ .

The final item in this subsection is the promised Lemma 4.2 which asserts conditions that guarantee the existence of but a single $\mathcal{O}(1)$ eigenvalue for $\mathbb{T}_{\ddagger}$.

**Lemma 4.2**: *There exists* $\kappa > 1$ *that depends only on the geometric data and given* $\mu \in (0, \frac{1}{4}\,]$*, there exists* $\kappa_{\mu} > \kappa$ *that depends only on* $\mu$ *and the geometric data; and these numbers having the following significance: Fix* $r > 1$ *and a pair* $(A, a)$ *that obeys* (2.11) *for some* $\tau \in [0,1]$. *Fix* $c > \kappa$ *and given* $p \in U$*, define* $r_{c\wedge}$ *and* $r_{cF}$ *using* $p$ *and* $(A, a)$. *Suppose that* $r_c$ *is less than* $r_{c\wedge}$, $r_{cF}$ *and* $c^{-1}$. *If* $N(r_c) < c^{-2}$ *and if* $r_c\, r_{\ddagger} \ge \kappa_{\mu}$, *then* $\mathbb{T}_{\ddagger}$ *has only one eigenvalue greater than* $\kappa_{\mu} c^{-1}$ *at each point of* $B_{(1-\mu/16)\, r_c}$ *and the corresponding eigenspace is 1-dimensional.*

***Proof of Lemma 4.2***: Assume that $c$ and $r_c$ are such that the Proposition 3.2 can be invoked. If $z > 1$ and if $r_c\, r_{\ddagger} \ge z^{1/2}$, then (3.1) requires

$$r_c^{-4} \int_{B_{r_c}} |a_\ddagger \wedge a_\ddagger|^2 \le z^{-2}\, c^{-2} ,$$

(4.8)

and so there exists a point in $B_{(1-\mu/16)r_c}$ where $|a_\ddagger \wedge a_\ddagger| \le c_0\, z^{-1}\, c^{-1}$. Let p´ denote such a point. The bound $|a_\ddagger \wedge a_\ddagger| \le c_0\, z^{-1}\, c^{-1}$ requires that $\mathbb{T}_\ddagger$ at p´ to have but one eigenvalue greater than $c_0\, z^{-2}\, c^{-2}$ because $\mathfrak{S}$ is a rank 1 Lie algebra. Granted this observation, then the second bullet of Proposition 3.2 implies the lemma's assertion when $r_c\, r_\ddagger > c_0$.

**c) Reducible connections**

Fix p ∈ U and a pair $(A, a)$ of connection on P and section of $T^*X \otimes (P \times_{SO(3)} \mathfrak{S})$ that obeys some $\tau \in [0, 1]$ version of (2.11). Fix $\mu \in (0, \frac{1}{4}]$ and supposing that $c > c_\mu$, then Propositions 3.1 and 3.2 can be invoked given a point p ∈ U and a positive number $r_c$ no greater than the minimum of $r_{c\wedge}$, $r_{cF}$ and $c^{-1}$ with the property that $N(r_c) \le c^{-2}$. Take $c$ to be 1000, and assume in addition that the corresponding endomorphism $\mathbb{T}_\ddagger$ has only one eigenvalue at each point in $B_{(1-\mu/16)r_c}$ that is greater than $c^{-1}$.

A connection to be denoted by $\hat{A}$ is defined on P over $B_{(1-\mu/16)r_c}$ by the formula

$$\hat{A} = A - \tfrac{1}{4}\,[\sigma_\ddagger, \nabla_A \sigma_\ddagger] .$$

(4.9)

This connection $\hat{A}$ is designed so that $\nabla_{\hat{A}}\sigma_\ddagger = 0$. The curvature of $\hat{A}$ is

$$F_{\hat{A}} = \sigma_\ddagger \langle \sigma_\ddagger F_A \rangle - \tfrac{1}{4}\, \nabla_A \sigma_\ddagger \wedge \nabla_A \sigma_\ddagger .$$

(4.10)

The $L^2$ norm of $F_{\hat{A}}$ on $B_{(1-\mu/16)r_c}$ is such that

$$\int_{B_{(1-\mu/16)r_c}} |F_{\hat{A}}|^2 \le c_\mu\, c^{-2} .$$

(4.11)

This follows directly from the formula in (4.10) because the $L^2$ norm of $F_A$ on $B_{r_c}$ is bounded by $c^{-1}$ as is the $L^4$ norm of $|\nabla_A \sigma_\ddagger|$. To prove this $L^4$ norm bound, use the top bullet of (4.7) to bound $|\nabla_A \sigma_\ddagger|$ by $c_0 |\nabla_A a_\ddagger|$ and then use the top bullet of (2.1) to bound the $L^4$ norm of $|\nabla_A a_\ddagger|$ on $B_{(1-\mu/16)r_c}$ by the sum of the $L^2$ norm on $B_{(1-\mu/16)r_c}$ of $|\nabla_A(\nabla_A a_\ddagger)|$ and $r_c^{-1}$ times the $L^2$ norm of $|\nabla_A a_\ddagger|$ on $B_{(1-\mu/16)r_c}$. The former norm is bounded courtesy of Proposition 3.1 by $c_\mu\, c^{-1}$ if both $r_c \le c^{-1}$ and $N(r_c) \le c^{-2}$. Meanwhile, the product of $r_c^{-1}$ times the latter norm is no greater than $c_0\, N(r_c)^{1/2}$.

As explained next, the bound in (4.11) implies that $\hat{A}$ is close in an $L^2_1$ sense to a flat connection on P. In particular, there is a principal SO(3) bundle isomophism over $B_{(1-\mu/16)r_c}$ from the product principal SO(3) to P that pulls $\hat{A}$ back as $\theta_0 + a_\diamond \sigma_\diamond$ with $\theta_0$ denoting the product connection, $\sigma_\diamond$ denoting the corresponding pull-back of $\sigma_\ddagger$ and $a_\diamond$ being an $\mathbb{R}$-valued 1-form. Moreover, $\sigma_\diamond$ is $\theta_0$-covariantly constant and $a_\diamond$ is coclosed and such that its norm and that of its derivative obey

$$\int_{B_{(1-\mu/16)r_c}} (|\nabla a_\diamond|^2 + r_c^{-2}\, |a_\diamond|^2) \le c_\mu\, c^{-2}\,. \tag{4.12}$$

This isomorphism is denoted by $h_\diamond$ and its construction is described in the next paragraph.

Having fixed an isomorphism between $P|_p$ and SO(3), use parallel transport by $\hat{A}$ outward on the geodesic arcs from p to define an isomorphism over $B_{(1-\mu/16)r_c}$ between product principal SO(3) bundle and P. This isomorphism maps $\sigma_\ddagger$ to a $\theta_0$-covariantly constant section of the product $\mathfrak{S}$ bundle, the latter being $\sigma_\diamond$. Meanwhile, this same isomorphism pulls $\hat{A}$ back as $\theta_0 + a_{\hat{A}}\tau$ with $a_{\hat{A}}$ being a real valued 1-form whose exterior derivative is the 2-form $\langle \sigma_\ddagger F_{\hat{A}} \rangle$.

Use the Neumann Green's function for Laplacian on $B_{(1-\mu/16)r_c}$ to construct a function on $B_{(1-\mu/16)r_c}$ to be denoted by $f$ that solves the equation $d{*}df = d{*}a_{\hat{A}}$ with radial derivative on $\partial B_{(1-\mu/16)r_c}$ equal to the radial component of $a_{\hat{A}}$ on $\partial B_{(1-\mu/16)r_c}$. Set $a_\diamond = a_{\hat{A}} - df$. This 1-form obeys (4.12) and it is the pull-back of $\theta_0 + a_{\hat{A}}\sigma_\diamond$ by the automorphism $e^{-f\sigma_\diamond}$ of the product principal SO(3) bundle.

**d) The proof of Proposition 4.1 when $r_c r_\ddagger$ is large**

The upcoming Lemma 4.3 is used to prove what is asserted by Proposition 4.1 when $r_c r_\ddagger$ is large. This lemma is invoked a second time in Section 6.

**Lemma 4.3**: *Given $\mu \in (0, \frac{1}{4}]$, there exists $\kappa_\mu > 1$ depending only on $\mu$ and the geometric data, and with the following significance: Assume that $c$ is greater than $\kappa_\mu$ and that $c$,, $r_c$ and $N(r_c)$ obey the bounds that are required to invoke the versions of Proposition 3.1-3.3 that has $\mu$ replaced by $\mu/16$. Assume in addition that (4.3) is true so as to define $\sigma_\ddagger$ and $\mathfrak{a}$ as described in (4.4). If, in addition, $r_c r_\ddagger > \kappa_\mu$, then $\sigma_\ddagger$ and $\mathfrak{a}$ are such that*

$$\int_{B_{(1-\mu/4)r_c}} (r_\ddagger^2|\nabla_A \mathfrak{a}|^2 + |\nabla_A(\nabla_A\sigma_\ddagger)|^2 + r_\ddagger^2(r_\ddagger^2|\mathfrak{a}|^2 + |\nabla_A\sigma_\ddagger|^2)) \le \kappa_{*\mu}\,(r_c r_\ddagger)^{-2}\, c^{-1}\,.$$

Lemma 4.3 is proved momentarily. Accept it as true for the time being.

***Proof of Proposition 4.1***: If $r_c r_\ddagger$ is greater than $c_\mu$, then Lemma 4.2 can be used to guarantee that what is said in (4.3) holds. This understood, use (4.4) to write $a_\ddagger = \nu\sigma_\ddagger + \mathfrak{a}$. It follows from the fourth bullet of (4.5) that the integral of $r^4|a\wedge a|^2$ over $B_{(1-\mu)r_c}$ is no greater than $c_\mu$ times the integral over $B_{(1-\mu)r_c}$ of $r_\ddagger^{\,4}|\mathfrak{a}|^2$. If $r_c r_\ddagger$ is larger than yet another version of $c_\mu$, then Lemma 4.3 bounds the latter integral by $c_\mu (r_c r_\ddagger)^{-2} c^{-1}$. This will be less than $\varepsilon$ when $r_c r_\ddagger$ is greater than $c_\mu c^{-1/2}\varepsilon^{-1/2}$. By way of a reminder, if the number $r_c r_\ddagger$ is less than $c_\mu c^{-1/2}\varepsilon^{-1/2}$ and if $N(r_c) \le c_\mu^{-1}\varepsilon^2$, then what is said in Section 4a can be invoked to bound the integral $r^4|a\wedge a|^4$ over $B_{(1-\mu)r_c}$ by $\varepsilon$

The rest of this subsection has five parts, the proof of Lemma 4.3 being in Part 5. Parts 1-4 set the stage for the proof of Lemma 4.3 in Part 5. These five parts of the subsection assume that $c$, $r_c$ and $N(r_c)$ obey the bounds that are needed to invoke the versions of Proposition 3.1-3.3 that has $\mu$ replaced by $\mu/16$. The values of $c$, $r_c$ and $N(r_c)$ are taken so that Lemma 4.2 can be used when $r_c r_\ddagger \ge c_\mu$ so as to invoke (4.3), write $a_\ddagger$ as in (4.4) and having done so, invoke what is written in (4.5) and (4.7).

*Part 1*: The discussions in the subsequent parts of this subsection invoke a convention with regards to metrics on the $|x| < c_0^{-1}$ ball in $\mathbb{R}^4$. The convention is to use the Euclidean metric to define the notions of self duality and anti-self duality and to define the Hodge star operator. For example, if $\omega$ is any given 2-form on this ball, then $\omega^+$ and $\omega^-$ refer to the respective parts of $\omega$ that are self dual and anti-self dual with the latter notions defined by the Euclidean metric. By way of a second relevant example, let $\mathcal{D}$ denote for the moment a given differential operator mapping sections of one tensor bundle over $\mathbb{R}^4$ to sections of another. The symbol $\mathcal{D}^\dagger$ is used below to denote the formal, Euclidean $L^2$ adjoint of $\mathcal{D}$, this defined by using Euclidean inner products and the Euclidean metric's volume form. A third example concerns covariant derivatives on tensor bundles. Unless stated to the contrary, these are defined using the Euclidean metric's Levi-Civita connection. The symbol $\nabla$ is used to denote these derivatives.

This convention with regards to metrics is needed because a Gaussian coordinate chart for the radius $c_0^{-1}$ ball centered at p is used momentarily to write the equations in (2.11) as equations on the $|x| < c_0^{-1}$ ball in $\mathbb{R}^4$. Such a coordinate chart pulls back the metric from X so as to define a metric on the $|x| < c_0^{-1}$ ball. The latter metric is denoted by $\mathfrak{m}$. This metric can differ from the Euclidean metric at any given point by at most $c_0|x|^2$. Meanwhile, the norms of its derivatives to first order are bounded at x by $c_0|x|$ and the norms of its second derivatives are bounded by $c_0$.

A second convention used below concerns the product connection on the product principal SO(3) bundle $\mathbb{R}^4 \times SO(3)$. This connection is denoted by $\theta_0$. Covariant derivatives on associated bundles that are defined using this connection are denoted by $\nabla$; and the corresponding exterior derivative for differential forms with values in an associated vector bundle is denoted by d.

*Part 2*: This part of the subsection introduces a certain model system of linear equations on $\mathbb{R}^4$. The definition of this model system requires three choices, these being a unit norm element in $\mathfrak{S}$ to be denote by $\sigma_\diamond$, a constant unit length section of $T^*\mathbb{R}^4$ to be denoted by $e$, and a number greater than 1 to be denoted by $m$.

Use $\mathfrak{h}$ to denote the orthogonal complement of $\sigma_\diamond$ in $\mathfrak{S}$, this being the kernel of the linear function $\langle\sigma_\diamond(\cdot)\rangle$. Use $\mathbb{V}$ to denote $T^*\mathbb{R}^4 \otimes \mathfrak{h}$ and use $\mathbb{W}^+$ and $\mathbb{W}^-$ to denote the respective vector spaces of $\mathfrak{h}$ valued, self dual and anti-self dual 2-forms on $\mathbb{R}^4$. Use $\mathbb{W}$ to denote the vector space $(\mathbb{W}^+ \oplus \mathfrak{h}) \oplus (\mathbb{W}^- \oplus \mathfrak{h})$ with $\mathfrak{h}$ used here to denote the product bundle with the kernel of $\langle\sigma_\diamond(\cdot)\rangle$ being the fiber.

The model operator maps $C^\infty(\mathbb{R}^4; \mathbb{V} \oplus \mathbb{V})$ to $C^\infty(\mathbb{R}^4; \mathbb{W})$; it is denoted by $\mathcal{L}$. This operator $\mathcal{L}$ is defined by the rule whereby a given pair $x = (\mathfrak{p}, \mathfrak{q}) \in C^\infty(\mathbb{R}^4; \mathbb{V} \oplus \mathbb{V})$ is sent to the element $\mathcal{L}x \in C^\infty(\mathbb{R}^4; \mathbb{W})$ whose respective $(\mathbb{W}^+ \oplus \mathfrak{h})$ and $(\mathbb{W}^- \oplus \mathfrak{h})$ components are the top and bottom bullets below.

- $\big((d\mathfrak{q} - m e\wedge[\sigma_\diamond, \mathfrak{p}])^+, *(d{*\mathfrak{q}} + m e\wedge[\sigma_\diamond, *\mathfrak{p}])\big)$ ,
- $\big((d\mathfrak{p} + m e\wedge[\sigma_\diamond, \mathfrak{q}])^-, *(d{*\mathfrak{p}} - m e\wedge[\sigma_\diamond, *\mathfrak{q}])\big)$ .

(4.13)

A look at the symbol of $\mathcal{L}$ proves $\mathcal{L}$ to be elliptic; and a computation finds

$$\mathcal{L}^\dagger\mathcal{L}x = \nabla^\dagger\nabla x + 4m^2 x \,. \tag{4.14}$$

The formula in (4.14) has the following integral version: If $x$ is an $L^2_1$ section of the bundle $\mathbb{V} \oplus \mathbb{V}$, then

$$\int_{\mathbb{R}^4} |\mathcal{L}x|^2 = \int_{\mathbb{R}^4} (\,|\nabla x|^2 + 4m^2|x|^2) \tag{4.15}$$

This holds for $L^2_1$ sections of $\mathbb{V} \oplus \mathbb{V}$ if it holds for smooth sections with compact support, and it is proved to hold for the latter sort by first taking the inner product of both sides of the identity in (4.14) with $x$, integrating the result over $\mathbb{R}^4$ and then integrating by parts.

*Part 3*: The operator $\mathcal{L}$ is used momentarily to prove Lemma 4.3. This part of the subsection sets the stage for this application. The stage setting starts with the choice of a Gaussian coordinate chart centered at p. This chart identifies the radius $c_0^{-1}$ ball in X centered at p with the same radius ball in $\mathbb{R}^4$ centered on the origin. As such, it writes tensors on the radius $c_0^{-1}$ ball centered at p as tensors on the same radius ball about the origin in $\mathbb{R}^4$. Subsequent discussions for the most part invoke implicitly this identification of the two balls and the associated identification of corresponding tensors on these balls.

The next part of the stage setting uses (4.4) to define $\sigma_‡$ and uses $\sigma_‡$ to define the connection $\hat{A}$ as done in (4.9). By way of a reminder, $\hat{A}$ is a connection over $B_{(1-\mu/16)r_c}$ of the bundle P with the property that $\nabla_{\hat{A}}\sigma_‡ = 0$. Use the isomorphism $h_\diamond$ that is described at the end of Section 4c to identify over $B_{(1-\mu/16)r_c}$ with the product principal SO(3) bundle. This isomorphism identifies $\sigma_‡$ with a unit length, constant element in $\mathfrak{S}$. The latter is denoted again denoted by $\sigma_\diamond$ because it plays the role in what follows that is played by its namesake in (4.13).

With $\sigma_\diamond$ as just described, the final part of the stage setting introduces the versions $e$ and $m$ that are needed to define $\mathcal{L}$ for use in the proof of Proposition 3.1. To define the constant 1-form $e$ on $\mathbb{R}^4$, use the Gaussian coordinate chart map to view the 1-form $\nu$ in (4.4) as a 1-form on the radius $(1 - \frac{1}{16}\mu)r_c$ ball in $\mathbb{R}^4$ centered on the origin. The desired version of $e$ is the constant 1-form on $\mathbb{R}^4$ that equals $|\nu|^{-1}\nu$ at the origin in $\mathbb{R}^4$. Take the number $m$ to be $\frac{1}{\sqrt{2}\,\pi}\, r_‡$.

*Part 4*: The four steps that follow use $\mathcal{L}$ to write the relevant $\tau \in [0, 1]$ version of the equations in (2.11) on the ball $B_{(1-\mu/4)r_c}$ as equations on the $|x| \le (1 - \frac{1}{4}\mu)r_c$ ball in $\mathbb{R}^4$.

<u>Step 1</u>: Use $\mathfrak{a}$ henceforth to denote the pull-back by the isomorphism $h_\diamond$ of its namesake in (4.4). This new incarnation of $\mathfrak{a}$ is viewed using the Gaussian coordinate chart map as an $\mathfrak{h}$-valued 1-form on the $|x| \le (1 - \frac{1}{16}\mu)r_c$ ball in $\mathbb{R}^4$.

The pull-back of $\hat{A}$ by $h_\diamond$ was written in Section 4c as $\theta_0 + a_\diamond\sigma_\diamond$ with $a_\diamond$ obeying (4.12). The notation used henceforth does not distinguish $a_\diamond$ from its pull-back via the Gaussian coordinate chart map to the $|x| \le (1 - \frac{1}{16}\mu)r_c$ ball in $\mathbb{R}^4$.

The pull-back of the connection A first by the isomorphism $h_\diamond$ and then by the Gaussian coordinate chart map is a connection on the product principal SO(3) bundle over the $|x| \le (1 - \frac{1}{16}\mu)r_c$ ball in $\mathbb{R}^4$. This pull-back of A is still denoted by A. It can be written as $\theta_0 + a_\diamond\sigma_\diamond + m\,\mathfrak{b}$ with $\mathfrak{b} = \frac{1}{4}\, m^{-1}[\sigma_\diamond, \nabla_A\sigma_\diamond]$ being another $\mathfrak{h}$ valued 1-form on the $|x| \le (1 - \frac{1}{16}\mu)r_c$ ball.

The pair $(\mathfrak{p}, \mathfrak{q})$ with $\mathfrak{p} = \tau\,\mathfrak{a} + (1-\tau)\,\mathfrak{b}$ and with $\mathfrak{q} = \tau\,\mathfrak{b} - (1-\tau)\,\mathfrak{a}$ is a section over the $|x| \le (1 - \frac{1}{16}\mu)\,r_c$ ball of the product $\mathbb{V}\oplus\mathbb{V}$ bundle. This section is denoted by $x_{(\mathfrak{a},\mathfrak{b})}$.

<u>Step 2</u>: The respective projections to $(\wedge^2T^*X)^+\otimes\mathfrak{h}$ and to $(\wedge^2T^*X)^-\otimes\mathfrak{h}$ of the equations in (2.11) over $B_{(1-\mu/4)r_c}$ are equivalent to equations on the $|x| \le (1 - \frac{1}{16}\mu)\,r_c$ ball in $\mathbb{R}^4$ that can be written schematically as

- $(d\mathfrak{q} - m\,e\wedge[\sigma_\Diamond, \mathfrak{p}])^+ + \mathcal{R}^+x_{(\mathfrak{a},\mathfrak{b})} = 0.$
- $(d\mathfrak{p} + m\,e\wedge[\sigma_\Diamond, \mathfrak{q}])^- + \mathcal{R}^-x_{(\mathfrak{a},\mathfrak{b})} = 0.$

(4.16)

with $\mathcal{R}^+$ and $\mathcal{R}^-$ being linear operators that are described in Part 1 of Section 4e. The upcoming Lemma 4.4 says that both $\mathcal{R}^+$ and $\mathcal{R}^-$ are small in a suitable sense. Look at (4.12) to see that the version of (4.16) with both $\mathcal{R}^+$ and $\mathcal{R}^-$ absent defines the respective parts of $\mathcal{L}x_{(\mathfrak{a},\mathfrak{b})}$ in the $(\wedge^2T^*X)^+\otimes\mathfrak{h}$ and $(\wedge^2T^*X)^-\otimes\mathfrak{h}$ of $\mathbb{W}$.

<u>Step 3</u>: The projection to $\mathfrak{h}$ of the equation $*d_A*a = 0$ on $B_{(1-\mu/4)r_c}$ can be written schematically on the $|x| \le (1 - \frac{1}{16}\mu)\,r_c$ ball in $\mathbb{R}^4$ as

$$*(d*\mathfrak{a} - m\,e\wedge[\sigma_\Diamond, *\mathfrak{b}]) + \mathfrak{R}^0{}_{\mathfrak{a}}(\mathfrak{a},\mathfrak{b}) = 0$$

(4.17)

with $\mathfrak{R}^0{}_{\mathfrak{a}}$ being a certain small linear operator. Part 1 of Section 4e says more about $\mathfrak{R}^0{}_{\mathfrak{a}}$. Meanwhile, Part 2 of Section 4e explains why $\mathfrak{a}$ and $\mathfrak{b}$ also obey the equation

$$*(d*\mathfrak{b} + m\,e\wedge[\sigma_\Diamond, *\mathfrak{a}]) + \mathfrak{R}^0{}_{\mathfrak{b}}(\mathfrak{a},\mathfrak{b}) = 0$$

(4.18)

with $\mathfrak{R}^0{}_{\mathfrak{b}}$ being yet another small linear operator.

The equations in (4.17) and (4.18) when written in terms of $\mathfrak{p} =$ $\mathfrak{p} = \tau\,\mathfrak{a} + (1-\tau)\,\mathfrak{b}$ and with $\mathfrak{q} = \tau\,\mathfrak{b} - (1-\tau)\,\mathfrak{a}$ read

- $*(d*\mathfrak{p} - m\,e\wedge[\sigma_\Diamond, *\mathfrak{q}]) + \mathcal{R}^0{}_{\mathfrak{p}}\,x_{(\mathfrak{a},\mathfrak{b})} = 0$ ,
- $*(d*\mathfrak{q} + m\,e\wedge[\sigma_\Diamond, *\mathfrak{a}]) + \mathcal{R}^0{}_{\mathfrak{q}}\,x_{(\mathfrak{a},\mathfrak{b})} = 0$ ,

(4.19)

with $\mathcal{R}^0{}_{\mathfrak{p}}$ and $\mathcal{R}^0{}_{\mathfrak{q}}$ being constant coefficient, linear combinations of the operators $\mathfrak{R}^0{}_{\mathfrak{a}}$ and $\mathfrak{R}^0{}_{\mathfrak{b}}$. The upcoming Lemma 4.4 explains the sense in which these operators are small.

Look at (4.13) to see that the version of (4.19) with both $\mathcal{R}^0{}_p$ and $\mathcal{R}^0{}_q$ absent defines the respective parts of $\mathcal{L}x_{(\mathfrak{a},\mathfrak{b})}$ in the two $\mathfrak{h}$ summands of $\mathbb{W}$.

*Part 5*: Lemma 4.4 follows directly. The lemma uses $\mathcal{R}$ to denote the linear operator from $C^\infty(\mathbb{R}^4; \mathbb{V}\oplus\mathbb{V})$ to $C^\infty(\mathbb{R}^4; \mathbb{W})$ given by the ordered set $((\mathcal{R}^+, \mathcal{R}^0{}_q), (\mathcal{R}^-, \mathcal{R}^0{}_p))$. As the chosen Gaussian coordinate chart for p identifies balls in X of radius less than $c_0^{-1}$ centered at p with balls of the same radius centered at the origin in $\mathbb{R}^4$, the lemma uses $B_r$ for $r \in (0, r_c]$ to denote the radius r ball about the origin in $\mathbb{R}^4$. This notation is used in subsequent arguments as well.

**Lemma 4.4**: *Fix* $\mu \in (0, \frac{1}{2}]$ *and there exists* $\kappa > 1$ *that depends only on* $\mu$ *and the geometric data whose significance is as follows: Suppose that* $r > 1$ *and and that* $(A, a)$ *obeys some* $\tau \in [0, 1]$ *version of* (2.11). *Fix* $p \in U$ *and given* $\delta \in (0,1)$, *take* $c > \kappa\delta^{-1}$ *and then define* $r_{c*}$ *using* p *and* $(A, a)$. *Suppose that* $r_c$ *is positive but no greater than the minimum of* $r_{c\wedge}$, $r_{cF}$ *and* $c^{-1}$. *Assume in addition that* $N(r_c) \le c^{-2}$ *and that* $r_c r_\ddagger \ge 1$. *The operator* $\mathcal{R}$ *is a first order differential operator whose principal symbol is bounded by* $\kappa$ *and whose operator norm is* $\delta$ *small relative to* $\mathcal{L}$ *in the following sense: If* $x$ *is a section of* $\mathbb{V}\oplus\mathbb{V}$ *with compact support on* $B_{(1-\mu/4)r_c}$, *then*

$$\int_{B_{(1-\mu/4)r_c}} |\mathcal{R}x|^2 \le \delta^2 \int_{B_{(1-\mu/4)r_c}} |\mathcal{L}x|^2 \quad .$$

This lemma is proved in Part 3 of Section 4e.

Lemma 4.4 is assumed true so as to complete Lemma 4.3's proof.

***Proof of Lemma 4.3***: The proof has five steps

<u>Step 1</u>: Assume that $c$ and $r_c$ are such that the conclusions of Lemma 4.4 are true for the version with $\delta = \frac{1}{4}$. The bound given by Lemma 4.4 and the identity in (4.15) lead directly to the following observation: If $x$ is a section of $\mathbb{V}\oplus\mathbb{V}$ with compact support in $B_{(1-\mu/4)r_c}$, then

$$\int_{B_{(1-\mu/4)r_c}} (|\nabla x|^2 + 4m^2|x|^2) \le \tfrac{4}{3} \int_{B_{(1-\mu/4)r_c}} |(\mathcal{L} + \mathcal{R})x|^2 \quad . \tag{4.20}$$

Step 2: By way of a reminder, the function $\chi_\mu$ equals 1 where the distance to p is less than $(1 - \frac{3}{4}\mu)\mathrm{r}_c$ and equals 0 where the distance to p is greater than $(1 - \frac{1}{4}\mu)\mathrm{r}_c$. In addition, $|d\chi_\mu| \le c_0\mu^{-1}\mathrm{r}_c^{-1}$. Let $x$ denote the section $(\chi_\mu\mathfrak{p}, \chi_\mu\mathfrak{q})$ of $\mathbb{V}\oplus\mathbb{V}$, this being a section with compact support in the radius $(1 - \frac{3}{4}\mu)\mathrm{r}_c$ ball centered at p. Use what is said in Part 4 to see that $x$ obeys an equation that can be written schematically as

$$\mathcal{L}x + \mathcal{R}x = \mathfrak{s}(d\chi_\mu)\, x_{(\mathfrak{a},\mathfrak{b})} \tag{4.21}$$

with $\mathfrak{s}(d\chi_\mu)$ being a homomorphism from $\mathbb{V}\oplus\mathbb{V}$ to $\mathbb{W}$ with norm bounded by $c_\mu\mathrm{r}_c^{-1}$ and with support in the annulus where the distance to p is between $(1 - \frac{3}{4}\mu)\mathrm{r}_c$ and $(1 - \frac{1}{4}\mu)\mathrm{r}_c$. By way of an explanation, what is denoted by $\mathfrak{s}$ is the principal symbol homomrphism of the operator $\mathcal{L} + \mathcal{R}$. It follows from (4.13) and what is said in Lemma 4.4 that $|\mathfrak{s}| \le c_\mu$. This implies in turn the bound $|\mathfrak{s}(d\chi_\mu)| \le c_\mu\mu^{-1}\mathrm{r}_c^{-1}$.

Step 3: Use (4.20) with (4.21) to conclude that

$$\int_{B_{(1-\mu/4)\mathrm{r}_c}} (|\nabla x|^2 + 4m^2|x|^2) \le c_\mu\mathrm{r}_c^{-2}( \int_{B_{(1-\mu/8)\mathrm{r}_c}} (|\mathfrak{a}|^2 + |\mathfrak{b}|^2) \ . \tag{4.22}$$

To bound the right hand side of (4.22), use the bound $|\mathfrak{a}| \le c_\mu c^{-1/2}$ from the second bullet of (4.5) to bound the integral of $|\mathfrak{a}|^2$ by $c_0 c^{-2}\mathrm{r}_c^4$. To bound the integral of $|\mathfrak{b}|^2$ on the right hand side of (4.22), first use the definition of $\mathfrak{b}$ to bound its absolute value by $c_0 m^{-1}|\nabla_A\sigma_\ddagger|$, and then invoke the fourth bullet of (4.7) to bound the latter by $c_0 m^{-1}|\nabla_A a_\ddagger|$. Use the definition of $\mathrm{N}(\mathrm{r}_c)$ to bound the integral of $|\nabla_A a_\ddagger|^2$ over $B_{\mathrm{r}_c}$ by $\mathrm{N}(\mathrm{r}_c)\mathrm{r}_c^2$. Use this bound to bound the integral of $|\mathfrak{b}|^2$ in (4.22) by $c_0(\mathrm{r}_c m)^{-2}\mathrm{N}(\mathrm{r}_c)\mathrm{r}_c^4$.

Since $m = \frac{1}{\sqrt{2\pi}} r_\ddagger$ and $\mathrm{r}_c r_\ddagger$ is assumed to be greater than 1 and $\mathrm{N}(\mathrm{r}_c) \le c^{-2}$, use of the bounds from the preceding paragraph in (4.22) leads to the bound

$$\int_{B_{(1-\mu/4)\mathrm{r}_c}} (|\nabla x|^2 + 4m^2|x|^2) \le c_\mu\, c^{-1}\mathrm{r}_c^2. \tag{4.23}$$

Step 4 will use (4.23) to bound the integral of $|x|^2$ by $c_\mu c^{-1}m^{-2}\mathrm{r}_c^2$.

Step 4: Use the function $\chi$ to construct a smooth function on $B_{(1-\mu/4)\mathrm{r}_c}$ that equals 1 where the distance to p is less than $(1 - \frac{7}{8}\mu)\mathrm{r}_c$ and equals zero where the distance to p is greater than $(1 - \frac{13}{16}\mu)r_c$. Denote this function by $\chi_{\mu 1}$. This function $\chi_{\mu 1}$ can and should be constructed so that $|d\chi_{\mu 1}| \le c_0 \mu^{-1} r_c^{-1}$.

Use $x_1$ to denote $\chi_{\mu 1} x$, this being a section of $\mathbb{V} \oplus \mathbb{V}$. The section $x_1$ obeys the equation $\mathcal{L} x_1 + \mathcal{R} x_1 = \mathfrak{s}(d\chi_\mu) x$ because the supports of $\chi_{\mu 1}$ and $d\chi_\mu$ are disjoint. This equation for $x_1$ with (4.20) lead to the bound

$$\int_{B_{(1-\mu/4)r_c}} (|\nabla x_1|^2 + 4m^2 |x_1|^2) \le c_0 \mu^{-2} r_c^{-2} \int_{B_{(1-\mu/4)r_c}} |x|^2 \ . \tag{4.24}$$

Use (4.23) to bound the $L^2$ norm of $|x|$ so as to conclude that

$$\int_{B_{(1-\mu/4)r_c}} (|\nabla x_1|^2 + 4m^2 |x_1|^2) \le c_\mu m^{-2} c^{-1} . \tag{4.25}$$

Step 5 will invoke (4.25) to bound the integral of $|x_1|^2$ by $c_0 m^{-4} c^{-1}$.

<u>Step 5</u>: Use the function $\chi$ to construct yet another smooth function on $B_{(1-\mu/4)r_c}$. This next function should equal 1 where the distance to p is less than $(1 - \frac{15}{16}\mu)r_c$ and equal to zero where the distance to p is greater than $(1 - \frac{29}{32}\mu)r_c$. This function is denoted by $\chi_{\mu 2}$. It can and should be constructed so that $|d\chi_{\mu 2}| \le c_0 \mu^{-1} r_c^{-1}$. Note that the supports of $\chi_{\mu 2}$ and $d\chi_{\mu 1}$ are disjoint. Use $x_2$ to denote $\chi_{\mu 2} x_1$. This section of $\mathbb{V} \oplus \mathbb{V}$ obeys the equation $\mathcal{L} x_2 + \mathcal{R} x_2 = \mathfrak{s}(d\chi_{\mu 2}) x_1$ because the supports of $\chi_{\mu 2}$ and $d\chi_{\mu 1}$ are disjoint. The same reasoning that lead to (4.24) leads to the bound

$$\int_{B_{(1-\mu/4)r_c}} (|\nabla x_2|^2 + 4m^2 |x_2|^2) \le c_0 \mu^{-2} r_c^{-2} \int_{B_{(1-\mu/4)r_c}} |x_1|^2 \ ; \tag{4.26}$$

and the latter with (4.25) leads in turn to the bound

$$\int_{B_{(1-\mu/4)r_c}} (|\nabla x_2|^2 + 4m^2 |x_2|^2) \le c_\mu (r_c m)^{-2} m^{-2} c^{-1} . \tag{4.27}$$

The bound that is asserted by Lemma 4.3 follows from (4.27). This is because the top bullet in (2.1) can be used with (4.27) and (4.12) to bound the $L^2$ norm of $\nabla_A x_2$ on $B_{(1-\mu/4)r_c}$ bounded by $c_0$ times the integral on the left hand side of (4.27).

**e) The operator $\mathcal{R}$ and the proof of Lemma 4.3**

The first part of this subsection describes the operators $\mathcal{R}^+$, $\mathcal{R}^-$ and $\mathfrak{R}^0_{\mathfrak{a}}$ that help define the operator $\mathcal{R}$. The second part of the subsection describes the operator $\mathfrak{R}^0_{\mathfrak{b}}$, this being the fourth contributor to $\mathcal{R}$. The third and final part of the subsection proves Lemma 4.4.

*Part 1*: This part of the subsection describes what are denoted by $\mathcal{R}^+$ and $\mathcal{R}^-$ in (4.16) and what is denoted by $\mathfrak{R}^0_{\mathfrak{a}}$ in (4.17). These descriptions use $x$ to denote a given section of $\mathbb{V}\oplus\mathbb{V}$ over $B_{(1-\mu/16)r_c}$ and write $x$ as $(\alpha, \beta)$.

The terms denoted by $\mathcal{R}^+$ and $\mathcal{R}^-$ have three sorts of contributions. There are contributions that account for the fact that the Euclidean self-dual projection in $\wedge^2 T^*\mathbb{R}^4$ differs from that of the pull-back via the Gaussian coordinate chart map of the metric on X. These contributions lead to respective term in $\mathcal{R}^+x$ and $\mathcal{R}^-x$ of the form

$$\mathfrak{t}_+(\nabla\beta - m[\sigma_\lozenge, \alpha]\otimes e) \quad \textit{and} \quad \mathfrak{t}_-(\nabla\alpha + m[\sigma_\lozenge, \beta]\otimes e) \tag{4.28}$$

with $\mathfrak{t}_+$ and $\mathfrak{t}_-$ being endomorphisms on the $|x| < c_0^{-1}$ ball with norms bounded by $c_0|x|^2$.

There are also contributions that account for the fact that $\nu$ differs from the constant 1-form $e$. These lead to respective terms in $\mathcal{R}^+x$ and $\mathcal{R}^-x$ of the form

$$-m\, j_+([\tau, \alpha]\otimes(\nu - e)) \quad \textit{and} \quad m\, j_-([\tau, \beta]\otimes(\nu - e))\ , \tag{4.29}$$

with $j_+$ and $j_-$ being endomorphisms on $B_{(1-\mu/16)r_c}$ with norms bounded by $c_0$.

There are contributions to $\mathcal{R}^+$ and $\mathcal{R}^-$ that account for the fact that the connection $\hat{A}$ need not be flat on $B_{(1-\mu/16)r_c}$. These contributions lead to terms in $\mathcal{R}^+$ and $\mathcal{R}^-$ that have the respective forms

$$j_+([\tau, \beta]\otimes a_\lozenge) \quad \textit{and} \quad j_-([\tau, \alpha]\otimes a_\lozenge) \tag{4.30}$$

with $a_\lozenge$ again denoting the 1-form on $B_{(1-\mu/16)r_c}$ that is defined by using the product structure for P from Part 3 to write $\hat{A}$ as $\theta_0 + a_\lozenge\sigma_\lozenge$ on $B_{(1-\mu/16)r_c}$.

What is denoted by $\mathfrak{R}^0_{\mathfrak{a}}$ in (4.17) also has three different contributions. There are contributions to account for the fact that the Euclidean metric and the metric pulled back from X via the Gaussian coordinate chart map differ. These lead to a term of the form

$$\mathfrak{t}_{01}(\nabla\alpha - m e\otimes[\tau,\beta]) + \mathfrak{t}_{00}(\alpha)$$

(4.31)

with $\mathfrak{t}_{01}$ and $\mathfrak{t}_{00}$ being tensor homomorphisms on the $|x| < c_0^{-1}$ ball whose norms are bounded respectively by $c_0|x|^2$ and by $c_0|x|$. Note also that $\mathfrak{t}_{01}$ when viewed as a tensor is symmetric. There is also a contribution to $\mathfrak{R}^0{}_{\mathfrak{a}}$ from the fact that $\nu$ and $e$ differ. This contribution is accounted for by a term of the form

$$m\, \mathfrak{j}_{0\mathfrak{a}}([\tau,\beta]\otimes(\nu - e))$$

(4.32)

with $\mathfrak{j}_{0\mathfrak{a}}$ being an endomorphism on $B_{(1-\mu/16)r_c}$ with norm bounded by $c_0$. There is also a term in $\mathfrak{R}^0{}_{\mathfrak{a}}$ that has the form

$$\mathfrak{j}_{0\mathfrak{a}}([\tau,\alpha]\otimes a_\Diamond)\ .$$

(4.33)

This third contribution to $\mathfrak{R}^0{}_{\mathfrak{a}}$ accounts for the fact that the connection $\hat{A}$ need not be flat.

*Part* 2: The operator $\mathfrak{R}^0{}_{\mathfrak{b}}$ that appears in (4.18) is given in the upcoming (4.38) at the very end of this part of the subsection. The intervening discussion defines the terms that appear in (4.38). By way of a look ahead, these definitions are obtained by rewriting $d{*}\mathfrak{b}$ and $m e\wedge[\tau,*\mathfrak{a}]$ using (2.5), (4.4), the identity in the middle bullet of (4.5) and (4.7).

The rewriting of the $d{*}\mathfrak{b}$ starts with the introduction of the endomorphism of the product bundle $B_{(1-\mu/16)r_c}\times\mathfrak{h}$ that is defined by the rule whereby a given element $\theta$ is sent to $2\pi^2(|\nu|^2\theta - \frac{1}{4}\langle[\theta,\sigma_\Diamond]\mathfrak{a}_\alpha\rangle[\sigma_\Diamond,\mathfrak{a}_\alpha])$. This endomorphism is denoted by $\mathfrak{D}$. It follows from the third bullet of (4.5) that $\mathfrak{D}$ is invertible if $c > c_0$; and it follows from Proposition 3.2 that both it and its inverse differ from the identity by a term with norm at most $c_0 c^{-1}$. Keeping in mind that $\mathfrak{b} = \frac{1}{4} m^{-1}[\sigma_\Diamond,\nabla_A\sigma_\Diamond]$, use the second bullet in (4.7) to obtain a formula for $\mathfrak{b}$ that has the schematic form

$$\mathfrak{b} = \tfrac{1}{2}\pi^2 m^{-1}\mathfrak{D}^{-1}\big([\sigma_\Diamond,\nu_\alpha(\nabla_A a_\ddagger)_\alpha] + \langle\sigma_\Diamond(\nabla_A a_\ddagger)_\alpha\rangle[\sigma_\Diamond,a_{\ddagger\alpha}]\big)$$

(4.34)

with the subscript indices indicating the components of $\mathfrak{b}$, $\nu$, $a_\ddagger$ and $\nabla_A$ when written using an orthonormal frame that is defined by the pull-back of the metric from X via the Gaussian coordinate chart map. Repeated indices are to be summed. This convention about indices is also used in the rest of Part 2. Note that the notation in (4.34) uses $a_\ddagger$ to denote the pull-back via $h_\Diamond$ and the Gaussian coordinate chart map of its namesake $(P\times_{SO(3)}\mathfrak{S})$-valued 1-form on $B_{(1-\mu/16)r_c}$. The incarnation of $a_\ddagger$ in (4.34) is a $\mathfrak{S}$-valued 1-form on the $|x| \le (1 - \frac{1}{16}\mu)\,r_c$ ball in $\mathbb{R}^4$. This same incarnation is used (4.35) and (4.36).

Introduce by way of notation $\nabla_A{}^m$ to denote the formal $L^2$ adjoint of the covariant derivative operator $\nabla_A$ with the pull-back of the metric from X used to define inner products and the volume 4-form. Differentiate (4.34) to obtain a formula for $\nabla_A{}^m\mathfrak{b}$ that has the schematic form

$$\nabla_A{}^m\mathfrak{b} = \tfrac{1}{4}\,m^{-1}\,\mathfrak{Q}^{-1}\big([\sigma_\lozenge, \nu_\alpha(\nabla_A{}^m\nabla_A a_\ddagger)_\alpha] + \langle\sigma_\lozenge\,(\nabla_A{}^m\nabla_A a_\ddagger)_\alpha\rangle\,[\sigma_\lozenge, a_{\ddagger\alpha}]\big) + \mathfrak{E}_0(\mathfrak{b}) + m^{-1}\,\mathfrak{E}_1(\nabla_A\mathfrak{a}) \tag{4.35}$$

with $\mathfrak{E}_0$ and $\mathfrak{E}_1$ being homomorphisms that are described directly. The action of the homomorphism $\mathfrak{E}_0$ on a given element $\beta$ is that of $\mathfrak{Q}^{-1}$ on

$$-(\nabla_{A,\mu}\mathfrak{Q})\,\mathfrak{D}^{-1}\,\beta_\mu\; - \;[[\beta_\mu,\sigma_\lozenge], \nu_\alpha(\nabla_{A,\mu}a_\ddagger)_\alpha] - \langle[\beta_\mu,\sigma_\lozenge]\,(\nabla_{A,\mu}a_\ddagger)_\alpha\rangle\,[\sigma_\lozenge, a_{\ddagger\alpha}] \\ - \langle\sigma_\lozenge\,(\nabla_{A,\mu}a_\ddagger)_\alpha\rangle\,[[\beta_\mu,\sigma_\lozenge], a_{\ddagger\alpha}]. \tag{4.36}$$

Note with regards to (4.36) that $[\beta,\sigma_\lozenge] = -\nabla_A\sigma_\lozenge$ when $\beta = \frac{1}{4}\,[\sigma_\lozenge, \nabla_A\sigma_\lozenge]$. Meanwhile, the action of the homomorphism $\mathfrak{E}_1$ on a given element $\mathfrak{q}$ of $(\otimes_2 T^*\mathbb{R}^4)\otimes\mathfrak{S}$ is that of the endomorphism $\mathfrak{Q}^{-1}$ on $-2\,(\nabla_\mu\nu)_\beta[\sigma_\lozenge, \mathfrak{q}_{\mu\beta}]$. Note in particular that the appearance in $\mathfrak{E}_1$ of the commutator with $\sigma_\lozenge$ implies that $\mathfrak{E}_1$ can be viewed as a homomorphism with domain $(\otimes_2 T^*\mathbb{R}^4)\otimes\mathfrak{S}$. This last fact has the following consequence: The term $m^{-1}\,\mathfrak{E}_1(\nabla_A\mathfrak{a})$ in (4.35) can be written using the identifications in Part 3 as $m^{-1}\,\mathfrak{E}_1(\nabla\mathfrak{a} + a_\lozenge\otimes[\sigma_\lozenge, \mathfrak{a}])$.

The next observation concerns the terms on the right hand side of (4.35) that contain $\nabla_A{}^m\nabla_A a_\ddagger$. These terms can be rewritten using (4.5) as $m^{-1}\,(r_\ddagger{}^2\,\mathfrak{E}_2(\mathfrak{a}\otimes\mathfrak{a}) + \mathfrak{E}_3(\mathfrak{a}))$ with $\mathfrak{E}_2$ and $\mathfrak{E}_3$ being endomorphisms with norms bounded by $c_0$.

The left hand side of (4.35) is $\nabla_A{}^m\mathfrak{b}$, and this can be written as

$$\nabla_A{}^m\mathfrak{b} = -{*d}({*\mathfrak{b}}) - \mathfrak{t}_{01}(\nabla\mathfrak{b}) - \mathfrak{t}_{00}(\mathfrak{b}) - \mathfrak{j}_{0\mathfrak{b}}(a_\lozenge\otimes[\tau, \mathfrak{b}]) \tag{4.37}$$

with $\mathfrak{t}_{01}$ and $\mathfrak{t}_{00}$ as in (4.31) and with $\mathfrak{j}_{0\mathfrak{b}}$ being a homomorphism with norm bounded by $c_0$. The fact that $\mathfrak{t}_{01}$ is a symmetric tensor homomorphism accounts for the appearance of only the $\nabla\mathfrak{b}$ and $a_\ddagger\otimes[\tau, \mathfrak{b}]$ parts of the covariant derivative when writing $\nabla_A{}^m\mathfrak{b}$.

A final observation concerns ${*}(m\,e\wedge[\sigma_\lozenge, {*\mathfrak{a}}])$: This term can be rewritten using the fact that $\nu_\alpha\mathfrak{a}_\alpha = 0$ as $-m\,\mathfrak{t}_{01}(e\otimes[\sigma_\lozenge, \mathfrak{a}]) - m\,\mathfrak{j}_{0\mathfrak{a}}([\sigma_\lozenge, \mathfrak{a}]\otimes(\nu - e))$ with $\mathfrak{t}_{01}$ as in (4.30) and with $\mathfrak{j}_{0\mathfrak{a}}$ as in (4.32).

With this admittedly large body of new notation in hand, define $\mathfrak{N}^0{}_\mathfrak{b}$ by the rule

$$\mathfrak{R}^0{}_{\mathfrak{b}}x = m^{-1}r_{\ddagger}{}^2\mathfrak{E}_2(\mathfrak{a}\otimes\alpha) + m^{-1}\mathfrak{E}_3(\alpha) + m^{-1}\mathfrak{E}_1(\nabla\alpha + \mathrm{a}_{\ddagger}\otimes[\tau,\alpha]) +$$
$$m\,\mathfrak{t}_{01}(e\otimes[\tau,\alpha]) + m\,\mathfrak{j}_{0\mathfrak{a}}([\tau,\alpha]\otimes(\nu - e)) + \mathfrak{E}_0(\beta) + \mathfrak{t}_{01}(\nabla\beta) + \mathfrak{t}_{00}(\beta) + \mathfrak{j}_{0\mathfrak{b}}(\mathrm{a}_{\diamond}\otimes[\tau,\beta]). \qquad (4.38)$$

The fact that (4.18) holds is a tautological consequence of the definitions of the terms on the right hand side of (4.38).

*Part 3*: This last part of the subsection addresses the assertions in Lemma 4.4.

***Proof of Lemma 4.4***: The proof has seven steps. The notation is such that $x$ denotes a section of $\mathbb{V}\oplus\mathbb{V}$ with compact support in $B_{(1-\mu/4)r_c}$. Keep in mind that the identity in (4.15) in the case when $x$ has compact support in $B_{(1-\mu/4)r_c}$ can be written as

$$\int_{B_{(1-\mu/4)r_c}} |\mathcal{L}x|^2 = \int_{B_{(1-\mu/4)r_c}} (|\nabla x|^2 + m^2|x|^2)\ . \qquad (4.39)$$

As before, $x$ is written at times as $(\alpha,\beta)$ with $\alpha$ and $\beta$ being sections of $\mathbb{V}$.

<u>Step 1</u>: It follows from the formulas in Parts 1 and 2 for $\mathcal{R}^+$, $\mathcal{R}^-$, $\mathfrak{R}^0{}_{\mathfrak{a}}$ and $\mathfrak{R}^0{}_{\mathfrak{b}}$ that $\mathcal{R}$ is a first order differential operator. Its principal symbol comes from the terms in (4.28), (4.31), the term that involves $\mathfrak{t}_{01}$ in (4.38) and the $m^{-1}\mathfrak{E}_1$ term (4.38). The terms in (4.38), (4.31) and the $\mathfrak{t}_{01}$ term in (4.38) account for the fact that the metric on $B_{r_c}$ is not the flat metric. It follows as a consequence that the norm of the contribution to the principal symbol of $\mathcal{R}$ from these terms is bounded by $c_0 r_c$.

Use what is said subsequent to (4.36) about $\mathfrak{E}_1$ to see that the contribution to the principal symbol of $\mathcal{R}$ from $m^{-1}\mathfrak{E}_1$ is bounded by $c_0 m^{-1}|\nabla\nu|$. This in turn is bounded by $c_0 m^{-1}|\nabla_A a_{\ddagger}|$ because of the bound in the fourth bullet of (4.7). Meanwhile, the third bullet of Proposition 3.3 bounds $|\nabla_A a_{\ddagger}|$ by $c_\mu r_{\ddagger}$ when $r_c r_{\ddagger}$ is greater than 1; and so a $c_\mu$ bound on the norm of the $m^{-1}\mathfrak{E}_1$ contribution to $\mathcal{R}$'s principal symbol follows because $m = \frac{1}{\sqrt{2\pi}} r_{\ddagger}$.

<u>Step 2</u>: This step derives bounds for the $L^2$ norms of $\mathcal{R}^+x$, $\mathcal{R}^-x$ and $\mathfrak{R}^0{}_{\mathfrak{a}}x$. Keep in mind that $\mathcal{R}^0{}_{\mathrm{p}}$ and $\mathcal{R}^0{}_{\mathrm{q}}$ are constant coefficient, linear combinations of $\mathfrak{R}^0{}_{\mathfrak{a}}$ and $\mathfrak{R}^0{}_{\mathfrak{b}}$ so bounds on the $L^2$ norms of $\mathfrak{R}^0{}_{\mathfrak{a}}x$ and $\mathfrak{R}^0{}_{\mathfrak{b}}x$ lead directly to bounds on the $L^2$ norms of $\mathcal{R}^0{}_{\mathrm{p}}x$ and $\mathcal{R}^0{}_{\mathrm{q}}x$. To start the derivation, use what is said about the pair of homomorphisms $\mathfrak{t}_+$ and $\mathfrak{t}_-$ in (4.28) and the pair $\mathfrak{t}_{01}$ and $\mathfrak{t}_{00}$ in (4.31) to see that the square of the pointwise norm of what the expressions in (4.28) and (4.31) are at most $c_\mu r_c{}^2(|\nabla x|^2 + m^2|x|^2)$. Meanwhile, the square of the pointwise norm of the homomorphisms in (4.29) and (4.32) are at most $c_\mu c^{-2} m^2 |x|^2$ because of the bound $|\nu - e| < c_0 c^{-1}$ from the fifth bullet of (4.5). It follows from the preceding that these parts of $\mathcal{R}^+x$, $\mathcal{R}^-x$ and $\mathfrak{R}^0{}_{\mathfrak{a}}x$ contribute at most

$$c_\mu (c^{-2} + r_c^2) \int_{B_{(1-\mu/4) r_c}} (|\nabla x|^2 + m^2 |x|^2) \tag{4.40}$$

to the square of the $L^2$ norms of $\mathcal{R}x$.

Step 3: The contributions to the $L^2$ norms of $\mathcal{R}^+x$ and $\mathcal{R}^-x$ from the terms in (4.30) and to that of $\mathfrak{R}^0{}_{\mathfrak{a}}x$ from (4.33) are no greater than

$$c_\mu \Big( \int_{B_{(1-\mu/4) r_c}} |a_\diamond|^4 \Big)^{1/2} \Big( \int_{B_{(1-\mu/4) r_c}} |x|^4 \Big)^{1/2} . \tag{4.41}$$

Meanwhile, (4.12) with two appeals to the top bullet in (2.1) bounds (4.41) by

$$c_\mu c^{-2} \int_{B_{(1-\mu/4) r_c}} (|\nabla x|^2 + r_c^{-2} |x|^2) . \tag{4.42}$$

To complete the story for the (4.30) and (4.33) contributions, invoke the definition of $m$ as $\frac{1}{\sqrt{2}\,\pi} r_\ddagger$ to see that $r_c^{-1} < \sqrt{2}\,\pi m$ when $r_c r_\ddagger \ge 1$. It follows as a consquence that (4.42) is at most $c_\mu c^{-2}$ times the expression on the right hand side of (4.39) and thus at most $c_\mu c^{-2}$ times $\int_{B_{(1-\mu/4) r_c}} |\mathcal{L}x|^2$ if $r_c r_\ddagger$ is greater than or equal to 1.

Step 4: This step and the next two supply bounds for the square of the $L^2$ norm of the various terms in the depiction by (4.38) of $\mathfrak{R}^0{}_{\mathfrak{b}}x$. The terms on the right hand side of (4.38) are discussed in turn starting from the leftmost and moving to the right.

The endomorphism $\mathfrak{E}_2$ that appears in (4.38) is bounded by $c_0$, and this leads to the pointwise bound on $m^{-1} r_\ddagger^2 \mathfrak{E}_2(\mathfrak{a} \otimes \alpha)$ by $c_0 m |\mathfrak{a}| |\alpha|$. Keeping in mind that $|\mathfrak{a}| \le c_0 c^{-1/2}$, this is no greater than $c_0 c^{-1/2} m |x|$. Use the latter bound to bound the square of the $L^2$ norm of the term $m^{-1} r_\ddagger^2 \mathfrak{E}_2(\mathfrak{a} \otimes \alpha)$ by $c_0 c^{-1} m^2$ times the integral of $|x|^2$. What with (4.39), the latter is no larger than $c_0 c^{-1}$ times the square of the $L^2$ norm of $\mathcal{L}x$.

The norm of the homomorphism $\mathfrak{E}_3$ is bounded by $c_0$ so the square of the $L^2$ norm of $m^{-1} \mathfrak{E}_3(\alpha)$ is at most $c_0 m^{-2}$ times the square $L^2$ norm of $|x|$, this being at most $c_0 m^{-4}$ times the square of the $L^2$ norm of $\mathcal{L}x$. As the assumption $r_c r_\ddagger \ge 1$ implies that $m^{-1} < r_c$, the square of the $L^2$ norm of $m^{-1}\mathfrak{R}_3(\alpha)$ is at most $c_0 r_c^4$ times that of $\mathcal{L}x$, this being at most $c_0 c^{-4}$ times the square of the $L^2$ norm of $\mathcal{L}x$ because $r_c$ is taken to be less than $c^{-1}$.

<u>Step 5</u>: The norm of the homomorphism $\mathfrak{E}_1$ is at most $c_0|\nabla\nu|$; and the fourth bullet of (4.7) asserts that this is at most $c_0|\nabla_A a_\ddagger|$. Since $m = \frac{1}{\sqrt{2}\,\pi} r_\ddagger$, the version of Proposition 3.3's second bullet with $\mu$ replaced by $\mu/16$ asserts that $|\nabla_A a_\ddagger|$ is no greater than $c_\mu c^{-1} m$. It follows from this that the square of the $L^2$ norm of $m^{-1}\mathfrak{E}_1(\nabla\alpha)$ is at most

$$c_\mu c^{-2} \int_{B_{(1-\mu/4)r_c}} |\nabla x|^2 \tag{4.43}$$

and that the square of the $L^2$ norm of $m^{-1}\mathfrak{E}_1(a_\diamond\otimes[\sigma_\diamond,\alpha])$ is at most $c_\mu c^{-2}$ times the expression in (4.40). Granted what is said subsequent to (4.41), this in turn is no less than $c_\mu c^{-2}$ times the sum of the squares of the $L^2$ norms of $|\nabla x|$ and $m|x|$.

Use the preceding observations with (4.39) to see that the square of the $L^2$ norm of $m^{-1}\mathfrak{E}_1(\nabla\alpha + a_\diamond\otimes[\sigma_\diamond,\alpha])$ is at most $c_\mu c^{-2}$ times the square of the $L^2$ norm of $\mathcal{L}x$.

<u>Step 6</u>: Since $|\mathfrak{t}_{01}| < c_0 r_c^2$, the square of the $L^2$ norm of $m\,\mathfrak{t}_{01}(e\otimes[\sigma_\diamond,\alpha])$ is at most $c_0 r_c^2$ times that of $m|x|$; and this is at most $c_0 c^{-2}$ times the square of the $L^2$ norm of $\mathcal{L}x$ because of (4.39) and because $r_c \le c^{-1}$. Since $|\nu - e| < c_0 c^{-1}$ and $|\mathfrak{j}_{0\mathfrak{a}}| \le c_0$, the square of the $L^2$ norm of $m\,\mathfrak{j}_{0\mathfrak{a}}([\sigma_\diamond,\alpha]\otimes(\nu - e))$ is at most $c_\mu c^{-2}$ times that $m|x|$ and thus at most $c_\mu c^{-2}$ times that of $\mathcal{L}x$.

A look at (4.36) finds the norm of the homomorphism $\mathfrak{E}_0$ less than $c_0|\nabla_A a|$. With this understood, then the argument in Step 5 can be repeated with at most cosmetic changes to prove that the square of the $L^2$ norm of $\mathfrak{E}_0(\beta)$ is at most $c_\mu c^{-2}$ times the square of the $L^2$ norm of $\mathcal{L}x$.

But for cosmetic changes, the arguments in Step 2 that concern the $\mathfrak{t}_{01}(\nabla\alpha)$ and $\mathfrak{t}_{00}(\alpha)$ parts of (4.38) can be repeated to prove that the square of the $L^2$ norms of $\mathfrak{t}_{01}(\nabla\beta)$ and $\mathfrak{t}_{00}(\beta)$ are no greater than $c_0 c^{-2}$ times the square of the $L^2$ norm of $\mathcal{L}x$. By the same token, the arguments in Step 2 that concern the endomorphism in (4.33) can be repeated in an almost verbatim fashion to prove that the square of the $L^2$ norm of $\mathfrak{j}_{0\mathfrak{b}}(a_\diamond\otimes[\sigma_\ddagger,\beta])$ is at most $c_\mu c^{-2}$ times the square of the $L^2$ norm of $\mathcal{L}x$.

<u>Step 7</u>: Fix $\delta < c_0^{-1}$. If $c^{-1}$ is less than $c_\mu^{-1}\delta$ and if $r_c r_\ddagger$ is greater than 1, then the bounds in Steps 2-6 lead directly to the assertion made by Lemma 4.4 about the small size of $\mathcal{R}$ relative to $\mathcal{L}$.

## 5. Monotonicity for the $(A, a)$ versions of N

Take X in this section to be an oriented, Riemannian 4-manifold and suppose that $U \subset X$ is a specified open set with compact closure. Use P to denote a given principal SO(3) bundle over X. As in Sections 3 and 4, the term *geometric data* refers to the data set that consists of the manifold X, the set U, the bundle P.

Suppose that $r > 1$ and $(A, a)$ are a pair of connection on P and section of the bundle $T^*X \otimes (P \times_{SO(3)} \mathfrak{G})$ that obey a given $\tau \in [0, 1]$ version of (2.11) on X. These equations imply among other things that $(A, a)$ obey (2.5) and (2.6). Some of the subsequent assertions in this section make the additional assumption that

$$\int_X |F_A - r^2 a \wedge a|^2 \le \mathrm{E}^2 \tag{5.1}$$

with E being a given positive number. If X is compact, then (5.1) holds if E is greater than $c_0\, |p_1(P \times_{SO(3)} \mathfrak{G})|$, this being a consequence of (2.17). By way of another relevant example, suppose that $\mathbb{I}$ and M and P are as described in Proposition 2.2. Let X denote the interior of $\mathbb{I} \times M$ and set $U = I \times M$ with $I \subset \mathbb{I}$ being an open interval with compact closure in the interior of $\mathbb{I}$. Then (3.1) holds for a suitable E if $(A, a)$ are boundary convergent. A lower bound for E is determined in this case by the boundary limits of $(A, a)$, this being a consequence of (2.19) and (2.26).

Fix a point $p \in U$ and use $(A, a)$ to define the functions K and N using the formulae in Section 3a. The upcoming Proposition 5.1 makes an assertion to the effect that N can not be small at large r if it is large at small r. By way of a heads up, the proposition defines a number from $(0, c_0^{-1}]$ it denotes by $r_{*p}$; and the definition uses implicitly the fact that the function $r \to r\,K(r)\,r$ on $[0, c_0^{-1}]$ is 1-1 onto its image in $[0, \infty)$. The function $r \to r\,K(r)\,r$ is 1-1 because K is increasing.

**Proposition 5.1**: *There exists $z_U > 1$, $\kappa > 1$ and given $\varepsilon \in (0, 1]$, there exists $\kappa_\varepsilon > \kappa$ with all three depending on the geometric data and with $\kappa_\varepsilon$ depending also on $\varepsilon$. These numbers have the following significance: Suppose that $r > \kappa$ has been specified and that $(A, a)$ is a pair of connection on P and section of $T^*X \otimes (P \times_{SO(3)} \mathfrak{G})$ that obeys the equations in (2.11) for a given $\tau \in [0, 1]$. Assume in addition that $(A, a)$ obeys the bounds in (5.1) with E being less than $\varepsilon^{-1}$. Fix $p \in U$. Assume that there exists $r \in (0, \kappa^{-1}]$ where $r\,K(r)\,r = z_U^{-1}$ and denote this value of r by $r_{*p}$. Assume in addition that*

$$\int_{\mathrm{dist}(\cdot,p)<\frac{1}{2}r_{*p}} |F_A|^2 \leq z_U^{-2}.$$

*Use* p *and* (A, *a*) *to define* N *and* K. *If* $r_0 \in (0, \kappa^{-1})$ *is such that* $N(r_0) \geq \varepsilon$ *and if* $r_1 \in (0, \kappa^{-1}]$ *is greater than the maximum of* $r_0$ *and* $r_{*p}$, *then* $N(r_1) \geq \kappa_\varepsilon^{-1}$.

Section 5b states and proves a lemma that leads directly to Proposition 5.1 when $r_0 K(r_0)\, r$ is large. The lemma in Section 5b is used in Section 5d to complete the proof of Proposition 5.1. Section 5a supplies a formula for the derivative of N that is used in Section 5b. Section 5c talks about Karen Uhlenbeck's compactness theorem and it defines the number $z_U$. Section 5e supplies the proof of a lemma in Section 5d

## a) The derivative of N

A formula for the derivative of N is given below in (5.2). The notation in this equation uses $\nabla_{A,r}a$ to denote the section of $P \times_{SO(3)} \mathfrak{S}$ over $B_r$–p that gives the pairing between $\nabla_A a$ and the vector field of unit length that points outward along the the geodesic arcs that start at p. The notation also uses $E_A$ to denote the $(P \times_{SO(3)} \mathfrak{S})$-valued 1-form on $B_r$–p that gives the pairing between $F_A$ and the same unit length vector field. Keeping in mind that a ball of radius $c_0^{-1}$ centered at p is well a Gaussian coordinate chart centered at p for a neighborhood of p in X, the upcoming formula in (5.2) holds for $r \in (0, c_0^{-1}]$.

The promised formula for N's derivative at r is as follows:

$$\frac{d}{dr} N = \frac{2}{r^2 K^2} \int_{\partial B_r} |\nabla_{A,r} a - \tfrac{1}{r} N a|^2 + \frac{r^2}{r^2 K^2} \int_{\partial B_r} |a \wedge a|^2 + \frac{1}{r^2 K^2 r^2} \int_{\partial B_r} (2\,|E_A|^2 - |F_A|^2) + \mathfrak{q}\,, \tag{5.2}$$

with $\mathfrak{q}$ being a function on $[0, c_0^{-1}]$ whose absolute value obeys

$$|\mathfrak{q}| \leq c_0 r + c_0 \frac{1}{r K^2} \int_{B_r} (|\nabla_A a|^2 + 2r^2\, |a \wedge a|^2 + r^{-2}\, |F_A|^2)\,. \tag{5.3}$$

Note by way of a parenthetical remark, the term denoted by $\mathfrak{q}$ is zero if the metric on $B_r$ is flat. The six steps that follow derive (5.2) and (5.3). The derivation owes allegiance to the derivations in [Al], [HHL] and [DF] of formulas for the derivatives of ancestors of N that were used to study singular level sets of solutions to elliptic equations.

<u>Step 1</u>: Differentiate the formula for N in (3.5) and use (3.6) to write

$$\frac{d}{dr} N = \frac{1}{r^2 K^2} \int_{\partial B_r} (|\nabla_A a|^2 + 2r^2 |a \wedge a|^2) - \frac{2}{r} N(1+N) \,.$$

(5.4)

This formula starts the derivation of (5.2) but it is useful in its own right.

<u>Step 2</u>: Fix an orthonormal frame for $T^*X$ on $B_r$. Use this frame to write the 4 directional covariant derivatives of $\nabla_A a$ along the dual vector fields as $\{\nabla_{A,\alpha} a\}_{\alpha \in \{1,2,3,4\}}$ and for each $\alpha \in \{1, 2, 3, 4\}$, write the 4 components of $\nabla_{A,\alpha} a$ as $\{(\nabla_{A,\alpha} a)_\beta\}_{\beta \in \{1,2,3,4\}}$. Meanwhile, write the components of the curvature $F_A$ with respect to the chosen frame as $\{F_{A,\alpha\beta}\}_{\alpha,\beta \in \{1,2,3,4\}}$. Granted this notation, let T denote the symmetric section over $B_r$ of $T^*X \otimes T^*X$ with components $\{T_{\alpha\beta}\}_{\alpha.\beta \in \{1.2,3,4\}}$ that are given by the formula

$$T_{\alpha\beta} = \langle (\nabla_{A,\alpha} a)_\nu (\nabla_{A,\beta} a)_\nu \rangle + r^{-2} \langle F_{A,\alpha\nu} F_{A,\beta\nu} \rangle - \tfrac{1}{2} \delta_{\alpha\beta} (|\nabla_A a|^2 + r^2 |a \wedge a|^2 + r^{-2} |F_A|^2) \,.$$

(5.5)

The formal, $L^2$ adjoint of the covariant derivative maps sections of $T^*X \otimes T^*X$ to $T^*X$. This operator is denoted by $\nabla^\dagger$. Step 7 proves that $\nabla^\dagger T = 0$ when the metric on $B_r$ is flat. When this is not the case, Step 6 proves that $\nabla^\dagger T$ can be written as

$$\nabla^\dagger T = \mathfrak{r}$$

(5.6)

where the notation has $\mathfrak{r}$ denoting a 1-form whose norm obeys $|\mathfrak{r}| \le c_0 |a| |\nabla_A a|$. A norm bound of this sort exists because $\mathfrak{r}$ can be written as $\mathfrak{r}_{\alpha\beta\sigma} \langle a_\alpha \nabla_{A,\beta} a_\sigma \rangle$ with each 1-form from the set $\{\mathfrak{r}_{\alpha\beta\sigma}\}_{\alpha,\beta,\sigma \in \{1,2,3,4\}}$ coming from the Riemann curvature tensor.

<u>Step 3</u>: Use $\hat{x}$ to denote the differential of the function $\frac{1}{2} \mathrm{dist}(p, \cdot)^2$. The norm of $\hat{x}$ is $\mathrm{dist}(p,\cdot)$ and its dual is a vector field that is tangent to the geodesic rays from the point p. The covariant derivative of $\hat{x}$ can be written as $\mathfrak{m} + \mathfrak{z}$ with $\mathfrak{m}$ being the metric tensor and $\mathfrak{z}$ being a tensor with norm bounded by $c_0 \mathrm{dist}(p,\cdot)^2$.

With the preceding understood, take the inner product of both sides of (5.6) with $\hat{x}$ and integrate the resulting identity over $B_r$. Having done so, integrate by parts to remove derivatives from T and $\mathfrak{R}$ to derive the identity

$$\frac{1}{2} \int_{\partial B_r} (|\nabla_A a|^2 + r^2 |a \wedge a|^2) = \int_{\partial B_r} (|\nabla_{A,r} a|^2 + \frac{1}{r^2} |E_A|^2 - \frac{1}{2r^2} |F_A|^2) + \frac{1}{r} \int_{B_r} (|\nabla_A a|^2 + 2r^2 |a \wedge a|^2)$$
$$+ \frac{1}{r} \int_{B_r} (\mathfrak{m}(\mathfrak{z}, T) - \mathfrak{m}(\hat{x}, \mathfrak{r}))$$

(5.7)

with $\mathfrak{m}(\cdot\,,\cdot)$ denoting the metric inner product on both $T^*X$ and $T^*X \otimes T^*X$.

Step 4: Use (5.7) with the definition of N to rewrite (5.4) so as to read

$$\frac{d}{dr}\mathrm{N} = \frac{2}{r^2K^2}\int_{\partial B_r} |\nabla_{A,r} a|^2 - \frac{2}{r}\mathrm{N}^2 + \frac{1}{r^2K^2}\int_{\partial B_r} (r^2|a\wedge a|^2 + \frac{2}{r^2}|E_A|^2 - \frac{1}{r^2}|F_A|^2)$$
$$- \frac{2}{r^3K^2}\int_{B_r} \mathfrak{m}(\hat{x},\mathfrak{r}) + \frac{2}{r^3K^2}\int_{B_r}\mathfrak{m}(\mathfrak{z},T)\,. \tag{5.8}$$

The next task is to rewrite the term $-\frac{2}{r}\mathrm{N}^2$ in (5.8). This is done by using Stokes theorem and the identity in (2.6) to write N as

$$\mathrm{N} = \frac{1}{r^2K^2}\int_{\partial B_r}\langle a\nabla_{A,r}a\rangle - \frac{1}{r^2K^2}\int_{B_r}\mathrm{Ric}(\langle a\otimes a\rangle)\ . \tag{5.9}$$

Use this identity to replace the two left most terms on the right hand side of (5.8) by

$$\frac{2}{r^2K^2}\int_{\partial B_r}|\nabla_{A,r}a - \frac{1}{r}\mathrm{N}a|^2 + \mathfrak{e}_1 + \mathfrak{e}_2 \tag{5.10}$$

with $\mathfrak{e}_1$ and $\mathfrak{e}_2$ being functions with absolute value obeying

$$|\mathfrak{e}_1| \le c_0\frac{1}{r^5K^2}(\int_{B_r}|a|^2)^2 \quad \textit{and} \quad |\mathfrak{e}_2| \le c_0\frac{1}{r^5K^2}(\int_{B_r}|a|^2)\int_{B_r}(|\nabla_A a|^2 + 2r^2|a\wedge a|^2)\,. \tag{5.11}$$

Replacing the left most two terms on the right hand side of (5.8) with (5.11) results in an equation for N's derivative that is identical to (5.2) with $\mathfrak{q}$ in (5.2) being the sum of $\mathfrak{e}_1$ and $\mathfrak{e}_2$ plus the terms in (5.8) that are proportional to the integrals over $B_r$ of $\mathfrak{m}(\hat{x},\mathfrak{r})$ and $\mathfrak{m}(\mathfrak{z},T)$.

Step 5: This step proves that $\mathfrak{q}$'s absolute value obeys (5.3). To this end, consider first the contribution to $\mathfrak{q}$ of the term

$$\frac{2}{r^3K^2}\int_{B_r}\mathfrak{m}(\mathfrak{z},T) \tag{5.12}$$

from (5.8). The integral of $\mathfrak{m}(\mathfrak{z},T)$ over $B_r$ is no greater than $c_0r^2$ times that of the sum $|\nabla_A a|^2 + 2r^2|a\wedge a|^2 + r^{-2}|F_A|^2$ and so the absolute value of what is written in (5.12) obeys the asserted bound for $|\mathfrak{q}|$ in (5.3). Such is also the case for the absolute value of what is denoted by $\mathfrak{e}_2$ in (5.10) as can be seen using (5.11) and (3.12). Meanwhile, the absolute value of what is denoted by $\mathfrak{e}_1$ in (5.10) is no greater than $c_0r$ as can be seen using (5.11) and (3.12) also. The last contribution to $\mathfrak{q}$ comes from the term in (5.8) with the integral of $\mathfrak{m}(\hat{x}, \mathfrak{r})$. Since $|\hat{x}| \le c_0 r$ and $|\mathfrak{r}| \le c_0 |a||\nabla_A a|$, this contribution is no greater than

$$c_0 \frac{1}{r^2 \kappa^2} \Big( \int_{B_r} |a|^2 \Big)^{1/2} \Big( \int_{B_r} |\nabla_A a|^2 \Big)^{1/2} . \tag{5.13}$$

It is a consequence of (3.12) that the expression in (5.13) is less than what is written on the right hand side of (5.3).

<u>Step 6</u>: This step derives (5.6). To start, introduce S to denote the symmetric section of $T^*X \otimes T^*X$ with components $\{S_{\alpha\beta}\}_{\alpha,\beta\in\{1,2,3,4\}}$ that are given by

$$S_{\alpha\beta} = \langle (\nabla_{A,\alpha} a)_\nu (\nabla_{A,\beta} a)_\nu \rangle - \tfrac{1}{2}\, \delta_{\alpha\beta} |\nabla_A a|^2 . \tag{5.14}$$

Commute derivatives to see that the components $\{(\nabla^\dagger S)_\beta\}_{\beta\in\{1,2,3,4\}}$ of the 1-form $\nabla^\dagger S$ are given by the formula

$$(\nabla^\dagger S)_\beta = \langle (\nabla^\dagger \nabla a)_\nu (\nabla_{A,\beta} a)_\nu \rangle + \langle (\nabla_{A,\alpha} a)_\nu [F_{A,\alpha\beta}, a_\nu] \rangle + \mathcal{R}_{\sigma\nu\alpha\beta} \langle (\nabla_{A,\alpha} a)_\nu\, a_\sigma \rangle \tag{5.15}$$

with $\{\mathcal{R}_{\sigma\nu\alpha\beta}\}_{\sigma,\nu,\alpha,\beta\in\{1,2,3,4]}$ being the components of the Riemann curvature tensor.

Use (2.5) to write the term $\langle (\nabla^\dagger \nabla a)_\nu (\nabla_{A,\beta} a)_\nu \rangle$ that appears in (5.15) as

$$\tfrac{1}{2}\, r^2\, \nabla_\beta |a \wedge a|^2 + \mathrm{Ric}(\langle a \otimes (\nabla_{A,\beta} a) \rangle) . \tag{5.16}$$

The term that is denoted by $\mathfrak{r}$ in (5.6) is the sum of the term $\mathcal{R}_{\sigma\nu\alpha\beta} \langle (\nabla_{A,\alpha} a)_\nu\, a_\sigma \rangle$ from (5.15) and the term $\mathrm{Ric}(\langle a \otimes (\nabla_{A,\beta} a) \rangle)$ from (5.16).

Use (3.28) to write the term $\langle (\nabla_{A,\alpha} a)_\nu [F_{A,\alpha\beta}, a_\nu] \rangle$ from (5.15) as $r^{-2} \langle F_{A,\alpha\beta} (\nabla_{A,\nu} F)_{\nu\alpha} \rangle$ and then use the Bianchi identity to write the latter as $\tfrac{1}{2}\, r^{-2} \nabla_\beta |F_A|^2 - r^{-2} \nabla_\alpha \langle F_{A,\alpha\nu} F_{A,\beta\nu} \rangle$. The formula in (5.6) follows directly from (5.15), (5.16), the preceding definition of $\mathfrak{q}$, and this rewriting of $\langle (\nabla_{A,\alpha} a)_\nu [F_{A,\alpha\beta}, a_\nu] \rangle$.

**b) Approximate monotonicity where r $\kappa$(r) $r$ is large**

The lemma that follows implies that $N(r)$ can not be small where r is greater than a given number $r_0$ if $N(r_0)$ is large at $r_0$ and if $r_0 \kappa(r_0)\, r$ is not too small.

**Lemma 5.2**: *There exists* $\kappa > 1$ *that depends only on the geometric data and has the following significance: Suppose that* $r > \kappa$ *has been specified and that* $(A, a)$ *is a pair of* *connection on* P *and section of* $T^*X \otimes (P \times_{SO(3)} \mathfrak{S})$ *that obeys the equations in (2.11) for a given* $\tau \in [0, 1]$. *Fix* E *so that* $(A, a)$ *obey the bounds in (5.1). Given* $p \in U$, *define* N *and* K *using* p *and the pair* $(A, a)$. *Fix* $\mathrm{r}_0 \in (0, \kappa^{-1}]$ *and* $\mathrm{r}_1 \in [\mathrm{r}_0, \kappa^{-1}]$. *Then*

$$\mathrm{N}(\mathrm{r}_1) \geq \mathrm{N}(\mathrm{r}_0) - \kappa \mathrm{r}_1{}^2 \Big(1 + \frac{1}{\mathrm{r}_0^2\, \mathrm{K}(\mathrm{r}_0)^2 \mathit{r}^2} \mathrm{E}^2\Big) - \kappa \frac{1}{\mathrm{r}_0\, \mathrm{K}(\mathrm{r}_0) \mathit{r}} \mathrm{E} \sqrt{\mathrm{N}(\mathrm{r}_0)}\ .$$

Lemma 5.2 implies what is asserted by Proposition 5.1 when $\mathrm{r}_0 \mathrm{K}(\mathrm{r}_0) \mathit{r} \geq c_0 \mathrm{E}\, (1 + \mathrm{N}(\mathrm{r}_0)^{-1})$.

***Proof of Lemma 5.2***: Introduce by way of notation W to denote $F_A - \mathit{r}^2 a \wedge a$. It follows from (5.2) and 5.3) that N's derivative has the lower bound

$$\frac{d}{d\mathrm{r}} \mathrm{N} \geq - c_0 \frac{1}{\mathrm{r}^2 \mathrm{K}^2 \mathit{r}^2} \int_{\partial B_{\mathrm{r}}} (|\mathrm{W}|^2 + \mathit{r}^2 |a \wedge a||\mathrm{W}|) - c_0 \frac{1}{\mathrm{r}^2 \mathrm{K}^2 \mathit{r}^2}\, \mathrm{r} \int_{B_{\mathrm{r}}} |\mathrm{W}|^2 - c_0 \mathrm{r} (\mathrm{N} + 1)\ . \tag{5.17}$$

No generality is lost with regards to the proof by assuming that N is no greater than $\mathrm{N}(\mathrm{r}_0)$ on $[\mathrm{r}_0, \mathrm{r}_1]$. With this understood, (5.17) will be used only where N's derivative is negative. If $\mathrm{r} \in [\mathrm{r}_0, \mathrm{r}_1]$ and if N's derivative at r is negative then

$$\frac{1}{\mathrm{r}^2 \mathrm{K}^2 \mathit{r}^2} \int_{\partial B_{\mathrm{r}}} \mathit{r}^4\, |a \wedge a|^2 \leq \frac{1}{\mathrm{r}} \mathrm{N} (1 + \mathrm{N})\ , \tag{5.18}$$

this being a consequence of the identity in (5.4). Use this bound in (5.17) to see that N's derivative obeys

$$\frac{d}{d\mathrm{r}} \mathrm{N} \geq - c_0 \frac{1}{\mathrm{r}^2 \mathrm{K}^2 \mathit{r}^2} \int_{\partial B_{\mathrm{r}}} |\mathrm{W}|^2 - c_0 \frac{1}{\mathrm{r}^{3/2} \mathrm{K} \mathit{r}} (\mathrm{N}(1 + \mathrm{N}))^{1/2} \Big(\int_{\partial B_{\mathrm{r}}} |\mathrm{W}|^2\Big)^{1/2} - c_0 \frac{1}{\mathrm{r}^2 \mathrm{K}^2 \mathit{r}^2}\, \mathrm{r} \int_{B_{\mathrm{r}}} |\mathrm{W}|^2 - c_0 \mathrm{r} (\mathrm{N} + 1)\ . \tag{5.19}$$

Keep in mind that there is an $\mathrm{E}^2$ upper bound for the integral of $|\mathrm{W}|^2$ over any ball in U because (5.1) asserts an $\mathrm{E}^2$ upper bound for the integral of $|\mathrm{W}|^2$ over X. Also keep in mind that K is increasing and so $\mathrm{K}(\mathrm{r}) > \mathrm{K}(\mathrm{r}_0)$ if r is greater than $\mathrm{r}_0$. Use the preceding facts when integrating (5.19) to see that

$$\mathrm{N}(\mathrm{r}_1) \geq (1 - c_0 \mathrm{r}_1{}^2) \mathrm{N}(\mathrm{r}_0) - c_0 \mathrm{E} \Big(\int_{\mathrm{r}_0}^{\mathrm{r}_1} \frac{1}{\mathrm{r}^3 \mathrm{K}(\mathrm{r})^2 \mathit{r}^2} \mathrm{N}(\mathrm{r})(1 + \mathrm{N}(\mathrm{r}))\, d\mathrm{r}\Big)^{1/2} - c_0 \mathrm{r}_1{}^2 \Big(1 + \frac{1}{\mathrm{r}_0^2\, \mathrm{K}(\mathrm{r}_0)^2 \mathit{r}^2} \mathrm{E}^2\Big). \tag{5.20}$$

The explicit integral on the right hand side of (5.20) can be bounded using the identity in (3.6) with the assumption that $\mathrm{N}(\mathrm{r})$ is less than $\mathrm{N}(\mathrm{r}_0)$ on $[\mathrm{r}_0, \mathrm{r}_1]$. In particular, it follows from the latter bound and (3.6) that the integrand for this integral is no greater than $-\mathrm{N}(r_0)$ times the derivative of the function $r \to \frac{1}{r^2\kappa^2 r^2}$. This being the case,

$$\int_{r_0}^{r_1} \frac{1}{r^3\kappa(r)^2 r^2}\mathrm{N}(r)dr \le c_0 \frac{1}{r_0^2\kappa(r_0)^2 r^2}\mathrm{N}(r_0). \tag{5.21}$$

Use this last bound with (5.20) to obtain the bound that is asserted by the lemma.

### c) Uhlenbeck's theorem and elliptic regularity

Uhlenbeck's compactness theorem plays a part in Proposition 5.1's proof. Appeals to Uhlenbeck's theorem are also made Section 6. The version used here is stated below in (5.22). The assertions in (5.22) summarize parts of Corollary 2.2 and Lemma 2.5 in [U]. The notation takes $\theta_0$ to be the product connection on a given product principal SO(3) bundle.

UHLENBECK'S THEOREM: *There exists* $\kappa_U > 1$ *with the following significance: Let* $B \subset U$ *denote a given ball of radius at most* $\kappa_U^{-1}$. *Suppose that* A *is a connection on* $P|_B$ *with the integral of* $|F_A|^2$ *over* B *being at most* $\kappa_U^{-2}$. *There exists an isomorphism from* $B \times SO(3)$ *to* P *that pulls* A *back as* $\theta_0 + \hat{a}_A$ *with* $\hat{a}_A$ *denoting here a* $\mathfrak{S}$ *valued 1-form on* B *with the properties listed below.*

- *The 1-form* $\hat{a}_A$ *is coclosed, thus* $d{*}\hat{a}_A = 0$.
- *The 3-form* ${*}\hat{a}_A$ *pulls back as zero to the boundary of* B.
- *The Sobolev* $L^2_1$ *norm of* $\hat{a}_A$ *on* B *is bounded by* $\kappa_U$ *times the* $L^2$ *norm of* $F_A$ *on* B.

*Conversely, if* $A = \theta_0 + \hat{a}_A$ *is a connection on* $B \times G$ *such that* $\hat{a}_A$ *obeys the first two bullets and has* $L^2_1$ *norm on* B *bounded by* $\frac{1}{2}\kappa_U^{-1}$, *then the third bullet is also obeyed.*

(5.22)

Proposition 5.1 refers to a number $z_U$. This $z_U$ is defined to be $100\,\kappa_U$.

Uhlenbeck's theorem is used subsequently to extract convergent subsequences from sequences of solutions to $r = 1$ versions of (2.11) on balls in $\mathbb{R}^4$ when certain conditions are met; one condition being an a priori $L^2$ norm bound on the corresponding sequence of curvatures.

### d) Proof of Proposition 5.1

The proof that follows uses the assumption that Proposition 5.1 is false to derive nonsense. To this end, suppose that there exists $\varepsilon \in (0, 1]$ with no $\kappa_\varepsilon$ as described by the proposition. If this is so, then there exists a sequence $\{r_n\}_{n=1,2,\ldots} \subset (1, \infty)$, a convergent sequence $\{\tau_n\}_{n=1,2,...} \subset [0, 1]$, a sequence $\{(A_n, a_n)\}_{n=1,2,...}$ of pairs consisting of a connection on P and section of $T^*X \otimes (P \times_{SO(3)} \mathfrak{G})$, a sequence of points $\{p_n\}_{n=1,2,...} \subset U$ and two sequences $\{r_{0n}\}_{n=1,2,...} \subset (0, c_0^{-1}]$ and $\{r_{1n}\} \subset (0, c_0^{-1}]$ with the former having limit zero and the latter described below. The members of these sequences have the following properties: Fix $n \in \{1, 2, \ldots\}$. The pair $(A_n, a_n)$ obeys the version of (2.11) defined by $r_n$ and $\tau_n$ and it obeys (5.1) with an n-dependent constant $E$ that is less than $\varepsilon^{-1}$. The number $r_{0n}$ is such that $N_n(r_{0n}) \geq \varepsilon$ with $N_n$ denoting the version of $N$ that is defined by $p_n$ and $(A_n, a_n)$. Meanwhile, the number $r_{1n}$ is greater than $r_{0n}$ and it obeys $N_n(r_{1n}) \leq \frac{1}{n}\varepsilon$ and $r_{1n} K_n(r_{1n}) r \geq z_U^{-2}$ with $K_n$ being the version of the function $K$ given by $p_n$ and $(A_n, a_n)$. One other condition is imposed, this stated directly in (5.23). The notation used in (5.23) and elsewhere has $r_{*n}$ denoting the version of the number $r_{*p}$ that is defined by $p_n$.

$$\int_{\mathrm{dist}(\cdot, p_n) < \frac{1}{2} r_{*n}} |F_{A_n}|^2 \leq \frac{1}{10{,}000} \kappa_U^{-2} , \tag{5.23}$$

with $\kappa_U$ being the number that appears in Uhlenbeck's theorem.

Nonsense follows directly from Lemma 5.2 if there is a subsequence of positive integers with each member of the corresponding subsequence from $\{r_{0n} K_n(r_{0n}) r_n\}_{n=1,2,...}$ being greater than $c_0 \varepsilon^{-3}$. This understood, nothing is lost by assuming that $r_{0n} K_n(r_{0n}) r_n$ is bounded by $c_0 \varepsilon^{-3}$ for each index n. The third bullet of the next lemma generates the desired nonsense when this bound is assumed.

**Lemma 5.3**: *Assume that* $\{r_n, \tau_n, (A_n, a_n), r_{0n}, r_{1n}\}_{n \in \{1,2,3...\}}$ *obey the conditions stated in the preceding paragraphs except for the condition in (5.23). Then there exists a subsequence* $\Xi \subset \{1, 2, \ldots\}$ *such that*

- $\lim_{n \to \infty} N_n(r_{*n}) = 0$ .
- $\lim_{n \in \Xi} r_{*n}^{-1} r_{0n} = 0$ .
- $\lim_{s \to 0} (\lim_{n \in \Xi} \int_{B_{s r_{*n}}} |F_{A_n}|^2 ) > \frac{1}{10} \kappa_U^{-2}$ .

Either the third bullet of this lemma holds or (5.23) hold, but they can't both hold. The assertion that both are true is the nonsense that proves Proposition 5.1.

***Proof of Lemma 5.3***: The proof has four parts.

*Part 1*: Fix $n \in \{1, 2, \ldots\}$, choose a Gaussian coordinate chart centered at $p_n$ and let $\phi_n$ denote the map from the $|x| \leq c_0^{-1} r_{*n}^{-1}$ ball in $\mathbb{R}^4$ to $B_{1/c_0}$ that is obtained by composing first the rescaling $x \to r_{*n} x$ of $\mathbb{R}^4$ and then the Gaussian coordinate chart map.

Define $A_{*n}$ to be $\phi_n{}^*A_n$ and $a_{*n}$ to be $r_{*n}{}^{-1}K_n(r_{*n})^{-1}\phi_n{}^*a_n$. The definitions are such that the pair $(A_{*n}, a_{*n})$ obey the $r = 1$ and $\tau = \tau_{*n}$ version of (2.11) on the $c_0{}^{-1}r_{*n}{}^{-1}$ ball about the origin in $\mathbb{R}^4$. Note however that the Riemannian metric that is used in this version of (2.11) is not the Euclidean metric, it is the product of $r_{*n}{}^{-2}$ times the metric that is pulled back from X by $\phi_n$. Even so, the latter metric differs from the Euclidean metric at any given $x \in \mathbb{R}^4$ with norm less than $c_0{}^{-1}r_{*n}{}^{-1}$ by $c_0|x|^2$, its derivative has norm at most $c_0|x|$ and the norms of its derivatives to any given higher order $k \in \{2, 3, \ldots\}$ are bounded by $c_k r_{*n}{}^{k-2}$ with $c_k$ being a number that depends only on the geometric data.

The notation used below has $K_{*n}$ and $N_{*n}$ denoting the versions of $K$ and $N$ that are defined by the origin and the pair $(A_{*n}, a_{*n})$. The definitions are such that $K_{*n}$ and $N_{*n}$ at any $s \in (0, c_0{}^{-1}r_{*n}{}^{-1})$ obey $K_{*n}(s) = K_n(r_{*n})^{-1}K_n(s r_{*n})$ and $N_{*n}(s) = N_n(s r_{*n})$. Note in particular that $K_{*n}(1) = 1$.

Fix $T \in [1, \infty)$ for the moment and suppose that there exists $s \in [1, c_0{}^{-1}r_{*n}{}^{-1})$ such that $N_{*n}(s) \geq (s K_{*n})^{-2}T$. Keeping in mind the identity $N_{*n}(s) = N_n(s r_{*n})$, an appeal to Lemma 5.2 proves that

$$N_n(r) \geq c_0{}^{-1}(s K_{*n})^{-2}T - c_0(r^2 r_{*n}{}^2 + (s K_{*n})^{-2}E^2) \quad \textit{for } r \in [sr_{*n}, c_0{}^{-1}] \,. \tag{5.24}$$

With (5.24) understood, let $s_{*n}$ denote the smallest of the number $s \in [1, c_0{}^{-1}r_{*n}{}^{-1})$ where $N_{*n}(s) \geq (s K_{*n})^{-2}T$. Since $s^2 K_{*n}(s)^2 N_{*n}(s) \leq T$ for all $s \in [1, s_{*n}]$ and since $K_{*n}(1) = 1$, integrating the $K_{*n}$ and $N_{*n}$ version of (3.6) leads to the inequality

$$K_{*n}(s) \leq 1 + \tfrac{1}{2}(1 - s^{-2})T \quad \textit{for } s \in [1, s_{*n}]. \tag{5.25}$$

Note in particular that the $T = c_0 E^2$ versions of (5.24) and (5.25) imply that

$$N_n(r) \geq c_0{}^{-1}s_{*n}{}^{-2}E^2(1 + E^2)^{-1} \quad \textit{when} \;\; r \geq s_{*n}r_{*n}. \tag{5.26}$$

To make something of (5.26), write $r_{1n}$ as $s_{1n}r_{*n}$. Note that $s_{1n} \geq 1$ because $K_n$ is increasing and the assumptions of the lemma posit that $r_{1n}K_n(r_{1n})\, r \geq 1$. If $s_{*n} \leq s_{1n}$, then (5.26) requires $s_{*n}{}^2 \geq \varepsilon^{-1} n E^2(1 + E^2)^{-1}$ because $N_{*n}(s_n)$ is $N_n(r_{1n})$ and the latter is at most $\frac{1}{n}\,\varepsilon$.

*Part 2*: It proves useful to distinguish between the cases when $\{r_{*n}{}^{-1}r_{1n}\}_{n\in\{1,2,\ldots\}}$ has a bounded subsequence and when it does not. Assume in this Part 2 that there is a bounded subsequence. The six steps that follow prove the assertion in Lemma 5.3 granted such a bound. The notation is such that this subsequence is henceforth relabled by consecutive integers starting at 1. The notation in these steps uses $s_*$ to denote an upper bound for $\{r_{*n}{}^{-1}r_{1n}\}_{n\in\{1,2,\ldots\}}$.

Step 1: Fix n and let $s_{1n}$ again denote $r_{*n}^{-1} r_{1n}$. Use what is said in the final paragraph of Part 2 to see that $N_n(s)$ for $s \in [1, s_{1n}]$ is bounded by $c_0 (s K)^{-2} E^2 (1+E^2)^{-1}$. Given this bound, the derivation of (5.25) can be repeated to bound $K_{*n}(s_{1n})$ by $1+c_0 E^2$. Use this bound on $K_{*n}(s_{1n})$ with the definition of $N_{*n}$ to conclude that

$$\int_{|x|\le s_{1n}} (|\nabla_{A_{*n}} a_{*n}|^2 + 2|a_{*n} \wedge a_{*n}|^2) \le c_0 (1+E^2) s_*^2 \tfrac{1}{n} \varepsilon \ . \tag{5.27}$$

Let $\Xi$ denote for the moment any given subsequence of positive integers. The fact that (5.27) holds for each $n \in \{1, 2,\ldots\}$ and the fact that $K_{*n}(1) = 1$ for each $n \in \{1, 2, \ldots\}$ leads directly to the top bullet in Lemma 5.2, this being that $\lim_{n\in\Xi} N_n(r_{*n}) = 0$.

Step 2: The bound in (5.27) implies in particular a $c_0 (1+E^2) s_*^2 \frac{1}{n} \varepsilon$ bound for the integral of $r_n^4 |a_n \wedge a_n|^2$ over the ball of radius $r_{1n}$ centered on p. The latter with (5.1) implies in turn an n-independent upper bound on the integral of $|F_{A_n}|^2$ over this same ball which when written using $F_{A_{*n}}$ asserts $\int_{|x|\le s_{1n}} |F_{A_{*n}}|^2 \le c_0 (1+E^2)$ .

Step 3: The $r = 1$ version of (2.11) are uniformly elliptic on any given ball in $\mathbb{R}^4$. With this understood, Uhlenbeck's theorem in (5.23) with the a priori bound in (5.27), the a priori bound in Step 2 and the $1+c_0 E^2$ bound on $\{K_{*n}(s_{1n})\}_{n\in\{1,2,\ldots\}}$ can be used to obtain a limit of sorts on the $|x| < 1$ ball from the sequence $\{(A_{*n}, a_{*n})\}_{n\in\{1,2,\ldots\}}$. The limit is a pair of connection on the product principal SO(3) bundle over the $|x| < 1$ ball and $\mathfrak{G}$-valued 1-form on this ball that obey some $\tau \in [0, 1]$ version of (2.11) with the Hodge star coming from the Euclidean metric. This pair is denoted by $(A_*, a_*)$. The manner of convergence to this limit is as follows: There exists a finite subset of the $|x| < 1$ ball to be denoted by $\Theta$, a subsequence of positive integers to be denoted by $\Xi$, and a corresponding sequence $\{g_n\}_{n\in\Xi}$ with n'th member being an isomorphism from the product principal bundle over the complement of $\Theta$ in the $|x| < 1$ ball to $\phi_n{}^*P$. This data is such that the sequence $\{|a_{*n}|\}_{n\in\Xi}$ converges weakly to $|a_*|$ in the $L^2_1$ topology on the $|x| \le 1$ ball; and it is such that $\{(g_n{}^*A_{*n}, g_n{}^*a_{*n})\}_{n\in\Xi}$ converges strongly to $(A_*, a_*)$ in the $C^\infty$ topology on compact subsets in the $|x| < 1$ ball that do not contain points in $\Theta$. The set $\Theta$ is characterized as follows: If $x \in \Theta$, then

$$\lim_{s\to 0} (\lim_{n\in\Xi} \int_{|x-(\cdot)|\le s} |F_{A_{*n}}|^2 ) \ge \tfrac{1}{10} \kappa_U^{-2} \tag{5.28}$$

with $\kappa_U$ being the number from Uhlenbeck's theorem in (5.22).

Step 4: The $\mathfrak{G}$-valued 1-form $a_*$ is $A_*$ covariantly constant and $a_* \wedge a_* = 0$. This follows from the $n \in \Xi$ versions of (5.27) and the fact that $\{|\nabla_{A_{*n}} a_{*n}|^2\}_{n\in\Xi}$ and $\{|a_n \wedge a_n|^2\}_{n\in\Xi}$ converge in the $C^\infty$ topology on compact sets in the complement of $\Theta$. Note that $a_*$ is not zero. In fact, its norm is $\frac{1}{\sqrt{2\pi}}$. The proof that this is so invokes three facts, the first being that $\kappa_{*n}(1) = 1$ and the second being that $\{|a_{*n}|\}_{n\in\Xi}$ converges weakly in $L^2{}_1$ to $|a_*|$. The third fact is that the restriction map from the space of smooth functions on the $|x| \le 1$ ball to space of smooth functions on the $|x| = 1$ sphere extends as a compact map from the Sobolev space of $L^2{}_1$ functions on the $|x| < 1$ ball to the space of $L^2$ functions on the $|x| = 1$ sphere.

The fact that $|a_*| = \frac{1}{\sqrt{2\pi}}$ and the aforementioned weak $L^2{}_1$ convergence of $\{|a_n|\}_{n\in\Xi}$ to $|a_*|$ imply in turn that $\lim_{n\in\Xi} \kappa_{*n}(s) = 1$ for all $s \in (0, 1]$.

Step 5: This step proves the second bullet in Lemma 5.3, this being the assertion that $\lim_{n\in\Theta} r_{*n}^{-1} r_{0n} = 0$. To this end, fix $n \in \Xi$ and write $r_{0n}$ as $s_{0n} r_{*n}$. It follows from (5.27) that $N_{*n} \le c_0 s_*^2 \frac{1}{n} E^2$ on the interval $[1, s_{1n}]$, and since $r_{0n} < r_{1n}$ and $N_n(r_{1n}) \ge \varepsilon$, this implies that $s_{0n} < 1$ when n is large.

To see that $\{s_{0n}\}_{n\in\Xi}$ must have limit zero, use the fact that $N_n(r_{0n}) \ge \varepsilon$ with (5.4) to see that $N_n(2r_{0n})$ is greater than $c_0^{-1}\varepsilon$ when n is large. This fact with (3.6) implies that $\kappa_{*n}(2s_{0n})$ is greater than $e^{\varepsilon/c_0} \kappa_{*n}(s_{0n})$ when n is large. With the preceding understood, suppose for the sake of argument that the lim-sup of the sequence $\{s_{0n}\}_{n\in\Xi}$ is positive. If such is the case, then it follows from what is said at the end of the preceding step that there exists large values of n in $\Xi$ with both $\kappa_{*n}(s_{0n})$ and $\kappa_*(2s_{0n})$ being very close to 1. But both can't be nearly 1 because $\kappa_{*n}(2s_{0n}) \ge e^{\varepsilon/c_0} \kappa_{*n}(s_{0n})$ when n is large

Step 6: This final step proves the third bullet in Lemma 5.3. This assertion follows from (5.28) if the origin is a point of $\Theta$. To prove that this is so, assume to the contrary that the origin is not in $\Theta$. If such is the case, then $\{|a_{*n}|\}_{n\in\Xi}$ converges in the $C^\infty$ topology on a neighborhood of 0 to the constant function $\frac{1}{\sqrt{2\pi}}$. This convergence implies in particular that the function $\kappa_{*n}(\cdot)$ is very close to 1 on some interval containing 0 when n is large. But the latter conclusion is nonsense because it runs afoul of the fact that $\kappa_{*n}(2s_{0n}) \ge e^{\varepsilon/c_0} \kappa_{*n}(s_{0n})$ when n is large.

*Part 3*: Suppose in this Part 3 that $\{r_{*n}^{-1} r_{1n}\}_{n\in\{1,2\ldots\}}$ has no bounded subsequences. It follows as a consequence of what is said in the last paragraph of Part 1 that there exists an increasing sequence $\{R_n\}_{n\in\{1,2,\ldots\}} \subset [1, \infty)$ with no bounded subsequences such that

$$\lim_{n\to\infty} \sup_{s\in[1,R_n]} s^2 \kappa_{*n}(s)^2 N_{*n}(s) \le c_0 (1+E^2) \; . \tag{5.29}$$

This with the definition of N implies that

$$\lim_{n\to\infty} \sup_{s\in[1,R_n]} \int_{|x|\le s} (|\nabla_{A_{*n}} a_{*n}|^2 + 2|a_{*n} \wedge a_{*n}|^2) \le c_0 E^2 \; . \tag{5.30}$$

What with (5.1), the bound from (5.30) on the integral of $|a_{*n} \wedge a_{*n}|^2$ implies in turn that

$$\lim_{n\to\infty} \sup_{s\in[1,R_n]} \int_{|x|\le s} |F_{A_{*n}}|^2 \le c_0 E^2 \; . \tag{5.31}$$

Meanwhile, (5.25) leads to the bound

$$\lim_{n\to\infty} \sup_{s\in[1,R_n]} \kappa_{*n} \le 1 + c_0 E^2 \; . \tag{5.32}$$

The $r = 1$ versions of (2.11) are uniformly elliptic on the $|x| \le c_0 r_{*n}^{-1}$ ball and so the bounds in (5.30)-(5.32) can be used in conjunction with Uhlenbeck's compactness theorem in (5.22) to obtain a limit from the sequence $\{(A_{*n}, a_{*n})\}_{n\in\{1,2,\ldots\}}$. The limit in this case is a pair of connection on the product principal SO(3) bundle on the whole of $\mathbb{R}^4$ and $\mathfrak{S}$-valued 1-form on the whole of $\mathbb{R}^4$ that obey some $\tau \in [0, 1]$ version of (2.11) with the Hodge star coming from the Euclidean metric. This pair is again denoted by $(A_*, a_*)$. The manner of convergence is essentially identical to that in Step 3 of Part 2. As before, there is a finite set $\Theta \subset \mathbb{R}^4$, a subsequence in $\{1, 2, \ldots\}$ to be denoted again by $\Xi$, and a corresponding sequence $\{g_n\}_{n\in\Xi}$ with n'th member being an isomorphism from the product principal bundle over the complement of $\Theta$ over the $|x| \le c_0^{-1} r_{*n}$ ball to $\phi_n{}^*P$. This data is such that $\{|a_{*n}|\}_{n\in\Xi}$ converges weakly in the $L^2{}_1$ topology on compact subsets of $\mathbb{R}^4$ to $|a_*|$; and it is such that $\{(g_n{}^*A_{*n}, g_n{}^*a_{*n})\}_{n\in\Xi}$ converges strongly in the $C^\infty$ topology on compact subsets of $\mathbb{R}^4 - \Theta$ to $(A_*, a_*)$. The set $\Theta$ is again characterized by (5.28).

The a priori bounds in (5.30) and (5.31) imply that the three functions $|F_{A_*}|^2$ and $|\nabla_{A_*} a_*|^2$ and $|a_* \wedge a_*|^2$ have finite integral on $\mathbb{R}^4$. The fact that $\kappa_{*n}(1) = 1$ for all $n \in \Xi$ implies that $|a_*|^2$ has integral equal to 1 on the $|x| = 1$ sphere and so $a_*$ is not identically zero. In any event, the bound in (5.32) leads to an a priori bound for $\sup_{s\in(0,\infty)} s^{-3} \int_{|x|=s} |a_*|^2$ .

As it turns out, the preceding conditions can be satisfied only in the case when the connection $A_*$ is flat and $a_*$ is $A_*$-covariantly constant. This is stated formally by **Lemma 5.4**: *Let* (A, $a$) *denote a pair of connection on the product principal* SO(3) *bundle over* $\mathbb{R}^4$ *and section of the associated product* $\mathfrak{G}$ *bundle that obeys some* $\tau \in [0, 1]$ *version of (2.11). If* $|F_A|^2$ *and* $|\nabla_A a|^2$ *and* $|a_* \wedge a_*|^2$ *have finite integral on* $\mathbb{R}^4$ *and if the average value of* $|a|^2$ *on the* $|x| = 1$ *sphere is bounded, then the connection* A *is flat and* $a$ *is* A-*covariantly constant and* $a \wedge a = 0$.

This lemma is proved in Section 5e. Accept it as true for now.

The fact that the sequence $\{(A_{*n}, a_{*n})\}_{n\in\Xi}$ converges in the $C^\infty$ topology on compact subsets of $\mathbb{R}^4-\Theta$ to $(A_*, a_*)$ and the fact that the norm of $a_*$ is constant implies that this norm is equal to $\frac{1}{\sqrt{2\pi}}$ . The argument is the same as that used in Step 4 of Part 2. The fact that $|a_*| = \frac{1}{\sqrt{2\pi}}$ and the $L^2_1$ weak convergence of $\{|a_{n*}|\}_{n\in\Xi}$ to $|a_*|$ implies convergence of the sequence of functions $\{\kappa_{*n}\}_{n\in\Xi}$ to 1 in the $C^\infty$ topology on compact subsets of $(0, \infty)$. This in turn implies that $\lim_{n\in\Xi} N_{*n}(1)$ must be zero. This is the assertion in the top bullet of Lemma 5.2. The argument is identical but for cosmetics to the argument used Step 5 of Part 3. The argument in Step 5 of Part 2 can be repeated verbatim to prove the second bullet in Lemma 5.2 and the argument in Step 6 of Part 2 can be repeated verbatim to prove the third bullet in Lemma 5.2.

### e) Proof of Lemma 5.4

The proof of this lemma has eleven parts. The proof uses $\kappa$ to denote the non-negative function on $[0, \infty)$ whose square is defined by the rule

$$r \to \kappa^2(r) = r^{-3} \int_{|x|=r} |a|^2 . \tag{5.33}$$

The function N on $(0, \infty)$ is defined by (3.5). The pair $\kappa$ and N obey (3.6).

*Part 1*: Let H denote the integral of $|\nabla_A a|^2 + 2|a\wedge a|^2$ over $\mathbb{R}^4$. Since $\kappa$ is increasing, the function N where $r \geq 1$ is no greater than $\frac{1}{\kappa(1)^2 r^2}$ H. Granted this bound, integrate (3.6) from $r = 1$ to values of r greater than 1 to see that $\kappa$ is bounded on $(0, \infty)$. Since $\kappa$ is also increasing, it has an $r \to \infty$ limit. Denote this limit by $\kappa_\ddagger$. Since the integral of $|\nabla_A a|^2 + 2|a\wedge a|^2$ over a large radius ball centered at the origin is nearly equal to H and since $\kappa \leq \kappa_\ddagger$, it follows that $N \geq c_0^{-1} H \frac{1}{\kappa_\ddagger^2 r^2}$ where r is large. Granted that this is so, integrate (3.6) again to see that

$$\kappa(r) \leq (1 - \tfrac{1}{2} c_0^{-1} H \frac{1}{\kappa_\ddagger^2 r^2})\kappa_\ddagger \quad \textit{when } r \textit{ is large}.$$

(5.34)

The subsequent parts of the proof derive a converse bound asserting that

$$\kappa(r) \geq (1 - c\,\tfrac{1}{r^3}\,)\,\kappa_{\ddagger}$$

(5.35)

with c being a positive constant. The preceding bound is not compatible with the bound in (5.34) unless $H = 0$. This is what Lemma 5.4 asserts.

*Part 2*: The function $|a|$ has square integrable differential on $\mathbb{R}^4$ because $|\nabla_A a|^2$ is integrable on $\mathbb{R}^4$. It follows as a consequence that $|a| - \frac{1}{\sqrt{2\pi}}\kappa_{\ddagger}$ is an $L^4$ function on $\mathbb{R}^4$ with $L^4$ norm bounded by $c_0 H^{1/2}$. By way of a reminder, if $f$ is any given function on $\mathbb{R}^4$ with square integrable derivatives, then there exists a number, $f_\infty$, such that $f - f_\infty$ is an $L^4$ function. Note by the way that $|a|^2$ obeys (2.6) and so the maximum principle can be invoked to conclude that $|a| \leq \frac{1}{\sqrt{2\pi}}\kappa_{\ddagger}$ on the whole of $\mathbb{R}^4$ with equality if and only if $|a|$ is constant (and thus $H = 0$). The preceding a priori bound for the norm of $a$ justifies in part the introduction and subsequent use of $a_{\ddagger} = \kappa_{\ddagger}^{-1} a$.

The $r = 1$ version of the equations in (2.11) with Uhlenbeck's theorem in (5.22) lead to pointwise bounds on the covariant derivatives of $a$ and on the norm of $F_A$. In particular, the fact that $|\nabla_A a_{\ddagger}|$ and $|F_A|$ are $L^2$ function on $\mathbb{R}^4$ and the fact that $|a|$ is a priori bounded lead via standard elliptic regularity theorems to bounds of the form

$$\sup_{x\in\mathbb{R}^4} |\nabla_A a_{\ddagger}| \leq z \quad \textit{and} \quad \sup_{x\in\mathbb{R}^4} |F_A| \leq z$$

(5.36)

with $z$ being a positive number. The number $z$ is determined by a positive lower bound for the $x \in \mathbb{R}^4$ versions of the number $r_{cF}$ as defined in (5.2) using $c = c_0\,\kappa_U$ with $\kappa_U$ being the constant that appears in (5.22). The set of $x \in \mathbb{R}^4$ versions of $r_{cF}$ has a positive lower bound because the integral over $\mathbb{R}^4$ of $|F_A|$ is finite.

*Part 3*: It is also the case that $\nabla_A(\nabla_A a_{\ddagger})$ has finite $L^2$ norm on $\mathbb{R}^4$, its square being less than $c_0(z+1)\kappa_{\ddagger}^{-2} H$. To see that this is so, repeat the derivation of (3.16) to see that $\nabla_A a_{\ddagger}$ obeys the $r = 1$ version of the equation therein with the endomorphisms $\mathcal{R}^1$ and $\mathcal{R}^0$ set equal to zero. Take the inner product of both sides of the latter version of (3.16) with $\nabla_A a_{\ddagger}$ to and use the $r = 1$ version of (3.28) to obtain the inequality

$$\tfrac{1}{2}\, d^{\dagger} d\, |\nabla_A a_{\ddagger}|^2 + |\nabla_A(\nabla_A a_{\ddagger})|^2 \leq c_0(|F_A| + 1)\,|\nabla_A a_{\ddagger}|^2 \,.$$

(5.37)

Fix R >> 1 and multiply both sides of (5.37) by the function $\chi(1 - R^{-1}|x|)$. Integrate the resulting inequality and then integrate by parts to remove derivatives from $|\nabla_A a_\ddagger|^2$. Use the bound for $|F_A|$ in (5.36) to then obtain an inequality that bounds the integral of $|\nabla_A(\nabla_A a_\ddagger)|^2$ over the $|x| \le R$ sphere by $c_0(z+1)$ times that of $|\nabla_A a_\ddagger|^2$ over the $|x| \le 2R$ ball and thus by $c_0(z+1)\kappa_\ddagger^{-2}H$. Granted this last bound, then the fact that $|\nabla_A(\nabla_A a_\ddagger)|^2$ is integrable on the whole of $\mathbb{R}^4$ with integral being less than $c_0(z+1)\kappa_\ddagger^{-2}H$ follows directly.

*Part 4*: The fact that $|a_\ddagger| - \frac{1}{\sqrt{2\pi}}$ is an $L^4$ function and the bound in (5.36) for $|\nabla_A a_\ddagger|$ imply that the $|a_\ddagger|$ converges pointwise to $\frac{1}{\sqrt{2\pi}}$ as $|x|$ gets ever larger. This is to say that

$$\lim_{|x|\to\infty} (|a_\ddagger| - \tfrac{1}{\sqrt{2\pi}}) = 0 . \tag{5.38}$$

To prove (5.38), fix $m \ge 1$ and suppose that $x \in \mathbb{R}^4$ is a point where $|a_\ddagger| \le (1 - m^{-1})\frac{1}{\sqrt{2\pi}}$. Then $|a_\ddagger| \le (1 - \frac{1}{2}m^{-1})\frac{1}{\sqrt{2\pi}}$ on the radius $c_0 z^{-1} m^{-1}$ ball centered at x and so the integral of $(|a_\ddagger| - (1 - m^{-1})\frac{1}{\sqrt{2\pi}})^4$ over this ball is greater than $c_0(z^{-1}m^{-2}\kappa_\ddagger)^4$. It follows from the latter observation that there can be at most $c_0(zm^2)^4H$ disjoint balls in $\mathbb{R}^4$ of radius $z^{-1}m^{-1}$ in $\mathbb{R}^4$ that contain points where $|a_\ddagger|$ is less than $(1 - m^{-1})\frac{1}{\sqrt{2\pi}}\kappa_\ddagger$.

Given that $|a \wedge a|$ is an $L^2$ function on $\mathbb{R}^4$ and that $|a_\ddagger|$ is bounded, much the same argument proves that

$$\lim_{|x|\to\infty} |a_\ddagger \wedge a_\ddagger| = 0. \tag{5.39}$$

*Part 5*: It is also the case that

$$\lim_{|x|\to\infty} |\nabla_A a_\ddagger| = 0 . \tag{5.40}$$

To prove this last claim, let x denote a given point in $\mathbb{R}^4$ and let $\chi_{(x)}$ denote the function $\chi(|x - (\cdot)| - 2)$. Multiply both sides of (5.37) by the function $\frac{1}{2\pi^2}\frac{1}{|x-(\cdot)|^2}\chi_{(x)}$. The latter function has compact support in the radius 2 ball centered at x and is the Green's function for the operator $d^\dagger d$ on the ball of radius 1 centered at x. Integrate the resulting inequality over the ball of radius 2 centered at x. Then integrate by parts to see that

$$|\nabla_A a_\ddagger|^2(x) \le c_0(z+1) \int_{|x-(\cdot)|<2} \frac{1}{|x-(\cdot)|^2} |\nabla_A a_\ddagger|^2 . \tag{5.41}$$

Fix $\rho \in (0, 1)$ for the moment and then break the integral in (5.41) into the contribution from the complement of the radius $\rho$ ball centered at x and the interior of this ball. The former is bounded by $c_0\rho^{-2}$ times the integral of $|\nabla_A a_{\ddagger}|^2$ over the radius 2 ball centered at x; and it follows from (5.36) that the latter is bounded by $c_0\rho^2 z^2$. This understood, take $\rho$ to be the minimum of $\frac{1}{2}$ and the $L^2$ norm of $|\nabla_A a_{\ddagger}|^2$ on the radius 2 ball centered at x. The resulting inequality asserts that

$$|\nabla_A a_{\ddagger}|^2(x) \le c_0 (z+1)^3 \Big( \int_{|x-(\cdot)|<2} |\nabla_A a_{\ddagger}|^2 \Big)^{1/2} . \tag{5.42}$$

The assertion in (5.40) follows from (5.42) because $|\nabla_A a_{\ddagger}|^2$ is square integrable.

*Part 6*: Fix $c > 1000$ and it follows from (5.38) and (5.39) that there exists $t_c > 1$ such that $|a_{\ddagger}| > (1 - \frac{1}{10} c^{-1}) \frac{1}{\sqrt{2\pi}}$ and $|a_{\ddagger} \wedge a_{\ddagger}| \le \frac{1}{1000} c^{-2}$ on the $|x| \ge t_c$ part of $\mathbb{R}^4$. This implies that the endomorphism $\mathbb{T}$ as defined in (4.2) has a unique eigenvalue greater than $c^{-2}$ at each point where $|x| \ge t_c$ . Moreover, the one large eigenvalue is larger than $\frac{1}{2\pi^2}(1 - c^{-1})$. Let $\sigma_{\ddagger}$ denote the norm 1 section of the product $\mathfrak{S}$ bundle over the $|x| \ge t_c$ part of $\mathbb{R}^4$ that restricts to each point as an eigenvector of $\mathbb{T}$ with largest eigenvalue. With $\sigma_{\ddagger}$ so defined, decompose $a_{\ddagger}$ as done in (4.4) on the $|x| \ge t_c$ part of $\mathbb{R}^4$. The 1-form $\nu$ and the $\mathfrak{S}$-valued 1-form $\mathfrak{a}$ obey the following analog of (4.5) where $|x|$ is very much larger than $t_c$:

- $(1 - c^{-1}) \frac{1}{\sqrt{2\pi}} \le |\nu| \le \frac{1}{\sqrt{2\pi}}$ .
- $|\mathfrak{a}| \le c^{-1}$ *and* $\lim_{|x|\to\infty} |\mathfrak{a}| = 0$.
- *The metric pairing between* $\nu$ *and* $\mathfrak{a}$ *is zero.*
- $|a_{\ddagger} \wedge a_{\ddagger}|^2 = 4|\nu|^2|\mathfrak{a}|^2 + |\mathfrak{a} \wedge \mathfrak{a}|^2$.
- *There is a constant 1-form, this denoted by e, with norm* $\frac{1}{\sqrt{2\pi}}$ *and such that* $|\nu - e| \le c^{-1}$ *and* $\lim_{|x|\to\infty} |\nu - e| = 0$.

(5.43)

The paragraph that follows explains why the assertions in (5.43) are true.

The bounds in the first two bullets follow from (5.38) and (5.39). The arguments for the third and fourth bullets are the same as those for the third and fourth bullets of (4.5). The proof of the fifth bullet uses the fact that the derivatives of $\nu$, $\sigma_{\ddagger}$ and $\mathfrak{a}$ obey (4.7) and in particular, the fourth bullet of (4.7). This fourth bullet implies that $\nabla\nu$ is square integrable on the complement of a large radius sphere in $\mathbb{R}^4$. This being the case, then $\nu$ differs from a constant 1-form by an $L^4$ function. This constant 1-form is $e$. The argument for (5.38) can be repeated with only notational changes to prove that $\lim_{|x|\to\infty} |\nu - e| = 0$. This implies what is asserted by the fifth bullet of (5.43).

Fix $c > 1000\, c_0$ and then $t_{c\ddagger} > c_0\, t_c$ so that (5.43) holds where $|x| \ge t_{c\ddagger}$.

*Part 7*: Use $\sigma_\ddagger$ as in (4.9) to define the connection $\hat{A}$ on the $|x| \ge t_{c\ddagger}$ part of $\mathbb{R}^4$. Keep in mind that $\nabla_{\hat{A}}\sigma_\ddagger = 0$. The curvature of $\hat{A}$ is depicted in (4.10). Since $F_A$ has finite $L^2$ norm on $\mathbb{R}^4$, the fact that $\nabla_A a$ is has pointwise bounded norm and has finite $L^2$ norm on $\mathbb{R}^4$ implies that $F_{\hat{A}}$ has finite $L^2$ norm on the $|x| \ge t_{c\ddagger}$ part of $\mathbb{R}^4$. This fact is used below to construct an automorphism of the product principal SO(3) bundle over $|x| \ge 4t_{c\ddagger}$ part of $\mathbb{R}^4$ that writes $\sigma_\ddagger$ as a constant element in $\mathfrak{S}$ to be denoted by $\sigma_\diamond$ and writes the connection $\hat{A}$ as $\theta_0 + a_\diamond \sigma_\diamond$ with $a_\ddagger$ being a coclosed 1-form on the $|x| \ge 4t_{c\ddagger}$ part of $\mathbb{R}^4$ that obeys

$$\int_{|x|>4t_{c\ddagger}} |\nabla a_\diamond|^2 \; + \; \Big( \int_{|x|>4t_{c\ddagger}} |a_\diamond|^4 \Big)^{1/2} \le c_0 \int_{|x|>t_{c\ddagger}} |F_{\hat{A}}|^2 \;. \tag{5.44}$$

Note that the exterior derivative of $a_\diamond$ obeys $da_\diamond = \langle \sigma_\ddagger F_{\hat{A}} \rangle$. The desired automorphism is constructed in the Part 8.

The norm of the 2-form $\langle \sigma_\ddagger F_{\hat{A}} \rangle$ is a priori bounded by $c_0\, z(1+\kappa_\ddagger^{\,2})$ with $z$ coming from (5.36). This bound follows directly from the bounds in (5.36) given the formula in (4.10) for $F_{\hat{A}}$. The a priori bound for $\langle \sigma_\ddagger F_{\hat{A}} \rangle$ leads to an a priori bound on the norm of the 1-form $a_\diamond$ where $|x| \ge 6t_{c\ddagger}$. This bound is written as $|a_\diamond| \le z_a$ with $z_a$ being the sum of the number on the right hand side of (5.44) and $c_0$ times $z(1+\kappa_\ddagger^{\,2})$. Moreover,

$$\lim_{|x|\to\infty} |a_\diamond| = 0 \;. \tag{5.45}$$

To prove (5.45) and the bound for $|a_\diamond|$, fix $x \in \mathbb{R}^4$ and reintroduce the function $\chi^{(x)}$, this being $\chi(|x-(\cdot)|-2)$. The 1-form $\chi_{(x)} a_\diamond$ obeys the system of equations

$$d(\chi_{(x)}\, a_\diamond) = \chi_{(x)} \langle \sigma_\ddagger F_{\hat{A}} \rangle + d\chi_{(x)} \wedge a_\diamond \quad \textit{and} \quad d{*}(\chi_{(x)}\, a_\diamond) = d\chi_{(x)} \wedge {*}a_\diamond \;. \tag{5.46}$$

Let $\mathcal{G}_x$ denote for the moment the Green's function for the operator $d + d^\dagger$ with pole at x. Keeping in mind that the norm of $\mathcal{G}_x$ is bounded by $c_0 \frac{1}{|x-(\cdot)|^3}$, it follows from (5.46) that

$$|a_\diamond|(x) \le c_0 \int_{|x-(\cdot)|<2} \frac{1}{|x-(\cdot)|^3} |F_{\hat{A}}| \; + \; c_0 \Big( \int_{|x-(\cdot)|<2} |a_\diamond|^4 \Big)^{1/4} \;. \tag{5.47}$$

As noted previously, $|F_{\hat{A}}|$ is at most $c_0 z(1+\kappa_{\ddagger}^2)$. Keeping this in mind, fix $\rho \in (0, 1)$ for the moment and write the integral on the right hand side of (5.47) with $|F_{\hat{A}}|$ as a sum of the contribution from where the distance to x is greater than $\rho$ and from where the distance to x is less than $\rho$. Their respective contributions to the right hand side of (5.47) are no greater than

$$c_0 \rho^{-1} \Big( \int_{|x-(\cdot)|<2} |F_{\hat{A}}|^2 \Big)^{1/2} \quad \textit{and} \quad c_0 \rho z(1+\kappa_{\ddagger}^2) . \tag{5.48}$$

Take $\rho$ to be $\big( \int_{|x-(\cdot)|<2} |F_{\hat{A}}|^2 \big)^{1/4}$ to see that

$$|a_\diamond|(x) \le c_0 \Big( \int_{|x-(\cdot)|<2} |F_{\hat{A}}|^2 \Big) + c_0 \Big( \int_{|x-(\cdot)|<2} |a_\diamond|^4 \Big)^{1/4} . \tag{5.49}$$

The sup norm bound for $|a_\diamond|$ and (5.45) follow from (5.49) because $|F_{\hat{A}}|$ is an $L^2$ function on $\mathbb{R}^4$ and because $|a_\diamond|$ is an $L^4$ function on $\mathbb{R}^4$.

*Part 8*: This part of the proof constructs the automorphism that pulls back the connection $\hat{A}$ as $\theta_0 + a_\diamond \sigma_\diamond$ with $a_\ddagger$ being coclosed and obeying (5.40). The desired automorphism is denoted subsequently by $h_\diamond$. It is constructed here as the composition of two automorphisms, the first is described in Step 1 and the second is described in Step 6.

<u>Step 1</u>: The norm 1 sphere in $\mathfrak{S}$ being a 2-sphere, there exists an automorphism of the product principal SO(3) bundle over the $|x| \ge t_{c\ddagger}$ part of $\mathbb{R}^4$ that pulls $\sigma_\ddagger$ as a constant map to $\mathfrak{S}$. The existence of a map of this sort follows from the fact that the Hopf fibration map from $S^3$ to $S^2$ defines an isomorphism between $\pi_3(S^3)$ and $\pi_3(S^2)$. Use such an automorphism of the product principal SO(3) bundle to identify $\sigma_\ddagger$ with a constant element $\sigma_\diamond \in \mathfrak{S}$. Since $\sigma_\ddagger$ is $\hat{A}$-covariantly constant, this automorphism must pull $\hat{A}$ back as $\theta_0 + a_{\hat{A}} \sigma_\diamond$ with $a_{\hat{A}}$ being a $\mathbb{R}$-valued 1-form with exterior derivative equal to $\langle \sigma_\ddagger F_{\hat{A}} \rangle$.

<u>Step 2</u>: Suppose that $B \subset \mathbb{R}^4$ is a ball with center where $|x| = 3t_{c\ddagger}$ and radius $t_{c\ddagger}$. The Neumann Green's function for the Laplacian on B can be used to find a smooth function $f_B$: $B \to \mathbb{R}$ such that the 1-form $a_{\hat{A},B} = a_{\hat{A}} + df_B$ is coclosed on B and is such that its norm and that of its derivative obeys

$$\int_B (|\nabla a_{\hat{A},B}|^2 + t_{c\ddagger}^{-2} |a_{\hat{A},B}|^2) \le c_0 \int_B |F_{\hat{A}}|^2 \tag{5.50}$$

Suppose that $B'$ is a second ball with center where $|x| = 3t_{c_\ddagger}$ and radius $t_{c_\ddagger}$. If $B'$ intersects $B$, then the 1-forms $a_{\hat{A},B}$ and $a_{\hat{A},B'}$ differ on $B \cap B'$ according to the rule that says $a_{\hat{A},B'} = a_{\hat{A},B} + df_{B',B}$ with $f_{B,B'}$ being a function on $B \cap B'$. It follows from (5.50) that

$$\int_{B\cap B'} (|\nabla df_{B,B'}|^2 + t_{c_\ddagger}^{-2}|df_{B,B'}|^2) \leq c_0 \int_{B\cup B'} |F_{\hat{A}}|^2 \ . \tag{5.51}$$

Any constant can be added to $f_{B,B'}$ and adding one to make its integral over $B\cap B'$ equal zero gives a version that obeys

$$t_{c_\ddagger}^{-4} \int_{B\cap B'} f_{B,B'}^{2} \leq c \int_{B\cup B'} |F_{\hat{A}}|^2 \tag{5.52}$$

with c having an upper bound that is determined a priori by a lower bound for the product of $t_c^{-4}$ times the volume of $B \cap B'$. The bound that is asserted in (5.52) follows from the fact that (5.51) bounds the $L^2$ norm of $df_{B,B'}$.

<u>Step 3</u>: Let $\mathcal{A}$ denote the spherical annulus where $\frac{5}{2} t_{c_\ddagger} \leq |x| \leq \frac{7}{2} t_{c_\ddagger}$. This annulus has a cover by $c_0$ balls of radius $t_{c_\ddagger}$ with centers on the $|x| = 3t_{c_\ddagger}$ sphere such that the intersection of any two balls in this cover is either empty or is a set with volume greater than $c_0^{-1} t_{c_\ddagger}^{4}$. Let $\mathfrak{U}$ denote this cover. The corresponding versions of the functions $\{f_{B,B'}\}_{B,B'\in\mathfrak{U}}$ can be glued together using a subbordinate partition of unity to construct a 1-form on this annulus, this denoted by $\hat{a}$, that has the following properties:

- $d\hat{a} = \langle \sigma F_{\hat{A}} \rangle$ .
- $\int_{\mathcal{A}} (|\nabla \hat{a}|^2 + t_{c_\ddagger}^{-2} |\hat{a}|^2) \leq c_0 \int_{|x|>t_{c_\ddagger}} |F_{\hat{A}}|^2$ .

(5.53)

The properties in (5.53) follow from (5.51) and (5.52) if the partition of unity is chosen so that the norm of the derivative of each of its member functions is bounded by $c_0 t_{c_\ddagger}^{-1}$.

The 1-forms $\hat{a}$ and $a_{\hat{A}}$ are related on $\mathcal{A}$ by the rule $\hat{a} = a_{\hat{A}} + df$ with $f$ being a function that is defined on a neighborhood of $\mathcal{A}$. Extend $f$ as a smooth function with compact support to the $|x| \geq \frac{7}{2} t_{c_\ddagger}$ part of $\mathbb{R}^4$. Also denote this extension by $f$ and also use $\hat{a}$ to denote the 1-form $a_{\hat{A}} + df$ over the rest of the $|x| \geq \frac{5}{2} t_{c_\ddagger}$ part of $\mathbb{R}^4$ .

<u>Step 4</u>: Use the function $\chi$ to construct a function that is equal to 1 where the distance to the origin is greater than $\frac{13}{4} t_{c_\ddagger}$ and is equal to zero where the distance to the origin is less than $\frac{11}{4} t_{c_\ddagger}$. This function can and should be constructed so that norm of its differential is bounded by $c_0 t_{c\ddagger}^{-1}$. Denote this function by $\chi_c$. Given $k \geq 2$, there exists a smooth function on the $|x| \leq 2^k t_{c\ddagger}$ ball in $\mathbb{R}^4$ to be denoted by $u_k$ that obeys the equation

$$d^\dagger du_k + d^\dagger(\chi_c \hat{a}) = 0 \tag{5.54}$$

where $|x| < 2^k t_c$ ball and is such that $*du_k + *\hat{a}$ pulls back as zero on the boundary of this ball. Let $\hat{a}_k$ denote $\chi_c\hat{a} + du_k$. The exterior derivative of $\hat{a}_k$ is zero where $|x| \leq \frac{11}{4} t_{c\ddagger}$ and it is equal to $\langle\sigma_\ddagger F_{\hat{A}}\rangle$ on the spherical annulus where $\frac{13}{4} t_{c\ddagger} \leq |x| \geq 2^k t_{c\ddagger}$. Its exterior derivative on $\mathcal{A}$ is $d\chi_c \wedge \hat{a} + \chi_c\langle\sigma_\ddagger F_{\hat{A}}\rangle$. These facts and (5.52) imply that the derivatives of $\hat{a}_k$ obey

$$\int_{|x|\leq 2^k t_{c\ddagger}} |\nabla\hat{a}_k|^2 \leq c_0 \int_{|x|>t_{c\ddagger}} |F_{\hat{A}}|^2 . \tag{5.55}$$

Add a constant 1-form to $\hat{a}_k$ if necessary so that the resulting 1-form has components with integral equal to zero on the $|x| \leq 2^k t_{c\ddagger}$ ball. Use $\tilde{a}_k$ to denote the latter. This 1-form obeys $d{*}\tilde{a}_k = 0$, it obeys $d\tilde{a}_k = \langle\sigma_\ddagger F_{\hat{A}}\rangle$ where $\frac{13}{4} t_{c\ddagger} \leq |x| \leq 2^k t_{c\ddagger}$ , and its derivatives obey the bound in (5.55). Since each component of $\tilde{a}_k$ has integral zero, the first and third bullets in (2.1) can be invoked to see that

$$\Big( \int_{|x|\leq 2^k t_{c\ddagger}} |\tilde{a}_k|^4 \Big)^{1/2} \leq c_0 \int_{|x|>t_{c\ddagger}} |F_{\hat{A}}|^2 . \tag{5.56}$$

<u>Step 5</u>: The sequence $\{\tilde{a}_k\}_{k\in\{2,3,\ldots\}}$ is such that if n is an integer greater than 2 and if $k \geq n$, then $\tilde{a}_k$ is defined on the spherical annulus where $\frac{13}{4} t_{c\ddagger} \leq |x| \leq 2^n t_{c\ddagger}$ and each has bounded $L^2_1$ norm on this annulus with the $L^2$ norm of $\nabla\tilde{a}_k$ and the $L^4$ norm of $\tilde{a}_k$ being bounded by $c_0$ times the $L^2$ norm of $\langle\sigma_\ddagger F_{\hat{A}}\rangle$ on the part of $\mathbb{R}^4$ where $|x| \geq t_{c\ddagger}$. It follows as a consequence that the sequence $\{\tilde{a}_k\}_{k\in\{2,3,\ldots\}}$ has a subsequence that converges in the $C^\infty$ topology on compact subsets of the $|x| \geq 4t_{c\ddagger}$ part of $\mathbb{R}^4$ to a smooth, coclosed 1-form whose differential is $\langle\sigma_\ddagger F_{\hat{A}}\rangle$. This is the 1-form $a_\diamond$. The convergence is in the $C^\infty$ topology because $\langle\sigma_\ddagger F_{\hat{A}}\rangle$ is smooth and because each $k \in \{2, 3, \ldots\}$ version of $\tilde{a}_k$ obeys the equations $d{*}\tilde{a}_k = 0$ and $d\tilde{a}_k = \langle\sigma_\ddagger F_{\hat{A}}\rangle$ on the spherical annulus where $\frac{13}{4} t_{c\ddagger} \leq |x| \leq 2^k t_{c\ddagger}$.

The integral bounds in (5.44) are consequences of the fact that the $L^4$ norm of any $k \in \{2, 3, \ldots\}$ version $\tilde{a}_k$ on its domain of definition and the $L^2$ norm of $\nabla\tilde{a}_k$ on this same domain are both bounded a priori by $c_0 \int_{|x|>t_{c\ddagger}} |F_{\hat{A}}|^2$ .

Step 6: The 1-forms $a_\diamond$ and $a_{\hat{A}}$ have the same differential on the $|x| \ge 4t_{c\ddagger}$ part of $\mathbb{R}^4$, this being $\langle \sigma_\ddagger F_{\hat{A}} \rangle$. It follows as a consequence that $a_\diamond$ can be written as $a_{\hat{A}} + du$ with u being a smooth function on this same part of $\mathbb{R}^4$. This implies that $e^{u\sigma_\diamond}$ defines an automorphism of the product principal SO(3) bundle over the $|x| \ge 4t_{c\ddagger}$ part of $\mathbb{R}^4$ that pulls back the connection $\theta + a_{\hat{A}}\sigma_\diamond$ as $\theta + a_\diamond \sigma_\diamond$.

*Part 9*: By way of a look ahead, this part of the proof of Lemma 5.4 borrows heavily from the discussion in Section 4d. To set the stage for what is to come, fix $c > 1000\, c_0$ and then fix $t_{c\diamond} > 6\, t_{c\ddagger}$ so that

- $|\nu - e| \le c^{-1}$ ,
- $|\nabla_A a_\ddagger| \le c^{-1}$ ,
- $|a_\diamond| \le c^{-1}$ .

(5.57)

The automorphism of the product principal SO(3) bundle over the $|x| \ge t_{c\diamond}$ part of $\mathbb{R}^4$ that was constructed in Part 5 is denoted by $h_\diamond$. Keep in mind that $h_\diamond$ pulls back $\sigma_\ddagger$ as $\sigma_\diamond$ and it pulls back $\hat{A}$ as $\theta_0 + a_\diamond \sigma_\diamond$. The pull-back of $\mathfrak{a}$ by $h_\diamond$ on the $|x| \ge t_{c\diamond}$ part of $\mathbb{R}^4$ is still denoted by $\mathfrak{a}$, the latter incarnation being an $\mathfrak{h}$-valued 1-form with $\mathfrak{h} \subset \mathfrak{G}$ denoting the kernel of the linear map $\langle \sigma_\diamond (\cdot) \rangle$. The $h_\diamond$ pull-back of the connection A is denoted by A also. Setting $m$ to equal $\frac{1}{\sqrt{2\pi}} \kappa_\ddagger$, write this new incarnation of A as $\theta_0 + a_\diamond \sigma_\diamond + \mathfrak{b}$ with $\mathfrak{b} = \frac{1}{4} m^{-1} [\sigma_\diamond, \nabla_A \sigma_\diamond]$ being another $\mathfrak{h}$-valued 1-form on the $|x| \ge t_{c\diamond}$ part of $\mathbb{R}^4$. The lemma that follows makes a central observation about $\mathfrak{a}$ and $\mathfrak{b}$.

**Lemma 5.5:** *The sections* $\mathfrak{a}$ *and* $\mathfrak{b}$ *and their covariant derivatives are square integrable on the* $|x| \ge t_{c\diamond}$ *part of* $\mathbb{R}^4$. *Moreover, there exists* $c_* > 1$ *such that*

$$\int_{|(\cdot)|>R} \left( |\nabla \mathfrak{a}|^2 + |\nabla \mathfrak{b}|^2 + |\mathfrak{a}|^2 + |\mathfrak{b}|^2 \right) \le c_* e^{-R/c_*} \quad \textit{if}\ R \ge t_{c\diamond}.$$

***Proof of Lemma 5.5***: The proof of the lemma has four steps.

Step 1: The fact that $\mathfrak{a}$ is square integrable on the $|x| \ge t_{c\diamond}$ part of $\mathbb{R}^4$ follows from the formula in the fourth bullet of (4.5) and the fact that $|a_\ddagger \wedge a_\ddagger|$ is square integrable. The fact that $\mathfrak{b}$ is square integrable on this same domain follows from the fourth bullet of (4.7) and the fact that $|\nabla_A a_\ddagger|$ is square integrable. The covariant derviative of $\nabla \mathfrak{a}$ is also square integrable on the $|x| \geq t_{c\diamond}$ part of $\mathbb{R}^4$, this being another consequence of the fourth bullet of (4.7) and the fact that $|\nabla_A a_{\ddagger}|$ is square integrable. The proof of the $\nabla\mathfrak{b}$ is square integrable follows from the fact noted in Part 3 about $|\nabla_A(\nabla_A a_{\ddagger})|^2$ being an integrable function on $\mathbb{R}^3$.

<u>Step 2</u>: Define the vector bundles $\mathbb{V}$ and $\mathbb{W}$ over $\mathbb{R}^4$ as in Part 2 of Section 4d. Use the pair $(\mathfrak{a}, \mathfrak{b})$ to define the section $(\mathfrak{p}, \mathfrak{q})$ of $\mathbb{V} \oplus \mathbb{V}$ over the $|x| \geq t_{c\diamond}$ part of $\mathbb{R}^4$ by the rule whereby $\mathfrak{p} = \tau\,\mathfrak{a} + (1-\tau)\,\mathfrak{b}$ and $\mathfrak{q} = \tau\,\mathfrak{b} - (1-\tau)\,\mathfrak{a}$. Denote this section by $x_{(\mathfrak{a},\mathfrak{b})}$.

Define the operator $\mathcal{L}$ as in (4.13) using $m = \frac{1}{\sqrt{2\pi}} K_{\ddagger}$. The arguments used in Part 4 of Section 4d can be repeated in the present context to see that $x_{(\mathfrak{a},\mathfrak{b})}$ obeys an equation on the $|x| \geq t_{c\diamond}$ part of $\mathbb{R}^4$ that has the schematic form

$$\mathcal{L}x_{(\mathfrak{a},\mathfrak{b})} + \mathcal{R}x_{(\mathfrak{a},\mathfrak{b})} = 0 \tag{5.58}$$

with $\mathcal{R}$ being a certain linear operator from sections of the product $\mathbb{V} \oplus \mathbb{V}$ bundle over the $|x| \geq t_{c\diamond}$ part of $\mathbb{R}^4$ to sections of the product $\mathbb{W}$ bundle over this same part of $\mathbb{R}^4$. With $\mathbb{W}$ written as the sum $(\mathbb{W}^+ \oplus \mathfrak{h}) \oplus (\mathbb{W}^- \oplus \mathfrak{h})$, the respective summand components of $\mathcal{R}$ are written as $(\mathcal{R}^+, \mathcal{R}^0{}_{\mathfrak{q}}), (\mathcal{R}^-, \mathcal{R}^0{}_{\mathfrak{p}}))$. The next paragraph describes these components.

Let $x$ denote a given section of $\mathbb{V} \oplus \mathbb{V}$ over the $|x| \geq t_{c\diamond}$ part of $\mathbb{R}^4$ and write its components as $(\alpha, \beta)$ with $\alpha$ and $\beta$ being sections of $\mathbb{V}$. Granted this notation, then the endomorphisms $\mathcal{R}^+$ and $\mathcal{R}^-$ are the sum of the respective $\pm$ terms in (4.29) and the corresponding $\pm$ terms in (4.30). The terms in (4.28) are zero in this case because the metric is flat. As in Section 4d, the endomorphisms $\mathcal{R}^0{}_{\mathfrak{q}}$ and $\mathcal{R}^0{}_{\mathfrak{p}}$ are constant coefficient linear combinations of endomorphisms that are denoted by $\mathfrak{R}^0{}_{\mathfrak{a}}$ and $\mathfrak{R}^0{}_{\mathfrak{b}}$. The endomorphism $\mathfrak{R}^0{}_{\mathfrak{a}}$ is the sum of the endomorphisms that are depicted in (4.32) and (4.33). The endomorphism in (4.31) is zero on $\mathbb{R}^4$. The endomorphism $\mathfrak{R}^0{}_{\mathfrak{b}}$ has the form that is depicted in (4.38) but with $\mathfrak{C}_3$, $t_{01}$ and $t_{02}$ set equal to 0.

<u>Step 3</u>: It follows from the preceding description of $\mathcal{R}$ where $|x| \geq t_{c\diamond}$ can be written as

$$\mathcal{R}x = \mathcal{T}x + \mathcal{C}(\nabla x)\ , \tag{5.59}$$

with $\mathcal{T}$ and $\mathcal{C}$ being homomorphisms of the following sort: The norm of the homomorphism $\mathcal{T}$ is no greater than $c_0(|\nu - e| + |\nabla_A a_{\ddagger}| + |a_{\diamond}|)$ and thus no greater than $c_0\, c^{-1}$. The components of $\mathcal{C}$ are constant coefficient, linear combinations of the components of $\langle \sigma_{\diamond} \nabla_A a_{\ddagger} \rangle$ and so $|\mathcal{C}| \leq c_0\, c^{-1}$ also.

The preceding remarks about $\mathcal{R}$ lead to the following observation: If $c > c_0$ and if $x$ is an $L^2_1$ section of $\mathbb{V} \oplus \mathbb{V}$ with support where $|x| \geq t_{c\lozenge}$, then

$$\int_{|(\cdot)|>t_{c_\lozenge}} (|\nabla x|^2 + 4m^2|x|^2) \leq \tfrac{4}{3} \int_{|(\cdot)|>t_{c_\lozenge}} |(\mathcal{L}+\mathcal{R})x|^2 \quad . \tag{5.60}$$

Note that the bounds for $C$ also lead to a $c_0\, c^{-1}$ bound on the norm of the principal symbol of the operator $\mathcal{R}$ where $|x| \geq t_{c\lozenge}$.

<u>Step</u> <u>4</u>: Fix an integer, k that is greater than 1 and let $\chi_k$ denote the function x on $\mathbb{R}^4$ given by the rule $x \to \chi(k - t_{c\lozenge}^{-1}|x|)$. This function equals 1 where $|x| \geq k\, t_{c\lozenge}$ and it equals 0 where $|x| \leq (k-1)\, t_c$. Use $x_k$ to denote $\chi_k\, x_{(\mathfrak{a},\mathfrak{b})}$. If $k > 2$, then $x_k$ obeys the equation

$$(\mathcal{L}+\mathcal{R})\, x_k = \mathfrak{s}(d\chi_{k-1})\, x_{k-1} \tag{5.61}$$

because $\chi_{k-1} d\chi_k = d\chi_k$ when $k > 2$. It follows as a consequence of (5.60) and what is said subsquently about the principal symbol of $\mathcal{R}$ that

$$\int_{\mathbb{R}^4} (|\nabla x_k|^2 + 4m^2|x_k|^2) \leq c_0\, t_{c\lozenge}^{-2} \int_{\mathbb{R}^4} |x_{k-1}|^2 \quad . \tag{5.62}$$

Keeping in mind that $L^2$ norm of $|\mathfrak{a}|^2 + |\mathfrak{b}|^2$ is in any event no greater than $c_0 \mathrm{H} \mathrm{K}_\ddagger^{\,2}$, it follows that the bound in (5.62) with its versions that have k replaced in turn by the integers from the set {3, …, k-1} imply the bound

$$\int_{\mathbb{R}^4} (|\nabla x_k|^2 + 4m^2|x_k|^2) \leq c_0\, (m\, t_{c\lozenge})^{-2k+3}\, \mathrm{H}\mathrm{K}_\ddagger^{\,2} \,. \tag{5.63}$$

The latter bound leads directly to Lemma 5.5's bound.

*Part 10*: What with (4.4) and (4.10) and the fact that $a_\lozenge$ is coclosed, the equations in (2.11) on the $|x| \geq t_{c\lozenge}$ part of $\mathbb{R}^4$ lead to equations for $\nu$ and $a_\lozenge$ that can be written as

- $(\tau\, da_\lozenge - (1-\tau)\, d\nu)^+ = -\tau \langle \sigma_\lozenge (m_\ddagger^{\,2}\, \mathfrak{b} \wedge \mathfrak{b})^- + \mathfrak{a} \wedge \mathfrak{a})^+ \rangle + (1-\tau) m \langle \sigma_\lozenge (\mathfrak{b} \wedge \mathfrak{a} + \mathfrak{a} \wedge \mathfrak{b})^+ \rangle$.
- $((1-\tau)\, da_\lozenge + \tau\, d\nu)^- = -(1-\tau)(\langle \sigma_\lozenge (m^2\, \mathfrak{b} \wedge \mathfrak{b}) + \mathfrak{a} \wedge \mathfrak{a})^- \rangle - \tau\, m \langle \sigma_\lozenge (\mathfrak{b} \wedge \mathfrak{a} + \mathfrak{a} \wedge \mathfrak{b})^- \rangle$.
- $d{*}\nu = -\, m \langle \sigma_\lozenge (\mathfrak{b} \wedge {*}\mathfrak{a} - \mathfrak{a} \wedge {*}\mathfrak{b}) \rangle$ *and* $d{*}a_\lozenge = 0$.

(5.64)

Let $w_+ = \tau\, a_\diamond - (1-\tau)(\nu - \frac{1}{\sqrt{2}\pi} e)$ and let $w_- = (1-\tau)\, a_\diamond + \tau(\nu - \frac{1}{\sqrt{2}\pi} e)$. The equations in (5.64) can be written schematically in terms of $w_+$ and $w_-$ as

$$(dw_+)^+ + d^\dagger w_+ = \mathfrak{q}_+ \quad \textit{and} \quad (dw)^- + d^\dagger w_- = \mathfrak{q}_- \,. \tag{5.65}$$

To make something of these equations, fix $x \in \mathbb{R}^4$ for the moment, and let $G^+_x$ denote the Green's for the operator for $(d(\cdot))^+ + d^\dagger(\cdot)$. Use $G^-_x$ to denote the Green's function with pole at x for the operator $(d(\cdot))^- + d^\dagger(\cdot)$. Both $|G^+_x|$ and $|G^-_x|$ are bounded by $c_0 \frac{1}{|x-(\cdot)|^3}$. Let $\mathfrak{s}_+$ and $\mathfrak{s}_-$ denote the the principal symbols of the operators $(d(\cdot))^+ + d^\dagger(\cdot)$ and $(d(\cdot))^- + d^\dagger(\cdot)$.

Since $|w_+|$ and $|w_-|$ are $L^4$ functions on the $|x| \geq t_{c\diamond}$ part of $\mathbb{R}^4$, it follows from what is said in Lemma 5.5 that $w_+$ and $w_-$ can be written at any $x \in \mathbb{R}^4$ with $|x| \geq 8\, t_c$ as

$$w_+|_x = \int_{\mathbb{R}^4} G^+_x(\chi_2 \mathfrak{q}_+ + \mathfrak{s}_+(d\chi_2) w_+) \quad \textit{and} \quad w_-|_x = \int_{\mathbb{R}^4} G^-_x(\chi_2 \mathfrak{q}_- + \mathfrak{s}_-(d\chi_2) w_-) \tag{5.66}$$

with $\chi_2$ denoting the function $\chi(2 - t_{c\diamond}^{-1}|x|)$. The next three paragraphs use (5.66) to bound the norms of $w_+$ and $w_-$ at $x \in \mathbb{R}^4$ with $|x|$ large by an x-independent multiple of $\frac{1}{|x|^3}$

The norms of the contributions to the right hand side of (5.66) from the integrals with $\mathfrak{s}_+$ and $\mathfrak{s}_-$ are bounded by a constant times $\frac{1}{|x|^3}$ because the integrands for these integrals have compact support.

To bound the norms of the integrals in (5.66) with $\mathfrak{q}_+$ and $\mathfrak{q}_-$, start by writing each integral as a sum of two integrals, the first where the domain of integration is the radius $\frac{1}{2}|x|$ ball about the origin and the second where the integration domain is the complement of this ball. The norm of the contribution to each integral from the radius $\frac{1}{2}|x|$ ball about the origin is at most $c_0 \frac{1}{|x|^3}$ times the sum of the squares of the $L^2$ norms of $\mathfrak{a}$ and $\mathfrak{b}$ because the distance from x to radius $\frac{1}{2}|x|$ ball about the origin is greater then $\frac{1}{2}|x|$ and the Green's functions are bounded by $c_0$ times the inverse of the third power of the distance to x.

Let $\mathcal{A}_x$ denote the region where the distance to the origin is greater than $\frac{1}{2}|x|$. Fix $\rho \in (0, 1]$ and break the contributions to the $\mathfrak{q}_+$ and $\mathfrak{q}_-$ integrals in (5.66) from $\mathcal{A}_x$ into the contributions from the complement in $\mathcal{A}_x$ of the radius $\rho$ ball centered on x and from the contribution from this same ball. The norm of former is no greater than $\rho^{-3}$ times the integral of $|\mathfrak{a}|^2$ and $|\mathfrak{b}|^2$ on the $\mathcal{A}_x$; and what is said in Lemma 5.5 implies that this is at most $c_0 \rho^{-3}\, e^{-|x|/2c_*}$. Meanwhile, the contribution from the ball of radius $\rho$ centered at x is at most $c_0 c^{-2}\rho$. This understood, take $\rho$ equal to $e^{-|x|/8c_*}$ to bound the norm of the contributions to the $\mathfrak{q}_+$ and $\mathfrak{q}_-$ integrals in (5.66) by $c_0\, e^{-|x|/8c_*}$.

*Part 11*: Since $|w_+|$ and $|w_-|$ are both bounded by a multiple of $\frac{1}{|x|^3}$, this must also be the case for $|\nu - \frac{1}{\sqrt{2\pi}} e|$. Fix $r > 100\, t_{c\Diamond}$ and it follows from the preceding observation and from Lemma 5.5 that

$$\frac{4}{r^4} \int_{r<|(\cdot)|<2r} |a_{\ddagger}|^2 \geq 1 - c\,\Big(\frac{1}{r^3} + e^{-r/c}\Big) \tag{5.67}$$

with c being a positive number. Since the function $r \to \kappa(r)$ is increasing, this bound implies in turn what is asserted in (5.35).

## 6. Convergence in the $C^0$ topology

Special cases of the upcoming Proposition 6.1 assert that the function $|\hat{a}_\Diamond|$ from either the second bullet of Proposition 2.1 or the second bullet of Proposition 2.2 is continuous on the complement of a finite set of points and that $\{|a_n|\}_{n\in\Lambda}$ converges to $|\hat{a}_\Diamond|$ in the $C^0$ topology on compact subsets in the complement of this finite set.

To set stage for Proposition 6.1, suppose that X is a smooth, oriented Riemannian manifold and that $P \to X$ is a principal SO(3) bundle. Let U denote a given open set in X with compact closure. As in the previous sections, the data X, P and U constitute the geometric data.

A set consisting of a number $r \geq 1$, a number $\tau \in [0, 1]$ and a pair $(A, a)$ of connection on P and section of $T^*X \otimes (P \times_{SO(3)} \mathfrak{G})$ that obeys (2.11) is said below to be a *solution to (2.11)*. A sequence $\{(r_n, \tau_n, (A_n, a_n))\}_{n\in\{1,2,\ldots\}}$ of solutions to (2.11) is said to *weakly converge on* U if the conditions in the list given momentarily in (6.1) are met.

One of the conditions in (6.1) refers to a certain Green's function for the operator $d^\dagger d + 1$. To say more, suppose first that X is compact. Given $p \in U$, the notation has $G_p$ denoting the Green's function with pole at p for $d^\dagger d + 1$ on X. Supposing that X is not compact, there is a compact subset of X with smooth boundary with all of its boundary points having positive distance from U. Fix such a manifold with boundary in X containing U and denoted it by $X_U$. If p is a given point in U, then $G_p$ denotes the Dirichelet Green's function on $X_U$ with pole at p for the operator $d^\dagger d + 1$.

- NUMERICAL CONSTRAINTS: *The sequence* $\{r_n\}_{n\in\{1,2,\ldots}$ *has no bounded subsequences and the sequence* $\{\tau_n\}_{n\in\{1,2,\ldots\}}$ *converges.*
- INTEGRAL CONSTRAINTS: *Use f to denote a* $C^2$ *function on* X.

a) *The sequences* $\{\int_X f\,|F_{A_n} - r^2 a_n \wedge a_n|^2\}_{n\in\{1,2,\ldots\}}$ *and* $\{\int_X f\,|d_{A_n} a_n|^2\}_{n\in\Lambda}$ *converge.*

b) *The sequences* $\{\int_X f\,|\nabla_{A_n} a_n|^2\}_{n\in\{1,2,\ldots\}}$ *and* $\{r_n{}^2 \int_X f\,|a_n \wedge a_n|^2\}_{n\in\{1,2,\ldots\}}$ *converge. The limit of the left most sequence is denoted in what follows by* $Q_{\nabla,f}$ *and that of the right most by* $Q_{\wedge,f}$.

c) *The sequence* $\{\int_U |a_n|^2\}_{n\in\{1,2,\ldots\}}$ *converges with non-zero limit.*

d) *The sequence* $\{\int_X |d|a_n||^2\}_{n\in\{1,2,\ldots\}}$ *converges.*

- POINTWISE CONSTRAINTS: *The sequence* $\{\sup_X |a_n|\}_{n\in\{1,2,\ldots\}}$ *converges.*
- BEHAVIOR OF THE LIMIT: *The sequence* $\{|a_n|\}_{n\in\{1,2,\ldots\}}$ *converges weakly in the* $L^2_1$ *topology on* X *and strongly in all* $p<\infty$ *versions of the* $L^p$ *topology. The limit function is denoted by* $|a_\Diamond|$.

a) *The function* $|a_\Diamond|$ *can be defined at each point in* X *by the rule*

$$p \to |\hat a_\Diamond|(p) = \lim\sup_{n\in\{1,2,\ldots\}} |a_n|(p).$$

b) *The sequence* $\{\langle a_n \otimes a_n\rangle\}_{n\in\{1,2,\ldots\}}$ *converges strongly in any* $q<\infty$ *version of the* $L^q$ *topology on the space of sections of* $T^*X \otimes T^*X$. *The limit section is denoted by* $\langle \hat a_\Diamond \otimes \hat a_\Diamond\rangle$ *and its trace is the function* $|\hat a_\Diamond|^2$.

c) *If* $f$ *is a* $C^2$ *function with compact support in* X, *then*

$$\tfrac{1}{2}\int_X d^* df\, |\hat a_\Diamond|^2 + Q_{\nabla,f} + Q_{\wedge,f} + \int_X f\, \mathrm{Ric}(\langle \hat a_\Diamond \otimes \hat a_\Diamond\rangle) = 0\,,$$

d) *The sequence that is indexed by the positive integers with n'th term being the integral of* $G_p(|\nabla_{A_n} a_n|^2 + 2r_n{}^2 |a_n \wedge a_n|^2)$ *is bounded. Let* $Q_{\Diamond,p}$ *denote the lim-inf of this sequence. The function* $|\hat a_\Diamond|^2$ *obeys the equation*

$$\tfrac{1}{2}|\hat a_\Diamond|^2(p) + Q_{\Diamond,p} = -\int_X G_p\,(\tfrac{1}{2}|\hat a_\Diamond|^2 - \mathrm{Ric}(\langle \hat a_\Diamond \otimes \hat a_\Diamond\rangle)) + Q_P$$

*with* $Q_{(\cdot)}$ *being a smooth function on* U.

(6.1)

Note that the condition in Item c) of the second bullet implies that the limit function $|\hat a_\Diamond|$ from the fourth bullet is not identically zero on U.

The subsequences that are described in the second bullet of Proposition 2.1 weakly converge on the compact manifold X when they are renumbered consecutively from 1. The sequences that are described in the second bullet of Proposition 2.2 weakly converge after renumbering if U is the interior of $I \times M$ and X is the interior of $I' \times M$ with $I' \subset \mathbb{I}$ being an open interval with compact closure that contains the closure of I.

**Proposition 6.1**: *Let* $\{(r_n, \tau_n, (A_n, a_n))\}_{n\in\{1,2,\ldots\}}$ *denote a weakly convergent sequence of solutions to (2.11) on* U *and let* $|\hat{a}_\diamond|$ *denote the corresponding limit function as described by the fourth bullet of (6.1). The function* $|\hat{a}_\diamond|$ *is continuous on the complement of a finite set in* U *and smooth at points in the complement of this set where it is positive. Moreover, the sequence* $\{|a_n|\}_{n\in\{1,2,\ldots\}}$ *converges to* $|\hat{a}_\diamond|$ *in the* $C^0$ *topology on compact subsets of the complement in* U *of this same finite set.*

By way of a parenthetical remark, the proof of Proposition 6.1 supplies a bound on the number of points in Proposition 6.1's finite set given the value of the limit of the sequence $\{\int_X |F_{A_n} - r^2 a_n \wedge a_n|^2\}_{n\in\{1,2,\ldots\}}$ from Item a) of the second bullet of (6.1).

The proof of Proposition 6.1 is in Section 6c. The intervening subsections establish results that are used in these proofs.

**a) The relations between of $\mathrm{r}_{c_\wedge}$ and $\mathrm{r}_{c\mathrm{F}}$**

The lemma stated momentarily supplies the key input to the proof of Proposition 6.1. It says implies among other things that most $p \in U$ versions of $\mathrm{r}_{c_\wedge}$ determine the corresponding version of $\mathrm{r}_{3c_\mathrm{F}}$. It also says that $\mathrm{r}_{c\mathrm{F}}$ in most cases determines $\mathrm{r}_{3c_\wedge}$.

**Lemma 6.2**: *There exists* $\kappa > 1$ *depending only on the geometric data with the following significance: Suppose that* $r > \kappa$ *and that* $(A, a)$ *is a pair of connection on* P *and section of* $T^*X \otimes (P \times_{SO(3)} \mathfrak{G})$ *that obeys (2.11) for some* $\tau \in [0,1]$. *Suppose in addition that* $\mathrm{E} \geq 1$ *and that* $(A, a)$ *obeys the bounds in (5.1). Fix* $c \geq 1$ *and there exists a set in* U *to be denoted by* $\Theta_c$ *with the properties listed below.*

- *The set* $\Theta_c$ *has at most* $\mathrm{E}^2 c^2$ *points*
- *If* $p \in U$ *and* $\mathrm{r} \in (0, \kappa^{-1}]$ *are such* p *has distance* $\kappa\mathrm{r}$ *or more from* $\Theta_c$, *then*

$$\int_{B_r} |F_A|^2 \leq 2 r^4 \int_{B_r} |a \wedge a|^2 + \tfrac{1}{4} c^{-2} .$$

Lemma 6.2 is proved momentarily. By way of a parenthetical remark, the upcoming proof of the second bullet of Lemma 6.2 can be repeated with only cosmetic changes to prove that

$$r^4 \int_{B_r} |a \wedge a|^2 \leq 2 \int_{B_r} |F_A|^2 + \tfrac{1}{4} c^{-2} \tag{6.2}$$

when $p \in U$ and $r \in (0, c_0^{-1}]$ are such that p has distance $c_0 r$ or more from $\Theta_c$. The inequality in (6.2) is not needed for subsequent arguments.

The next lemma is a corollary of sorts to Lemma 6.2. The lemma reintroduces the number $z_U$ from Proposition 5.1.

**Lemma 6.3**: *There exists* $\kappa > 1$ *with the following significance: Use* $c \geq z_U$ *to define the set* $\Theta_c$*. Fix* $p \in U$ *and supposed that there exists* $r \in (0, \kappa^{-1}]$ *such that* $r\kappa(r)\, r = z_U^{-1}$*. Denote this number by* $r_{*p}$*. If* $\mathrm{dist}(p, \Theta_c) > \kappa r_{*p}$*, then* $\int_{\mathrm{dist}(\cdot, p) < \frac{1}{2} r_{*p}} |F_A|^2 \leq \frac{1}{2} z_U^{-2}$.

By way of a look ahead, Lemma 6.3 is used to certify Proposition 5.1's curvature requirement for points that are uniformly far from any given $c > c_0 z_U$ version of $\Theta_c$.

***Proof of Lemma 6.2***: The proof uses $\mathfrak{w}$ to denote $F_A - r^2 a \wedge a$. The three step iterative algorithm constructs the set $\Theta_c$.

Step 1: Let $r_0$ denote the largest $r \in (0, c_0^{-1}]$ with the property that the integral of $|\mathfrak{w}|^2$ on every radius r ball centered in U is less than $\frac{1}{8} c^{-2}$. There exists a maximal set of points in U with the following two properties: The integral of $|\mathfrak{w}|^2$ over the radius $r_0$ ball about each such point from the set is equal to $\frac{1}{8} c^{-2}$ and distinct points in the set are separated by a distance no less than $2r_0$. Denote this set by $\Theta_{c0}$ and use $n_0$ to denote the set of points in $\Theta_{c0}$. Note that $n_0 \leq \mathrm{E}^2 c^{-2}$.

Step 2: Suppose that $k \in \{0, 1, 2, \ldots\}$ and that a set $\Theta_{ck}$ has been constructed with two properties, the first being that if $p \in \Theta_{ck}$, then there exists $r_p \in [r_0, 2^k r_0]$ such that the integal of $|\mathfrak{w}|^2$ on the radius $r_p$ ball centered at p is greater than or equal to $\frac{1}{8} c^{-2}$. The second property is as follows: If p and p´ are distinct points from $\Theta_{ck}$ then the ball of radius $r_p$ centered at p is disjoint from the ball of radius $r_{p'}$ centered at p´.

To construct $\Theta_{ck+1}$, fix a maximal set of points in $U - \Theta_{ck}$ with the following three properties: The integral of $|\mathfrak{w}|^2$ on the radius $2^{k+1} r_0$ ball centered at each point of this set is greater than or equal to $\frac{1}{8} c^{-2}$, distinct points in this set are separated by a distance no less than $2^{k+2} r_0$, and each point in this set has distance no less than $2^{k+2} r_0$ from $\Theta_{ck}$. Take $\Theta_{ck+1}$ to be the union of this set and $\Theta_{ck}$.

<u>Step 3</u>: The iterative construction in Step 2 must end after a finite number of iterations because the integral of $|w|^2$ over the whole of X is at most $E^2$. Supposing that N is the total number of iterations, take $\Theta_c$ to be $\Theta_{cN}$.

The set $\Theta_c$ obeys the constraint given by the lemma's first bullet because the construction provides a set of disjoint balls with centers at the points of $\Theta_c$ such that $|w|^2$ has integral equal to $\frac{1}{8}c^{-2}$ on each ball. The upcoming proof of the second bullet of the lemma invokes the triangle inequality in the guise

$$\int_{B_r} |F_A|^2 \le 2r^4 \int_{B_r} |a\wedge a|^2 + 2\int_{B_r} |w|^2 . \tag{6.3}$$

To start the proof of the lemma's second bullet, suppose that $p \in U$ and $r \in (0, c_0^{-1}]$ are such that p has distance $4r$ from each point of $\Theta_c$. If r is less than $r_0$, then the integral of $|w|^2$ over $B_r$ is less than $\frac{1}{8}c^{-2}$ and so the lemma's second bullet follows from (6.3). Suppose next that $k \in \{0, 1, \ldots, N-1\}$ and $r \in [2^k r_0, 2^{k+1} r_0)$. The point p has distance greater than 4r from $\Theta_c$ so it has distance greater than $2^{k+2} r_0$ from each point in $\Theta_{ck+1}$. Were the integral of $|w|^2$ over $B_r$ greater than $\frac{1}{8}c^{-2}$ then p would be in $\Theta_{ck+1}$ and this is not possible because $\Theta_{ck+1}$ is maximal. This observation with (6.3) imply what is asserted by the lemma's second bullet. Very much the same argument proves that the second bullet of the lemma is true if $r \ge 2^{N+1} r_0$.

***Proof of Lemma 6.3***: Let u denote the function on X that is defined as follows:

- $u = 1$ *where* $\mathrm{dist}(\cdot, p) \le \frac{1}{2} r_{*p}$ .
- $u = 2 - 2\,\mathrm{dist}(\cdot, p)$ *where* $\frac{1}{2} r_{*p} \le \mathrm{dist}(\cdot, p) \le r_{*p}$ .
- $u = 0$ *where* $r_{*p} \le \mathrm{dist}(\cdot, p)$ .

(6.4)

Multiply both sides of (2.6) by $u^2$ and then integrate the resulting identity over the ball of radius $r_{*p}$ centered at p. Integrate by parts and then use Holder's inequality

$$2r^2 \int_{\mathrm{dist}(\cdot,p) < \frac{1}{2} r_{*p}} |a\wedge a|^2 \le r_{*p}^{-2}(1 + c_0 r_{*p}^{2}) \int_{\mathrm{dist}(\cdot,p) < r_{*p}} |a|^2 . \tag{6.5}$$

Use the fact that the function $\kappa$ is increasing to bound the right hand side of (6.5) by $\frac{1}{4} r_{*p}^{2} \kappa(r_{*p})^2 (1 + c_0 r_{*p}^{2})$. Use the latter bound in (6.5) with the definition of $r_{*p}$ to see that

$$r^4 \int_{\text{dist}(\cdot,p)<\frac{1}{2}r_{*p}} |a\wedge a|^2 \le \tfrac{1}{8}\, z_U^{-2}(1 + c_0 r_{*p}^{\,2})\ . \tag{6.6}$$

The assertion made by Lemma 6.3 follows from (6.6) and Lemma 6.2's second bullet.

**b) A lower bound for $r_{c\wedge}$**

Suppose that $(r,\tau,(A,a))$ is a solution to (2.11) on X. Fix $\varepsilon$ so that (5.1) holds. Given $c > c_0$, the next lemma derives an a priori positive lower bound for the versions of the number $r_c$ that are defined by points $p \in U$ and the pair $(A,a)$ if $|a|(p)$ is positive and p is not in Lemma 6.2's set $\Theta_c$.

**Lemma 6.4**: *There exists* $\kappa > 1$ *and given* $\delta \in (0, 1]$, *there exists* $\kappa_\delta \ge \kappa\delta^{-1}$*; and these are such that the following is true: Suppose that* $r > \kappa_\delta$ *and that* $(A,a)$ *is a pair of connection on* P *and section of* $T^*X \otimes (P\times_{SO(3)} \mathfrak{S})$ *that obeys (2.11) for some* $\tau \in [0, 1]$. *Assume that (5.1) holds with* $\varepsilon \le \delta^{-1}$ *and that* $\sup_X |a| \le \delta^{-1}$. *Fix* $c > \kappa z_U$ *and let* $\Theta_c$ *denote the set from Lemma 6.2. Suppose that* $p \in U$ *is a point in* U *where* $|a|$ *is no smaller than* $\delta$ *and where the distance to* $\Theta_c$ *is no smaller than* $\delta$. *Then the numbers* $r_{c\wedge}$ *and* $r_{cF}$ *that are defined by* p *and* $(A,a)$ *are both greater than* $\kappa_\delta^{-1}$.

***Proof of Lemma 6.4***: The proof has six parts. Their arguments assume that $c$ is greater than 1000 and large enough so that the Propositions 3.1-3.3 and the $\mu = \frac{1}{10}$ version of Proposition 4.1 can be invoked when needed. The value of $c$ is assumed to be 1000 times greater than Proposition 5.1's version of $\kappa$ and also large enough to use in Lemma 6.3.

The proof distinguishes the versions of Section 3a's functions K and N that are defined by distinct points in U: With q being a point in U, the notation uses $K_{(q)}$ and $N_{(q)}$ to denote the versions that are defined by a q and the pair $(A,a)$. With $r \in (0, c_0^{-1}]$ fixed, the notation uses $N_{(\cdot)}(r)$ to denote the function on U that is defined by the rule $q \to N_{(q)}(r)$.

*Part 1*: Let p for the moment denote any given point in U. Define $r_{c\lozenge}$ as in (3.1) but with $c$ replaced by $c^2$. This is to say that $r_{c\lozenge}$ is the largest $r \in (0, c_0^{-1}]$ such that

$$r^4 \int_{B_r} |a \wedge a|^2 \le c^{-4}\ . \tag{6.7}$$

The number $r_{c\lozenge}$ is less than the corresponding version of $r_{c\wedge}$. It is also less than the corresponding version of $r_{cF}$ if $r_{c\lozenge} < \frac{1}{10}\text{dist}(p,\Theta_c)$, this being a consequence of Lemma 6.2.

In any event, if both $r_{c\wedge}$ and $r_{c\mathrm{F}}$ are greater than $r_{c\Diamond}$, then Proposition 4.1 can be invoked with $r_c = r_{c\Diamond}$.

The notation that follows uses $r_{c\Diamond p}$ to denote the version of $r_{c\Diamond}$ that is defined by $(A, a)$ and a given point $p \in U$ when the dependence on the chosen point p is germane to the discussion.

Suppose now that p is a point in U where the corresponding version of $r_{c\Diamond p}$ is less than $\frac{1}{1000}\,\mathrm{dist}(p, \Theta_c)$. If q is a point in U with distance less than $9\,r_{c\Diamond p}$ from p, then $r_{c\Diamond q}$ is no greater than $10 r_{c\Diamond p}$ because the ball of radius $10 r_{c\Diamond p}$ centered at q contains the ball of radius $r_{c\Diamond p}$ centered at p. This being the case, it follows that q's distance to $\Theta_c$ is no smaller than $99\,r_{c\Diamond q}$ and so Lemma 6.2 can be invoked to conclude that $r_{c\Diamond q}$ is no greater than the q and $(A, a)$ version of the number $r_{c\mathrm{F}}$. This implies in turn that Proposition 4.1 can be invoked using q in lieu of p and using $r_c = r_{c\Diamond q}$.

*Part 2*: Suppose again that p is a point in U where $r_{c\Diamond p} \le \frac{1}{1000}\,\mathrm{dist}(p, \Theta_c)$. As explained directly, there exists a point to be denoted by q with distance at most $3 r_{c\Diamond p}$ from p with the following property:

$$\textit{If } q' \textit{ has distance } 2 r_{c\Diamond q} \textit{ or less from } q, \textit{ then } r_{c\Diamond q'} \ge \tfrac{1}{100}\, r_{c\Diamond q}\,. \tag{6.8}$$

To see that this is so, suppose for the sake of argument that (6.8) does not hold for the case q = p. Let $p_1$ denote a point in the closure of radius $2 r_{c\Diamond p}$ ball centered at p where the function $r_{c\Diamond(\cdot)}$ takes on its smallest value. Keep in mind in this regard that the function $r_{c\Diamond(\cdot)}$ is continuous, this being a direct consequence of its definition in (6.7). Either (6.8) holds with $q = p_1$ or not. If not, take $p_2$ to be a point in the closure of the radius $2\,r_{c\Diamond p_1}$ ball centered on $p_1$ where the function $r_{c\Diamond(\cdot)}$ takes on its minimum. Either (6.8) holds for $q = p_2$ or not. If not, then continue in this vein to generate an in turn points $\{p_1, p_2, \ldots\}$ until a point that obeys (6.8) is found. Note that a point obeying (6.8) is reached sooner of later because the value of $r_{c\Diamond(\cdot)}$ for the point selected at any given iteration is less than $\frac{1}{100}$ times the value of $r_{c\Diamond(\cdot)}$ at the previous point. Meanwhile, $r_{c\Diamond(\cdot)}$ is in any event no smaller than $c_0 |a|_\infty^{-1} r^{-1} c^{-1}$ with $|a|_\infty$ denoting $\sup_X |a|$.

Let N denote the number of iterations that are needed to find a point that obeys (6.8). The latter point has distance at most $2\,(\sum_{k=0,1,\ldots N-1} (\frac{1}{100})^k)\, r_{c\Diamond p}$ from p because the distance from the k'th point in the sequence to the next point in the sequence for $k \le N-1$ is no greater than $2\,(\frac{1}{100})^k\, r_{c\Diamond p}$. Since $\sum_{k=0,1,\ldots N-1} (\frac{1}{100})^k$ is less than $\frac{100}{99}$ and since twice this is less than 3, it follows that there is a point that obeys (6.8) with distance less than $3r_{c\Diamond p}$ from the point p.

*Part 3*: Let q denote a point with distance at most $3r_{c\Diamond p}$ from p where (6.8) is obeyed. Given $\varepsilon \in (0, 1]$, reintroduce Proposition 4.1's number $\kappa_{*\varepsilon}$. The following observation plays a central role in what is to come:

*Fix $c > c_0$ and then fix $\varepsilon \le c_0^{-1} c^{-4}$. If $r_{c\Diamond q} \le c_0^{-1} \kappa_{*\varepsilon}^{-1}$, then there is a point with distance at most $\frac{3}{2} r_{c\Diamond q}$ from q where $N_{(\cdot)}(\frac{1}{100} r_{c\Diamond q}) \ge \min(\kappa_{*\varepsilon}^{-1}, c^{-2})$.*

(6.9)

The proof of (6.9) has three steps.

<u>Step 1</u>: Suppose that $\mu \in (0, \frac{1}{4}]$ has been chosen. The radius $r_{c\Diamond q}$ ball centered at q can be covered by $c_0$ balls of radius $\frac{1}{100}(1-2\mu) r_{c_q}$ centered at points in the radius $\frac{3}{2} r_{c\Diamond q}$ ball centered at q. Fix such a cover, use $\mathfrak{U}$ to denote the set of the center points of the balls in the cover and use $\mathfrak{N}$ to denote the number of points in $\mathfrak{U}$.

<u>Step 2</u>: Fix $\varepsilon \in (0, 1]$ for the moment. It follows from what is said in Part 1 that Proposition 4.1 can be invoked with $\mu$ set equal to $\frac{1}{10}$, with $c \ge c_0$, with the given value of $\varepsilon$, with the number $r_c$ set equal to $\frac{1}{100} r_{c\Diamond q}$ and with the point p replaced by any point from the set $\mathfrak{U}$ where the function $N_{(\cdot)}(\frac{1}{100} r_{c\Diamond q})$ is less than both $\kappa_{*\varepsilon}^{-1}$ and $c^{-2}$.

<u>Step 3</u>: If $N_{(\cdot)}(\frac{1}{100} r_{c\Diamond q})$ is less than both $\kappa_{*\varepsilon}^{-1}$ and $c^{-2}$ at every point in $\mathfrak{U}$, then Proposition 4.1 can be brought to bear to obtain an $\varepsilon$ upper bound for the integral of $r^4|a \wedge a|^2$ on the radius $(1-\mu)\frac{1}{100} r_{c\Diamond q}$ ball centered at each point of $\mathfrak{U}$. These bounds would lead in turn to an $\mathfrak{N}\varepsilon$ upper bound for the integral of $r^4|a \wedge a|^2$ over the radius $r_{c\Diamond q}$ ball centered at q. Since $\mathfrak{N} \le c_0$, such a bound leads to nonsense if $\varepsilon \le c_0^{-1} c^{-4}$ because the integral of $r^4|a \wedge a|^2$ over this ball is exactly $c^{-4}$.

*Part 4*: Use $r_{*p}$ to value of $r \in (0, c_0^{-1}]$ where $r K_{(p)}(r) r = z_U^{-1}$ if such r exists, or $c_0^{-1}$ if $r K_{(p)}(r) r < z_U^{-1}$ on $(0, c_0^{-1}]$. Fix $\nu \in (0, 1]$ for the moment and suppose that $r_0 \in (0, c_0^{-1}]$ is such that $N_{(p)}(r_0) \ge \nu$. Let $\hat{\kappa}$ denote the $\varepsilon = \nu$ version of Proposition 5.1's number $\kappa_\varepsilon$. Use Lemma 6.3 and Proposition 5.1 to see that $N_{(p)}(\cdot) \ge \hat{\kappa}^{-1}$ on the interval between the $\max(r_0, r_{*p})$ and $c_0^{-1}$ if $r_{*p} < \frac{1}{100}\mathrm{dist}(p, \Theta_c)$. Suppose that this is so. Fix $r \in [\max(r_0, r_{*p}), c_0^{-1}]$ and $s \in [r, c_0^{-1}]$. The stated lower bound for $N_{(p)}$ with the formula in (3.6) implies that

$$\mathrm{K}_{(p)}(s) \ge (\tfrac{s}{r})^{1/\hat{\kappa}}\, \mathrm{K}_{(p)}(r)\ . \tag{6.10}$$

Let $|a|_\infty = \sup_X |a|$. Since $\mathrm{K}_{(p)}(s)$ is in any event no greater than $c_0 |a|_\infty$, the bound in (6.10) when turned on its head implies that

$$\mathrm{K}_{(p)}(r) \le c_0 |a|_\infty\, r^{1/\hat{\kappa}} \quad \textit{if}\ \ r \in [\max(r_0, r_{*p}), c_0^{-1}) \tag{6.11}$$

With regards to the assumption $r_{*p} \le \frac{1}{100}\operatorname{dist}(p, \Theta_c)$, keep in mind that $\mathrm{K}_{(p)}$ is increasing and so if $|a|$ is positive at p, then

$$r_{*p} \le \frac{1}{\sqrt{2\pi}}\, \frac{1}{|a|(p)}\, z_U^{-1}\, r^{-1}. \tag{6.12}$$

The number $r_{*p}$ is upper semi-continuous in a crude sense. To say more, suppose that p is a given point in U. Fix $r \in (0, c_0^{-1}]$. The integral of $|a|^2$ over the radius 2r ball centered at p is no smaller than $\frac{15}{4}(1 - c_0 r^2)\,\mathrm{K}_p(r)^2 r^4$ because $\mathrm{K}_{(p)}$ is increasing. If q is in the radius r ball centered at p, then the integral of $|a|^2$ over the radius 3r ball centered at q is no smaller than this same number because the radius 2r ball centered at p is contained in the radius 3r ball centered at q. Meanwhile, the latter integral is no greater than $\frac{81}{4}(1 + c_0 r^2)\,\mathrm{K}_{(q)}(3r)\, r^4$ because $\mathrm{K}_{(q)}$ is increasing. Comparing this upper bound with the previous lower bound leads to the observation that

$$\mathrm{K}_{(q)}(3r) \ge \tfrac{5}{27}(1 - c_0 r^2)\,\mathrm{K}_{(p)}(r)\ . \tag{6.13}$$

This inequality implies that $(6r)\,\mathrm{K}_{(q)}(6r) > r\,\mathrm{K}_{(p)}(r)$ if $r \le c_0^{-1}$ because $\frac{5}{27}$ is less than $\frac{1}{6}$ and $\mathrm{K}_{(q)}$ is increasing. By way of an example, if r is greater than $r_{*p}$ then 6r is greater than $r_{*q}$.

*Part 5*: Fix $\delta \in (0, 1]$ and assume that $p \in U$ is such that $|a|(p) \ge \delta$ and that $\operatorname{dist}(p, \Theta_c) > \delta$ also. If $r > c_0 \delta^{-2}$, then (6.12) implies that $r_{*p} \le \frac{1}{1000}\operatorname{dist}(p, \Theta_c)$. It then follows from (6.13) that $r_{*p'} \le \frac{1}{100}\operatorname{dist}(p', \Theta_c)$ if $\operatorname{dist}(p', p) \le c_0^{-1}\delta$. In particular, if it is case that $r_{c\Diamond p} < c_0^{-1}\delta$ and if p´ is a point in the radius $10\, r_{c\Diamond p}$ ball centered at p, then $r_{*p'}$ will be less than $\frac{1}{100}\operatorname{dist}(p', \Theta_c)$.

With the preceding understood, assume in what follows that $r_{c\Diamond p}$ is less than $c_0^{-1}\delta$ and that p´ is a point in the radius $10\, r_{c\Diamond p}$ ball centered at p with with $\mathrm{N}_{p'}(r_{c\Diamond p'}) \ge \nu$. There is a p´ version of (6.11) because $r_{*p'}$ is less than $\frac{1}{100}\,\mathrm{dist}(p', \Theta_c)$ and the latter asserts the bound $K_{p'}(r) \le c_0 |a|_\infty r^{1/\hat{\kappa}}$ for any given $r \in [\max(r_{c\lozenge p'}, r_{*p'}), c_0^{-1}]$.

To see the significance of the preceding bound, suppose first that $r_{c\lozenge p}$ is greater than $6 r_{*p}$. It then follows from what is said at the end of Part 4 that $r_{c\lozenge p}$ is also greater than $r_{*p'}$. It is also the case that $r_{c\lozenge p'} \le 11\, r_{c\lozenge p}$ because the ball of radius $r_{c\lozenge p}$ centered at p is contained in the ball of radius $11 r_{c\lozenge p}$ centered at p´. It follows as a consequence that the bound $K_{(p')}(r) \le c_0 |a|_\infty r^{1/\hat{\kappa}}$ holds in particular at $r = 11\, r_{c\lozenge p}$. Since $K_{(p')}$ is increasing, the integral of $|a|^2$ over the radius $11\, r_{c\lozenge p}$ ball centered at p´ is therefore no greater than $c_0 |a|_\infty^2\, r_{c\lozenge p}^{4+2/\hat{\kappa}}$. As just noted, the radius $r_{c\lozenge p}$ ball centered at p is contained in the radius $11\, r_{c\lozenge p}$ ball centered at p´, and so this same upper bound holds for the integral of $|a|^2$ over the radius $r_{c\lozenge p}$ ball centered at p. Meanwhile, the latter integral can not be smaller than $2\pi^2 (1 - c_0 r_{c\lozenge p}^{\,2})\, r_{c\lozenge p}^{\,4} |a|^2(p)$ because $K_{(p)}$ is increasing. These upper and lower bounds are compatible only in the event that

$$r_{c\lozenge p} \ge c_0^{-1} \left(\frac{|a|(p)}{|a|_\infty}\right)^{\hat{\kappa}} . \tag{6.14}$$

The next paragraph explains why $r_{c\lozenge p}$ is always greater than $6r_{*p}$ if $r \ge c_0 |a_\infty|^{\hat{\kappa}}\, \delta^{-(1+\hat{\kappa})}$.

To prove the $r_{c\lozenge p} > 6r_{*p}$ bound, suppose for the sake of argument that $r_{c\lozenge p}$ is not greater than $6 r_{*p}$. Since p´ has distance at most $10 r_{c\lozenge p}$ from p, it has distance at most $60\, r_{*p}$ from p. The ball of radius $r_{*p}$ centered at p is contained in the ball of radius $61 r_{*p}$ centered at p´ and $61\, r_{*p}$ is greater than $r_{*p'}$ because 61 is greater than 6. With this understood, use the p´ version of (6.11) to see that $K_{(p')}(61\, r_{*p}) \le c_0 |a|_\infty\, r_{*p}^{\,1/\hat{\kappa}}$. This leads to a $c_0 |a|_\infty^2\, r_{*p}^{\,4+2/\hat{\kappa}}$ upper bound for the integral of $|a|^2$ over the radius $61\, r_{*p}$ ball centered at p´ and thus a $c_0 |a|_\infty^2\, r_{*p}^{\,4+2/\hat{\kappa}}$ upper bound for the integral of $|a|^2$ on the radius $r_{*p}$ ball centered at p. Since the latter is in no event less than $\frac{1}{4}(1 - c_0 r_{*p}^{\,2}) |a|^2(p)\, r_{*p}^{\,4}$, this upper bound holds only in the event that

$$r_{*p} \ge c_0^{-1} \left(\frac{|a|(p)}{|a|_\infty}\right)^{\hat{\kappa}} . \tag{6.15}$$

Given that $r_{*p} K_{(p)}(r_{*p})\, r$ is supposed to be $z_U^{-1}$ and given that $K_{(p)}(r_{*p}) \ge \sqrt{2}\,\pi |a|(p)$, the lower bound in (6.15) can hold only in the event that

$$r \le c_0 |a_\infty|^{\hat{\kappa}} \left(\frac{1}{|a|(p)}\right)^{-(1+\hat{\kappa})} . \tag{6.16}$$

This lower bound is nonsensical if $r$ is larger than $c_0|a_\infty|^{\hat{\kappa}}\,\delta^{-(1+\hat{\kappa})}$ because $|a|(p)$ is assumed to be greater than $\delta$.

*Part 6*: Fix $c > c_0$ and $\varepsilon < c_0 c^{-4}$ so that (6.9) is true if $r_{c\Diamond p}$ is less than both $c^{-2}$ and $c_0^{-1}\kappa_{*\varepsilon}^{-1}$. Take $\nu$ to be the minimum of $c^{-2}$ and what is denoted by $\kappa_{*\varepsilon}^{-1}$. Let q be as described in (6.9) and let p´ denote a point in the radius $\frac{3}{2}r_{c\Diamond q}$ ball centered at q with $N_{(p')}(\frac{1}{100}r_{c\Diamond q}) \geq \nu$. Since q is in the radius $3r_{c\Diamond p}$ ball centered at p, it follows from the constructions that p´ is in the radius $4r_{c\Diamond p}$ ball centered at p. If $r_{c\Diamond p}$ is also less than $c_0^{-1}\delta$ and if $r > c_0|a_\infty|^{\hat{\kappa}}\,\delta^{-(1+\hat{\kappa})}$, then the last paragraph in Part 5 asserts that $r_{c\Diamond p}$ is greater than $6r_{*p}$ and then (6.11) asserts a lower bound for $r_{c\Diamond p}$. A lower bound for $r_{c\Diamond p}$ is an automatic lower bound for p's version of $r_{c\wedge}$ and it follows from what is said in Part 1 that the minimum of $r_{c\Diamond p}$ and $\frac{1}{100}\mathrm{dist}(p,\Theta_c)$ is a lower bound for p's version of $r_{cF}$.

**c) Proof of Proposition 6.1**

The proof of this proposition has ten parts.

*Part 1*: Fix $c > 1$ to be larger than the number $\kappa$ that appears in Lemma 6.4. The number $c$ should also be chosen to be larger than the versions of $\kappa$ that are in Propositions 3.1 and 3.2 so that the latter can be invoked, and it should be larger than the $\mu = \frac{1}{1000}$ versions of the number $\kappa_\mu$ that are in Propositions 3.3, Lemma 4.2 and Lemma 4.3. It is also necessary to take $c$ greater than what is needed to invoke Lemma 6.3. Additional lower bounds on $c$ are given below, but all a priori determined by the geometric data.

Let $\{(r_n, \tau_n, (A_n, a_n))\}_{n\in\{1,2,\ldots\}}$ denote a weakly convergent sequence of solutions to (2.11). Given $n \in \{1, 2, \ldots\}$, let $\Theta_{cn}$ denote the version of Lemma 6.2's set $\Theta_c$ that is obtained using $r = r_n$, $\tau = \tau_n$ and $(A, a) = (A_n, a_n)$. This set has at most $E^2c^2$ elements with $E^2$ being an upper bound for the $f = 1$ version of the sequence in Item a) of the second bullet in (6.1). Item a) of the second bullet in (6.1) implies that the sequence $\{\Theta_{cn}\}_{n\in\{1,2,\ldots\}}$ converges in the following sense: Let $\partial U$ denote the boundary of the closure of U. There is a set $\Theta \subset U$ with at most $E^2c^2$ elements such that

$$\lim_{n\in\{1,2,\ldots\}}\left(\sup_{p\in\Theta_{c_n}}\mathrm{dist}(p,\Theta\cup\partial U) + \sup_{p\in\Theta}\mathrm{dist}(p,\Theta_{cn})\right) = 0. \tag{6.17}$$

Suppose that p is a point in $U-\Theta$ where the function $|\hat{a}_\Diamond|$ is positive. Fix a subsequence $\Xi_p \subset \{1, 2, \ldots\}$ such that

$$\lim_{n\in\Xi_p} |a_n|(p) = |\hat{a}_\Diamond|(p) .$$

(6.18)

Since p is disjoint from both $\Theta\cup\partial U$, there exists a positive number that is smaller than half of both $|\hat{a}_\Diamond|(p)$ and the distance from p to each point in $\Theta\cup\partial U$. Let $\delta$ denote such a number. With $\delta$ chosen, there exists $n_\delta \in \Xi_p$ such that $\text{dist}(p,\Theta_{c_n}) > \delta$ and $|a_n|(p) > \delta$ when $n \in \Xi_p$ is greater than $n_\delta$.

If $n \in \Xi_p$ is sufficiently large, then the preceding observations license an appeal to the version of Lemma 6.4 with the input data $r = r_n$, $\tau = \tau_n$ and $(A, a) = (A_n, a_n)$. The appeal to Lemma 6.4 results in a number greater than 1 to be denoted by $\kappa_\delta$ that depends on $\delta$ but not p and has the following signficance: If $n \in \Xi_p$ is greater than $n_\delta$ and $r_n$ is greater than $\kappa_\delta$, then the versions of $r_{c\wedge}$ and $r_{cF}$ when defined using p and the data set $(r_n, \tau_n, (A_n, a_n)$ are bounded from below $\kappa_\delta^{-1}$.

*Part 2*: Supposing that $n \in \Xi_p$, let $K_n$ and $N_n$ denote the version of the functions $K$ and $N$ from Section 3a that are defined using p and $(r_n, (A_n, a_n))$. Given $\varepsilon \in (0,]$, there exists $r_\varepsilon \in (0, c_0^{-1}]$ such that all sufficienctly large $n \in \Xi_p$ versions of $K_n$ and $N_n$ obey

- $K_n > \frac{1}{\sqrt{2\pi}} |\hat{a}_\Diamond|(p) - \varepsilon$ *on* $[0, c_0^{-1})$ .
- $N_n < \varepsilon$ *on* $(0, r_\varepsilon)$ .

(6.19)

To prove the top bullet, use the identity $K_n(0) = \frac{1}{\sqrt{2\pi}} |a_n|(p)$ with (6.18) to see that $K_n(0)$ differs from $\frac{1}{\sqrt{2\pi}} |\hat{a}_\Diamond|(p)$ by at most $\varepsilon$ when n is large. The top bullet in (6.19) follows directly from this observation because $K_n$ is an increasing function. To prove the lower bullet, let $r_{*p(n)}$ denote the $(r_n, (A_n, a_n))$ version of the number $r_{*p}$ that is defined in Part 4 of Section 6b. It follows from (6.12) that $r_{*p(n)} < c_0 \delta^{-1} r_n^{-1}$ when n is large. In particular, this will be less than $\frac{1}{1000} \text{dist}(p, \Theta_{c_n})$ when n is large. With the preceding understood, it follows that if n is large and $r_0$ is such that $N_n(r_0) \geq \varepsilon$, then Lemma 6.3 and Proposition 5.1 can be invoked. These supply a number $\kappa_\varepsilon > 1$ whose inverse is a lower bound for $N_n(r)$ if r is greater than both $r_0$ and $r_{*p(n)}$ but less than $c_0^{-1}$. The latter bound with (3.6) and the $|a_n|(p) > \delta$ bound when $n > n_\delta$ leads in turn to the bound

$$K_n(r) \geq \left(\frac{r}{r_\varepsilon}\right)^{1/\kappa_\varepsilon} \delta$$

(6.20)

in the event that r is greater than both $r_\varepsilon$ and $r_{*p}$ but less than $c_0^{-1}$. The bound in (6.20) leads to a lower bound for $r_\varepsilon$ that is independent of $n \in \Xi_p$ if n is large. This is because $r_{*p}$ is less than $c_0 \delta^{-1} r_n^{-1}$ and because $K_n$ is, in any event, no greater than $c_0$ times the upper bound for the sequence in the third bullet of (6.1).

*Part 3*: It follows from what is said in Part 2 that there exists $r_{p\diamond} \in (0, \frac{1}{100} c^{-1}]$ such that if $n \in \Xi_p$ is large, then $N_n(r_{p\diamond}) \le c^{-4}$. What is said in Part 1 allows $r_{p\diamond}$ to be chosen so as to be less than half the versions of $r_{c\wedge}$ and $r_{cF}$ that are defined by p and each sufficiently large n version of $(r_n, (A_n, a_n))$. It follows as a consequence that Propositions 3.1 and 3.2 can be invoked with $r_c = 2r_{p\diamond}$ when $n \in \Xi_p$ is large. The second bullet of Proposition 3.2 for such n implies that $|a_n| \ge \frac{3}{4}|a_n|(p)$ on the radius $r_{p\diamond}$ ball centered at p when n is large.

The preceeding conclusions lead to the following observation:

$$\textit{The points in } U-\Theta \textit{ where } |\hat{a}_\diamond| \textit{ is positive constitute an open set.} \tag{6.21}$$

To prove this, suppose that p is a point in $U-\Theta$ where $|\hat{a}_\diamond|$ is positive. Define $r_{p\diamond}$ as just instructed and let q denote any given point in the radius $r_{p\diamond}$ ball centered at p. The fact that $|a_n|(q) \ge \frac{3}{4}|a_n|(p)$ if $n \in \Xi_p$ is large implies that $\limsup_{n\in\{1,2,\ldots\}} |a_n|(q) > 0$, and this in turn implies that $|\hat{a}_\diamond|(q) > 0$.

*Part 4*: Supposing as before that p is a point in $U-\Theta$ where $|\hat{a}_\diamond|$ is positive, define the number $r_{p\diamond}$ as in Part 3. If $n \in \Xi_p$ is large, the number $\kappa_n(r_{p\diamond})$ will be greater than $\delta$ and therefore $2r_{p\diamond}\kappa_n(2r_{p\diamond})\, r_n$ will be very much larger than the $\mu = \frac{1}{100}$ version of the numbers $\kappa_\mu$ that appear in Proposition 3.3, Lemma 4.2 and Lemma 4.3. This being the case, define $r_c = 2r_{p\diamond}$ and $\kappa_{\ddagger n} = \kappa_n(r_c)$. With $\kappa_{\ddagger n}$ understood, set $a_{\ddagger n} = \kappa_{\ddagger n}^{-1} a_n$. Invoke Lemma 4.2 when n is large to write $a_{\ddagger n}$ as in (4.4) on $B_{(1-\mu/16)r_c}$ and denote by $\nu_n$, $\sigma_{\ddagger n}$ and $\mathfrak{a}_n$ the corresponding versions of $\nu$, $\sigma_\ddagger$ and $\mathfrak{a}$. To be sure of the notation, this is to say that

$$a_{\ddagger n} = \nu_n \sigma_{\ddagger n} + \mathfrak{a}_n \tag{6.22}$$

Let $\Xi'_p \subset \Xi_p$ denote the subsequence where the decomposition in (6.22) can be made. As explained momentarily,

$$\lim_{n\in\Xi'_p} |\mathfrak{a}_n| = 0 \;\; \textit{on radius } \tfrac{3}{4} r_c \textit{ ball centered at } p. \tag{6.23}$$

The proof of (6.23) has three steps. The notation in these steps has $r_{\ddagger n}$ denoting $\kappa_{\ddagger n} r_n$.

<u>Step 1</u>: The $r = r_n$ and $(A, a) = (A_n, a_n)$ version of (2.5) on $B_{(1-\mu/16)r_c}$ when written in terms of $\nu_n$, $\sigma_{\ddagger n}$ and $\mathfrak{a}_n$ leads to a differential equation for $|\mathfrak{a}_n|^2$ that reads

$$\tfrac{1}{2}\, d^{\dagger}d\,|\mathfrak{a}_n|^2 + |\nabla_{A_n}\mathfrak{a}_n|^2 + 4\,r_{\ddagger n}{}^2|\nu_n|^2|\mathfrak{a}_n|^2 + 2\,\nu_{n\beta}\,\langle \nabla_{A_n,\alpha}\sigma_{\ddagger n}(\nabla_{A_n,\alpha}\mathfrak{a}_n)_\beta\rangle + \mathrm{Ric}(\langle \mathfrak{a}_n \otimes \mathfrak{a}_n\rangle) = 0\ . \tag{6.24}$$

To explain the notation, an orthonormal basis for T*X on $B_{r_c}$ has been chosen and the components of $\nu_n$ and $\mathfrak{a}_n$ with respect to this basis are $\{\nu_{n\beta}\}_{\beta=1,2,3,4}$ and $\{\mathfrak{a}_{n\beta}\}_{\beta=1,2,3,4}$. This basis is also used to define the components of the covariant derivatives. As in previous sections, repeated indices are summed. Keep in mind when deriving (6.24) that both $\langle\sigma_{\ddagger n}\,\mathfrak{a}_n\rangle$ and $\nu_{n\beta}\,\mathfrak{a}_{n\beta}$ are everywhere zero.

<u>Step 2</u>: Use the function $\chi$ to construct a function on $B_{r_c}$ that equals 1 where the distance to p is less than $\frac{7}{8}\,r_c$ and equals 0 where the distance to p is greater than $\frac{15}{16}\,r_c$. Since $\mu = \frac{1}{100}$, this function has compact support in the radius $(1-\mu)\,r_c$ ball centered at p. Use $\chi_p$ to denote this function. Given a point q in the radius $\frac{7}{8}\,r_c$ ball centered at p, let $G_q$ denote the Dirichelet Green's function on $B_{r_c}$ with pole at q for the operator $d^{\dagger}d$. Multiply both sides of (6.24) by $\chi_p G_q$ and integrate the resulting inequality over $B_{r_c}$. Having done so, integrate by parts and use the fact that $G_q$ is positive but less than $c_0\,\frac{1}{\mathrm{dist}(q,\cdot)^2}$ and that $|\nabla G_q| \le c_0\,\frac{1}{\mathrm{dist}(q,\cdot)^3}$ to derive the inequality

$$|\mathfrak{a}_n|^2(q) \le c_0\, r_c^{-4} \int_{B_{(1-\mu)r_c}} |\mathfrak{a}_n|^2 \; + c_0 \int_{B_{(1-\mu)r_c}} \frac{1}{\mathrm{dist}(q,\cdot)^2}\,|\nabla_{A_n}\sigma_{\ddagger n}|^2\ . \tag{6.25}$$

<u>Step 3</u>: Invoke the $r = r_n$ and $(A_n, a_n)$ version of Lemma 4.3. This lemma bounds $L^2$ norm of $|\mathfrak{a}_n|$ by $c_\mu r_{\ddagger n}{}^{-4}\,(r_c\, r_{\ddagger n})^{-2}\, c^{-1}$ and thus the left most term in (6.22) by $c_\mu\,(r_c\, r_{\ddagger n})^{-6}\, c^{-1}$. Lemma 4.3 and the second bullet in (2.1) jointly bound the right most integral in (6.25) by $c_\mu\,(r_c\, r_{\ddagger n})^{-2}\, c^{-1}$. Use these bounds in (6.25) to see that $|\mathfrak{a}_n| \le c_\mu\,(r_c\, r_{\ddagger n})^{-1}\, c^{-1/2}$. Since $\kappa_{\ddagger n}$ is greater than $\delta$ when n is large, this bound implies what is asserted by (6.23).

*Part 5*: The $r = r_n$ and $(A, a) = (A_n, a_n)$ version of (2.5) on $B_{(1-\mu/16)\,r_c}$ also leads to an equation for $\nu_n$ that reads

$$\nabla^{\dagger}\nabla\nu_n + \nu_n\,|\nabla_{A_n}\sigma_{\ddagger n}|^2 + \langle \nabla_{A_n}{}^{\dagger}\nabla_{A_n}\sigma_{\ddagger n}\,\mathfrak{a}_n\rangle + 2\,\langle\nabla_{A_n,\alpha}\sigma_{\ddagger n}\nabla_{A_n,\alpha}\mathfrak{a}_n\rangle + \mathrm{Ric}(\nu_n) = 0\ . \tag{6.26}$$

The same version of (6.26) leads to an equation for $\langle\nabla_{A_n}{}^{\dagger}\nabla_{A_n}\sigma_{\ddagger n}\,\mathfrak{a}_n\rangle$ that writes it as

$$-\,\nabla^{\dagger}\nabla\nu_{n\alpha}\langle\mathfrak{a}_{n\alpha}\,\mathfrak{a}_n\rangle + \ \mathfrak{w}_n$$

(6.27)

with $\mathfrak{w}_n$ being a 1-form on $B_{(1-\mu/16)r_c}$ whose norm obeys the bound

$$|\mathfrak{w}_n| \le c_0 |\nabla \nu_n| |\mathfrak{a}_n| ( |\nabla_{A_n} \sigma_{\ddagger n}| + |\nabla_{A_n} \mathfrak{a}_n| ) \,.$$

(6.28)

What with (6.27) and the fact that $|\mathfrak{a}_n|$ is small, it follows that (6.26) can be written schematically as

$$\nabla^{\dagger}\nabla \nu_n + \mathrm{Ric}(\nu_n) = \mathfrak{e}_n \,,$$

(6.29)

with the norm of $\mathfrak{e}_n$ obeying

$$|\mathfrak{e}_n| \le c_0 |\nabla \nu_n| |\mathfrak{a}_n| ( |\nabla_{A_n} \sigma_{\ddagger n}| + |\nabla_{A_n} \mathfrak{a}_n| ) + c_0 |\nabla_{A_n} \sigma_{\ddagger n}| |\nabla_{A_n} \mathfrak{a}_n| + c_0 |\mathfrak{a}_n|^2 \,.$$

(6.30)

Note for future reference that the combination $|\nabla \nu_n||\mathfrak{a}_n|$ that appears in (6.28) and (6.30) is no greater than $c_0 c^{-3/2}$, this being a consequence of (4.7), what is said by the second bullet of the $r = r_n$ and $(A, a) = (A_n, a_n)$ version of Proposition 3.3 and the bound $|\mathfrak{a}_n| \le c_\mu (r_c r_{\ddagger n})^{-1} c^{-1/2}$ from Step 3 in Part 4.

*Part 6*: Write $\nu_n$ on $B_{(1-\mu)r_c}$ as $\nu_n = t_n + s_n$ where $s_n$ is the solution to the equation

$$\nabla^{\dagger}\nabla s_n + \mathrm{Ric}(s_n) = \mathfrak{e}_n$$

(6.31)

on $B_{(1-\mu)r_c}$ that equals zero on the $\partial B_{(1-\mu)r_c}$. There exists a unique solution if $r_c \le c_0^{-1}$ and since $r_c$ is in any event less than $2c^{-1}$, this will be the case if $c > c_0$. The Dirichelet Green's function for $\nabla^{\dagger}\nabla + \mathrm{Ric}(\cdot)$ on $B_{(1-\mu)r_c}$ can be used to see that the norm of $s_n$ is bounded by

$$|s_n| \le c_0 \int_{B_{(1-\mu)r_c}} \frac{1}{\mathrm{dist}(q,\cdot)^2} |\mathfrak{e}_n| \,.$$

(6.32)

To say something about the left hand side of (6.32), introduce by way of notation

$$I_n = \int_{B_{(1-\mu)r_c}} \frac{1}{\mathrm{dist}(q,\cdot)^2} ( |\nabla_{A_n} \sigma_{\ddagger n}|^2 + |\nabla_{A_n} \mathfrak{a}_n|^2 ) \,.$$

(6.33)

Use the bound $|\nabla \nu_n||\mathfrak{a}_n| \le c_0\, c^{-3/2}$ from the last paragraph of Part 5 with the $c_\mu(r_c r_{\ddagger n})^{-1} c^{-1/2}$ bound on $|\mathfrak{a}_n|$ from Step 3 in Part 4 in (6.30) to bound the integral that appears on the right hand side of (6.32) by $c_\mu\,(c^{-3/2}\sqrt{\mathrm{I}_n} + \mathrm{I}_n + (r_c r_{\ddagger n})^{-2} c^{-1})$. Meanwhile, the second bullet in (2.1) with the $r = r_n$ and $(A, a) = (A_n, a_n)$ version of Lemma 4.2 leads to a $c_\mu(r_c r_{\ddagger n})^{-2} c^{-1}$ bound for $\mathrm{I}_n$. Use this bound to see that $|s_n| \le c_\mu\, c^{-2}(r_c r_{\ddagger n})^{-1}$.

*Part 7*: Take $n \in \Xi'_p$ so that $s_n$ and $t_n$ are defined. What is denoted by $t_n$ obeys the homogeneous equation $\nabla^\dagger\nabla t_n + \mathrm{Ric}(t_n) = 0$ on $B_{(1-\mu)r_c}$ and it equals $\nu_n$ on $\partial B_{(1-\mu)r_c}$. As explained in the next paragraph, the sequence $\{t_n\}_{n\in\Xi'_p}$ converges in the $C^\infty$ topology on compact subsets of $B_{(1-\mu)r_c}$ as $n \in \Xi'_p$ gets ever larger. Let $t_{(p)}$ denote the limit 1-form, This is a nowhere zero solution on $B_{(1-\mu)r_c}$ to the equation $\nabla^\dagger\nabla t + \mathrm{Ric}(t) = 0$.

To prove what was just said about $\{t_n\}_{n\in\Xi'_p}$, note first that that the sequence $\{\kappa_{\ddagger n}\}_{n\in\Xi'_p}$ converges with the square of its limit being

$$r_c^{-3} \int_{\partial B_{r_c}} |\hat{a}_\Diamond|^2 \ . \tag{6.34}$$

This is so because the sequence $\{|a_n|\}_{n\in\{1,2,\ldots\}}$ converges weakly in the $L^2_1$ topology on U to $|\hat{a}_\Diamond|$ and the restriction map from $C^\infty(B_{r_c})$ to $C^\infty(\partial B_{r_c})$ extends to define a compact map from the space of $L^2_1$ functions on $B_{r_c}$ to the space of $L^2$ functions on $\partial B_{r_c}$. Use $\kappa_{\Diamond(P)}$ to denote the positive square root of what is written in (6.34). Since both $\{\mathfrak{a}_n\}_{n\in\Xi'_p}$ and $\{s_n\}_{n\in\Xi'_p}$ converge to zero, it follows that $\{|t_n|\}_{n\in\Xi'_p}$ converges weakly in the $L^2_1$ topology and strongly in the $L^2$ topology $\kappa_{\Diamond(P)}|\hat{a}_\Diamond|$ on $B_{(1-\mu)r_c}$ because the sequence $\{|a_n|\}_{n\in\{1,2,\ldots\}}$ converges weakly in the $L^2_1$ topology and strongly in the $L^2$ topology to $|\hat{a}_\Diamond|$. This strong $L^2$ convergence with the fact that each $n \in \Xi'_p$ version of the 1-form $t_n$ obeys the equation $\nabla^\dagger\nabla t_n + \mathrm{Ric}(t_n) = 0$ implies via standard elliptic regularity techniques that the sequence $\{t_n\}_{n\in\Xi'_p}$ converges in the $C^\infty$ topology on compact subsets of $B_{(1-\mu)r_c}$. The $C^\infty$ convergence implies in turn that the limit, $t_{(p)}$, obeys the equation $\nabla^\dagger\nabla t_{(p)} + \mathrm{Ric}(t_{(p)}) = 0$. The limit is nowhere zero on $B_{(1-\mu)r_c}$ because each large $n \in \Xi'_p$ version of $a_{\ddagger n}$ is described by (4.5) on the large ball $B_{(1-\mu/16)r_c}$.

Granted what was just said, it follows that the sequence $\{|a_n|\}_{n\in\Xi_p}$ converges in the $C^0$ topology on compact subsets of $B_{(1-\mu)r_c}$ to the smooth function $\kappa_{\Diamond(P)}|t_{(p)}|$.

*Part 8*: Suppose now that q is a point in U with distance less than $\frac{1}{8}\,r_{\lozenge p}$ from p. Since $r_{\lozenge p}$ is half of $r_c$, it follows that q's version of $r_{c\wedge}$ and $r_{cF}$ are no smaller than $\frac{15}{16}\,r_c$. This being the case, what with Lemma 6.4, the analysis in Parts 4-7 can be carried out with the point q replacing the point p. Doing so supplies a 1-form $t_{(q)}$ that is defined on a ball centered at q that contains the radius $\frac{1}{4}\,r_c$ ball centered at p, and a subsequence $\Xi_q \subset \{1, 2, \ldots\}$ such that $\{|a_n|\}_{n\in\Xi_q}$ converges in the $C^0$ topology on the aforementioned ball about q to $\kappa_{\lozenge q}|t_{(q)}|$.

Both $\kappa_{\lozenge(p)}|t_{(p)}|$ and $\kappa_{\lozenge(q)}|t_{(q)}|$ are smooth functions on the radius $\frac{1}{4}\,r_c$ ball centered at p. Since $\{|a_n|\}_{n\in\{1,2,\ldots\}}$ converges strongly in the $L^2$ topology to $|\hat{a}_\lozenge|$ it follows that both of these functions must equal $|\hat{a}_\lozenge|$ almost everywhere. Since both are a priori continuous, it follows that these $\kappa_{\lozenge(p)}|t_{(p)}| = \kappa_{\lozenge(q)}|t_{(q)}|$ at each point in the radius $\frac{1}{4}\,r_c$ ball centered at p.

This last observation can be made using any given point in the radius $\frac{1}{4}\,r_c$ ball centered at p, and doing so proves that the sequence $\{|a_n|\}_{n\in\Xi_p}$ converges to $|\hat{a}_\lozenge|$ in the $C^0$ topology on the radius $\frac{1}{4}\,r_c$ ball centered at p. Granted that this is so, the identification of $|\hat{a}_\lozenge|$ with $\kappa_{\lozenge(p)}|t_{(p)}|$ implies that $|\hat{a}_\lozenge|$ is a smooth function on the radius $\frac{1}{4}\,r_c$ ball centered at p.

*Part 9*: The subsequence $\Xi_p$ was specifically choosen so that $\lim_{n\in\Xi_p} |a_n|(p)$ and $\limsup_{n\in\{1,2,\ldots\}} |a_n|(p)$ are equal. However, what is said at the end of Part 8 implies that $\lim_{n\in\Xi_p} |a_n|(q) = \limsup_{n\in\{1,2,\ldots\}} |a_n|(q)$ for all points q in the radius $\frac{1}{4}\,r_c$ ball centered at p. It is actually the case that

$$\lim_{n\in\Xi_p} |a_n|(p) = \limsup_{n\in\{1,2,\ldots\}} |a_n|(q) \quad \textit{at all points } q \in U-\Theta\,. \tag{6.35}$$

That this is so follows by definition at points where $|\hat{a}_\lozenge|$ is zero. It is also true at a dense set of points where $|\hat{a}_\lozenge|$ is positive because $\{|a_n|\}_{n\in\{1,2,\ldots\}}$ converges strongly in the $L^2$ topology to $|\hat{a}_\lozenge|$. Since (6.35) holds on a dense set where $|\hat{a}_\lozenge| > 0$, it then follows from what is said in Part 8 that it must hold at all points where $|\hat{a}_\lozenge| > 0$.

Let $Z \subset U$ denote the set where $|\hat{a}_\lozenge|$ is zero. The identity in (6.35) with what is said in Part 8 implies that $|\hat{a}_\lozenge|$ is smooth on $U-(Z\cup\Theta)$ and that $\{|a_n|\}_{n\in\{1,2,\ldots\}}$ converges to $|\hat{a}_\lozenge|$ in the $C^0$ topology on compact subsets of $U-(Z\cup\Theta)$.

*Part 10*: Suppose next that $p \in Z$. This is to say that $\lim_n |a_n|(p) = 0$. This part of the proof explains why $|\hat{a}_\lozenge|$ is continuous near p and why $\{|a_n|\}_{n\in\{1,2,\ldots\}}$ converges in the $C^0$ topology to $|\hat{a}_\lozenge|$ on a neighborhood of p. That this is so follows from what is said in Parts 1 and 2. To elaborate, suppose for the sake of argument that there exists $\delta > 0$ and a subsequence $\Lambda \subset \{1, 2,\ldots\}$ and a sequence $\{p_n\}_{n\in\Lambda}$ that converges to p and is such that $|a_n|(p_n) > 2\delta$ for each $n \in \Lambda$. No generality is lost by taking $\delta$ to be less than half the distance from p to $\Theta \cup \partial U$ and then taking $\Lambda$ so that $\text{dist}(p_n, \Theta_{c_n}) > \delta$.

To generate nonsense from a sequence of this sort, use the fact that $|a_n|(p_n) > 2\delta$ and the fact that $\text{dist}(p_n, \Theta_{c_n}) > \delta$ with the analysis in Part 1 to see that each $n \in \Lambda$ version of $r_{c\wedge}$ and $r_{cF}$ is greater than $\kappa_\delta^{-1}$ with $\kappa_\delta > 1$ coming from Lemma 6.3. Let $\hat{N}_n$ now denote the version of the function N that is defined using the point $p_n$ and the data set $(r_n, (A_n, a_n))$. The analysis in Part 2 can be repeated for each sufficiently large $n \in \Lambda$ to see that there exists $r_\diamond \in (0, \frac{1}{100} c^{-1}]$ such that $\hat{N}_n(r_\diamond) \le c^{-4}$. This in turn allows for an appeal to Proposition 3.2 to see that $|a_n|$ must be greater than $\frac{3}{4}|a_n|(p_n)$ on the radius $r_\diamond$ ball centered at $p_n$ when $n \in \Lambda$ is large. This implies in turn that $|a_n|(p) \ge \delta$ when n is large, which is nonsense because $\limsup_{n\in\{1,2,\ldots\}} |a_n|(p) = 0$.

## 7. The data Z, $\mathcal{I}$ and $\nu$, the connection $A_\Delta$ and the homomorphism $\sigma_\Delta$

Supposing that $\{(r_n, \tau_n, (A_n, a_n)\}_{n\in\{1,2,\ldots\}}$ is a weakly convergent sequence of solutions to (2.11) on U, let $|\hat{a}_\diamond|$ denote the limit function from the fourth bullet of (6.1) and let $\Theta \subset U$ denote Proposition 6.1's finite set. Keep in mind in this regard that the set when $|\hat{a}_\diamond| = 0$ is a closed set and so the set where $|\hat{a}_\diamond| > 0$ is an open set. That this is so follows from Proposition 6.1. The lemma that follows directly implies that the zero locus in U of $|\hat{a}_\diamond|$ is nowhere dense. Meanwhile, the upcoming Proposition 7.2 describes the convergence of $\{(A_n, a_n)\}_{n\in\{1.2,\ldots\}}$ at the points in $U-\Theta$ where $|\hat{a}_\diamond|$ is positive.

**Lemma 7.1**: *There exists* $\kappa > 1$ *that depends only on the geometric data and has the following significance: Let* $\{(r_n, \tau_n, (A_n, a_n)\}_{n\in\{1,2,\ldots\}}$ *denote a weakly convergent sequence of solutions to (2.11). Let* $|a_\diamond|$ *denote the function from the fourth bullet of (6.1) and let* $\Theta$ *denote the finite set from Proposition 6.1. Fix* $p \in U-\Theta$ *and define the function* $h_\diamond$ *on* $[0, \kappa^{-1}]$ *by the rule* $r \to h_\diamond(r) = \int_{\partial B_r} |\hat{a}_\diamond|^2$ .

- *The function* $h_\diamond$ *is strictly positive on* $(0, \kappa^{-1}]$.
- *Iff* $s \in (0, \kappa^{-1}]$ *and* $r \in [s, \kappa^{-1}]$, *then* $h_\diamond(r) \ge (\frac{r}{s})^3 e^{-\kappa(r^2-s^2)} h_\diamond(s)$.

This lemma is proved in Section 7a.

The notation in the upcoming Proposition 7.2 uses Z to denote the part of the zero locus $|\hat{a}_\diamond|$ in $X-\Theta$. The notation also has $\tau$ denoting the $n \to \infty$ limit of the sequence $\{\tau_n\}_{n\in\{1,2,\ldots\}}$. By way of a reminder from Section 1c, the term *real line bundle* is used here to describe a vector bundle with fiber $\mathbb{R}$ that is associated to a principal $\mathbb{Z}/2\mathbb{Z}$ bundle.

Keep in mind that a real line bundle has a canonical fiber metric. In addition, the Levi-Civita connection canonically defines a metric compatible covariant derivative on tensor valued sections of a real line bundle. The symbol $\nabla$ is used to denote this covariant derivative; and d is used to denote the corresponding exterior covariant derivative.

**Proposition 7.2**: *Let* $\{(r_n, \tau_n, (A_n, a_n)\}_{n\in\{1,2,\ldots\}}$ *denote a weakly convergent sequence of solutions to (2.11). Let* $|a_\diamond|$ *denote the function from the fourth bullet of (6.1) and let* $\Theta$ *denote the finite set from Proposition 6.1. There exists a real line bundle over* $U-(Z\cup\Theta)$ *to be denoted by* $\mathcal{I}$*; a section of* $\mathcal{I}\otimes T^*X$ *over* $U-(Z\cup\Theta)$*, this denoted by* $\nu$*; a principal SO(3) bundle over* $U-(Z\cup\Theta)$ *to be denoted by* $P_\Delta$*; a connection on* $P_\Delta$ *to be denoted by* $A_\Delta$*; and an isometric bundle homomorphism* $\sigma_\Delta\colon \mathcal{I}\to P_\Delta\times_{SO(3)}\mathfrak{S}$*. Their properties are listed below.*

- $|\nu| = |\hat{a}_\diamond|$.
- *The section* $\nu$ *is harmonic in the sense that both* $d\nu = 0$ *and* $d{*}\nu = 0$.
- $|\nabla\nu|$ *is an* $L^2$ *function on* $U-(Z\cup\Theta)$ *that extends as an* $L^2$ *function to* $U$.
- *For any point* p *in* U*, the function* $|\nabla\nu|(\mathrm{dist}(p,\cdot)^{-1}$ *extends as an* $L^2$ *function on* U *with a bound on its* $L^2$ *norm that is independent of* p.
- *The curvature tensor of* $A_\Delta$ *obeys* $d_{A_\Delta} * F_{A_\Delta} = 0$.
- *The homomorphism* $\sigma_\Delta$ *is* $A_\Delta$*-covariantly constant.*

*There is, in addition, a subsequence* $\Xi \subset \{1, 2, \ldots\}$ *and a sequence* $\{g_n\}_{n\in\Xi}$ *of isomorphisms from* $P_\Delta$ *to* $P|_{U-(Z\cup\Theta)}$ *such that*

i) $\{g_n{}^*a_n\}_{n\in\Xi}$ *converges to* $\nu\sigma_\Delta$ *in the* $L^2_1$ *topology on compact subsets of* $U-(Z\cup\Theta)$ *and in the* $C^0$ *topology on compact subsets of* $U-\Theta$.

ii) $\{g_n{}^*A_n\}_{n\in\Xi}$ *converges on compact subsets of* $U-(Z\cup\Theta)$ *to* $A_\Delta$ *in the* $L^2{}_1$ *topology* .

With regards to what is said in the fifth bullet, the fact that the $d_{A_\Delta} * F_{A_\Delta}$ is zero makes $A_\Delta$ harmonic because $d_{A_\Delta}F_{A_\Delta}$ is already zero by virtue of the Bianchi identity.

Proposition 7.2 is proved in Section 7c using a local version of the proposition that is stated and proved in Section 7b.

## a) Proof of Lemma 7.1

The proof of the lemma has two parts.

*Part 1*: Fix $p \in U-\Theta$ so as to define $h_\diamond$. Suppose that $s \in (0, c_0{}^{-1}]$ is such that $h_\diamond(s) > 0$. This part of the proof explains why $h_\diamond(r) \ge (\frac{r}{s})^3\, e^{-c_0(r^2-s^2)}\, h_\diamond(s)$ when $r \in [s, c_0{}^{-1}]$.

Fix $n \in \Xi$ for the moment and introduce by way of notation $h_n$ to denote the version of the function h from (3.4) that is defined by the point p and the pair $(A_n, a_n)$. It follows from Items c) and d) of the second bullet of (6.1) that the sequence of functions $\{h_n\}_{n\in\Xi}$ converges to $h_\diamond$ in the $C^0$ topology on $[0, c_0^{-1}]$. This is so because the restriction map from the space of smooth functions on a closed ball to the space of functions on the ball's boundary extends as a continous map from the space of $L^2_1$ functions on the ball to the space of $L^2$ functions on the boundary.

Given that $h_\diamond = \lim_{n\to\infty} h_n$, the assertion that $h_\diamond(r) \geq (\frac{r}{s})^3 e^{-c_0(r^2-s^2)} h_\diamond(s)$ follows from the fact that (3.11) holds when $h \in \{h_n\}_{n\in\Xi}$.

*Part 2*: This part of the proof explains why the function $h_\diamond$ is positive on $(0, c_0^{-1}]$. To start the explanation, note that if $h(r) = 0$ for some $r \in (0, c_0^{-1}]$, then h must vanish on the whole of $[0, r]$ because $h_\diamond(r) \geq (\frac{r}{s})^3 e^{-c_0(r^2-s^2)} h_\diamond(s)$. With the preceding understood, suppose for the sake of argument that there is a point in $U-\Theta$ whose version of $h_\diamond$ is zero on an interval in $(0, c_0^{-1}]$. The three steps that follow get nonsense from this assumption.

Step 1: The function $|\hat{a}_\diamond|$ is not identically zero because of Item c) of the second bullet in (6.1). This implies in turn that there exists $p \in U-\Theta$ and numbers $r_\Delta \in (0, c_0^{-1})$ and $r_\diamond \in [0, r_\Delta)$ such that p's version of $h_\diamond$ is zero on $[0, r_\diamond]$ and non-zero on $(r_\diamond, r_\Delta]$. Use $r_\ddagger$ to denote $\frac{1}{2}(r_\diamond + r_\Delta)$ and use $\delta$ to denote the positive square root of $r_\ddagger^{-3} h_\diamond(r_\ddagger)$. Keep in mind that $\delta > 0$.

Step 2: Fix $\varepsilon \in (0, 1]$ for the moment. Since the function $h_\diamond$ is continuous, there exists $r_\varepsilon \in (r_\diamond, r_\ddagger)$ with the property that $r_\varepsilon^{-3} h_\diamond(r_\varepsilon) = \varepsilon^2 \delta^2$. For $n \in \Xi$, let $\kappa_n$ denote the version of Section 3a's function $\kappa$ that is defined by the point p and the pair $(A_n, a_n)$. Since the sequence $\{h_n\}_{n\in\Xi}$ converges to $h_\diamond$, there exists $n_\varepsilon \in \Xi$ such that if $n \in \Xi$ and n if is greater than $n_\varepsilon$, then

$$c_0^{-1} \varepsilon \delta \leq \kappa_n(r_\varepsilon) \leq c_0 \varepsilon \delta \quad \textit{and} \quad \kappa_n(r_\ddagger) > \tfrac{1}{2}\delta. \tag{7.1}$$

Integrate the $(A_n, a_n)$ version of (3.6) using the upper bound for $\kappa_n(r_\varepsilon)$ in (7.1) and the upper bound for $\kappa_n(r_\ddagger)$ in (7.1 to see that

$$1 \leq c_0 \exp\Big( \int_{r_\varepsilon}^{r_\ddagger} \frac{N_n(s)}{s} ds \Big) \varepsilon . \tag{7.2}$$

It follows as a consequence that there exists $r_n \in [r_\varepsilon, r_\ddagger]$ such that $N_n(r_n) \geq c_0^{-1} \frac{r_\diamond}{r_\ddagger} \ln(\varepsilon^{-1})$.

Step 3: The sequence $\{r_\varepsilon \kappa_n(r_\varepsilon)\, r_n\}_{n\in\Xi \text{ and } n > n_\varepsilon}$ has no bounded subsequences because $r_\varepsilon > r_\lozenge$ and because (7.1) bounds $\kappa_n(r_\varepsilon)$ from below by $c_0^{-1}\varepsilon\,\delta$ when n is large. Since Item a) of the second bullet in (6.1) is assumed, the fact that $\{r_\varepsilon \kappa_n(r_\varepsilon)\, r_n\}_{n\in\Xi \text{ and } n > n_\varepsilon}$ has no bounded subsequences is used directly to call on Lemma 5.2 for a lower bound on the large $n \in \Xi$ versions of $N_n$ at value of r greater than $r_n$. To say more, choose $n \in \Xi$ with $n > n_\varepsilon$ and invoke the $(A_n, a_n)$ version of Lemma 5.2 using $r_n$ for what is denoted in the lemma by $r_0$ to see that

$$N_n > c_0^{-1}\, \frac{r_\lozenge}{r_\ddagger}\, |\ln(\varepsilon)| - c_0 r_\Delta^{\,2} \quad \textit{on} \quad [r_n, r_\Delta]. \tag{7.3}$$

Integrate the $(A_n, a_n)$ version of (3.6) on $[r_\ddagger, r_\Delta]$ and then use the bounds in (7.3) and (7.5) to obtain a lower bound for $\kappa_n(r_\Delta)$ that reads

$$\kappa_n(r_\Delta) \geq c_0^{-1}\,\delta\, \Big(\frac{r_\Delta}{r_\ddagger}\Big)^{c_0^{-1}|\ln\varepsilon|\, r_\lozenge / r_\ddagger} \quad . \tag{7.4}$$

This last inequality is nonsensical when ε is small and n is large because $\kappa_n(r_\Delta)$ when n is large is no greater than $2\kappa_\lozenge(r_\Delta)$.

**b) Local convergence where $|\hat{a}_\lozenge|$ is positive**

The proposition that follows uses $\sigma_\lozenge$ to denote a fixed, unit normed element in $\mathfrak{S}$.

**Proposition 7.3**: *Let $\{(r_n, \tau_n, (A_n, a_n)\}_{n\in\{1,2,\ldots\}}$ denote a weakly convergent sequence of solutions to (2.11). Use $|\hat{a}_\lozenge|$ to denote the limit function from the fourth bullet of (6.1) and use Θ to denote the set from Proposition 6.1. If* p *is a point in* U−Θ *and if $|\hat{a}_\lozenge|$ is positive at* p*, then $\{(A_n, a_n)\}_{n\in\Lambda}$ converges on a neighborhood of* p *in the following sense: There is a data set consisting of a ball in* X *centered at* p *to be denoted by* B*, a smooth, harmonic 1-form on* B *to be denoted by $\nu_B$, and a second smooth, coclosed 1-form on* B*, this denoted by $a_{\lozenge B}$. Their properties are listed below.*

- $|\nu_B| = |\hat{a}_\lozenge|$
- $\nabla^\dagger\nabla a_{\lozenge B} + \mathrm{Ric}(a_{\lozenge B}) = 0.$

*There is, in addition, a sequence $\{h_{Bn}\}_{n\in\Lambda}$ of isomorphisms from the product principal* SO(3) *bundle over* B *to* P *with the property that $\{h_{Bn}{}^*a_n\}_{n\in\{1,2,\ldots\}}$ converges on B in the $C^0$ and $L^2_1$ topologies to $\nu_B\,\sigma_\lozenge$ and a subsequence of $\{h_{Bn}{}^*A_n\}_{n\in\{1,2,\ldots\}}$ converges on* B *in the $L^2_1$ topology to $\theta_0 + a_{\lozenge B}\sigma_\lozenge$.*

***Proof of Proposition 7.3***: The proof of this proposition has three parts.

*Part 1*: Let $p \in U-\Theta$ denote a point where $|\hat{a}_\lozenge|(p)$ is positive. Define the number $r_{p\lozenge}$ as in Part 3 of Section 6c and set $r_c = 2r_{p\lozenge}$. Set $\mu$ equal to $\frac{1}{100}$ in what follows. The ball B to be used in Proposition 7.3 can be any ball centered at p with radius less than $(1-\mu)r_c$.

Keeping in mind that $\{|a_n|\}_{n\in\{1,2,\ldots\}}$ converges in the $C^0$ topology to $|\hat{a}_\lozenge|$ on compact subsets of $U-\Theta$, what is said in Part 4 of Section 6c can be repeated to see that each sufficiently large n version of $a_{\ddagger n}$ can be written on $B_{(1-\mu/16)r_c}$ as in (6.22). Fix a positive integer to be denoted by $n_p$ with the property that the $n \geq n_p$ versions of $a_{\ddagger n}$ can be written on $B_{(1-\mu/16)r_c}$ as in (6.22), and in so doing, defining the data $\nu_n$, $\sigma_{\ddagger n}$ and $\mathfrak{a}_n$.

As noted in Part 4 of Section 6c, the norm of any given $n \geq n_p$ version of $|\mathfrak{a}_n|$ is bounded by $c_0(r_c r_{\ddagger n})^{-1} c^{-1/2}$. Meanwhile Part 6 of Section 6c writes $\nu_n$ when $n \geq n_p$ as $t_n + s_n$ with $s_n$ obeying $|s_n| \leq c_0(r_c r_\ddagger)^{-1} c^{-3/2}$ and with $t_n$ obeying $\nabla^\dagger\nabla t_n + \mathrm{Ric}(t_n) = 0$. The arguments in Part 7 of Section 6c prove that $\lim_{n\to\infty} t_n$ converges. Let $t_{(p)}$ again denote the limit, this being a nowhere zero solution on $B_{(1-\mu)r_c}$ to the equation $\nabla^\dagger\nabla t_{(p)} + \mathrm{Ric}(t_{(p)}) = 0$. Define $\kappa_{\lozenge(p)}$ as in Part 7 of Section 6c to be the positive square root of the number that is depicted in (6.34). What Proposition 6.1 denotes as $\nu_B$ is the 1-form $\kappa_{\lozenge(p)}\, t_{(p)}$.

*Part 2*: Fix $n \geq n_p$ and let $\hat{A}_n$ denote the $(r_n, (A_n, a_n))$ version of the connection on P's restriction to $B_{(1-\mu/16)r_c}$ that is depicted in (4.9). Use $h_{Bn}$ to denote the isomorphism between the product principal SO(3) bundle over $B_{(1-\mu/16)r_c}$ and P that writes $\sigma_{\ddagger n}$ as the constant element $\sigma_\lozenge \in \mathfrak{S}$ and pulls $\hat{A}_n$ back as $\theta_0 + a_{\lozenge n}\sigma_\lozenge$ with $a_{\lozenge n}$ being a coclosed 1-form on $B_{(1-\mu/16)r_c}$ that obeys the $a_\lozenge = a_{\lozenge n}$ version of (4.12). The notation used subsequently does not distinguish $A_n$ from its pull-back via $h_{Bn}$, nor is $\mathfrak{a}_n$ distinguished from its $h_{Bn}$ pull-back. What Proposition 7.2 denotes as $\sigma_B$ is $\sigma_\lozenge$.

Fix $n \geq n_p$. The pull-back to $B_{(1-\mu/16)r_c}$ of the connection $A_n$ by the isomorphism $h_{Bn}$ can be written as written as $\theta_0 + a_{\lozenge n}\sigma_\lozenge - \frac{1}{4}[\sigma_\lozenge, h_n{}^*(\nabla_{A_n}\sigma_{\lozenge n})]$. The bounds in the $r = r_n$ and $(A, a) = (A_n, a_n)$ version of (4.27) lead to a $c_\mu(r_c r_{\ddagger n})^{-1} c^{-1/2}$ bound for the $L^2{}_1$ norm of the $\mathfrak{S}$-valued 1-form $h_{Bn}{}^*A_n - \theta_0 - a_{\lozenge n}\sigma_\lozenge$ on the ball $B_{(1-\mu/4)r_c}$. The various $n \geq n_p$ versions of this bound imply that the sequence $\{h_{Bn}{}^*A_n - \theta_0 - a_{\lozenge n}\sigma_\lozenge\}_{n\geq n_p}$ converges to zero in the $L^2{}_1$ topology on $B_{(1-\mu/4)r_c}$.

To see about the convergence of the sequence $\{a_{\lozenge n}\}_{n\geq n_p}$, note first the $a_\lozenge = a_{\lozenge n}$ versions of (4.12) for $n \geq n_p$ supply an a priori $L^2{}_1$ bound for this sequence on $B_{(1-\mu/16)r_c}$. The two steps that follow explain why there is a subsequence of $\{a_{\lozenge n}\}_{n\geq n_p}$ that converges strongly in the $L^2{}_1$ topology on $B_{(1-\mu)r_c}$ and why the limit, $a_{\lozenge B}$, obeys $\nabla^\dagger\nabla a_{\lozenge B} + \mathrm{Ric}(a_{\lozenge B}) = 0$.

Step 1: Fix $n \geq n_p$. The $\sigma_\Diamond$ component of the pull-back via $h_{Bn}$ of the equation in (3.28) can be written on $B_{(1-\mu/16)r_c}$ as

$$*d*d\, a_{\Diamond n} = r_{\ddagger n}{}^2 \langle \sigma_{\Diamond n} [\mathfrak{a}_{n\alpha} \nabla_{A_n} \mathfrak{a}_{n\alpha}] \rangle . \tag{7.5}$$

Use $b_n$ to denote $\mathfrak{S}$-valued 1-form on $B_{(1-\mu)r_c}$ that vanishes on $\partial B_{(1-\mu)r_c}$ and obeys the equation

$$\nabla^\dagger \nabla b_n + \mathrm{Ric}(b_n) = r_{\ddagger n}{}^2 \langle \sigma_{\Diamond n} [\mathfrak{a}_{n\alpha} \nabla_{A_n} \mathfrak{a}_{n\alpha}] \rangle . \tag{7.6}$$

If $c > c_0$, then there exists a unique solution of (7.6) that vanishes on $\partial B_{(1-\mu)r_c}$. Moreover, the $L^2_1$ norm of this solution is bounded by $c_0 r_c$ times the $L^2$ norm on $\partial B_{(1-\mu)r_c}$ of the term on the right hand side of (7.6). Meanwhile, the latter $L^2$ norm is observedly no greater than $c_0 r_{\ddagger n}{}^2$ times the product of the sup-norm of $|\mathfrak{a}_n|$ on $\partial B_{(1-\mu)r_c}$ and the $L^2$ norm of $\nabla_{A_n} \mathfrak{a}_n$ on $B_{(1-\mu)r_c}$. Part 4 of Section 6c bounds $|\mathfrak{a}_n|$ by $c_\mu (r_c r_{\ddagger n})^{-1} c^{-1/2}$ and Lemma 4.3 bounds the $L^2$ norm of $\nabla_{A_n} \mathfrak{a}_n$ by $c_\mu r_{\ddagger n}{}^{-1} (r_c r_{\ddagger n})^{-1} c^{-1/2}$; and so it follows as a consequence that there is a $c_\mu (r_c r_{\ddagger n})^{-1} c^{-1}$ bound on the $L^2_1$ norm of $b_n$.

Step 2: Write $a_{\Diamond n}$ on $B_{(1-\mu)r_c}$ as $w_{\Diamond n} + b_n$. Since $a_{\Diamond n}$ is coclosed, the equations in (7.5) and (7.6) lead to an equation for $w_{\Diamond n}$ on $B_{(1-\mu)r_c}$ that reads

$$\nabla^\dagger \nabla w_{\Diamond n} + \mathrm{Ric}(w_{\Diamond n}) = 0 . \tag{7.7}$$

Since the sequence $\{b_n\}_{n \geq n_p}$ converges to zero on $B_{(1-\mu)r_c}$ in the $L^2_1$ topology and the sequence $\{a_{\Diamond n}\}_{n \geq n_p}$ is bounded on $B_{(1-\mu)r_c}$ in the $L^2_1$ topology, it follows that the sequence $\{w_{\Diamond n}\}_{n \geq n_p}$ is bounded on $B_{(1-\mu)r_c}$ in the $L^2_1$ topology. This being the case, it has a subsequence that converges weakly in the $L^2_1$ topology. The limit is what Proposition 7.3 refers to as $a_B$. Standard elliptic regularity theorems can be invoked using the fact that each $n \geq n_p$ version of $w_{\Diamond n}$ obeys (7.7) to see $a_B$ obeys the same equation as each $n \geq n_p$ version of $w_{\Diamond n}$ and that the subsequence of $\{w_{\Diamond n}\}_{n \geq n_p}$ that converges weakly in the $L^2_1$ topology to $a_B$ converges in the $C^\infty$ topology to $a_B$ on compact subsets of $B_{(1-\mu)r_c}$.

*Part 3*: As noted previously the sequences $\{\mathfrak{a}_n\}_{n\ge n_p}$ and $\{s_n\}_{n\ge n_p}$ converge pointwise to zero on $B_{(1-\mu)r_c}$. These sequences and the sequence $\{\nabla_{A_n}\sigma_n\}_{n\ge n_p}$ also converge to zero on $B_{(1-\mu)r_c}$ in the $L^2_1$ topology. That this is so for $\{\mathfrak{a}_n\}_{n\ge n_p}$ and $\{\nabla_{A_n}\sigma_n\}_{n\ge n_p}$ follows directly from the versions of (4.27) that use the elements from the sequence $\{(r = r_n, (A = A_n, a = a_n))\}_{n\ge n_p}$. The fact that $\{s_n\}_{n\ge n_p}$ converges to zero on $B_{(1-\mu)r_c}$ in the $L^2_1$ topology is proved by first taking the inner product of both sides of any given $n \ge n_p$ version of (6.31) with $s_n$. Integrate the resulting equation and integrate by parts on the left hand side to obtain the integral of $|\nabla s_n|^2$. The $c_\mu(r_c r_{\ddagger n})^{-1} c^{-3/2}$ pointwise bound for $|s_n|$ and the $c_\mu(r_c r_{\ddagger n})^{-2} c^{-1}$ bound in Part 6 of Section 6c for the number $\mathrm{I}_n$ leads to a $c_\mu (r_c r_{\ddagger n})^{-2} c^{-5/2}$ bound for the integral on the right hand side.

It follows from what was just said about $\{\mathfrak{a}_n\}_{n\ge n_p}$, $\{\nabla_{A_n}\sigma_n\}_{n\ge n_p}$ and $\{s_n\}_{n\ge n_p}$ that the sequence $\{h_n{}^*a_n\}_{n\in\{1,2,\ldots\}}$ converges on $B_{(1-\mu)r_c}$ in the $C^0$ and $L^2_1$ topologies to $\nu_B\sigma_B$. The fact that $\nu_B$ is harmonic follows Item a) of the second bullet in (6.1) because the sequences $\{\mathfrak{a}_n\}_{n\ge n_p}$, $\{\nabla_{A_n}\sigma_n\}_{n\ge n_p}$, $\{s_n\}_{n\ge n_p}$ and the sequence $\{h_{Bn}{}^*A_n - \theta_0 - a_{\lozenge n}\sigma_\lozenge\}_{n\ge n_p}$ from Part 2 all converges to zero on $B_{(1-\mu)r_c}$ in the $L^2_1$ topology,

**b) Proof of Proposition 7.2**

The proof of the proposition has five parts.

*Part 1*: Fix a locally finite cover of $U-(Z\cup\Theta)$ by open balls of the sort that are described in Proposition 7.2. Denote this cover by $\mathfrak{U}$. If B is a ball from $\mathfrak{U}$, then Proposition 7.3 assigns to B a pair $(\nu_B, a_{\lozenge B})$ of 1-forms on B having specified properties, one being that $|\nu_B| = |\hat{a}_\lozenge|$.

Proposition 7.3 also specifies a sequence $\{h_{Bn}\}_{n\in\Lambda}$ with each member being an isomorphism from the product principal SO(3) bundle over B to $P|_B$ such that $\{h_{Bn}{}^*a_n\}_{n\in\{1,2,\ldots\}}$ converges to $\nu_B\sigma_\lozenge$ and such that a subsequence of $\{h_{Bn}{}^*A_n\}_{n\in\{1,2,\ldots\}}$ converges to $\theta_0 + a_{\lozenge B}\sigma_\lozenge$. This subsequence can be taken to be independent of the ball B from $\mathfrak{U}$ because the cover of $\mathfrak{U}$ is necessarily countable. Let $\Lambda$ denote such a subsequence

Let B and B´ denote two balls from $\mathfrak{U}$ that intersect. Fix $n \in \Lambda$ and use $h_{BB'n}$ to denote $h_{Bn}(h_{B'n})^{-1}$, this being an automorphism of the product principal SO(3) bundle over $B \cap B'$. The method of convergence described in Proposition 7.3 implies that $\{h_{BBn}\sigma_\lozenge(h_{BBn})^{-1}\}_{n\in\Lambda}$ converges to either $\sigma_\lozenge$ or $-\sigma_\lozenge$ in the $L^2_1$ topology on compact subsets of $B\cap B'$. This implies in particular that either $\nu_B = \nu_{B'}$ or $\nu_B = -\nu_{B'}$ on $B \cap B'$. If it is the former, set $\iota_{BB'} = 1$ and if it is the latter, set $\iota_{BB'} = -1$.

The data $\{\iota_{BB'}\}_{B,B'\in\mathfrak{U} \text{ and } B\cap B'\neq\emptyset}$ is the cocycle data that defines a principal $\mathbb{Z}/2\mathbb{Z}$ bundle over U−(Z∪Θ). The associated real line bundle is the bundle $\mathcal{I}$ to use for Proposition 7.2. Meanwhile, the data $\{\nu_B\}_{B\in\mathfrak{U}}$ defines a section over U−(Z∪Θ) of $\mathcal{I}\otimes T^*X$. This is the section of $\mathcal{I}\otimes T^*X$ to use for Proposition 7.2 Denote it by $\nu$.

*Part 2*: What is said by the first bullet of Proposition 7.2 to the effect that $|\nu| = |\hat{a}_\diamond|$ follows directly from the first bullet of Proposition 7.3. The assertion of the second bullet of Proposition 7.2 follows from the third, fourth and fifth bullet of Proposition 7.3 and the fact that each $B \in \mathfrak{U}$ version of $\nu_B$ is coclosed.

The proofs of the third bullet and fourth bullets of Proposition 7.2 starts by returning to the milieu of the proof of Proposition 7.3. By way of a reminder of what is done their, let B denote a ball in $\mathfrak{U}$ and let p denote its center point. The radius of B is a number less than $(1-\mu)r_c$ with $\mu$ being $\frac{1}{100}$ and with $r_c$ being a number that is determined by $|\hat{a}_\diamond|(p)$ and a choice for $c$ that is greater than $c_0$. Part 1 of the proof of Proposition 7.3 defines a positive integer to be denoted by $n_p$ which is such that when $n \geq n_p$, then $a_{\ddagger n}$ on $B_{(1-\mu/16)r_c}$ can be written as in (6.22). The latter defines $\nu_n$. Since (4.7) holds, each $n \geq n_p$ version of $|\nabla\nu_n|$ is pointwise bounded on B by $c_0|\nabla_{A_n} a_{\ddagger n}|$. Granted that this is so, the fact that $\nu_B$ is the limit of the sequence $\{K_{\ddagger n}\nu_n\}_{n\geq n_p}$ has the following implication: Let $f$ denote a smooth, non-negative, $C^2$ function on U. Then

$$\int_{U-(Z\cup\Theta)} f\,|\nabla\nu|^2 \leq c_0 \lim_{n\to\infty} \int_{U-(Z\cup\Theta)} f\,|\nabla_{A_n} a_n|^2 \;. \tag{7.8}$$

Note that the left hand side of (7.8) is finite, this being a consequence of Item b) of the second bullet in (6.1). The assertion made by the third bullet of Proposition 7.2 follows from the $f = 1$ version in (7.8). The fact that $|\nabla\nu|$ extends to U−Z as an $L^2$ function follows from the fact that Θ is a finite set. Extend this function as 0 over Z to define the extension of $|\nabla\nu|$ to the whole of U

The fourth bullet of Proposition 7.2 also follows from (7.8) and the dominated convergence theorem. To how this comes about, fix $p \in U$ and then fix $k \in \{1, 2, \ldots\}$ so as to define the function $f_k$ on U by the rule

$$x \to f_k(x) = \frac{1}{\mathrm{dist}(p,x)^2}\,\chi(2 - k\,\mathrm{dist}(p,x)) \;. \tag{7.9}$$

The function $f_k$ is no greater than $\frac{1}{\mathrm{dist}(p,\cdot)^2}$ and it is equal to this function where the distance to p is greater than $2k^{-1}$. Meanwhile, the function $f_k$ is equal to zero where the distance to p is less than $k^{-1}$. Item d) of the fourth bullet in (6.1) implies among other things that there is a k independent bound for the limit on the right hand side of (7.8) when $f = f_k$. The fourth bullet of Proposition 7.2 follows from the preceding observation using the dominated convergence theorem.

*Part 3*: The principal bundle $P_\Delta$ is defined by cocycle data that is associated to the cover $\mathfrak{U}$ of $U-(Z\cup\Theta)$. By way of a reminder, the cocycle data associates to a given ordered pair (B, B´) of intersecting balls from U a map $h_{BB'}$: $B\cap B' \to SO(3)$. This assignment is such that $h_{BB'} = h_{B'B}^{-1}$ and such that if B, B´ and B´´ are three balls from U with points in common, then $h_{BB'}h_{B'B''}h_{B''B} = 1$ on $B\cap B'\cap B''$.

Let $\mathfrak{U}_2 \subset \mathfrak{U}\times\mathfrak{U}$ denote the set of pairs that have non-empty intersection. The four steps that follows prove that there is a subsequence in {1, 2, …}, this being $\Xi$, such that if (B,B´) is from $\mathfrak{U}_2$, then$\{h_{BB'n}\}_{n\in\Xi}$ converges in the $L^2_2$ and $C^0$ topology to a smooth map from $B \cap B'$ to SO(3). Denote this map by $h_{BB'}$. A fifth step explains why the set $\{h_{BB'}\}_{(B,B')\in\mathfrak{U}_2}$ defines cocycle data for a smooth principal SO(3) bundle over $U-(Z\cup\Theta)$. The latter bundle is $P_\Delta$.

<u>Step 1</u>: Fix $B \in \mathfrak{U}$ and for $n \in \Lambda$ and large, write $h_{Bn}{}^*A_n$ as

$$h_{Bn}{}^*A_n = \theta_0 + a_{\lozenge B}\,\sigma_\lozenge + A_{Bn} \tag{7.10}$$

with $\{A_{Bn}\}_{n\in\Lambda}$ being a sequence of $\mathfrak{S}$-valued , coclosed 1-forms that converges to zero in the $L^2_1$ topology on B. Supposing that $B' \in \mathfrak{U}$ and that $B \cap B'$ is non-empty, then comparing the B and B´ versions of (7.10) leads to an equation that identifies

$$dh_{BB'n} = - a_{\lozenge B}\,\sigma_\lozenge h_{BB'n} - A_{Bn}h_{BB'n} + a_{\lozenge B'}\,h_{BB'n}\,\sigma_\lozenge + h_{BB'n}A_{B'n}\,{}_{\cdot p}{}^{-2}(1 + c_0\,r_{*p}{}^2)\ . \tag{7.11}$$

Since $|h_{BB'n}| = 1$ for all n, and since $\{A_{Bn}\}_{n\in\Lambda}$ and $\{A_{B'n}\}_{n\in\Lambda}$ converge to zero in the $L^2_1$ topology on $B \cap B'$, the equation in (7.11) with (2.1) imply that there is a subsequence $\Xi_{BB'} \subset \Lambda$ such that $\{h_{BB'n}\}_{n\in\Xi_{BB'}}$ converges in the $L^2_2$ topology on $B \cap B'$.

<u>Step 2</u>: Fix $n \in \Xi_{BB'}$ and act on both sides of the corresponding version of (7.11) by $d^\dagger$. Having done so, use the fact that $a_{\lozenge B}$, $a_{\lozenge B'}$, $A_{Bn}$ and $A_{Bn'}$ are coclosed to see that $h_{BB'n}$ obeys an equation that has the schematic form

$$d^\dagger d\,h_{BB'n} = Q_{BB'n} \tag{7.12}$$

with the norm of $Q_{BB'n}$ such that $|Q_{BB'n}| \le c_0\,|dh_{BB'n}|(1+|A_{Bn}|+|A_{B'n}|)$. Fix $p \in B \cap B'$ and let $G_p$ denote the Dirichelet Green's function for B with pole at p. Use $\chi$ to construct a smooth, non-negative function with compact support on $B \cap B'$ that is equal to 1 near p.

Denote this function by $\chi_p$. Multiply both sides of (7.12) by $\chi_p G_p$ and then integrate. An integration by parts leads to the identity

$$h_{BB'n}(p) = \int_{B\cap B'} (\chi_p G_p Q_{BB'n} + (d^\dagger d\chi_p G_p + 2\mathfrak{m}(d\chi_p, dG_p)) h_{BB'n}) \,. \tag{7.13}$$

with $\mathfrak{m}(\cdot\,,\cdot)$ used here to denote the metric inner product.

<u>Step 3</u>: Let $\varepsilon_{BB'n} \in (0, \infty)$ denote an upper bound for the $L^2{}_1$ norms of $A_{Bn}$ and $A_{B'n}$ and let $C_{BB'}$ denote an upper bound for the $L^2{}_2$ norm of $h_{BB'n}$. Supposing that p´ is a second point in $B \cap B'$, subtract the p and p´ versions of (7.13) and use the Step 2's bounds for $Q_{BB'n}$ with what is said by the second bullet in (2.1) to see that

$$|h_{BB'n}(p) - h_{BB'n}(p')| \le c_0 (\mathrm{dist}(p,p')^{1/2} + \varepsilon_{BB'n}) C_{BB'} \,. \tag{7.14}$$

Since $\lim_{n\in\Xi_{BB'}} \varepsilon_{BB'n} = 0$, this last equation implies that the sequence $\{h_{BB'n}\}_{n\in\Xi_{BB'}}$ converges in the $C^0$ topology on compact subsets of $B \cap B'$.

<u>Step 4</u>: Since $\mathfrak{U}_2$ is countable, an enumeration of $\mathfrak{U}_2$ can be used to construct a single subsequence to be denoted by $\Xi$ such that if $(B, B') \in \mathfrak{U}_2$, then $\{h_{BB'n}\}_{n\in\Xi}$ converges in the $L^2{}_2$ topology on $B \cap B'$ and in the $C^0$ topology on compact subsets of $B \cap B'$ to an $L^2{}_2$ and $C^0$ map from $B \cap B'$ to SO(3). Denote the limit by $h_{BB'}$. The next paragraph explains why this map to SO(3) is smooth.

The manner of convergence of any $B \in \mathfrak{U}$ version of $\{h_{Bn}{}^*A_n\}_{n\in\Xi}$ to the connection $\theta_0 + a_{\Diamond B}\sigma_\Diamond$ by Proposition 7.3 and the manner of convergence of any $(B,B') \in \mathfrak{U}_2$ version of $\{h_{BB'n}\}_{n\in\Xi}$ to $h_{BB'}$ as described in Part 3 imply that

$$a_{\Diamond B}\sigma_\Diamond = a_{\Diamond B'} h_{BB'}\sigma_\Diamond (h_{BB'})^{-1} - dh_{BB'}\,(h_{BB'})^{-1} \,. \tag{7.15}$$

Since $a_{\Diamond B}$ and $a_{\Diamond B'}$ are smooth and $h_{BB'}$ is continuous, the equation in (7.15) implies that $h_{BB'}$ is smooth.

By way of a parenthetical remark, it follows from what is said in Part 1 that

$$h_{BB'}\sigma_\Diamond (h_{BB'})^{-1} = \iota_{BB'}\sigma_\Diamond. \tag{7.16}$$

This understood, then it follows from (7.15) that $dh_{BB'}(h_{BB'})^{-1}$ can be written as $df_{BB'}\sigma_\Diamond$ with $f_{BB'}$ being a harmonic function on $B \cap B'$.

Step 5: Since $h_{BB'n} = h_{B'Bn}^{-1}$ for all $n \in \Xi$, it follows that $h_{BB'} = h_{B'B}^{-1}$. Much the same argument proves that $h_{BB'}h_{B'B''}h_{B''B} = 1$ on $B \cap B' \cap B''$ when the latter set is not empty. It follows as a consequence that $\{h_{BB'}\}_{(B,B')\in\mathfrak{U}_2}$ define cocycle data for a smooth principal SO(3) bundle over $U-(Z\cup\Theta)$.

*Part 4*: If $B \in \mathfrak{U}$, then $a_{\lozenge B}\sigma_\lozenge$ is a $\mathfrak{S}$-valued 1-form on B; and it follows from (7.15) that the collection $\{a_{\lozenge B}\sigma_\lozenge\}_{B\in\mathfrak{U}}$ of $\mathfrak{S}$-valued 1-forms is the cocycle data for a connection on $P_\Delta$. This is the connection $A_\Delta$. The assertion in the fifth bullet of Proposition 7.2 to the effect that $d_{A_\Delta} * F_{A_\Delta} = 0$ follows from the second bullet in Proposition 7.3.

To define the homomorphism $\sigma_\Delta$ for the sixth bullet of Proposition 7.2, first fix $B \in \mathfrak{U}$ and define a homomorphism $\sigma_B: \mathbb{R} \to \mathfrak{S}$ by the rule $t \to t\sigma_\lozenge$. Since $\{\iota_{BB'}\}_{(B,B')\in\mathfrak{U}_2}$ is the cocycle data for $\mathcal{I}$ and $\{h_{BB'}\}_{(B,B')\in\mathfrak{U}_2}$ is the cocycle data for $P_\Delta$, it follows from (7.16) that the collection $\{\sigma_B\}_{B\in\mathfrak{U}}$ define a homomorphism from $\mathcal{I}$ to $P_\Delta\times_{SO(3)}\mathfrak{S}$. The fact that this homomorphism is $A_\Delta$-covariantly constant follows from the fact that $\{a_{\lozenge B}\sigma_\lozenge\}_{B\in\mathfrak{U}}$ defines the cocycle data for $A_\Delta$.

*Part 5*: The three steps that follow address the convergence assertion made by Items i) and ii) at the end of Proposition 7.2.

Step 1: The arguments in Section 3 of [U] (see Proposition 3.2 in [U]) can copied with only cosmetic changes to construct a sequence of $\{h_n\}_{n\in\Xi}$ from $P_\Delta$ to $P|_{U-(Z\cup\Theta)}$ such that $\{h_n{}^*A_n\}_{n\in\Xi}$ converges on compact subsets of $U-(Z\cup\Theta)$ to $A_\Delta$. Note in this regard that this appropriation of Uhlenbeck's arguments requires the $C^0$ convergence of each $(B, B') \in \mathfrak{U}_2$ version of the sequence $\{h_{BB'n}\}_{n\in\Xi}$ to $h_{BB'}$. The $L^2_2$ convergence is, of course, also needed as is the $L^2_1$ convergence of each $B \in \mathfrak{U}$ version of $\{h_{Bn}{}^*A_n\}_{n\in\Xi}$, but these are not sufficient.

By way of a quick summary, this appropriation of Uhlenbeck's arguments starts by fixing a cover of $U-(Z\cup\Theta)$ by open balls whose element are in 1-1 coorrespondence with the elements in $\mathfrak{U}$. The correspondence is such that a ball from this new cover has the same center but smaller radius than its partner in $\mathfrak{U}$. Let $\mathfrak{U}_*$ denote this new cover. If $B \in \mathfrak{U}$, then its partner in $\mathfrak{U}_*$ is denoted by $B_*$.

To continue, fix $n \in \Xi$ and write any given $(B, B') \in \mathfrak{U}_2$ version of $h_{BB'n}$ as $h_{BB'}\exp(u_{BB'n})$ with $u_{BB'n}$ being a map from $B\cap B'$ to $\mathfrak{S}$. This map has small $L^2_2$ norm when n is large, and pointwise small on $B_*\cap B_*'$. This is to say that the sequence $\{u_{BB'n}\}_{n\in\Xi \text{ and nis large}}$ converges to zero on $B_*\cap B_*'$ in the $C^0$ and $L^2_2$ topologies. Since the cover $\mathfrak{U}$ is locally finite, the arguments in [U] for Proposition 3.2 can be used to write $h_{BB'n}$ on $B_* \cap B_*'$ when n is sufficiently large as

$$h_{BB'n} = \exp(q_{Bn})\, h_{BB'} \exp(-q_{B'n}) \tag{7.17}$$

with $q_{Bn}$ being a smooth map from $B_*$ to $\mathfrak{G}$ and $q_{B'n}$ being a smooth map from $B_*'$ to $\mathfrak{G}$. Moreover, both are such that their $L^2_2$ norms and $C^0$ norms are bounded by those of $u_{BB'n}$ on $B_* \cap B_*'$.

Given $B \in \mathfrak{U}$, let $n_B$ be such that the map $q_{Bn}$ is defined for $n \ge n_B$. It follows from what was said at the end of the preceding paragraph that the sequence $\{q_{Bn}\}_{n\in\Xi \text{ and } n>n_B}$ converges to zero in the $C^0$ and $L^2_1$ topologies on $B_*$.

Step 2: Since a compact subset of $U-(Z\cup\Theta)$ is covered by at most a finite set of balls from $\mathfrak{U}$, the identity in (7.17) can be assumed to hold on all balls from the latter subset of $\mathfrak{U}$ when n is large. The corresponding versions of $\exp(q_{Bn})$ with B from this subset of $\mathfrak{U}$ and n large then define the desired bundle isomorphism $h_n$ between $P_\Delta$ and P over this compact set.

With the preceding understood, the fact that $\{h_n^*A_n\}_{n\in\Xi}$ converges on compact set to $A_\Delta$ follows from what is said by Proposition 7.3 about the $L^2_1$ convergence of each $B \in \mathfrak{U}$ version of $\{h_{Bn}^*A_n\}_{n\in\Xi}$ using the fact that each $B \in \mathfrak{U}$ version of the sequence $\{q_{Bn}\}_{n\in\Xi \text{ and } n>n_B}$ converges to zero on $B_*$ in both the $C^0$ and $L^2_2$ topologies.

Step 3: The fact that $\{h_n^*a_n\}_{n\in\Xi}$ converges to $\nu\sigma_\Delta$ in the $C^0$ topology on compact subsets of $U-(Z\cup\Theta)$ follows from two facts, the first being that each $B \in \mathfrak{U}$ version of $\{q_{Bn}\}_{n\in\Xi \text{ and } n>n_B}$ converges to zero in the $C^0$ topology on compact subsets of B. The second fact was noted by Proposition 7.3, this being that each $B \in \mathfrak{U}$ version of $\{h_{Bn}^*a_n\}_{n\in\{1,2,\ldots\}}$ converges to $\nu\sigma_\Diamond$ in the $C^0$ topology on B.

The assertion that $C^0$ convergence occurs on compact subsets of $U-\Theta$ follows from what was said in the preceding paragraph with two additional facts, the first being that $|\nu| = |\hat{a}_\Diamond|$; and the second being Proposition 6.1's statement that $\{|a_n|\}_{n\in\{1,2,\ldots\}}$ converges to $|\hat{a}_\Diamond|$ in the $C^0$ topology on compact subsets of $U-\Theta$.

The assertion that $\{h_n^*a_n\}_{n\in\Xi}$ converges to $\nu\sigma_\Delta$ in the $L^2_1$ topology on compact subsets of $U-(Z\cup\Theta)$ follows from two facts, the first being Proposition 7.3's statement that each $B \in \mathfrak{U}$ version of $\{h_{Bn}^*a_n\}_{n\in\{1,2,\ldots\}}$ converges to $\nu\sigma_\Diamond$ in the $L^2_1$ topology on B; and the second being that each $B \in \mathfrak{U}$ version of $\{q_{Bn}\}_{n\in\Xi \text{ and } n>n_B}$ converges to zero in the $L^2_1$ topology on $B_*$.

## 8. Hölder continuity of $|\nu|$

The proposition that follows implies that the the norm of the 1-form $\nu$ from Proposition 7.2 is Hölder continuous on compact subsets of $X-\Theta$.

**Proposition 8.1**: *There exists* $\kappa > 1$ *that depends only on the geometric data and with the following significance: Supposing that* $\{(r_n, \tau_n, (A_n, a_n)\}_{n\in\{1,2,\ldots\}}$ *is a weakly convergent sequence of solutions to (2.11), let* $|a_\diamond|$ *denote the function from the fourth bullet of (6.1), let* $\Theta$ *denote the finite set from Proposition 6.1, and let* $Z$, $\mathcal{I}$ *and* $\nu$ *denote corresponding data that from Proposition 7.2. If* $B \subset X-\Theta$ *is a ball with compact closure and if* $p, q$ *are any two points in* $B$*, then* $||\nu|(p) - |\nu|(q)| \le x_B \,\mathrm{dist}(p,q)^{1/\kappa}$ *with* $x_B$ *depending only on* $B$.

The proof of this proposition is in Section 8c. The intervening subsections introduce a version of Almgren's frequency function for $\nu$ that is used in the proof.

### a) The definition of the frequency function

To set the notation used below, fix for the moment a point $p \in X$. Given $p$, there exists a positive number to be denoted by $r_p$ such that the ball of radius $r_p$ centered at $X$ is well inside a Gaussian coordinate chart centered at $p$ for neighborhood of $p$ in $X$. This number $r_p$ is assumed to be less than 1.

With $p \in X$ chosen, two auxilliary functions on $[0, r_p]$ are needed in order to define the analog of Almgren's frequency function. The first of function is the analog for $\nu$ of the function in (3.3). This version is denoted by $h_\diamond$ and it is defined by the rule

$$r \to h_\diamond(r) = \int_{\partial B_r} |\nu|^2 \quad . \tag{8.1}$$

The lemma that follows states two relevant facts about the function $h_\diamond$. The first bullet of Lemma 7.1 says that $h_\diamond > 0$ on $(0, r_p)$ and the second bullet gives a lower bound for the ratio of values of $h_\diamond$. The second of the required function is denoted by $\mathfrak{d}_\diamond$ and it is defined in below in (8.2). It is the analog for $\nu$ of the function that is depicted in (3.4). As in (3.4), what (8.2) denotes by $\mathrm{M}$ is defined by writing the trace of the second fundamental form of a given $r \in (0, c_0^{-1}]$ version of $\partial B_r$ as $\frac{3}{r} + \mathrm{M}$. Note that $|\mathrm{M}(r)| \le c_0 r^2$. The rule that defines $\mathfrak{d}_\diamond$ is

$$r \to \mathfrak{d}_\diamond(r) = \int_0^r \Big(\frac{1}{h_\diamond(s)}\Big(\int_{B_s} \mathrm{Ric}(\nu\otimes\nu) + \tfrac{1}{2}\int_{\partial B_s} \mathrm{M}|\nu|^2\Big)\Big)\, ds \quad . \tag{8.2}$$

Note that $\mathfrak{d}_\diamond$ is continuous and it obeys $|\mathfrak{d}_\diamond(r)| \le c_0 r^2$. These assertions are consequences of Lemma 7.1 because that lemma implies that the integral of $\mathrm{Ric}(\nu\otimes\nu)$ over any given radius $s \in (0, r_p]$ ball is no greater than $c_0\, s\, h_\diamond(s)$.

With $h_\diamond$ and $\eth_\diamond$ in hand, define the function $\kappa_\diamond$ on $[0, r_p]$ to be the positive square root of $e^{-2\eth_\diamond} r^{-3} h_\diamond$. It follows from Lemma 7.1, what was said about $\eth_\diamond$ and the fact that $|\nu|$ is continuous that this function $\kappa_\diamond$ is continuous and bounded by $c_0$. If $p \in U-\Theta$, then

$$\kappa_\diamond(0) = \sqrt{2}\,\pi|\nu|(p)\ ; \tag{8.3}$$

this is because $|\nu|$ is a continuous function on U.

The desired analog of Almgren's frequency function is denoted by $N_\diamond$. It is the function on $(0, r_p]$ that is defined by the rule

$$r \to N_\diamond(r) = \frac{1}{r^2\kappa_\diamond(r)^2} \int_{B_r} |\nabla\nu|^2 \tag{8.4}$$

with it understood that $|\nabla\nu|$ is defined by $\nu$ on $U-(Z\cup\Theta)$ and then extended over Z and $\Theta$ by fiat as 0. The next lemma states some salient properties of $N_\diamond$.

**Lemma 8.2**: *Given an open set* $U \subset X$ *with compact closure, there exists* $\kappa > 1$ *with the following significance: Fix* $p \in U$ *and define the functions* $\kappa_\diamond$ *and* $N_\diamond$. *The function* $N_\diamond$ *is continuous on* $(0, \kappa^{-1}]$ *and it extends as a continuous function to* $[0, \kappa^{-1}]$. *Moreover,*

- $\frac{d}{dr}\kappa_\diamond = \frac{1}{r} N_\diamond \kappa_\diamond$ *on* $(0, \kappa^{-1}]$.
- *If* $s \in (0, \kappa^{-1}]$ *and* $r \in [s, \kappa^{-1}]$, *then* $N_\diamond(r) \geq N_\diamond(s) - \kappa(r^2 - s^2)$.
- *If* $p \in U-\Theta$ *and* $|\nu|(p) > 0$, *then* $\lim_{r\to 0} N_\diamond(r) = 0$.
- *If* $p \in U-\Theta$ *and* $|\nu|(p) = 0$, *then* $\lim_{r\to 0} N_\diamond(r) > \kappa^{-1}$.

Lemma 8.2 is proved in Section 2b. Lemma 8.2 implies that $N_\diamond$ extends as a continuous function to the closed interval $[0, c_0^{-1}]$.

**b) Proof of Lemma 8.2**

The proof of this lemma has six parts. The arguments use implicitly the fact that there is a positive upper bound for the norms of the Riemannian curvature tensor on any given compact set in X.

*Part 1*: Given $\rho \in (0, 1]$, let $Z_\rho$ denote the set of points in $U-\Theta$ where $|\nu| < \rho$. As explained directly, the function $|\nabla\nu|$ is such that

$$\lim_{\rho\to 0} \int_{Z_\rho - Z} |\nabla\nu|^2 = 0\ . \tag{8.5}$$

The limit in (8.5) is zero because $|\nabla\nu|$ is an $L^2$ function on $U-Z$ and because the volume of $Z_\rho - Z$ gets ever smaller with limit zero as $\rho$ limits to 0. This is also why

$$\lim_{\rho\to 0} \int_{\operatorname{dist}(\cdot,\Theta)<\rho} |\nabla\nu|^2 = 0 \,. \tag{8.6}$$

The fact that $|\nabla\nu|$ is an $L^2$ function on $U-(Z\cup\Theta)$ and the fact that the limits in (8.5) and (8.6) are zero imply that $N_\lozenge$ is continuous on the interval $(0, r_p]$. The fact that it extends as a continuous function to $[0, r_p]$ follows from the second bullet of the lemma and the fact that $N_\lozenge$ is non-negative.

*Part 2*: The five steps that follow derive the identity in the lemma's first bullet.

Step 1: The fact that $\nu$ is harmonic leads to the equation

$$\nabla^\dagger\nabla\nu + \operatorname{Ric}(\nu) = 0 \tag{8.7}$$

on $U-(Z\cup\Theta)$. This implies in turn that $|\nu|^2$ obeys the equation

$$d^\dagger d|\nu|^2 + |\nabla\nu|^2 + \operatorname{Ric}(\nu,\nu) = 0 \,. \tag{8.8}$$

Step 2: To exploit (8.8), fix $\rho \in (0, 1]$ and use $\chi_\rho$ in what follows to denote the function $\chi(2(1-\rho^{-1}|\nu|))$. This function is equal to 0 where $|\nu| < \frac{1}{2}\rho$ and it is equal to 1 where $|\nu| \geq \rho$. Use $\upsilon_\rho$ to denote the function $\prod_{q\in\Theta}\chi(2(1-\rho^{-1}\operatorname{dist}(q,\cdot))$. The latter function equals 1 where the distance to each point in $\Theta$ is greater than $\rho$ and it is equal to 0 where the distance to some point in $\Theta$ is less than $\frac{1}{2}\rho$.

Step 3: Fix $p \in U-\Theta$ and $r \in (0, r_p]$. Take the inner product of both sides of (8.8) with $(\chi_\rho\upsilon_\rho)^2$ and integrate both sides of the result over $B_r-(B_r\cap(Z\cup\Theta))$. Having done so, integrate by parts to obtain an identity that can be written as

$$\frac{1}{2}\int_{\partial B_r} \nabla_r(\chi_\rho^2\upsilon_\rho^2\,|\nu|^2) = \int_{B_r-(B_r\cap(Z\cup\Theta))} (|\nabla\nu|^2 + \operatorname{Ric}(\nu,\nu)) + Q_\rho \tag{8.9}$$

with $Q_\rho$ denoting an integral with absolute value bounded by

$$c_0 \int_{Z_\rho - Z} (|\nabla\nu|^2 + |\nu|^2) \; + c_0 \int_{\mathrm{dist}(\cdot,\Theta)<\rho} (|\nabla\nu|^2 + \rho^{-2}|\nu|^2) \;, \tag{8.10}$$

and with $\nabla_r$ denoting the directional derivative along the outward pointing unit length normal vector to $\partial B_r$.

Step 4: The left hand side of (8.10) can be written as

$$\tfrac{1}{2} r^3 \tfrac{d}{dr} (r^{-3} \int_{\partial B_r} \chi_\rho^2 \upsilon_\rho^2 \, |\nu|^2 ) - \tfrac{1}{2} \int_{\partial B_r} \chi_\rho^2 \upsilon_\rho^2 M \, |\nu|^2 \;. \tag{8.11}$$

To put this formula in perspective, let $h_{\diamond\rho}$ denote for the moment the function on $[0, r_p]$ whose value at any given $r \in [0, r_p]$ is the integral over $\partial B_r$ of the function $\chi_\rho{}^2 \upsilon_\rho{}^2 |\nu|^2$. Use (8.3) and (8.11) to write (8.9) as

$$\tfrac{1}{2} r^3 \tfrac{d}{dr} (r^{-3} h_{\diamond\rho}) \;=\; \int_{B_r - (B_r \cap (Z \cup \Theta))} |\nabla\nu|^2 \;+\; h_\diamond(r) \left( \tfrac{d}{dr} \mathfrak{d}_\diamond \right) + Q_\rho + \mathfrak{q}_\rho \tag{8.12}$$

with $\mathfrak{q}_\rho$ being the integral over $\partial B_r$ of $(1 - \chi_\rho{}^2 \upsilon_\rho{}^2) |\nu|^2$.

Step 5: The $\mathfrak{q}_\rho$ term in (8.12) is no greater than $c_0 r^3 \rho^2$ and so its $\rho \to 0$ limit is 0. Meanwhile, the fact that $\rho$ is bounded and the fact that the limits in (8.5) and (8.6) are zero imply that $\lim_{\rho\to 0} Q_\rho = 0$ also. Granted that these limits are zero, the identity in (8.12) with the fact that $|h_\diamond(r) - h_{\diamond\rho}(r)| \le c_0 r^{-3} \rho^2$ has two implications; the first being that the function $r^{-3} h_\diamond$ is differentiable, and the second being that its derivative is given by

$$\tfrac{1}{2} r^3 \tfrac{d}{dr} (r^{-3} h_\diamond) \;=\; \int_{B_r - (B_r \cap (Z \cup \Theta))} |\nabla\nu|^2 \;+\; h_\diamond(r) \left( \tfrac{d}{dr} \mathfrak{d}_\diamond \right) \tag{8.13}$$

This last identity when written in terms of $\mathrm{K}_\diamond$ is the assertion in the lemma's first bullet.

*Part 3*: Since $\mathrm{K}_\diamond$ is positive on $(0, r_p{}^{-1}]$, the assertion made by the second bullet of Lemma 8.2 in the case when $p \notin \Theta$ follows from Lemma 5.2 and the next lemma. The notation in this next lemma uses $\mathrm{N}_n$ to denote the version of the function $\mathrm{N}$ from Section 3a that is defined by the given point p and the data set $(r_n, (A_n, a_n))$

**Lemma 8.4**: *Fix* $p \in U - \Theta$ *so as to define the functions* $\mathrm{K}_\diamond$ *and* $\mathrm{N}_\diamond$. *If* $r \in (0, r_p]$, *then* $\lim_{n \in \Xi} \mathrm{N}_n(r) = \mathrm{N}_\diamond(r)$.

This lemma is proved momentarily.

The second bullet of Lemma 2.2 holds at any given point $p \in \Theta$ if it holds at all points in a neighborhood of p. This is because the function on U that assigns to a given point in U the value of the point's version of $N_\Diamond(r)$ for a *fixed* $r \in (0, c_0^{-1}]$ is a continuous function on U. To prove the latter claim, fix $r \in (0, c_0^{-1}]$ and define a function on U by the rule that assigns to each point the integral of $|\nabla\nu|^2$ over the radius r ball centered at that point. This is a continuous function because $|\nabla\nu|$ is an $L^2$ function. Meanwhile, the function that assigns to a given point the integral of $\mathrm{Ric}(\nu,\nu)$ over the radius r ball about point is also continuous. By the same token, the function that assigns to each point the integral of $|\nu|^2$ over the boundary of the closed, radius r ball centered at the point is also continuous as is the function that assigns the integral of $M|\nu|^2$. These observations imply that the terms that appear in the numerator and denominator in (8.4)'s definition of any fixed $r \in (0, c_0^{-1}]$ version of $N_\Diamond(r)$ vary continuously as functions of the chosen point in U.

The lemma that follows asserts a fact that is used in the proof of Lemma 8.4.

**Lemma 8.5**: *There exists* $\kappa > 1$ *that depends only on the geometric data and has the following significance: Suppose that* $(r, \tau, (A, a))$ *is a solution to (1.2). Given* $\delta \in (0, 1]$, *let* $U_\delta$ *denote the subset of* U *where* $|a| < \delta$. *Fix* $p \in U$ *and* $r \in (0, r_p]$. *Then*

$$\int_{B_r \cap U_\delta} (|\nabla_A a|^2 + r^2|a \wedge a|^2) \le \kappa\delta^\kappa \Big( \int_{B_{2r}} (|\nabla_A a|^2 + r^2|a \wedge a|^2) + r^2 \Big).$$

The proof of Lemma 8.5 is in Part 4 of this subsection. Accept it as true for now.

***Proof of Lemma 8.4***: Since r is positive, $K_\Diamond(r)$ is positive and so $\lim_{n\to\infty} K_n(r) = K_\Diamond(r)$. It follows as a consequence that Lemma 8.4's assertion is true if

$$\lim_{n\in\Xi} \int_{B_r} (|\nabla_{A_n} a_n|^2 + r_n^2 |a_n \wedge a_n|^2) = \int_{B_r - (B_r \cap (Z \cup \Theta))} |\nabla\nu|^2 \; . \tag{8.14}$$

To prove (8.14), fix for the moment $\rho \in (0, 1]$. Let $\Theta_\rho$ denote the union of the radius $\rho$ balls centered at the points in $\Theta$. Let $A_\Delta$ denote the connection from Proposition 7.2. It follows from what is said by Item ii) of Proposition 7.2 that

$$\lim_{n\in\Xi} \int_{B_r - (B_r \cap (Z_\rho \cap \Theta_\rho))} |F_{A_n}|^2 = \int_{B_r - (B_r \cap (Z_\rho \cap \Theta_\rho))} |F_{A_\Delta}|^2 \; . \tag{8.15}$$

The important point to take away from (8.15) is that the limit on the left hand side exists. Since Item a) of the second bullet in (6.1) is assumed, it follows from (8.15) that

$$\lim_{n\in\Xi} \int_{B_r-(B_r\cap(Z_\rho\cap\Theta_\rho))} r_n^2|a_n\wedge a_n|^2 = 0 \,. \tag{8.16}$$

Meanwhile, what is said by Items i) and ii) of Propositions 7.2 imply that

$$\lim_{n\in\Xi} \int_{B_r-(B_r\cap(Z_\rho\cap\Theta_\rho))} |\nabla_{A_n} a_n|^2 = \int_{B_r-(B_r\cap(Z_\rho\cup\Theta_\rho))} |\nabla v|^2 \,. \tag{8.17}$$

The identity in (8.14) follows from the various $\rho > 0$ versions of (8.16) and (8.17) given that the limits on the left hand side of (8.5) and (8.6) are zero, and given Lemma 8.5 and given Item d) of the fourth bullet of Proposition 2.2

*Part 4*: This part of the subsection contains the

***Proof of Lemma 8.5***: With $\delta > 0$ given, use $x_\delta$ to denote the function $\chi(\delta^{-1}|a| - 1)$. This function is equal to 1 where $|a| < \delta$ and it is equal to 0 where $|a| > 2\delta$. Use u in this proof to denote the function $\chi(r^{-1}\mathrm{dist}(\cdot, p) - 1)$. Note in particular that u is equal to 1 where the distance to p is less than r and equal to 0 where the distance to p is greater than 2r. Let $f$ denote the function on [0, 1] given by the rule

$$\delta \to f(\delta) = \int_{B_{2r}} x_\delta u^2(|\nabla_A a|^2 + r^2|a\wedge a|^2) \tag{8.18}$$

The identity in (2.6) asserts that if $(r, \tau, (A, a))$ obey (1.2), then $|a|^2$ obeys the equation

$$\tfrac{1}{2} d^\dagger d|a|^2 + |\nabla_A a|^2 + 2r^2|a\wedge a|^2 + \mathrm{Ric}(\langle a\otimes a\rangle) = 0. \tag{8.19}$$

Multiply both sides of this identity by $x_\delta u^2$ and integrate over the support of u. Integrate by parts and use the fact that the support of $dx_\delta$ lies where $\delta \le |a| \le 2\delta$ to derive an inequality involving $f_\delta$ and $f_{2\delta}$ asserting that

$$f_\delta \le \tfrac{c}{1+c} f_{2\delta} + c_0\delta^2 \tag{8.20}$$

with c being positive, independent of $\delta$ and less than $c_0$. Let $N_\delta$ denote the greatest integer that is less than $c_0^{-1}\,\frac{|\ln(\delta)|}{\ln(2)}$. The inequality in (8.20) can be iterated $N_\delta$ times, the first to express $f_{2\delta}$ in terms of $f_{4\delta}$, the next to express $f_{4\delta}$ in terms of $f_{8\delta}$, and so on with the result being an inequality that implies directly the one that is asserted by Lemma 2.4.

*Part 5*: To prove the third bullet of Lemma 8.2, first fix $r_\diamond \in (0, c_0^{-1}]$ and then fix $r \in (0, r_\diamond]$. Use the formula in the the first bullet of Lemma 8.2 to write

$$K_\diamond(r_\diamond) = \exp\Big(\int_r^{r_\diamond} \frac{N_\diamond(s)}{s}\, ds\Big)\, K_\diamond(r)\ . \tag{8.21}$$

With this formula understood, suppose that $\varepsilon \in (0, 1)$ and that $r_\varepsilon \in (0, r_\diamond)$ is such that $N_\diamond(r_\varepsilon) \geq \varepsilon$. The second bullet of Lemma 8.2 says that $N_\diamond(r) \geq \varepsilon - c_0\,(r^2 - r_\varepsilon^{\,2})$ for $r \geq r_\varepsilon$ and this implies in turn that

$$\int_{r_\varepsilon}^{r_\diamond} \frac{N_\diamond(s)}{s}\, ds \;\geq \varepsilon \ln\Big(\frac{r_\diamond}{r_\varepsilon}\Big) - c_0. \tag{8.22}$$

The function $K_\diamond$ at $r = r_\varepsilon$ obeys $K_\diamond(r_\varepsilon) \geq \sqrt{2}\,\pi|\nu|(p)$ because of (2.3) and because $K_\diamond$ is non-decreasing. Since $K_\diamond$ is in any event is no greater than $c_0\, \sup_X|\nu|$, these bounds with (2.22) used in (2.21) leads to the lower bound $r_\varepsilon \geq c_0^{-1}\,\big(\frac{|\nu|(p)}{\sup_X|\nu|}\big)^{1/\varepsilon}\, r_\diamond$ when $|\nu|(p) > 0$.

*Part 6*: The proof of Lemma 8.2’s fourth bullet has eight steps. Steps 1 and 2 prove that the fourth bullet of Lemma 8.2 is true if there is a subsequence in $\Xi$ with the criteria for membership being that $a_n(p) = 0$. Step 3 proves that the fourth bullet of Lemma 8.2 is true if $|a_n|(p) > 0$ for $n \in \Xi$ sufficiently large and if there is a subsequence in $\Xi$ with a positive lower bound for the corresponding versions of the number $r_{c\diamond}$ that is defined in (6.7) using a suitable choice for the number $c$. The remaining steps prove that the fourth bullet of Lemma 8.2 is true if $|a_n|(p) > 0$ for all large n from $\Xi$ and if there is no subsequence of the sort just described.

<u>Step 1</u>: Let $p \in Z\!-\!\Theta$, thus a point where $|\nu| = 0$. Keeping in mind that $|\nu| = |\hat{a}_\diamond|$, it follows that $\lim_{n\in\Xi} |a_n|(p) = 0$. For $n \in \Xi$, let $K_n$ and $N_n$ denote the versions of the functions K and N from Section 3 that are defined using p and $(A_n, a_n)$. If $a_n(p) = 0$, then the limit as $r \to 0$ of $N_n(r)$ is a positive integer. Assuming that this is the case, it then

follows from Lemma 6.3 and Proposition 5.1 that $\mathrm{N}_n(\mathrm{r}) \ge \mathrm{c}_0^{-1}$ at points $\mathrm{r} \in (0, \mathrm{c}_0^{-1}]$ where $\mathrm{r}\,\mathrm{K}_n(\mathrm{r})\, r_n \ge z_U^{-1}$.

Step 2: Fix $\mathrm{r} \in (0, \mathrm{c}_0^{-1}]$. If n is sufficiently large, then $\mathrm{K}_n(\mathrm{r})$ is greater than $\frac{1}{2}\mathrm{K}_\Diamond(\mathrm{r})$, and thus $\mathrm{r}\,\mathrm{K}_n(\mathrm{r})\, r_n$ will be larger than $\frac{1}{2}\mathrm{r}\,\mathrm{K}_\Diamond(\mathrm{r})\, r_n$. This in turn is greater than $z_U^{-1}$ when n is very large because $\mathrm{K}_\Diamond(\mathrm{r})$ is positive. If it is also the case that $a_n(\mathrm{p}) = 0$, then it follows from what is said in Step 1 that $\mathrm{N}_n(\mathrm{r})$ is greater than $\mathrm{c}_0^{-1}$. If there is a subsequence $\vartheta \subset \Xi$ with $a_n(\mathrm{p}) = 0$, then it follows from Lemma 8.4 that $\mathrm{N}_\Diamond(\mathrm{r})$ is greater than $\mathrm{c}_0^{-1}$ also.

The preceding observations prove the fourth bullet of Lemma 8.2 if there is a subsequence of $\Xi$ with any given integer n being a member if $a_n(\mathrm{p}) = 0$.

Step 3: Suppose henceforth that $|a_n|(\mathrm{p}) > 0$ if n is from $\Xi$ and sufficiently large. Fix $\mathrm{n} \in \Xi$ so that $|a_n|(\mathrm{p}) > 0$. Fix $c > \mathrm{c}_0$ as done in the proof of Proposition I.6.1 and define the number $\mathrm{r}_{c\Diamond}$ as in (I.6.7) using the data set $(r = r_n, (\mathrm{A} = \mathrm{A}_n, a = a_n))$. Denote this version of $\mathrm{r}_{c\Diamond}$ by $\mathrm{r}_{c\Diamond(n)}$. This step explains why $\lim_{n\to\infty} \mathrm{r}_{c\Diamond(n)} = 0$ unless $\lim_{\mathrm{r}\to 0} \mathrm{N}_\Diamond(\mathrm{r}) \ge \mathrm{c}_0^{-1} c^{-2}$.

The latter claim is proved by assuming it false so as to generate nonsense. To this end, fix $z > 1$ for the moment and suppose that $\lim_{\mathrm{r}\to 0} \mathrm{N}_\Diamond(\mathrm{r}) \le z^{-1} c^{-2}$. Suppose in addition that there is a subsequence $\vartheta \subset \Xi$ and $\delta > 0$ such that $\mathrm{r}_{c\Diamond(n)} > \delta$ when $\mathrm{n} \in \vartheta$. Fix $\mathrm{r}_c \in (0, \delta)$ so as to be less than $\frac{1}{100}\,\mathrm{dist}(\mathrm{p}, \Theta)$. Assume in addition that $\mathrm{N}_\Diamond(\mathrm{r}_c) \le 2 z^{-1} c^{-2}$.

Since $\mathrm{N}_\Diamond(\mathrm{r}_c) \le 2 z^{-1} c^{-2}$, it follows from Lemma 8.4 that $\mathrm{N}_n(\mathrm{r}_c) < 3 z^{-1} c^{-2}$ when $\mathrm{n} \in \Xi$ is large. Meanwhile, if $\mathrm{n} \in \vartheta$, then $\mathrm{r}_c$ is less than the versions of $\mathrm{r}_{c\wedge}$ and $\mathrm{r}_{cF}$ that are defined by p and the data set $(r = r_n, (\mathrm{A} = \mathrm{A}_n, a = a_n))$. If $z > \mathrm{c}_0$, then this $\mathrm{N}_n(\mathrm{r}_c) \le 3 z^{-1} c^{-2}$ bound enables the $(r = r_n, (\mathrm{A} = \mathrm{A}_n, a = a_n))$ version of Proposition 3.2. In particular, if $z > \mathrm{c}_0$, then Proposition 3.2 with $\mu = \frac{1}{100}$ asserts that $|a_n| \le 2\,|a_n|(\mathrm{p})$ on the whole of the radius $(1 - \mu)\,\mathrm{r}_c$ ball centered at p. Since $\lim_{n\in\Xi} |a_n|(\mathrm{p}) = 0$, it follows that $|\nu| = 0$ on the whole radius $(1 - \mu)\,\mathrm{r}_c$ ball centered at p. This last conclusion is the sought after nonsense because it runs afoul of Lemma 7.1.

The conclusions of the fourth bullet of the Lemma 8.2 hold if $\lim_{\mathrm{r}\to 0} \mathrm{N}_\Diamond(\mathrm{r}) \ge \mathrm{c}_0^{-1} c^{-2}$, this being a corollary of the lemma's second bullet.

Step 4: Suppose that $|a_n|(\mathrm{p}) > 0$ if $\mathrm{n} \in \Xi$ is large and that $\lim_{n\to 0} \mathrm{r}_{c\Diamond(n)} = 0$. If this is so, then $\mathrm{r}_{c\Diamond(n)} < \frac{1}{1000}\,\mathrm{dist}(\mathrm{p}, \Theta)$ when n is large. Fix $\mathrm{n} \in \Xi$ so that $|a_n|(\mathrm{p})$ is positive and so that the preceding bound holds. What is said in Part 2 of Section 6b can be repeated to find a point to be denoted by $\mathrm{q}_n$ with distance at most $3\mathrm{r}_{c\Diamond(n)}$ from p that obeys the upcoming (8.23), this being an analog of (6.5). The notation has $r_{c\Diamond(n)q}$ denoting the $q_n$ version of $r_{c\Diamond}$ that is defined using the data $r = r_n$ and $(A, a) = (A_n, a_n)$.

*If $q'$ has distance $2r_{c\Diamond(n)q}$ or less from $q_n$, then the version of $r_{c\Diamond}$ that is defined by $q'$ using the data set $(r = r_n, (A = A_n, a = a_n))$ is greater than $\frac{1}{100} r_{c\Diamond(n)q}$ .*

(8.23)

Keep in mind for what follows that $r_{c\Diamond(n)q}$ is in any event no greater than $4r_{c\Diamond(n)}$. This is so because $q_n$ has distance at most $3r_{c\Diamond(n)}$ from p.

With $q_n$ understood, the argument in Part 3 of Section 6b can be repeated with only notational changes to prove the next assertion.

*Fix $\varepsilon \in (0, c_0^{-1} c^{-4})$. If n is large, there is a point with distance at most $\frac{3}{2} r_{c\Diamond(n)q}$ from $q_n$ with the following property: Let $\hat{N}_n$ denote the version of the function N as defined using this point and the data set $(r = r_n, (A = A_n, a = a_n))$. The value of $\hat{N}_n$ at $r = \frac{1}{100} r_{c\Diamond(n)q}$ is greater than the minimum $\kappa_{*\varepsilon}^{-1}$ and $c^{-2}$.*

(8.24)

By way of a reminder, $\kappa_{*\varepsilon}$ is given by Proposition 4.1. Fix $z > c_0$ so that $\varepsilon_* = z^{-1} c^{-4}$ can be used in (8.24). Supposing that $n \in \Xi$ is large, use $p_n$ to denote a point that is described by (8.24) using $\varepsilon = \varepsilon_*$.

<u>Step 5</u>: Fix $n_*$ so that if $n \in \Xi$ is larger than $n_*$, then $p_n$ is defined. For $n > n_*$, define $r_{*p_n}$ to be the value of r where the $p_n$ and $(r = r_n, (A = A_n, a = a_n))$ version of the function $\kappa$ obeys $r\kappa(r) r_n = z_U^{-1}$. It is necessarily the case that $\lim_{n\in\Xi \text{ and } n > n_*} r_{*p_n} = 0$. To see that this is so, fix $n \in \Xi$ with $n > n_*$ and let $\hat{\kappa}_n$ denote the version of the function $\kappa$ that is defined by $p_n$ and the data set $(r = r_n, (A = A_n, a = a_n))$. If there exists $r \in (0, c_0^{-1}]$ and a subsequence $\vartheta \subset \Xi$ with the property that $r_{*p_n} > r$, then $\lim_{n\in\vartheta} \hat{\kappa}_n(r)$ would be equal to zero. Since the function $\kappa$ is non-decreasing, this event can occur only if

$$\lim_{n\in\vartheta} \int_{\mathrm{dist}(\cdot, p_n) < r} |a_n|^2 = 0$$

(8.25)

Since $\{p_n\}_{n\in\Xi \text{ and } n > n_*}$ converges to p, the latter event would run afoul of Lemma 7.1.

<u>Step 6</u>: Let E denote twice the limit of the $f = 1$ version of the left most sequence that appears in Item a) of the second bullet of (6.1). Use $\varepsilon$ in what follows to denote the minimum of $E^{-1}$, $\kappa_{*\varepsilon_*}^{-1}$ and $c^{-2}$. Fix $r \in (0, c_0^{-1}]$. It follows from the definition of $\{p_n\}_{n\in\Xi \text{ and } n > n_*}$ in Step 4 and from what is said in Step 5 that Lemma 6.3 and Proposition 5.1 can be invoked when n is large to prove that $\hat{N}_n(r) > \kappa_\varepsilon^{-1}$. Note in this regard that the data $p_n$ and ($r = r_n$, ($A = A_n$, $a = a_n$)) obeys the assumptions of Lemma 6.3 and Proposition 5.1 when n is large because the radius $c_0 r_{*p_n}$ ball centered on $p_n$ will be disjoint from $\Theta$ when n is sufficiently large.

<u>Step</u> <u>7</u>: The lower bound from Step 6 on $\hat{N}_n(r)$ can be written to say that

$$\int_{\mathrm{dist}(\cdot,p_n)<r} (|\nabla_{A_n} a_n|^2 + 2r_n^2 |a_n \wedge a_n|^2) \geq \kappa_\varepsilon^{-1} r^2 \hat{K}_n(r)^2 . \tag{8.26}$$

If n is large, then the radius r ball about $p_n$ is contained in the radius 2r ball about p. This being the case, (8.26) implies in turn that

$$\int_{\mathrm{dist}(\cdot,p)<2r} (|\nabla_{A_n} a_n|^2 + 2r_n^2 |a_n \wedge a_n|^2) \geq \kappa_\varepsilon^{-1} r^2 \hat{K}_n(r)^2 . \tag{8.27}$$

This last inequality can be rewritten to say that

$$N_n(2r) \geq c_0^{-1} \kappa_\varepsilon^{-1} \left(\frac{\hat{K}_n(r)}{K_n(2r)}\right)^2 . \tag{8.28}$$

The inequality in (8.28) is invoked momentarily.

<u>Step</u> <u>8</u>: To see about the ratio in (8.28), note first that

$$\int_{\frac{1}{2}r < \mathrm{dist}(\cdot,p_n) < r} |a_n|^2 \leq \tfrac{1}{4} r^4 \hat{K}_n(r)^2 (1 + c_0 r^2) \tag{8.29}$$

because the function $\hat{K}_n$ is increasing. As the function $K_n$ is also increasing,

$$\int_{\frac{5}{8}r < \mathrm{dist}(\cdot,p) < \frac{7}{8}r} |a_n|^2 \geq \tfrac{1}{64} r^4 K_n(\tfrac{1}{2} r)^2 (1 - c_0 r^2) . \tag{8.30}$$

If n is large, the integration domain for the integral in (8.29) contains that for the integral in (8.30). This understood, these inequalities imply that

$$\hat{K}_n(r) \geq c_0^{-1} K_n(\tfrac{1}{2} r) . \tag{8.31}$$

Meanwhile, the first and third bullets of Lemma 8.2 can be invoked to prove that

$$K_n(\tfrac{1}{2}r) \geq (1 - c_0 r^2)\, 4^{-N_n(2r)}\, K_n(2r) \ .$$

(8.32)

The inequalities (8.31) and (8.32) lead to a lower bound for the ratio that appears on the right hand side of (8.28). Use this lower bound to see that

$$2^{4N_n(2r)}\, N_n(2r) \ \geq c_0^{-1}\, \kappa_\varepsilon^{-1} \ .$$

(8.33)

The bound in (8.33) implies in turn a bound of the form $N_n(2r) \geq c_0^{-1}\, \kappa_\varepsilon^{-1}$.

Lemma 8.4 with the fact that $N_n(2r) \geq c_0^{-1}\, \kappa_\varepsilon^{-1}$ holds for all large $n \in \Xi$ implies that $N_\lozenge(2r) \geq c_0^{-1}\, \kappa_\varepsilon^{-1}$. This bound proves the fourth bullet of Lemma 8.2 when $|a_n|(p) > 0$ for all large $n \in \Xi$ and $\lim_{n\in\Xi} r_{c\lozenge(n)} = 0$ because it holds for any given $r \in (0, c_0^{-1}]$.

**c) Proof of Proposition 8.1**

Given the first and fourth bullets of Lemma 8.2, the argument is virtually identical to the argument in Section 3e of [T1]. For the sake of completeness, this two part argument is given below.

*Part 1*: Differentiating the equation in (8.7) leads to an equation that has the form

$$\nabla^\dagger\nabla(\nabla\nu) + \mathcal{R}_1 \nabla\nu + \mathcal{R}_0 \nu \ = 0$$

(8.34)

where $\mathcal{R}_1$ and $\mathcal{R}_0$ having bounded norms. This identity will be used momentarily.

Given $\delta \in (0,1]$ and let $\beta_\delta$ denote the function $\chi(2 - \delta^{-1}\mathrm{dist}(\cdot, Z))$. This function is equal to 1 where the distance to Z is greater than $2\delta$ and it is equal to zero where the distance to Z is less than $\delta$. Let $B_*$ denote for the moment a given ball with compact closure in $X-\Theta$. Use p to its center and let $r_*$ denote its radius. Define a function $\gamma$ to be $\chi(2r_*^{-1}\mathrm{dist}(\cdot, p) - 1)$. This function is equal to 1 where the distance to p is less than $\frac{1}{2}r_*$ and equal to 0 where the distance to p is greater than $r_*$. The equation in (3.18) leads to the $L^2$ norm bound

$$\int_X \beta_\delta^2 \gamma\, |(\nabla^\dagger\nabla)\nabla\nu|^2 \leq c_0 \, .$$

(8.35)

Let $B' \subset B_*$ denote the concentric ball with half the radius of $B_*$. An integration by parts then leads from (8.35) to an $L^2_3$ bound for $\nu$ of the form

$$\int_{B'} |\nabla^{\otimes 3}(\beta_\delta \nu)|^2 \ \leq c_*\, \delta^{-6} \, ,$$

(8.36)

with $c_*$ here and in what follows denoting a number that is greater than 1 and in any event bounded by $c_0(1+r_*)^{c_0}$. The precise value of $c_*$ can be assumed to increase between successive appearances.

The bound in (8.36) leads via a standard Sobolev inequality to a $c_*\delta^{-3}$ bound for the $L^8$ norm $\beta_\delta|\nabla\nu|$ on B´. The latter bound gives a corresponding $c_B\delta^{-3}$ bound for the $L^8$ norm of the 1-form $d|\beta_\delta\upsilon|$ on B´; and this last $c_*\delta^{-3}$ bound leads in turn (via another Sobolev inequality) to the exponent $\frac{1}{4}$ Hölder norm bound

$$||\nu(p)|-|\nu(q)|| \le c_*\,\delta^{-3}\,\mathrm{dist}(p,q)^{1/4} \tag{8.37}$$

if both p and q are points in B´ with distance 2δ or more from Z.

*Part 2*: Meanwhile, it follows from the second and fourth bullets of Lemma 8.2 that if $\delta \in (0, \frac{1}{100}r_0)$ and if p has distance less than δ from Z, then

$$|\nu(p)| \le c_0\,\delta^{\upsilon} \tag{8.38}$$

with υ positive and independent of δ and p and $B_*$. Of course, the analogous inequality holds for $|\nu(q)|$ if q has distance less than δ from Z.

Label p and q so that $\mathrm{dist}(p,Z) \le \mathrm{dist}(q,Z)$. If $\mathrm{dist}(p,q)^{1/16} \le \mathrm{dist}(p,Z)$, then it follows from (8.37), that

$$||\nu(p)|-|\nu(q)|| \le c_B|\mathrm{dist}(p,q)|^{1/16}. \tag{8.39}$$

Suppose on the other hand that $\mathrm{dist}(p,q)^{1/16} \ge \mathrm{dist}(p,Z)$. Assume in addition that $\mathrm{dist}(q,Z) > 10\,\mathrm{dist}(p,Z)$. This assumption implies that $\mathrm{dist}(q,Z) \le \frac{10}{9}\,\mathrm{dist}(p,q)$; and since

$$||\nu(p)|-|\nu(q)|| \le |\nu(p)|+|\nu(q)|, \tag{8.40}$$

it follows from (8.38) and its $|\nu(q)|$ analog that

$$||\nu(p)|-|\nu(q)|| \le c_0\,\mathrm{dist}(p,q)^{\upsilon}. \tag{8.41}$$

The last case to consider has $\mathrm{dist}(p,q)^{1/8} \ge \mathrm{dist}(p,Z)$ and $\mathrm{dist}(q,Z) \le 10\,\mathrm{dist}(p,Z)$. Granted these inequalities, it follows from (8.38), its $|\nu|(q)$ analog and from (8.41) that

$$||\nu(p)|-|\nu(q)|| \le c_0\,\mathrm{dist}(p,q)^{\upsilon/16}. \tag{8.42}$$

The claimed Hölder continuity of $|v|$ follows from (8.39), (8.41) and (8.42).